\ifx \labelsONmargin\undefined
  \ifx\RequirePackage\undefined
    \expandafter\ifx\csname amsart.sty\endcsname\relax
        \documentstyle[12pt,amssymb,amscd]{amsart}
        \def\mathbb{\Bbb}
        \def\mathfrak{\frak}

        \makeatletter
        \newenvironment{proof}{%
          \@ifnextchar[{%
                        \expandafter\let\expandafter\end@proof
                          \csname endpf*\endcsname
                        \my@proof
                       }{\let\end@proof\endpf\pf}%
        }{\end@proof}
        \def\my@proof[#1]{\@nameuse{pf*}{#1}}
        \makeatother

        \def\xrightarrow[#1]#2{@>{#2}>{#1}>}
        \def\xleftarrow[#1]#2{@<{#2}<{#1}<}
    \fi
  \else
      \ifx\usepackage\RequirePackage
      \else
        \documentclass[12pt]{amsart}

        \advance\oddsidemargin -0.1\textwidth
        \evensidemargin\oddsidemargin
        \usepackage{amssymb,amscd}
        \def\Sb#1\endSb{_{\substack{#1}}}
        \def\Sp#1\endSp{^{\substack{#1}}}
      \fi
  \fi
\makeatletter
\@ifundefined{DeclareFontShape} 
     {\@ifundefined{selectfont}
             {
                \message{We need AmS-LaTeX, this version of LaTeX does not 
                        support it}\@@end
             }{
                \def\mathcal{\cal}
                
                \def\pcyr{%
                        \def\default@family{UWCyr}%
                        \let\oldSl@\sl
                        \def\sl{\def\default@shape{it}\oldSl@}%
                        \cyracc
                        \language\Russian\family{UWCyr}\selectfont
                }
             }}{
                \DeclareFontEncoding{OT2}{}{\noaccents@}
                \DeclareFontSubstitution{OT2}{cmr}{m}{n}
                \DeclareFontFamily{OT2}{cmr}{\hyphenchar\font45 }
                \DeclareFontShape{OT2}{cmr}{m}{n}{%
                     <5><6><7><8><9><10>gen*wncyr %
                     <10.95><12><14.4><17.28><20.74><24.88> wncyr10 %
                }{}
                \DeclareFontShape{OT2}{cmr}{m}{it}{%
                     <5><6><7><8><9><10> gen * wncyi%
                     <10.95><12><14.4><17.28><20.74><24.88> wncyi10%
                }{}
                \DeclareFontShape{OT2}{cmr}{bx}{n}{%
                     <5><6><7><8><9><10> gen * wncyb%
                     <10.95><12><14.4><17.28><20.74><24.88> wncyb10%
                }{}
                \DeclareFontShape{OT2}{cmr}{m}{sl}{%
                     <-> ssub * cmr/m/it%
                }{}
                \DeclareFontShape{OT2}{cmr}{m}{sc}{%
                     <5><6><7><8><9><10>%
                     <10.95><12><14.4><17.28><20.74><24.88> wncysc10%
                }{}
                \DeclareFontFamily{OT2}{cmss}{\hyphenchar\font45 }
                \DeclareFontShape{OT2}{cmss}{m}{n}{%
                     <8><9><10> gen * wncyss%
                     <10.95><12><14.4><17.28><20.74><24.88> wncyss10%
                }{}
                
                \def\cyrencodingdefault{OT2}
                \def\pcyr{%
                        \cyracc
                        \let\encodingdefault\cyrencodingdefault
                        \language\Russian\fontencoding{OT2}\selectfont
                }
             }   

\@ifundefined{theorembodyfont} 
     {
       \let\theorembodyfont\thm@bodyfont
       \let\normalshape\relax
     }{}   

\makeatletter
  \def\@sect@my#1#2#3#4#5#6[#7]#8{%
\ifnum #2>\c@secnumdepth
   \let\@svsec\@empty
 \else
   \refstepcounter{#1}%
\edef\@svsec{\ifnum#2<\@m
             \@ifundefined{#1name}{}{\csname #1name\endcsname\ }\fi
\noexpand\rom{\csname the#1\endcsname.}\enspace}\fi
 \@tempskipa #5\relax
 \ifdim \@tempskipa>\z@ 
   \begingroup #6\relax
   \@hangfrom{\hskip #3\relax\@svsec}{\interlinepenalty\@M #8\par}%
   \endgroup
   \if@article\else\csname #1mark\endcsname{%
        \ifnum \c@secnumdepth >#2\relax\csname the#1\endcsname. \fi#7}\fi
\ifnum#2>\@m \else
       \let\@tempf\\ \def\\{\protect\\}\addcontentsline{toc}{#1}%
{\ifnum #2>\c@secnumdepth \else
             \protect\numberline{%
               \ifnum#2<\@m
               \@ifundefined{#1name}{}{\csname #1name\endcsname\ }\fi
               \csname the#1\endcsname.}\fi
           #8}\let\\\@tempf
     \fi
 \else
  \def\@svsechd{#6\hskip #3\@svsec
    \@ifnotempty{#8}{\ignorespaces#8\unskip
       \ifnum\spacefactor<1001.\fi}%
        \ifnum#2>\@m \else
          \let\@tempf\\ \def\\{\protect\\}\addcontentsline{toc}{#1}%
            {\ifnum #2>\c@secnumdepth \else
              \protect\numberline{%
                \ifnum#2<\@m
                \@ifundefined{#1name}{}{\csname #1name\endcsname\ }\fi
                \csname the#1\endcsname.}\fi
             #8}\let\\\@tempf\fi}%
 \fi
\@xsect{#5}}
  \ifx\RequirePackage\undefined
  \let\@sect\@sect@my             
  \fi
  \def\th@remark@my{\theorempreskipamount6\p@\@plus6\p@
    \theorempostskipamount\theorempreskipamount
    \def\theorem@headerfont{\it}\normalshape}
  \ifx\RequirePackage\undefined
  \let\th@remark\th@remark@my
  \fi
\let\myLabel\@gobble
\def\labelsONmargin{\@mparswitchfalse\def\myLabel##1{\@bsphack\marginpar
                                  {\normalshape\tiny\rm Label ##1}\@esphack}}
\def\cyracc{\def\u##1{
                \if \i##1\char"1A%
                \else \if I##1\char"12%
                \else \accent"24 ##1\fi\fi }%
\def\"##1{\if e##1{\char"1B}%
                \else \if E##1{\char"13}%
                \else \accent"7F ##1\fi\fi }%
\def\9##1{\if##1z\char"19 
\else\if##1Z\char"11 
\else\if##1E\char"03 
\else\if##1e\char"0B 
\else\if##1u\char"18 
\else\if##1U\char"10 
\else\if##1A\char"17 
\else\if##1a\char"1F 
\else\if##1p\char"7E 
\else\if##1P\char"5E 
\else\if##1Q\char"5F 
\else\if##1q\char"7F 
\else\if##1i\char"1A 
\else\if##1I\char"12 
\else\if##1N\char"7D 
\fi
\fi
\fi
\fi
\fi
\fi
\fi
\fi
\fi
\fi
\fi
\fi
\fi
\fi
\fi
}%
\def\cydot{{\kern0pt}}}%
\def\cydot{$\cdot$}

\ifx\Russian\undefined
        \def\Russian{0\relax
    \message{Don't know the hyphenation rules for Russian^^J
                        Please do INITeX with `input  russhyph' in the 
                        command line}%
                \gdef\Russian{0\relax}%
        }
\fi
\def\@putname#1#2#3#4{\def\@@ref{#3}\let\old@bf\bf
        \def\bf##1{\old@bf\if?\noexpand##1?{#4}\else##1\fi}#1{#2}%
        \let\bf\old@bf}
\let\my@ref=\ref
\def\ref#1{\@putname\my@ref{#1}{#1}{\tiny\rm\@@ref}}
\let\my@pageref=\pageref
\def\pageref#1{\@putname\my@pageref{#1}{#1}{\tiny\rm\@@ref}}
\let\my@cite=\cite
\def\cite#1{\@putname\my@cite{#1}{\@citeb}{\tiny\rm\@@ref}}
\makeatother
\fi
\relax 
\theoremstyle{plain} 
\theorembodyfont{\sl}

\newtheorem{nwthrmi}{Fact } 
\theoremstyle{definition}
\newtheorem{nwthrmii}{Condition A } 
\theoremstyle{plain}
\newtheorem{nwthrmiii}{Fact }

\numberwithin{equation}{section}
\theorembodyfont{\rm}
\theoremstyle{definition}

\newtheorem{definition}{Definition}[section]
\newtheorem{conjecture}[definition]{Conjecture}
\newtheorem{example}[definition]{Example}

\theoremstyle{remark}

\newtheorem{remark}[definition]{Remark} 

\theoremstyle{plain} 
\theorembodyfont{\sl}

\newtheorem{theorem}[definition]{Theorem}
\newtheorem{lemma}[definition]{Lemma}
\newtheorem{corollary}[definition]{Corollary}
\newtheorem{proposition}[definition]{Proposition}
\newtheorem{amplification}[definition]{Amplification}

\begin{document}
\bibliographystyle{amsplain}
\relax 

\title[Quasi-algebraic geometry]{ Quasi-algebraic geometry of curves I \\
Riemann--Roch theorem and Jacobian }

\author{ Ilya Zakharevich }

\address{ Department of Mathematics, Ohio State University, 231 W.~18~Ave,
Columbus, OH, 43210 }

\dedicatory{Dedicated to Yu.~I.~Manin at his $ 60 $th birthday, to a Teacher
from a pupil}

\email {ilya@@math.ohio-state.edu}

\date{ August 1997\quad Printed: \today }

\setcounter{section}{-1}

\maketitle
\begin{abstract}
We discuss an analogue of Riemann--Roch theorem for curves with
an infinite number of handles. We represent such a curve $ X $ by its Shottki
model, which is an open subset $ U $ of $ {\mathbb C}P^{1} $ with infinite union of circles as
a boundary. An {\em appropriate\/} bundle on $ X $ is $ \omega^{1/2}\otimes{\mathcal L} $, $ {\mathcal L} $ being a bundle with
(say) constants as gluing conditions on the circles. An {\em admissible\/}
section of an appropriate bundle on $ X $ is a holomorphic half-form on $ U $
with given gluing conditions and $ H^{1/2} $-smoothness condition.

We study the restrictions on the mutual position of the circles and
the gluing constants which guarantee the finite dimension of the space of
appropriate sections of admissible bundles, and make the Riemann--Roch
theorem hold. The resulting Jacobian variety is described as an
infinite-dimension analogue of a torus.

\end{abstract}
\tableofcontents

\section{Introduction }

The need to carve out a set of curves of infinite genus for which
``most'' theorems of algebraic geometry are true comes out from the
following observations:
\begin{enumerate}
\item
The existence of algebro-geometric description of solutions of
infinite-dimensional integrable systems;
\item
The ability to describe the series of perturbation theory for string
amplitudes as integrals over moduli spaces;
\item
The hope that the union of compactifications of moduli spaces may
have a simpler geometry than the moduli spaces themselves.
\end{enumerate}

Different approaches which would result in different sets of curves
are possible (as in papers of Feldman, Kn\"orrer and Trubowitz
cf.~\cite{FelKnoTru96Inf}), thus we first motivate our choice of
tools
(Shottki model, $ H^{1/2} $-topology, capacities and half-forms) as (probably) a
best one to fulfill the expectations of the above origins of the theory.
Until Section~\ref{s0.10} we discuss how the above topics motivate the choice
of the Shottki model as a way to describe a curve of infinite genus.
After this (up to section~\ref{s0.40}) we describe motivations for the
choice of half-forms, $ H^{1/2} $-topology, and generalized-Sobolev-spaces (or
capacities) to describe admissible sections on a given curve. In the
remaining part of the introduction we do a walk-through the methods and
results one can find in this paper, as well as some historic remarks.

\subsection{Integrable systems and algebro-geometric solutions } A great
break-through in the first topic came with the work \cite{McKTru76Hil}. A
hyperelliptic curve of genus $ g $ with real branching points and a divisor of
degree $ g $ on this curve allows one to construct a so-called
{\em algebro-geometric\/} solution of KdV equation
\begin{equation}
u_{xx x}+12uu_{x}-4u_{t}=0.
\notag\end{equation}
Such solutions are called $ g $-{\em gap potentials}. McKean and Trubowitz studied
what substitutes algebraic geometry for solutions of KdV which are {\em not\/}
finite gap potentials. To such a solution they associated a curve (i.e.,
a complex manifold of dimension 1) which was hyperelliptic of {\em infinite
genus}, i.e., had infinitely many branching points, and some substitution
for the notion of a divisor of degree $ g $. It was shown that for the curves
which are related to KdV equation one can construct a well-behaved
analogue of $ 1 $-dimensional algebraic geometry. One should consider this
analogue as a generalization of algebraic geometry to hyperelliptic
curves of infinite genus (in fact only to some special curves of this
type).

Note that other integrable systems lead to different classes of
curves which appear in algebro-geometric solutions for the systems. Thus
one may expect that infinite-dimensional integrable systems may
lead to generalizations of algebraic geometry to different classes of
curves of infinite genus.

For this approach one of the richest systems is so called KP system
\begin{equation}
\left(u_{x xx}+12uu_{x}-4u_{t}\right)_{x} +3u_{yy}=0.
\notag\end{equation}
To describe an algebro-geometric solution of a KP system one starts with
an arbitrary algebraic curve, and an arbitrary linear bundle on this
curve. Thus to generalize the approach of \cite{McKTru76Hil} to the KP
equation, one needs to study curves of infinite genus of {\em generic\/} form (as
opposed to a hyperelliptic curve) and bundles on them.

Suppose for a minute that we have such a generalization, i.e., a
collection of curves and bundles on them. Call the members of these
families {\em admissible}. We assume that to any admissible curve and a bundle
on it we can associate some solution of KP equation. Let us investigate
what can we deduce about this collection from the known properties of KP
system (the properties below are applicable at least to algebro-geometric
solutions).

The dynamics of KP leaves the curve the same, but changes the
bundle. Thus the collection of admissible bundles should be reach enough
to include all the bundles obtained by time-flow of KP. In fact KP can be
generalized to include an infinite collection of commuting flows (with
different time variables), and they (taken together) can transform any
bundle to any other one (at least in finite-genus case). Thus we should
expect that we need our collection to include {\em all\/} the possible bundles,
thus the whole Jacobian.

The dynamics of KP is described in terms of locally affine structure
on the Jacobian, hence one should be able to describe the Jacobian as a
quotient of a vector space by a lattice. Algebraic geometry identifies
the Jacobian with a quotient of the space of global holomorphic forms by
the forms with integer periods. Thus one needs something similar to this
description.

There is an alternative description of a solution of KP equation in
terms of so called $ \tau $-{\em functions}. The relation of algebro-geometric
description with the description in terms of $ \tau $-functions needs a
construction of a Laurent series of a meromorphic global section of an
admissible bundle. This ceases to be trivial in infinite-genus case,
since already the description in \cite{McKTru76Hil} shows that one may need to
apply this operation at the {\em infinity\/} of the curve, i.e., at the points
added to the curve to compactify it. These points {\em are not\/} smooth points
of the curve, the curve is not even a topological manifold near these
points.

Finally, to get somewhat {\em explicit\/} description of solutions of KP as
functions of several variables, one needs a way to explicitly describe
global sections of linear bundles. Moreover, the bundles we need to
consider should have a finite number (preferably one!) of independent
global sections. A tool to construct such bundles in the case of finite
genus is the Riemann--Roch theorem, which gives an estimate on the
dimension of sections of the given bundles, and this estimate is precise
in the case of bundles in generic position. These sections are described
in terms of $ \theta $-functions, thus we will also need to describe $ \theta $-functions.

Collecting all this together, we see that we need to describe curves
$ X $ which can be equipped with a linear bundle $ {\mathcal L} $ such that the space of
global sections is finite-dimensional. One should be able to define what
does it mean that two bundles $ {\mathcal L} $ and $ {\mathcal L}' $ are equivalent, describe the
equivalence classes of bundles in terms of global holomorphic $ 1 $-forms,
give a local description of sections of the bundles, and give a global
description of sections of the bundles in terms of $ \theta $-functions.

In this paper we do not complete this program. However, we describe
all the ingredients but the last one.

\subsection{Universal Grassmannian }\label{s0.4}\myLabel{s0.4}\relax  The other two origins, the string theory
and geometry of moduli spaces, come into play if we consider the
algebro-geometric methods of solving integrable systems in the other
direction, as a way to find information of algebro-geometric type from
solutions of integrable systems (similar to solution of Shottki problem
in \cite{Shio86Char}).

The key idea is that in the investigated cases the set of
algebro-geometric solutions is {\em dense\/} in the set of all solutions, thus
the set of the solutions is a {\em completion\/} of the set of algebro-geometric
solutions. Since the the set of algebro-geometric solutions is a moduli
space of appropriate structures, and the set of solutions is a linear
space (due to the possibility to solve the Cauchy problem), we see that
the moduli space has a completion which has a topology of a linear space.

While the moduli spaces carry remarkable measures \cite{BeiMan86Mum} used
in the integrals of the string theory, this linear space has a symplectic
structure, so one may expect that this simplectic structure may have a
relationship to the measure on the subsets corresponding to
algebro-geometric solutions.

To carry out this program one needs to investigate how an arbitrary
solution may be approximated by algebro-geometric solutions. For KdV
equation this question was answered by \cite{McKTru76Hil} (for a fixed
spectrum), ameliorated to describe the inclusion of the subset of finite
gap potentials up to a diffeomorphism in \cite{ZakhFinGap2,ZakhFinGap3}, up
to a symplectomorphism in \cite{BatBloGuil95Sym}. The answer is that the above
inclusion is isomorphic to inclusion of trigonometric polynomials into
the space of all functions (for a very wide range of spaces of
functions). Since trigonometric polynomials fill finite-dimensional
coordinate subspaces in spaces of functions (using Fourier coefficients
as coordinates), we see that this picture is very similar to one of a
divisor with normal crossings.

Much simpler problem is how an algebro-geometric solution may be
approximated by (simpler) algebro-geometric solutions. A remarkable fact
is that such problems for different classes of integrable systems may be
solved via a uniform approach: compactification of moduli spaces via
{\em Universal Grassmannian}.

The reason is that the standard way to associate a solution of an
integrable system to an algebro-geometric data comes from consideration
of the mapping to {\em Universal Grassmannian}. The Universal Grassmannian
gives a convenient way to say that the space of sections of a linear
bundle on one curve is close to the space of sections on another linear
bundle on another curve.

Consider an algebraic curve $ X $ with a fixed point $ P\in X $ and a local
coordinate system $ x $ in a neighborhood of $ P $ which maps $ P $ to $ 0\in{\mathbb C} $. Let $ {\mathcal L} $ be
a linear bundle on $ X $, fix a trivialization of this bundle in a
neighborhood of $ P $. One can associate a Laurent series $ {\mathfrak l}\left(\varphi\right) $ to any section
$ \varphi $ of $ {\mathcal L} $ in a punctured neighborhood of $ P $: using a coordinate system on $ X $
and the trivialization of $ {\mathcal L} $, one can write this section as a function of
$ z\approx0 $, $ z\not=0 $. Let $ V $ be the space of meromorphic sections of $ {\mathcal L} $ which are
holomorphic outside of $ P $. (Abusing divisor notations, one can write $ V $ as
$ \Gamma\left(X,{\mathcal L}\left(\infty\cdot P\right)\right)=\bigcup_{k}\Gamma\left(X,{\mathcal L}\left(k\cdot P\right)\right) $.) Then $ {\mathfrak l}\left(V\right) $ is a subspace of the space $ {\mathfrak L} $ of
Laurent series. Let the universal Grassmannian $ \operatorname{Gr}\left({\mathfrak L}\right) $ be the Grassmannian
of subspaces of $ {\mathfrak L} $. Then $ {\mathfrak l}\left(V\right)\in\operatorname{Gr}\left(V\right) $ depends on $ X $, $ {\mathcal L} $, the coordinate system
near $ P $, and on the trivialization of $ {\mathcal L} $ near $ P $.

On $ {\mathfrak L} $ there are natural actions of the group of formal
diffeomorphisms\footnote{I.e., invertible $ \infty $-jets of mappings $ {\mathbb C} \to {\mathbb C} $ with $ 0\in{\mathbb C} $ being a fixed point.} of $ \left({\mathbb C},0\right) $ and the group of multiplication by invertible
Taylor series. Let $ G $ be the semidirect product of these groups. The
action of $ G $ on $ {\mathfrak L} $ corresponds (via $ {\mathfrak l} $) to changes of coordinate system on
$ X $, and a change of the trivialization of $ {\mathcal L} $. We see that to a triple $ \left(X,P,{\mathcal L}\right) $
we can naturally associate a point $ {\mathfrak l}\left(X,P,{\mathcal L}\right)\in\operatorname{Gr}\left({\mathfrak L}\right)/G $.

Note that while only the mapping to $ \operatorname{Gr}\left({\mathfrak L}\right)/G $ is invariantly defined,
one can get a canonical lifting to $ \operatorname{Gr}\left({\mathfrak L}\right) $ as far as $ X\not={\mathbb C}P^{1} $. The description
below does not behave well w.r.t. deformations of the curve/bundle, but
can be easily modified to do so. Given a bundle $ {\mathcal L} $, we can find $ k\in{\mathbb Z} $ such
that $ {\mathcal L}\left(k\cdot P\right) $ has only one independent section, and $ {\mathcal L}\left(\left(k+1\right)\cdot P\right) $ has another
one. The ratio $ z $ of these two sections identifies a neighborhood of $ P $
with a neighborhood of $ \infty\in{\mathbb C}P^{1} $, and this identification is defined up to an
affine transformation. Now the only section of $ {\mathcal L}\left(k\cdot P\right) $ gives a local
identification of $ {\mathcal L} $ with $ {\mathcal O}\left(-k\cdot P\right) $. Together with a coordinate system in a
neighborhood of $ P $ it gives a local section of $ {\mathcal L} $ defined up to
multiplication by a constant.

What remains is to pick up a coordinate system from those which
differ by an affine transformation. If $ X\not={\mathbb C}P^{1} $, then for some $ l>k+1 $ there
is a section of $ {\mathcal L}\left(l\cdot P\right) $ of the form $ z^{l}+p_{l-1}\left(z\right)+az^{-m}+bz^{-m-1}+O\left(z^{-m-2}\right) $ with
$ a\not=0 $, $ p_{l-1} $ being a polynomial of degree $ l-1 $, and $ m>0 $. Taking the minimal
possible $ l $, we see that the condition that $ a=1 $, $ b=0 $ picks up a coordinate
system out of the above class, unique up to multiplication by a root of
degree $ l+m $ of 1. (If $ X={\mathbb C}P^{1} $, then of course there is no canonically
defined coordinate system, since $ \left(X,P,{\mathcal L}\right) $ has automorphisms.) Now using
this coordinate system and the corresponding trivialization of $ {\mathcal L} $ one gets
a canonically defined image of $ \left(X,P,{\mathcal L}\right) $ in $ \operatorname{Gr}\left({\mathfrak L}\right) $.

Let the {\em Teichmuller\/}--{\em Jacoby\/} space $ {\mathcal N}_{g,1,d} $ be the moduli space (i.e.,
the ``set'' of equivalence classes) of triples $ \left(X,P,{\mathcal L}\right) $ with $ g\left(X\right)=g $, and
$ \deg \left({\mathcal L}\right)=d $. What is of primary importance to us is the fact that $ {\mathfrak l}|_{{\mathcal N}_{g,1,d}} $
is an injection (even in the case when $ X $ may have double points). Indeed,
\begin{enumerate}
\item
Meromorphic functions on $ X $ may be described as ratios of elements in
$ {\mathfrak l}\left(X,P,{\mathcal L}\right) $, thus normalization $ \operatorname{Norm}\left(X\right) $ of $ X $ may be reconstructed basing on
$ {\mathfrak l}\left(X,P,{\mathcal L}\right) $;
\item
The point $ P $ can be reconstructed since we know all the meromorphic
functions on $ X $, and know the order of pole at $ P $;
\item
To reconstruct the lifting $ \bar{{\mathcal L}} $ of $ {\mathcal L} $ to $ \operatorname{Norm}\left(X\right) $ note that if $ \varphi\in V $, then
the divisor (of zeros) of $ \varphi $ can be reconstructed as poles of
functions in $ {\mathfrak l}\left(\varphi\right)^{-1}{\mathfrak l}\left(V\right) $.
\item
Finally, to describe the gluings one needs to perform to get $ X $ from
$ \operatorname{Norm}\left(X\right) $ it is sufficient to consider a complement to $ \Gamma\left({\mathcal L}\left(\infty\cdot P\right)\right) $ in
$ \Gamma\left(\bar{{\mathcal L}}\left(\infty\cdot P\right)\right) $.
\end{enumerate}

As a corollary, the mapping $ {\mathfrak l} $ defines an inclusion of $ {\mathcal N}_{g,1,d} $ into
$ \operatorname{Gr}\left({\mathfrak L}\right)/G $. Since the image has a natural structure of a topological space,
this inclusion defines some natural {\em compactification\/} $ \bar{{\mathcal N}}_{g,1,d} $ of the
Teichm\"uller--Jacoby space. Indeed, consider the closure of $ {\mathfrak l}\left({\mathcal N}_{g,1,d}\right) $ in
$ \operatorname{Gr}\left({\mathfrak L}\right)/G $. As we will see it shortly, this closure is an image of a compact
smooth manifold\footnote{In fact we will see that already the {\em normalization\/} of the image is
smooth.}. The description of this compactification is very similar
to the description of Deligne--Mumford compactification, it is carried
out by adding to $ {\mathcal N}_{g,1,d} $ so called {\em semistable\/} objects. We will see that
these object are non-smooth curves and sheaves on them.

The mapping to Universal Grassmannian is a way to collect all the
moduli spaces for different genera into ``one big heap''. It is relevant to
curves of infinite genus since on a back yard of this heap one expects to
find moduli spaces of curves of infinite genus.

\subsection{Semistable elliptic curves }\label{s0.5}\myLabel{s0.5}\relax  In this section we consider the case
$ g=1 $, where one can calculate the mapping to the Universal Grassmannian
explicitly. We will see that in this case the closure of the image is
smooth\footnote{As we will see it later, $ g=1 $ is last case when one does not need
any normalization of the image.}.

Fix an elliptic curve $ E $. Then the set of classes of equivalence of
linear bundles of degree 0 (i.e., the Jacobian variety) is isomorphic to
the curve itself (after we fix a point $ P $ on the curve). Given any point
$ Q\in E $, the divisor $ Q-P $ determines a linear bundle $ {\mathcal O}\left(Q-P\right) $ on $ E $, and any
bundle of degree 0 is isomorphic to exactly one bundle of this form.

In the case of $ g=1 $ it is not instructive to consider the
Teichm\"uller--Jacoby manifold {\em literally}, since it is tainted by the fact
that an involution $ \sigma $ of an elliptic curve $ E $ sends any bundle $ {\mathcal L} $ of degree
0 to $ {\mathcal L}^{-1} $, thus the moduli space of bundles of degree 0 on an elliptic
curve fixed {\em up to an isomorphism\/} is $ E/\sigma $, i.e., a rational curve. To fix
the problem one can add some harmless discrete parameter which would
prohibit $ \sigma $ to be an automorphism, say, consider collections $ \left(X,P,{\mathcal L},\alpha\right) $, $ \alpha $
being a homology class modulo $ {\mathbb Z}/3Z $. However, knowing that we can fix this
problem, we are going to ignore it whatsoever, since it can also be
avoided by considering curves close to the given one, what we are going
to do anyway.

For every elliptic curve $ E $ there is a unique number $ j\in{\mathbb C} $ such that $ E $
is isomorphic to the curve
\begin{equation}
y^{2}=\left(500j-1\right)x^{3}-15jx-j.
\notag\end{equation}
The number $ j $ is called the $ j $-{\em invariant\/} of the curve. The moduli
(Teichm\"uller) space of curves of genus 1 may be identified with $ {\mathbb C} $ via
$ j $-invariant. The Teichm\"uller--Jacoby space is fibered over the
Teichm\"uller space with the fiber being the Jacobian, i.e., the elliptic
curve itself (i.e., in the standard notations this space coincides with
$ {\mathcal M}_{1,2} $). To compactify the Teichm\"uller space one adds a point with $ j=\infty $,
obtaining $ {\mathbb C}P^{1} $ as the compactified Teichm\"uller space.

To visualize the compactification of the Teichm\"uller space it is
more convenient to consider the family of elliptic curves $ X_{\varepsilon} $ given by
$ y^{2}=x^{2}-x^{3}-\varepsilon $. When $ \varepsilon \to 0 $ (so $ j=\frac{1}{15^{3}\varepsilon\left(4-27\varepsilon\right)} \to \infty $) one gets a rational
curve with a double point. Below we discuss how the corresponding
compactification of Teichm\"uller--Jacoby space looks like. We will see
that instead of adding one point, we need to add a rational curve with a
selfintersection. The compactified Teichm\"uller--Jacoby space maps to the
compactified Teichm\"uller space, with the fiber over $ \infty $ being the above
singular curve.

Consider an arbitrary elliptic curve $ E $ with large a $ j $-invariant, and
a bundle of degree 0 on $ E $. One can identify $ E $ with $ X_{\varepsilon} $ for an appropriate
small $ \varepsilon $. Fix a point $ P_{\varepsilon} $ on $ X_{\varepsilon} $, then for any bundle $ {\mathcal L} $ of degree 0 on $ X_{\varepsilon} $
one can identify $ {\mathcal L} $ with $ {\mathcal O}\left(Q_{\varepsilon}-P_{\varepsilon}\right) $, $ Q_{\varepsilon} $ being an appropriate point on $ X_{\varepsilon} $.
Fix $ P_{\varepsilon} $ to be the point of $ X $ which has $ y $-coordinate being 0, and is close
to $ \left(1,0\right)=P_{0} $. Suppose that $ {\mathcal L}\not={\mathcal O} $, thus $ Q_{\varepsilon}\not=P_{\varepsilon} $.

Let us describe meromorphic sections of $ {\mathcal L}={\mathcal O}\left(Q_{\varepsilon}-P_{\varepsilon}\right) $ with the only
pole being at $ P_{\varepsilon} $. It is the same as to describe sections of $ {\mathcal O}\left(Q_{\varepsilon}+kP_{\varepsilon}\right) $,
$ k\gg0 $, i.e., meromorphic functions on $ X_{\varepsilon} $ with a (possible) simple pole at
$ Q_{\varepsilon} $, and any pole at $ P_{\varepsilon} $. These sections are uniquely determined by their
singular part at $ P_{\varepsilon} $ (since there is no meromorphic function on an
elliptic curve $ X_{\varepsilon} $ which has only a simple pole at $ Q_{\varepsilon} $). Thus to describe
the space of these sections, it is sufficient to describe functions $ f_{k} $ on $ X_{\varepsilon} $
which have a pole of prescribed order $ k $ at $ P_{\varepsilon} $ and possibly an additional
pole at $ Q_{\varepsilon} $. (Note that $ f_{0}\equiv 1 $ should be considered as a section of $ {\mathcal O}\left(Q_{\varepsilon}-P_{\varepsilon}\right) $
with a pole at $ P_{\varepsilon} $, since holomorphic sections of $ {\mathcal O}\left(Q_{\varepsilon}-P_{\varepsilon}\right) $ are identified
with functions having a zero at $ P_{\varepsilon} $.)

Note that $ y $ is a meromorphic function on $ X_{\varepsilon} $ with a triple pole at
infinity of $ X_{\varepsilon} $, $ x $ has a double pole at infinity of $ X_{\varepsilon} $. Let $ P_{\varepsilon} $ have
coordinates $ \left(a_{\varepsilon},0\right) $, the line through $ P_{\varepsilon} $ and $ Q_{\varepsilon} $ have the equation
$ c_{\varepsilon}x+y=d_{\varepsilon} $, the third point $ R_{\varepsilon} $ of intersection of this line and $ X_{\varepsilon} $ has
$ x $-coordinate $ b_{\varepsilon} $. If $ k=2k_{1} $ is even, one can take $ f_{k}=\frac{1}{\left(x-a_{\varepsilon}\right)^{k_{1}}} $. If
$ k=2k_{1}-1 $, and $ k_{1}>1 $, then $ yf_{k+1} $ has no pole at infinity, thus can be
considered as $ f_{k} $. What remains is to describe $ f_{1} $. Note that $ \frac{1}{c_{\varepsilon}x+y-d_{\varepsilon}} $
has a pole at $ P_{\varepsilon} $, $ Q_{\varepsilon} $, and at $ R_{\varepsilon} $. Moreover, it has a zero of third order
at infinity. Multiplying this function by $ x-b_{\varepsilon} $, we kill the pole at $ R_{\varepsilon} $,
thus get $ f_{1}=\frac{x-b_{\varepsilon}}{c_{\varepsilon}x+y-d_{\varepsilon}} $.

Let us investigate how the space spanned by $ f_{i} $, $ i\geq0 $, depends on $ X_{\varepsilon} $
and $ Q_{\varepsilon} $. First of all, $ \varepsilon $ depends smoothly on $ Q_{\varepsilon} $ (if coordinates of $ Q_{\varepsilon} $ are
$ q_{1},q_{2} $, then $ \varepsilon=q_{2}^{2}-q_{1}^{2}-q_{1}^{3} $). From this moment on, we may denote $ Q_{\varepsilon} $ by just
$ Q $, since $ \varepsilon $ is a function of $ Q $. Second, $ x $ depends smoothly on $ y $ and $ \varepsilon $ in a
neighborhood of $ P_{\varepsilon} $ on $ X_{\varepsilon} $, {\em including\/} the case $ \varepsilon=0 $. Thus $ P_{\varepsilon} $ depends
smoothly on $ Q $, {\em including\/} cases $ \varepsilon=0 $ and $ Q=\left(0,0\right)\in{\mathbb C}^{2} $. Hence $ c_{\varepsilon} $, $ d_{\varepsilon} $ depend
smoothly on $ Q $. Substituting $ y=-c_{\varepsilon}x+d_{\varepsilon} $ into $ y^{2}=x^{2}-x^{3}-\varepsilon $, we get a cubic
equation in $ x $, which has $ q_{1} $ as a root, and two well-separated roots at $ a_{\varepsilon} $
and $ b_{\varepsilon} $. In particular, $ b_{\varepsilon} $ depends smoothly on $ Q $. Note that one can
consider $ y $ as a local coordinate near $ P_{\varepsilon} $, thus the functions $ f_{k} $ can be
written as function of $ y $ if $ y $ is small. We see that all the functions
$ f_{k}\left(y\right) $ depend smoothly on $ \left(q_{1},q_{2}\right) $ (at least as far as $ q_{1} $, $ q_{2} $ are small, $ y $
is small).

If $ q_{2}^{2}-q_{1}^{2}-q_{1}^{3}\not=0 $, the described above space $ V_{q_{1},q_{2}}=\left<f_{i}\left(y\right)\right>_{i\geq0}\subset{\mathfrak L} $ is
the image of $ \left(X_{\varepsilon},P_{\varepsilon},{\mathcal L}\right) $ in the Universal Grassmannian, i.e., it is in
$ {\mathfrak l}\left({\mathcal N}_{1,1,0}\right) $. We see that the closure of $ {\mathfrak l}\left({\mathcal N}_{1,1,0}\right) $ contains at least the
subspaces $ V_{q_{1},q_{2}} $ with $ q_{2}^{2}-q_{1}^{2}-q_{1}^{3}=0 $. What remains is to give a
description of these subspaces as spaces of sections of sheaves on
non-smooth curves. We will see that if $ \left(q_{1},q_{2}\right)\not=\left(0,0\right) $, it is sufficient to
consider {\em linear bundles\/} on non-smooth curves, but if $ \left(q_{1},q_{2}\right)=\left(0,0\right) $, a
consideration of a sheaf is unavoidable.

It is obvious {\em which\/} non-smooth curve we need to consider: the limit
$ X_{0} $ of $ X_{\varepsilon} $ when $ \varepsilon \to $ 0, with the equation $ y^{2}=x^{2}-x^{3} $. There is a
parameterization mapping $ {\mathbb C}P^{1} \to X_{0} $, $ \lambda \mapsto \left(x,y\right) $ if $ y=\lambda x $. It is a bijection
outside of $ \lambda=\pm1 $, these two points are both mapped to (0,0). Thus $ X_{0} $ is a
rational curve with two points glued together. To describe a linear
bundle $ {\mathcal L} $ on $ X_{0} $ one needs to describe its lifting to $ {\mathbb C}P^{1} $, and the
identification of two fibers over $ \lambda=\pm1 $. If the lifting has degree 0, it
is a trivial bundle, thus one can describe the gluing by a number $ \theta\in{\mathbb C}^{*} $.
Denote this bundle $ {\mathcal L}_{\theta} $, it has a global holomorphic section iff $ \theta=1 $.

Note that if $ P_{0}=\left(1,0\right)\in X_{0} $, $ Q\in X_{0} $, $ Q\not=\left(0,0\right) $, then the mapping $ f\left(\lambda\right) \to
\frac{\lambda-\lambda_{0}}{\lambda}f\left(\lambda\right) $ gives an isomorphism of $ {\mathcal O}\left(Q-P_{0}\right) $ to $ {\mathcal L}_{\theta} $ if $ \theta=\frac{1-\lambda_{0}}{1+\lambda_{0}} $, $ \lambda_{0} $
being the coordinate of $ Q $. The defined above functions $ f_{k}\left(y\right) $ form a basis
of meromorphic sections of $ {\mathcal O}\left(Q-P_{0}\right) $, thus we see that for $ Q\in X_{0} $ the defined
above subspace $ V_{q_{1},q_{2}}\subset{\mathfrak L} $ may be interpreted as $ {\mathfrak l}\left(X_{0},P_{0},{\mathcal L}_{\theta}\right) $.

Now investigate what happens if $ Q \to \left(0,0\right) $ along $ X_{0} $, i.e., $ \theta \to 0 $ or
$ \theta \to \infty $. Instead of $ y $, consider the coordinate $ \lambda $ in a neighborhood of $ P_{0} $
on $ X_{0} $. In this coordinate the space $ V\left(X_{0},P_{0},{\mathcal L}_{\theta}\right) $ is the set of polynomials
$ p\left(\mu\right) $ in $ \mu=\lambda^{-1} $ which satisfy $ p\left(1\right)=\theta p\left(-1\right) $. Denote this subspace of $ {\mathbb C}\left[\mu\right] $ by
$ V_{\theta} $. Obviously, the subspace $ V_{\theta} $ depends smoothly on $ \theta\in{\mathbb C}P^{1} $, $ V_{0} $ consists of
polynomials with a zero at 1, $ V_{\infty} $ consists of polynomials with a zero at $ -1 $.
Moreover, images of $ V_{0} $ and $ V_{\infty} $ in $ \operatorname{Gr}\left({\mathfrak L}\right)/G $ coincide, since $ V_{\infty}=\frac{\mu+1}{\mu-1}V_{0} $,
and $ \frac{\mu+1}{\mu-1} $ is smooth near $ \mu=\infty $, therefore $ f\left(\mu\right) \to \frac{\mu+1}{\mu-1}f\left(\mu\right) $
corresponds to an element of $ G $.

Returning back to the representation of $ V_{q_{1},q_{2}} $, $ q_{1}=q_{2}=0 $, in terms of
$ f_{k} $, we see that $ f_{1}=x/y=\lambda^{-1} $, $ f_{2}=\left(x-1\right)^{-1}=-\lambda^{-2} $, thus $ f_{2k}=\pm\lambda^{-2k} $; $ y=\lambda\left(1-\lambda^{2}\right) $,
thus $ f_{2k+1}=\frac{1-\lambda^{2}}{\lambda^{2k-1}}=\lambda^{-2k-1}-\lambda^{-2k+1} $. We see that this space is $ {\mathbb C}\left[\lambda^{-1}\right] $.
Hence the limit of spaces of sections of $ {\mathcal O}\left(Q-P\right) $ when $ Q \to \left(0,0\right) $ is
isomorphic to the space of meromorphic functions on the normalization $ {\mathbb C}P^{1} $
of $ X_{0} $ with the only pole at $ \lambda=0 $, i.e., with $ {\mathfrak l}\left({\mathbb C}P^{1},0,{\mathcal O}\right) $. If we want to
describe this space in the same terms as we described the mapping $ {\mathfrak l} $, we
should consider sections of $ {\mathcal O} $ on $ {\mathbb C}P^{1} $ as sections of the direct image $ \pi_{*}{\mathcal O} $
on $ X_{0} $, $ \pi $ being the projection of $ {\mathbb C}P^{1} $ to $ X_{0} $. Since $ \pi_{*}{\mathcal O} $ is no longer a
sheaf of sections of a linear bundle, we conclude that we {\em need\/} to
consider sheaves instead of bundles.

The only thing which remains to prove is that the completion of the
image $ {\mathfrak l}\left({\mathcal N}_{1,1,0}\right) $ is smooth. We had already shown the part of the
completion we consider here is an image of a neighborhood of 0 in $ {\mathbb C}^{2} $.
What remains to prove is that the mapping $ {\mathbb C}^{2} \xrightarrow[]{V_{\bullet}} \operatorname{Gr}\left({\mathfrak L}\right)/G $ has a
non-degenerate derivative at (0,0). It is sufficient to pick up a section
of the projection $ \operatorname{Gr}\left({\mathfrak L}\right) \to \operatorname{Gr}\left({\mathfrak L}\right)/G $, and show that a lifting of $ V_{\bullet} $ to this
section has a non-degenerate derivative at (0,0). Again, it is sufficient to
restrict this mapping to $ X_{0} $ and show that the derivatives along two
branches of $ X $ at (0,0) are independent.

One branch of $ X_{0} $ gives a family of subspaces $ V_{\theta} $, $ \theta\approx0 $, $ V_{\theta}\subset{\mathbb C}\left[\lambda^{-1}\right] $. It
is a simple but tedious calculation to show that the lifting defined in
Section~\ref{s0.4} results in a smooth lifting of $ V_{\bullet} $, and that the smooth
liftings resulting from two branches are non-degenerate, and have
different derivatives indeed.

We see that to describe the elements in the completion of $ {\mathfrak l}\left({\mathcal N}_{1,1,0}\right) $
in a way similar to the description of $ {\mathfrak l} $ one needs to consider both
degenerated curves and degenerated bundles (i.e., sheaves) on these
curves. Using the fact that the space of sections of $ \pi_{*}{\mathcal O} $ is close to the
set of sections of $ {\mathcal O}\left(Q\right) $ if $ Q $ is close to the double point (0,0) of $ X_{0} $, it
is reasonable to abuse notations and write $ \pi_{*}{\mathcal O} $ as $ {\mathcal O}\left(O\right) $, here $ O=\left(0,0\right) $.

\subsection{Compactified moduli spaces } We start with a description of elements of
the would-be compactified moduli spaces. As we have seen it in Section
~\ref{s0.5}, we need to allow the curves to have double points, and allow a
generalization of a notion of bundle. We call such bundles {\em non-smooth\/}
bundles, in fact they are sheaves which are (locally) direct images of
bundles on the normalization\footnote{Recall that in the case $ \dim =1 $ {\em normalization\/} coincides with ungluing
double points and straightening out cusps.} of the curve.

\begin{definition} Let $ Y $ be a connected curve with only singularities $ \operatorname{Sing}\left(Y\right) $
being double points, let $ \widetilde{Y} $ be the normalization of $ Y $. Define
$ g\left(Y\right)=g\left(\widetilde{Y}\right)+\operatorname{card}\left(\operatorname{Sing}\left(Y\right)\right) $, for a subset $ S_{0}\subset\operatorname{Sing}\left(Y\right) $ let $ Y_{S_{0}} $ be the curve
obtained by ungluing points in $ S_{0} $. A {\em smooth linear bundle\/} $ {\mathcal L} $ on $ Y $ is a
linear bundle $ \widetilde{{\mathcal L}} $ on $ \widetilde{Y} $ with a fixed identification of fibers at points of $ \widetilde{Y} $
over the same point on $ Y $. Define {\em degree\/} of a smooth linear bundle as the
$ \deg \left(\widetilde{{\mathcal L}}\right) $. A {\em section\/} of $ {\mathcal L} $ is a section of $ \widetilde{{\mathcal L}} $ compatible with identifications
of fibers.

A (not necessarily smooth) {\em linear bundle\/} $ {\mathcal L} $ on $ Y $ is a collection
$ \left(S_{0},\bar{{\mathcal L}}\right) $ of $ S_{0}\subset\operatorname{Sing}\left(Y\right) $ and a smooth linear bundle $ \bar{{\mathcal L}} $ on $ Y_{S_{0}} $. A {\em section\/} of $ {\mathcal L} $
is a section of $ \bar{{\mathcal L}} $. Let $ \deg \left({\mathcal L}\right)=\deg \left(\bar{{\mathcal L}}\right)+\operatorname{card}\left(S_{0}\right) $, call $ S_{0} $ the set of {\em double
points\/} of $ {\mathcal L} $. Let $ P\subset\operatorname{Smooth}\left(Y\right) $ be a finite collection of points. Call the
collection $ \left(Y,P,{\mathcal L}\right) $ {\em semistable\/} if $ \left(Y,P,{\mathcal L}\right) $ has no infinitesimal
automorphisms, i.e., if any rational connected component of $ \widetilde{Y} $ has at
least three points in $ P $ or over double points of $ Y $, and if $ Y $ is a
(smooth) elliptic curve, then either $ P\not=\varnothing $, or $ \deg {\mathcal L}\not=0 $. \end{definition}

In other words, for a ``non-smooth linear bundle'' $ {\mathcal L} $ the double points
of the curve are broken into two subsets: for one (smooth points of $ {\mathcal L} $) we
fix identifications of fibers of the lifting to normalizations, for
another one (as we have seen in Section 0.5, one should consider them as
{\em poles\/} of $ {\mathcal L} $) local sections of $ {\mathcal L} $ may have ``different'' values at two
branches of $ Y $ which meat at a given point.

Note that the smooth structures on the Teichm\"uller space and on the
Jacobian of a smooth curve provide an atlas on the set of equivalence
classes of smooth curves/bundles. The explicit deformation which we are
going to describe now will provide an atlas in a neighborhood of a class
of a semistable curve.

Consider a curve $ Y $ as in the definition, let $ Q\in Y $ be a double point,
and $ y_{1} $, $ y_{2} $ be coordinates on the branches of $ Y $ near $ Q $, so $ Y $ is locally
isomorphic to $ y_{1}y_{2}=0 $. Let $ {\mathcal L} $ be a linear bundle over $ Y $ with a fixed
trivialization near $ Q $. Suppose that $ Q $ is a double point of $ {\mathcal L} $, thus local
sections of $ {\mathcal L} $ can be written as (unrelated) functions $ f_{1}\left(y_{1}\right) $, $ f_{2}\left(y_{2}\right) $.
Note that trivialization of $ {\mathcal L} $ induces identification of fibers of $ \widetilde{{\mathcal L}} $ over
$ Q $. Let $ \overset{\,\,{}_\circ}{{\mathcal L}} $ be the bundle on $ Y $ (smooth near $ Q $) obtained by addition of
this identification to $ {\mathcal L} $. Extend the bundle $ \overset{\,\,{}_\circ}{{\mathcal L}} $ from $ Y $ to $ Y\cup\left\{\left(y_{1},y_{2}\right) |
|y_{1}|,|y_{2}|\ll 1\right\} $ trivially using the trivialization of $ {\mathcal L} $ near $ Q $, denote this
extension as $ \bar{{\mathcal L}} $.

Define a deformation $ \left(Y_{\varepsilon_{1},\varepsilon_{2}},{\mathcal L}_{\varepsilon_{1},\varepsilon_{2}}\right) $ of $ \left(Y,{\mathcal L}\right) $ (here $ \varepsilon_{1},\varepsilon_{2}\in{\mathbb C} $ are
fixed small numbers) by the following recipe:
\begin{enumerate}
\item
if $ \varepsilon_{1}\varepsilon_{2}=0 $, then $ Y_{\varepsilon_{1},\varepsilon_{2}}=Y $;
\item
if $ \varepsilon_{1}=\varepsilon_{2}=0 $, then $ {\mathcal L}_{\varepsilon_{1},\varepsilon_{2}}={\mathcal L} $;
\item
if $ \varepsilon_{1}\not=0 $, but $ \varepsilon_{2}=0 $, then $ {\mathcal L}_{\varepsilon_{1},\varepsilon_{2}} $ is the bundle $ \overset{\,\,{}_\circ}{{\mathcal L}}\left(P_{\varepsilon_{1}}\right) $, here $ P_{\varepsilon_{1}} $ is the
point on $ Y $ with $ y_{1}=\varepsilon_{1} $; similarly for $ \varepsilon_{1}=0 $, but $ \varepsilon_{2}\not=0 $;
\item
If $ \varepsilon_{1}\not=0 $, $ \varepsilon_{2}\not=0 $, then $ Y_{\varepsilon_{1},\varepsilon_{2}} $ is obtained by gluing the hyperbola
$ y_{1}y_{2}=\varepsilon_{1}\varepsilon_{2} $ to $ Y $ (with a small neighborhood of $ Q $ removed) via coordinate
projections, and $ {\mathcal L}_{\varepsilon_{1},\varepsilon_{2}}=\overset{\,\,{}_\circ}{{\mathcal L}}\left(P_{\varepsilon_{1},\varepsilon_{2}}\right) $, here $ \overset{\,\,{}_\circ}{{\mathcal L}}_{\varepsilon_{1},\varepsilon_{2}} $ is the restriction of $ \bar{{\mathcal L}} $ to
$ Y_{\varepsilon_{1},\varepsilon_{2}} $, the point $ P_{\varepsilon_{1},\varepsilon_{2}} $ is the point on hyperbola with coordinates
$ \left(\varepsilon_{1},\varepsilon_{2}\right) $.
\end{enumerate}
If $ Y $ had some marked points $ \left\{P_{i}\right\}\subset\operatorname{Smooth}\left(Y\right) $, then $ Y_{\varepsilon_{1},\varepsilon_{2}} $ has the same
marked points (correctly defined since $ Y_{\varepsilon_{1},\varepsilon_{2}} $ is identified with $ Y $
outside a small neighborhood of $ Q $).

Note that if $ P\in\operatorname{Smooth}\left(Y\right) $, then $ {\mathfrak l}\left(Y,P,{\mathcal L}\right) $ is well-defined in $ \operatorname{Gr}\left({\mathfrak L}\right)/G $.

\begin{lemma} $ g\left(Y_{\varepsilon_{1},\varepsilon_{2}}\right)=g\left(Y\right) $, $ \deg {\mathcal L}_{\varepsilon_{1},\varepsilon_{2}}=\deg {\mathcal L} $. If $ \left(Y,P,{\mathcal L}\right) $ is semistable, then
$ \left(Y_{\varepsilon_{1},\varepsilon_{2}},P,{\mathcal L}_{\varepsilon_{1},\varepsilon_{2}}\right) $ is semistable too. If $ P\in\operatorname{Smooth}\left(Y\right) $, then
$ {\mathfrak l}\left(Y_{\varepsilon_{1},\varepsilon_{2}},P,{\mathcal L}_{\varepsilon_{1},\varepsilon_{2}}\right)\in\operatorname{Gr}\left({\mathfrak L}\right)/G $ depends smoothly on $ \varepsilon_{1},\varepsilon_{2} $. \end{lemma}

Similarly, if $ Q $ is a double point of $ Y $, but not a double point of $ {\mathcal L} $,
define a one-parametric deformation $ \left(Y_{\varepsilon},{\mathcal L}_{\varepsilon}\right) $ by gluing to $ Y $ a hyperbola
$ y_{1}y_{2}=\varepsilon $ without changing $ {\mathcal L} $ (far from $ Q $) and the trivialization of $ {\mathcal L} $. A
statement similar to the above lemma continues to be true.

Let $ \bar{{\mathcal N}}_{g,1,d} $ be the set of equivalence classes of semistable
collections $ \left(Y,P,{\mathcal L}\right) $, let $ \bar{{\mathcal N}}_{g,1,d}^{\left(m,k\right)}\subset\bar{{\mathcal N}}_{g,1,d} $ be the subset of $ \bar{{\mathcal N}}_{g,1,d} $
consisting of curves with exactly $ m $ double points and bundles with
exactly $ k $ double points $ \left(k\leq m\right) $. Note that $ \bar{{\mathcal N}}_{g,1,d}^{\left(0,0\right)}={\mathcal N}_{g,1,d} $, and the
complement to $ {\mathcal N}_{g,1,d} $ in $ \bar{{\mathcal N}}_{g,1,d} $ is a disjoint union of $ \bar{{\mathcal N}}_{g,1,d}^{\left(m,k\right)} $. Call
these subsets {\em strata\/} of $ \bar{{\mathcal N}}_{g,1,d} $. Moreover, note that each $ \bar{{\mathcal N}}_{g,1,d}^{\left(m,k\right)} $
carries a natural smooth structure (inherited from the smooth structures
on Teichm\"uller spaces and Jacobians). (In fact the above smooth
structures can be refined to structures of orbifolds, but we ignore this
refinement here.)

Let $ \left(Y,y_{0},{\mathcal L}\right)\in\bar{{\mathcal N}}_{g,1,d}^{\left(m,k\right)} $, $ \varepsilon\in{\mathbb C}^{m+k} $ is small. Taking coordinate systems
in neighborhoods of double points, and trivializations of $ {\mathcal L} $ in these
neighborhoods, we obtain a deformation $ \left(Y_{\varepsilon},P,{\mathcal L}_{\varepsilon}\right)\in\bar{{\mathcal N}}_{g,1,d} $. Since
$ \bar{{\mathcal N}}_{g,1,d}^{\left(m,k\right)} $ has a natural smooth structure, we see that {\em a piece\/} of $ \bar{{\mathcal N}}_{g,1,d} $
is fibered over $ \bar{{\mathcal N}}_{g,1,d}^{\left(m,k\right)} $ with fibers being small balls in $ {\mathbb C}^{m+k} $.

Define a structure of an manifold on $ \bar{{\mathcal N}}_{g,1,d} $ using the above
fibration as an atlas in a neighborhood of a point in $ \bar{{\mathcal N}}_{g,1,d}^{\left(m,k\right)} $.

\begin{proposition} $ \bar{{\mathcal N}}_{g,1,d} $ is a smooth compact manifold (orbifold). $ {\mathfrak l}\left(\bar{{\mathcal N}}_{g,1,d}\right) $
coincides with the closure of $ {\mathfrak l}\left({\mathcal N}_{g,1,d}\right) $. \end{proposition}

\begin{proposition} \label{prop0.33}\myLabel{prop0.33}\relax  $ {\mathfrak l}\left(X,P,{\mathcal L}\right)={\mathfrak l}\left(X',P',{\mathcal L}'\right) $ iff $ \left(X,P,{\mathcal L}\right) $ becomes isomorphic
to $ \left(X',P',{\mathcal L}'\right) $ after ungluing of double points of $ {\mathcal L} $ and $ {\mathcal L}' $. \end{proposition}

This means that $ \bar{{\mathcal N}}_{g,1,d} $ is a natural smooth compactification of
$ {\mathcal N}_{g,1,d} $. Moreover, the explicit coordinates in a neighborhood of a stratum
show that the {\em boundary\/} $ \bar{{\mathcal N}}_{g,1,d}\smallsetminus{\mathcal N}_{g,1,d} $ is a divisor with normal
intersections. As it follows from Proposition~\ref{prop0.33}, the image
$ {\mathfrak l}\left(\bar{{\mathcal N}}_{g,1,d}\right) $ can be obtained from $ \bar{{\mathcal N}}_{g,1,d} $ by contracting smooth
submanifolds into points, thus the normalization of the image coincides
with $ \bar{{\mathcal N}}_{g,1,d} $.

\subsection{Adjacency of moduli spaces and Shottki model }\label{s0.8}\myLabel{s0.8}\relax  We have seen that
the consideration of the mapping to the universal Grassmannian leads to a
remarkable compactification of Teichm\"uller--Jacoby space. Moreover, one
can use this compactification to define inclusions of compactified
Teichm\"uller--Jacoby spaces for different $ g $, $ d $ one into another. These
inclusions are going to be compatible with the mapping $ {\mathfrak l} $ to Universal
Grassmannian, and may be constructed studying the properties of the
mapping $ {\mathfrak l} $.

Consider a curve $ X $ with a bundle $ {\mathcal L} $ (we allow $ X $ and $ {\mathcal L} $ to have double
points). Given two smooth points $ Q_{1} $, $ Q_{2} $ on $ X $ one can construct a curve
$ X_{Q_{1}Q_{2}} $ obtained by gluing $ Q_{1} $ with $ Q_{2} $. This curve has one more double
point, and genus $ g\left(X\right)+1 $. One can consider the bundle $ {\mathcal L}_{Q_{1}Q_{2}} $ on $ X_{Q_{1}Q_{2}} $ which
has an extra double point at $ Q_{1}\sim Q_{2} $, otherwise coincides with $ {\mathcal L} $. Note that
sections of $ {\mathcal L} $ can be considered as sections of $ {\mathcal L}_{Q_{1}Q_{2}} $ and visa versa,
hence $ {\mathfrak l}\left(X,P,{\mathcal L}\right) = {\mathfrak l}\left(X_{Q_{1}Q_{2}},P,{\mathcal L}_{Q_{1}Q_{2}}\right) $ for any $ P $.

Using this remark one can easily include the compactification of one
moduli space into the boundary of another one. Indeed, one can generalize
the construction of $ \bar{{\mathcal N}}_{g,1,d} $ to the case of $ n $ marked points $ \left\{P_{i}\right\} $ (instead
of one $ P $), the only change being that one needs to consider $ \operatorname{Gr}\left({\mathfrak L}^{n}\right) $
instead of $ \operatorname{Gr}\left({\mathfrak L}\right) $. Now associate an element of $ \bar{{\mathcal N}}_{g+l,n-2l,d+l} $ to an element
$ \left(Y,\left\{P_{i}\right\},{\mathcal L}\right)\in\bar{{\mathcal N}}_{g,n,d} $ as follows: take first $ 2l $ points out of $ \left\{P_{i}\right\} $, and glue
them pairwise. Since we do not glue fibers of $ {\mathcal L} $ at these points, the
resulting double points on the resulting curve $ Y_{1} $ are double points of a
bundle $ {\mathcal L}_{1} $. Clearly, $ \left(Y_{1},\left\{P_{2l+i}\right\},{\mathcal L}_{1}\right)\in\bar{{\mathcal N}}_{g+l,n-2l,d+l} $. Thus we obtain a
mapping $ \bar{{\mathcal N}}_{g,n,d} \to \bar{{\mathcal N}}_{g+l,n-2l,d+l} $. Note that this correspondence does not
change the space of global sections of $ {\mathcal L} $, thus the image under the map $ {\mathfrak l} $.

Since we are especially interested in $ \bar{{\mathcal N}}_{g,1,d} $, let us modify the
above construction to get a mapping $ \bar{{\mathcal N}}_{g,1,d} \to \bar{{\mathcal N}}_{g+l,1,d+l} $. To do this it
is sufficient to define a mapping $ \bar{{\mathcal N}}_{g,1,d} \to \bar{{\mathcal N}}_{g,1+2l,d} $. We will define such a
mapping for every {\em most degenerate\/} curve, which is a semistable curve with
the class on a stratum of $ \dim =0 $ on the moduli space. The mapping will be
given by {\em gluing\/} this curve to the given curve. Since there are finitely
many most degenerate curves (they are enumerated by appropriately colored
trees), we will obtain a finite number of mappings $ \bar{{\mathcal N}}_{g,1,d} \to \bar{{\mathcal N}}_{g,1+2l,d} $.

A most degenerate curve $ \left(Z,\left\{z_{j}\right\}\right) $ of genus 0 with $ 2l+2 $ marked points
$ z_{j} $, $ 0\leq j\leq2l+1 $ has $ 2l-1 $ double points, and its normalization has $ 2l $
connected components. Any bundle of degree 0 over $ Z $ is trivial. Given a
curve $ \left(Y,P_{0},{\mathcal L}\right) $, glue $ Z $ to $ Y $ by identifying $ P_{0} $ and $ z_{0} $ and identify $ {\mathcal L}|_{P_{0}} $
with $ {\mathcal O}_{Z}|_{z_{0}} $ arbitrarily, denote the resulting bundle $ {\mathcal L}' $. The resulting
collection $ \left(Y\cup Z,\left\{z_{j}\right\},{\mathcal L}'\right) $ is obviously an element of $ \bar{{\mathcal N}}_{g,1+2l,d} $, and
gluing together $ z_{j} $, $ 1\leq j\leq2l $, we obtain an element of $ \bar{{\mathcal N}}_{g+l,1,d+l} $.

We see that $ {\mathfrak l}\left(X,P,{\mathcal L}\right) $ can be approximated by $ {\mathfrak l}\left(X',P',{\mathcal L}'\right) $ (here $ X $
and $ X' $ are smooth curves) with $ g\left(X'\right)<g\left(X\right) $ if $ \left(X,P,{\mathcal L}\right) $ is close to a
subset of codimension 2 which consists of semistable curves/bundles such
that a bundle has a double point.

To visualize better the above surgery we improve the description in
two ways. First, let us revisit the process of $ \varepsilon $-deformation.

We glue a hyperbola $ y_{1}y_{2}=\varepsilon $ to the coordinate lines $ y_{1}=0 $, $ y_{2}=0 $ via
the coordinate projections. One can momentarily see that this is
equivalent to gluing together two regions $ \left\{|y_{1}|\geq\sqrt{|\varepsilon|}\right\} $, $ \left\{|y_{2}|\geq\sqrt{|\varepsilon|}\right\} $
along the boundary via $ y_{2}=\varepsilon/y_{1} $. Given a bundle $ {\mathcal L} $ with trivializations
near $ Q_{1} $ and $ Q_{2} $ one obtains a bundle $ \overset{\,\,{}_\circ}{{\mathcal L}} $ on the resulting curves (its
sections are sections of $ {\mathcal L} $ outside of the disks with compatible
restrictions to the disk boundaries). The deformed bundle is $ \overset{\,\,{}_\circ}{{\mathcal L}}\left(Q\right) $, here $ Q $
is close to the circle $ |y_{1}|=|y_{2}|=\sqrt{|\varepsilon|} $ on the deformed curve.

There is an alternative description of the deformed curve. Above we
were removing two disks of radii $ \varepsilon_{1} $ and $ \varepsilon_{2} $ near points $ y_{1} $, $ y_{2} $ and gluing
the boundaries. The radii of disks were variable, so it was hard to
visualize the picture. Remove instead bigger disks of (fixed) radii $ \delta_{1} $,
$ \delta_{2} $. What remains is to glue two annuli $ \left\{\varepsilon_{i}<|z|<\delta_{i}\right\} $, $ i=1,2 $, to the
resulting manifold with a boundary, and glue inner boundaries of annuli
together. Each annulus is conformally equivalent to a cylinder with ratio
length/radius being $ \log \frac{\delta_{i}}{\varepsilon_{i}} $. Thus we need to glue in a cylinder
$ S^{1}\times\left(0,L\right) $ of {\em conformal length\/} $ L=\log \frac{\delta_{1}\delta_{2}}{\varepsilon_{1}\varepsilon_{2}} $. If $ \varepsilon_{1,2} $ decrease, it is
equivalent to gluing longer and longer {\em handles\/} between circles of radii
$ \delta_{1,2} $. Note also that we modify $ {\mathcal L} $ by adding a pole inside this handle.

Second, note that the above gluing mapping $ y_{2}=\varepsilon/y_{1} $ is defined not
only on the boundary of the above regions, but in the regions themselves.
Thus continuing the identification of boundaries, one can identify the
region $ \left\{|y_{2}|\geq\sqrt{|\varepsilon|}\right\} $ with $ \left\{|y_{1}|\leq\sqrt{|\varepsilon|}\right\} $. If a part of curve corresponds to
a subset of $ \left\{|y_{2}|\geq\sqrt{|\varepsilon|}\right\} $, then it is identified with a subset of
$ \left\{|y_{1}|\leq\sqrt{|\varepsilon|}\right\} $. Suppose that the curve we deform has rational smooth
components. If one component $ Y_{i} $ contains $ l_{i} $ double points, then to
perform the above deformation we remove $ l_{i} $ disks of radii $ \varepsilon_{ik} $, $ k=1,2 $,
around these points. What remains is a sphere without $ l_{i} $ disks, and we
need to glue several such spheres together identifying boundaries by
fraction-linear mappings. Say, if $ i\not=j $ and $ Y_{j} $ is glued to $ Y_{i} $ via gluing
$ \varphi_{i} $, then $ \varphi\left(Y_{j}\right) $ is a subset of a removed disk for $ Y_{i} $. Consider now the
union $ Y_{i}\cup\varphi_{i}\left(Y_{j}\right) $. It is again a sphere with several disks removed, and the
boundaries of this disks are identified with the boundaries of other
components (or different circles on $ Y_{i}\cup\varphi_{i}\left(Y_{j}\right) $) via fraction-linear
mappings.

Continuing this process as long as we can, the result is {\em one\/} sphere
with several disks removed, and boundaries of these disks are identified
pairwise via fraction-linear mappings. This is so called {\em Shottki model\/} of
the curve, and we see that the analogue of Deligne--Mumford
compactification we consider here leads naturally to this model of
deformed curves. One can do a similar thing in the case when the initial
curve contains non-rational components as well. In this case one should
restrict attention to the ``connected part'' of the curve which consists of
rational components.

We conclude that in the simplest case the above process of
deformation may be described as this (for one particular choice of
degenerate rational curve with 3 marked points): Take two points $ Q_{1} $, $ Q_{2} $
which are close to the marked point $ P $ and much closer to each other than
to $ P $. Now remove two non-intersecting disks around $ Q_{1} $ and $ Q_{2} $, and glue a
very long handle between boundaries of these disks. Here a disk which
contains $ P $, $ Q_{1} $ and $ Q_{2} $ plays the r\^ole of the ``old'' marked point, we assume
that $ z_{0} $ and $ z_{4} $ (using the notations from the beginning of this section)
are on the same smooth component of the degenerate rational curve.

\subsection{Second compactification } The previous section shows that compactified
moduli spaces are included one into another, thus one can consider the
direct limit, i.e., the {\em union\/} of these spaces. The mapping $ {\mathfrak l} $ identifies
this union $ {\mathfrak N}_{d-g}=\bigcup_{l}\bar{{\mathcal N}}_{g+l,1,d+l} $, with a subset of the Universal
Grassmannian, and (at last!) we have the ingredients necessary for the
discussion what the closure of this subset may look like.

Consider an element $ \left(Y,P,{\mathcal L}\right) $ of $ \bar{{\mathcal N}}_{g+l,1,d+l} $ for a very big $ l $. Suppose
that it is close to an element $ \left(Y',P',{\mathcal L}'\right) $ of $ \bar{{\mathcal N}}_{g,1,d} $. To get $ Y $ we glue a
most degenerate rational curve to $ Y' $, the result is a curve with $ 2l $
components and $ 2l+2 $ marked points. Next we glue marked points pairwise,
and $ \varepsilon $-deform the resulting curve at the double points. However, the
double points are naturally broken into two categories: double points on
the attached rational curve, and doubled marked points. Let us change the
order of gluing and deformations: glue in the degenerate curve, $ \varepsilon $-deform
the double points on this curve, then glue marked points together and
$ \varepsilon $-deforming them. Note that the former double points are not double
points of the linear bundle, but the latter ones are double points for
the bundle.

Deformation of double points on a degenerate rational curve leads
to a Shottki model for some curve, i.e., a sphere with several disks
removed. However, since the genus is 0, there is no removed disks at all,
thus we get a sphere with $ 2l+2 $ marked points. The condition that the
parameters of deformation are small is translated into the fact that
double ratios of marked points on the sphere are very big, i.e., that
{\em conformally\/} the centers are {\em well separated}. Deform now the double point
where the degenerated curve is glued to $ Y' $. This leads to the part of the
above sphere being identified (by a fraction-linear mapping) with a disk
around $ P' $, and the marked points go to points of this disk. We conclude
that this part of deformation corresponds to picking up a collection of
points which are close to $ P' $ and conformally well separated.

The second part of the deformation is the removing of small disks
around the marked points, and gluing together the boundaries of these
disks (or, what is the same, gluing long handles to boundaries of bigger
fixed disks). In the most important case $ g=0 $, thus the description of a
neighborhood of a point on $ {\mathfrak N}_{k} $ is related to studying curves of genus $ g $
obtained by gluing together $ 2g $ small conformally well separated disks on
$ {\mathbb C}P^{1} $. Note that the notion of being {\em well separated\/} depends on the
combinatorics of degenerated curves we glue in to get the limit point of
$ {\mathfrak N}_{k} $.

Since we expect that $ {\mathfrak l}\left({\mathfrak N}_{k}\right) $ is dense in the space of solutions of an
integrable system, one can approximate a point in this space by a
sequence of curves and bundles. Assuming the best case scenario, we can
obtain a term of this sequence by a small deformation of the previous
term (here a small deformation is taken in the sense of algebraic
geometry, i.e., one takes a couple of points on the curve, glues them
together and deforms the double point slightly). We see that such a point
of the space can be described as a {\em curve of potentially infinite\/} genus,
i.e., a sequence of the curves where the next one is obtained by gluing
long handles to the previous one.

In this paper we show that one can consider a {\em curve of actually
infinite\/} genus instead. Such a curve is obtained in the same manner as
the above sequence of finite-genus curves: one takes an infinite
collection of ``well-separated'' points on $ {\mathbb C}P^{1} $, removes a disk of a small
radius around each point, and glues the boundaries pairwise.\footnote{The exact meaning of the above terms ``well-separated collection'' and
circles being ``of small radii'' is the main topic of this paper. The
answers vary a little bit depending on what problem one considers, for
the list of results see Section~\ref{s0.120}.} We show that
under suitable conditions the standard theorems of algebraic geometry
hold for the resulting curves as well. This shows that the {\em completed\/}
moduli space $ \bar{{\mathfrak N}}_{-1} $ of such curves$ + $bundles may be important in algebraic
geometry. {\em If\/} this moduli space coincides with the whole space of
solutions (i.e., the best case scenario has place indeed), one can see
that the completed moduli space has a very simple topology. Indeed, by
Cauchy--Kovalevskaya theorem the space of solutions is identified with
the space of initial data for the solutions, which is a topological
linear space.

One can see that under assumption that the moduli space described in
this paper coincides with the set of solutions of integrable systems, the
study of compactified moduli spaces simplifies a lot by a transition to
the case of infinite genus. This simplification may suggest additional
approaches to the problem of studying the moduli spaces in finite genus
as well.

Unfortunately, at this moment it is unclear whether the best case
scenario is applicable to the moduli space as a whole. As we noted above,
the results of \cite{McKTru76Hil} show that this is true for the real part of
the hyperelliptic subset of the moduli space.

\subsection{Growth conditions and divisors }\label{s0.10}\myLabel{s0.10}\relax  The main target of this
paper is to describe how to quantify the conditions on the collection of
circles to be ``conformally well separated'' and ``of small radii'' so that
the main theorems of algebraic geometry are still valid for the resulting
curves of infinite genus. However, before we discuss the conditions on
the nonlinear data (curves and bundles), we should discuss the
simpler conditions on the linear data (sections of above bundles).

The basic theorem of $ 1 $-dimensional algebraic geometry states that
any linear bundle on a compact analytic curve has a finite-dimensional
space of global sections. The moment we try to drop the compactness
condition this theorem breaks, since any non-compact analytic curve is a
Stein manifold, thus any sheaf on it has a giant space of global
sections.

One way to fix this situation is to consider growth conditions. Let
$ \bar{X} $ be a compact curve with a linear bundle $ \bar{{\mathcal L}} $. Suppose that $ \bar{X} $ and $ \bar{{\mathcal L}} $ are
provided with Hermitian metrics. Pick up a point $ \infty\in\bar{X} $, let $ X=\bar{X}\smallsetminus\left\{\infty\right\} $, $ {\mathcal L}=\bar{{\mathcal L}}|_{X} $.
Then bounded sections of $ {\mathcal L} $ can be uniquely extended to analytic sections
of $ \bar{{\mathcal L}} $, thus $ {\mathcal L} $ has a finite-dimensional space of bounded sections.
Similarly, if we consider sections of $ {\mathcal L} $ with magnitude going to 0 when $ x
\to \infty $, it is the same as to consider sections of $ \bar{{\mathcal L}}\left(-1\cdot\infty\right) $ on $ \bar{X} $, here $ -1\cdot\infty $
is a divisor on $ \bar{X} $.

Going in a different direction, consider $ L_{2} $-sections of $ {\mathcal L} $. It is
easy to see that this is also equivalent to consideration of sections of
$ \bar{{\mathcal L}} $. If we consider a different metric on $ X $, say, $ \frac{dx^{2}}{\operatorname{dist}\left(x,\infty\right)^{2\alpha}} $, then
to consider holomorphic on $ X L_{2} $-sections is the same as to consider
$ \bar{{\mathcal L}}\left(\left[\alpha\right]\cdot\infty\right) $. Similarly, different choice of metric on $ {\mathcal L} $ will lead to
different shift of $ \bar{{\mathcal L}} $ by a divisor at $ \infty $.

We see that different growth conditions applied to sections of $ {\mathcal L} $
lead to different ``effective'' continuations of $ {\mathcal L} $ to $ \bar{X} $. Obviously, the
situation becomes more complicated when $ \bar{X}\smallsetminus X $ consists of more than
one point. If $ \bar{X}\smallsetminus X $ is discrete, different possible choices of growth
conditions lead to a lattice which is isomorphic to the lattice of
divisors on $ \bar{X} $ concentrated on $ \bar{X}\smallsetminus X $. Situation goes out of control if $ \bar{X}\smallsetminus X $
is a ``massive'' set, or $ \bar{X} $ is not a smooth manifold at all. In such cases
different choices of metrics on $ X $ and $ {\mathcal L} $ (or some other data controlling
the growth of sections) form a very complicated lattice. However, there
is a way to select ``good'' elements of this lattice.

\subsection{Riemann--Roch and growth control }\label{s0.20}\myLabel{s0.20}\relax  The Riemann--Roch theorem says
that on a compact analytic curve there is a relationship between a {\em degree
of\/} $ {\mathcal L} $, which is an easily calculable geometric characteristics of $ {\mathcal L} $, and
the dimensions of spaces of global sections of $ {\mathcal L} $ and $ \omega\otimes{\mathcal L}^{-1} $, $ \omega $ being the
linear bundle of holomorphic forms:
\begin{equation}
\dim  \Gamma\left({\mathcal L}\right)-\dim \Gamma\left(\omega\otimes{\mathcal L}^{-1}\right) = \deg {\mathcal L}-g+1.
\notag\end{equation}
Let us show how this theorem might be applied to picking up ``correct''
growth conditions on non-compact analytic curves.

Consider some fixed growth conditions for $ X $. They select a new sense
for the functor $ \Gamma\left(X,\bullet\right) $. Suppose again that $ \bar{X}\smallsetminus X $ is discrete, $ {\mathcal L}=\bar{{\mathcal L}}|_{X} $. We
have seen above that consideration of these growth conditions is
equivalent to consideration of $ \bar{{\mathcal L}}\left(D\right) $, here $ D $ is some divisor at infinity,
i.e., $ \Gamma\left(X,{\mathcal L}\right)=\Gamma\left(\bar{X},\bar{{\mathcal L}}\left(D\right)\right) $. One should expect that $ D $ does not depend on $ {\mathcal L} $ (at
least when $ {\mathcal L} $ does not change a lot). Suppose that there is a way to
determine $ \deg \bar{{\mathcal L}} $ basing on $ {\mathcal L} $, say, we describe bundles by
divisors on $ X $.

If the dimensions of sections of linear bundles on $ X $ with given
growth conditions satisfy the Riemann--Roch theorem, then an easy
calculation shows that $ 2D\sim0 $. In particular, suppose that $ {\mathcal L} $ was a
restriction of a bundle $ \bar{{\mathcal L}} $ of half-forms on $ \bar{X} $, i.e., $ \bar{{\mathcal L}}\otimes\bar{{\mathcal L}}\simeq\omega $. Then we see
that the bundle on $ \bar{X} $ which corresponds to $ \bar{{\mathcal L}}|_{X} $ with given growth
conditions is $ \bar{{\mathcal L}}\left(D\right) $. Note that $ \bar{{\mathcal L}}\left(D\right)\otimes\bar{{\mathcal L}}\left(D\right)=\omega $. In particular, the growth
conditions which satisfy Riemann--Roch theorem ``preserve'' the set of
bundles of half-forms (or $ \theta $-{\em characteristics\/}): the sections on $ X $ of one
$ \theta $-characteristic $ \bar{{\mathcal L}} $ which satisfy the growth conditions can be naturally
identified with sections on $ \bar{X} $ of some other $ \theta $-characteristic $ \bar{{\mathcal L}}\left(D\right) $.

Suppose now that the growth conditions on sections of $ {\mathcal L} $ are picked
up in some invariant way---whatever it means. Because of that we expect
that $ D $ is a linear combination of points on $ \bar{X}\smallsetminus X $ with the same
coefficients. Since $ \deg  D=0 $, thus $ D $ is 0.

We see that Riemann--Roch provides a selection criterion for picking
up growth conditions which ``do not add'' points at infinity to a given
divisor. By analogy, one can apply the same criterion in cases when $ \bar{X}\smallsetminus X $
is massive or $ \bar{X} $ does not exist: If there is an invariant way to describe
the growth conditions, and the Riemann--Roch theorem holds for a class of
bundles, then one may describe the part at infinity of the divisor of a
section of such a bundle: a section has no ``poles'' at infinity if it
satisfies the growth conditions.

\subsection{Handles, Sobolev spaces, and representations of $ \protect \operatorname{SL}_{2}\left({\mathbb C}\right) $ }\label{s0.30}\myLabel{s0.30}\relax  Let us
apply heuristics from the previous section to the case of a curve $ X $ of
infinite genus, i.e., a Riemannian sphere with an infinite number of
disks removed, and infinite number of handles glued in along the cut
lines (or just any surface glued in from ``pants'' in such a way that the
graph of gluing is of infinite genus). This is not a compact Riemannian
surface, so one needs to consider ``infinities'' of this surface, and fix
some growth conditions near these infinities.

Moreover, if the genus is infinite, one needs to be especially
careful with Riemann--Roch theorem, since the right-hand side contains
genus $ g\left(X\right) $ (i.e., number of handles) of the curve, which is not a number
any more, but infinity. To save the theorem, let us consider the quantity
$ \deg {\mathcal L}-g+1 $ as a unity. One way to do it is to consider some fixed bundle $ {\mathcal M} $
on $ X $ which is ``naturally constructed'' and satisfies the condition
$ \deg {\mathcal M}=g-1 $ for finite-genus $ X $, and consider sections of $ {\mathcal L}\otimes{\mathcal M} $ instead of
sections $ {\mathcal L} $. Then the Riemann--Roch theorem may be rewritten as
\begin{equation}
\dim  \Gamma\left({\mathcal L}\otimes{\mathcal M}\right)-\dim \Gamma\left(\omega\otimes{\mathcal M}^{-1}\otimes{\mathcal L}^{-1}\right) = \deg {\mathcal L}.
\notag\end{equation}
Note that $ \deg \omega\otimes{\mathcal M}^{-1}=\deg {\mathcal M} $, so $ {\mathcal N}=\omega\otimes{\mathcal M}^{\otimes-2} $ is a ``naturally constructed'' bundle
of degree 0. There are $ 2^{g} $ different square roots of a given bundle of
degree 0, and fixing one solution $ {\mathcal R} $ to $ {\mathcal R}^{2}=\omega\otimes{\mathcal M}^{\otimes-2} $ we may change $ {\mathcal M} $ to
$ {\mathcal M}\otimes{\mathcal R}^{-1} $. After such a change the formula simplifies to
\begin{equation}
\dim  \Gamma\left({\mathcal L}\otimes{\mathcal M}\right)-\dim \Gamma\left({\mathcal L}\otimes{\mathcal M}^{-1}\right) = \deg {\mathcal L},
\notag\end{equation}
which has an additional advantage of being symmetric w.r.t. $ {\mathcal L} \mapsto {\mathcal L}^{-1} $.
Moreover, $ {\mathcal M} $ is now a solution of $ {\mathcal M}^{2}=\omega $.

We come to the following heuristic: to consider Riemann--Roch
theorem for a curve $ X $ with infinitely many handles one should fix a
bundle of half-forms $ {\mathcal M}=\omega^{1/2} $ on $ X $, consider a linear bundle $ {\mathcal L} $ of finite
degree, and apply appropriate growth conditions to sections of $ {\mathcal L}\otimes{\mathcal M} $. We
expect that the choice of growth conditions is very restricted, since any
possible change is equivalent to consideration of different square root
of $ \omega $ (see Section~\ref{s0.20}).

One should expect that thus obtained growth conditions are in some
way ``invariant'' (since so rigid), thus geometrically defined. It is still
meaningful to apply these conditions in the case $ {\mathcal L}={\mathcal O} $, when $ {\mathcal L}\otimes{\mathcal M} $ is the
sheaf of half-forms. We come to the following conclusions:

Conditions of boundness (i.e., $ L_{\infty} $-topology) cannot be applied (since
there is no geometrically-defined norm on the fibers of the bundle of
half-forms), as well as $ L_{2} $-type restrictions (since a square of a
holomorphic half-form is a half-form on $ X_{{\mathbb R}} $, thus cannot be invariantly
integrated). One can easily see that the fourth degree of a half form $ \alpha $
(more precise, $ \alpha^{2}\otimes\bar{\alpha}^{2} $) is a top-degree-form on $ X_{{\mathbb R}} $, thus a good candidate
is the condition of integrability of fourth degree, i.e., $ L_{4} $-topology.

It turns out that generic Banach spaces (like $ L_{4} $) are not convenient
for cut-and-glue operations we are going to perform, but there is a way
to get a Hilbert structure instead of the Banach one. Note that $ s $-th
derivative of $ \alpha $-form changes as a $ \left(\alpha+s\right) $-form under coordinate
transformations (in the main term), thus its square is a top-degree form
on $ X_{{\mathbb R}} $ (in the main term) if $ \alpha+s=1 $. We conclude that the notion that
$ 1/2 ${\em -th derivative of a half-form is square-integrable\/} has a chance to be
geometrizable. This notion leads us to consideration of $ H^{1/2} $-Sobolev
spaces with values in (holomorphic) half-forms. (We discuss the basics of
the theory of Sobolev spaces in Section~\ref{h2}) In fact taking different $ s $
one gets an entire hierarchy of Banach spaces $ W_{\frac{2}{s+1/2}}^{s} $ which
interpolate between $ L_{4} $ and $ H^{1/2} $. However, only one of these spaces is a
Hilbert one, and having a Hilbert space is very convenient for our method
of divide and conquer.

As we will see it in Section~\ref{s2.70}, the spaces of functions we
obtain in such a way are closest possible analogues of Hardy spaces. Note
also that they are very small modifications of a particular
representation of $ \operatorname{SL}_{2}\left({\mathbb C}\right) $ from a supplementary series (see Section
~\ref{s2.60}). Indeed, the supplementary series is realized in $ s $-forms on
$ \left({\mathbb C}P^{1}\right)_{{\mathbb R}} $, and the Hilbert structure in these spaces is equivalent to
$ H^{1-2s} $-structure. Taking $ s=1/4 $, we see that the $ H^{1/2} $-structure on sections
of $ \omega^{1/4}\otimes\bar{\omega}^{1/4} $ may be defined\footnote{Recall that $ H^{s} $-structure on sections of bundles on manifolds is defined
only up to equivalence of topologies.} to be invariant w.r.t. fraction-linear
transformations. One should compare this result with the above heuristic
about $ H^{1/2} $-structure on sections of $ \omega^{1/2} $ being ``invariant in the main
term'' w.r.t. diffeomorphisms (we use an additional heuristic that
$ \omega^{1/4}\otimes\bar{\omega}^{1/4} $ is ``close'' to $ \omega^{1/2} $).

\subsection{Thick infinity }\label{s0.40}\myLabel{s0.40}\relax  We conclude that a working definition of {\em admissible\/}
sections on a sphere with some disks removed is that the section is a
half-form on the domain which can be extended as a $ H^{1/2} $-section of $ \omega^{1/2} $
to the whole sphere (though a usual heuristic about $ H^{1/2} $ is that it is a
smoothness class, it restricts the growth as well).

The usual definition of $ H^{s} $ on a manifold with a boundary (see
Section~\ref{s2.20}) allows one to work with function within the interior of
the manifold only. In our case the smooth part of the boundary (union of
circles) is not closed (if the number of circles is infinite), so one
loses the info about what happens on the {\em dust}, i.e., in the accumulation
points of the disks. This has some very inconvenient consequences. Let $ S $
be the sphere with interiors of disks removed, $ V $ be the closure of $ {\mathbb C}P^{1}\smallsetminus S $.
Then the sections we consider are elements of $ H^{1/2}\left({\mathbb C}P^{1}\smallsetminus V\right) $. The condition
of such a $ \varphi $ being holomorphic can be written as $ \bar{\partial}\varphi=0\in H^{-1/2}\left({\mathbb C}P^{1}\smallsetminus V\right) $ (since
taking derivative moves one notch down on Sobolev scale). However, if the
dust $ S\cap V $ is massive enough, there are functions of smoothness $ H^{-1/2} $ with
support in $ S\cap V $, thus there functions $ \varphi $ of smoothness $ H^{1/2} $ such that $ \bar{\partial}\varphi $
has support on the dust only. Since $ \bar{\partial}\varphi=0\in H^{-1/2}\left({\mathbb C}P^{1}\smallsetminus V\right) $, such functions
look like holomorphic ones, but in fact the dust is the part of infinity
on our curve of infinite genus, thus one should consider the divisor of
these functions as having components at infinity (since $ \bar{\partial}\varphi $ is not zero
there). This ruins our preparations to count the divisor at infinity
correctly.

One possible way out of this deadlock is to prohibit collections of
circles with a massive dust, say, restrict the Hausdorff dimension of the
dust to be smaller than 1. We pick up a different approach: we massage
the definition of the Sobolev space for a subset of a manifold in such a
way that it now takes into account the behaviour on the boundary even if
the boundary is massive. These are so called {\em generalized\/} Sobolev spaces.
As a corollary, we can consider pretty monstrous curves of infinite
genus.

For example, one can take a Serpinsky carpet\footnote{Which is a two-dimensional analogue of the Cantor set, obtained by
repeated removing of a smaller rectangle inside a bigger one (or removing
a triangle with vertices on the sides of a bigger triangle).} on a complex plane,
and take one disk inside each thrown away triangle/rectangle. Gluing the
boundaries of disk pairwise gives a curve we can deal with (if radii of
disks decrease quick enough). In particular, the dust can have a positive
measure (in fact in our theory the only condition on the dust is that it
is nowhere dense, so can be a set of accumulation points).

One can see that smooth points on the resulting curve have
infinitely many connected components, and it is the dust which keeps
these components together. The components are tubes, and the boundaries
of these tubes are glued to the dust, which is connected. We see that
points of the dust need to be considered as legitimate points on the
curve of infinite genus, and the curve is {\em not a topological manifold\/} in
neighborhoods of these points.

The striking fact is that it is possible to strengthen the
restriction on the radii of the circles in such a way that any global
holomorphic function has an {\em asymptotic\/} Taylor expansion near any point of
the dust (see Section~\ref{s4.95})! Thus the dust would consist of points
which one has a full right to call {\em smooth points}.

\subsection{Fight with $ H^{1/2} $ } While consideration of half-forms of smoothness $ H^{1/2} $
{\em enormously\/} simplifies the work with infinities on the curve, it is one of
the worst choices when we consider gluing conditions on sections on the
glued together circles on the complex sphere. The reason is that Sobolev
spaces with half-integer indices do not satisfy a lot of properties of
other Sobolev spaces when restrictions to hypersurfaces is considered
(this is the case when one usually considers Besov spaces instead of
Sobolev ones).

To fight with this, we define a notion of {\em mollified restriction\/} to a
hypersurface. Applying this restriction to the case of one removed disk
$ \left\{|z|>1\right\} $ one can momentarily see that the space $ H^{1/2} $ we consider coincides
with the Hardy space $ {\mathcal H}^{2} $ for the circle $ \left\{|z|=1\right\}! $ We see that the space of
``global holomorphic functions'' we consider is a generalization of the
Hardy space to the case when the boundary consists of many circles. What
is more, the Hardy space $ {\mathcal H}^{\infty} $ (of multiplicators in $ {\mathcal H}^{2} $) also plays an
important r\^ole when we consider equivalence of bundles in Section~\ref{s8.7}.

\subsection{Main results }\label{s0.120}\myLabel{s0.120}\relax  Since technical details take a lot of place
in the
course of discussion, we make a concise list of main results (without
discussing what the different conditions on the curve/bundle mean).

In Section~\ref{s35.40} we prove (an analogue of) Riemann--Roch theorem
for $ \bar{\partial}\colon {\mathcal O} \to \bar{\omega} $ in assumption that the matrix $ {\bold M}_{2}=\left(e^{-l_{ij}}-\delta_{ij}\right) $ (formed basing
on conformal distances $ l_{ij} $ between disks) gives a compact mapping $ l_{2} \to
l_{2} $. In Section~\ref{s7.90} we show that the mapping $ \bar{\partial}\colon {\mathcal O} \to \bar{\omega} $ has a maximal
rank compatible with its index. In Section~\ref{s9.20} we show an analogues
result for $ \bar{\partial}\colon \omega \to \omega\otimes\bar{\omega} $ under an additional assumption that disks have a
{\em thickening\/} (in fact we show that the mapping of taking $ A $-periods of
global homomorphic forms is an isomorphisms). In Section~\ref{s9.40} we
describe the image of the $ \left(A,B\right) $-period mapping. In all these cases the
mappings in question have an infinite index, so we massage these mappings
to get a mapping of index 0 (by extending the domain or target spaces)
which is in fact an isomorphism.

In Section~\ref{s7.30} we show that the space of global holomorphic
sections of $ \omega^{1/2}\otimes{\mathcal L} $ is finite-dimensional provided the matrix
$ {\bold M}_{1}=\left(a_{i}e^{-l_{ij}/2}-a_{i}\delta_{ij}\right) $: give a compact mapping $ l_{2} \to l_{2} $. Here $ \left(a_{i}\right) $ is some
sequence associated to the bundle $ {\mathcal L} $ (it consists of norms of gluing
cocycle for $ {\mathcal L} $). In the section~\ref{s7.40} we show that the duality
\begin{equation}
\operatorname{Coker}\left(\bar{\partial}\colon \omega^{1/2}\otimes{\mathcal L} \to \omega^{1/2}\otimes{\mathcal L}\otimes\bar{\omega}\right) \simeq \operatorname{Ker}\left(\bar{\partial}\colon \omega^{1/2}\otimes{\mathcal L}^{-1} \to \omega^{1/2}\otimes{\mathcal L}^{-1}\otimes\bar{\omega}\right)
\notag\end{equation}
holds as far as $ {\bold M}_{1} $ gives a bounded operator $ l_{2} \to l_{2} $.

In Section~\ref{s5.60} we show that the above mappings satisfy the
Riemann--Roch theorem as far as both $ \left(a_{i}e^{-l_{ij}/2}-a_{i}\delta_{ij}\right) $ and
$ \left(a'_{i}e^{-l_{ij}/2}-a'_{i}\delta_{ij}\right) $ give compact operators $ l_{2} \to l_{2} $ (here $ a'_{i} $ is the
sequence associated to $ {\mathcal L}^{-1} $). In the same section we show that these
mappings are Fredholm under an additional condition that $ \left(a_{i}\right) $, $ \left(a_{i}'\right) $ are
bounded.

In Section~\ref{s4.50} we show that the above conditions give no
restrictions on the {\em dust\/} of the corresponding M\"obius group. In Section
~\ref{s7.70} we discuss heuristics which show that our version of the
Riemann--Roch theorem is close to the strongest possible.

In Section~\ref{s8.40} we show that an existence of so called
{\em Hilbert\/}--{\em Schmidt\/} bundle on a curve shows that $ \omega^{1/2} $ is also
Hilbert--Schmidt. In Section~\ref{s8.50} we show that the set of
Hilbert--Schmidt bundles is convex (in the sense of multiplication of
bundles), thus it is possible to define a {\em group\/} of {\em strong
Hilbert\/}--{\em Schmidt\/} bundles. In Section~\ref{s8.7} we define the notion of two
bundles being {\em equivalent}, and define the Jacobian. In Section~\ref{s8.80} we
define divisors, and in Section~\ref{s8.90} the mapping to the Universal
Grassmannian (using asymptotic series at points of the dust defined in
Section~\ref{s4.95}).

In Section~\ref{s7.90} we show that if $ \sum_{i\not=j}e^{-l_{ij}}<\infty $, then any admissible
bundle of degree 0 can be described by locally-constant cocycles (and
some amplifications of this result). In Section~\ref{s9.60} we show that the
period matrix is symmetric, and the imaginary part is positive. In
Section~\ref{s9.70} we show that the Jacobian coincides with the quotient of
the space of possible periods by the lattice of {\em integer\/} periods, as far
as $ \sum_{i\not=j}e^{-l_{ij}}<\infty $. Finally, in Section~\ref{s9.80} we show how to modify the
notion of Hodge structure to describe the curves we study in this paper.

\subsection{Main tools } We already noted that to deal with massive infinities we
need the notion of generalized Sobolev space (introduced in Section
~\ref{s2.40}), and to deal with gluing conditions along the boundaries of
disks we need mollified restriction/extension mappings (defined in
Section~\ref{s5.31}).

We discuss the relationship of the ``invariant'' growth conditions (we
need them to have a good enumeration of divisor at infinity, see Section
~\ref{s0.20}) with the supplementary series of representations of $ \operatorname{SL}_{2}\left({\mathbb C}\right) $ in
Section~\ref{s2.60}. In Section~\ref{s2.70} we demonstrate that the classical Hardy
space $ {\mathcal H}^{2} $ is a particular case of spaces we consider here (for ``curves of
genus 1/2'', when the fundamental domain is a disk). In Section~\ref{s8.7} we
show that the space of multiplicators (needed for defining equivalence of
bundles) gives $ {\mathcal H}^{\infty} $ in the particular case of a disk as a fundamental
domain.

To handle divide-and-conquer strategy of dealing with curves of
infinite genus we consider {\em almost perpendicular\/} families of subspaces
(studied in Section~\ref{s5.61}). To study duality we employ spaces of
{\em strong sections\/} of bundles (see Section~\ref{s6.50}), to study $ \bar{\partial} $ as a
Fredholm operator we employ spaces of {\em weak sections\/} (see Section
~\ref{s5.61}). (In most important cases these spaces coincide.)

To define $ B $-periods on a curve we apply a method of averaging
(Section~\ref{s9.41}) of paths connecting two sides of an $ A $-cycle.

\subsection{Historic remarks } In the spring of 1981 Yu.~I.~Manin was overwhelmed
by unusual amount of sophomore students who asked him to be their
undergraduate advisor. Contrary to his habits of the time, Yuri Ivanovich
presented a list of topics he thought would be interesting and
instructive for us to work on. One of these topics was a generalization
of algebro-geometric description of solutions of KP to a curve of
infinite genus represented by its Shottki model. (It was R.~Ismailov who
started to work on this topic.)

Five years later, when integrable systems reached the peak of their
popularity, we obtained first results on the Riemann--Roch theorem which
were similar to the results of this paper. The curves were represented in
the way that was in fashion that time: as collections of pants glued
along boundaries. The restrictions on the curves/bundles were
jaw-breaking and very $ \operatorname{ad} $-hockish. Since Manin's suggestions were
completely forgotten already, there was nothing similar to the Shottki
model approach. However, half-forms and $ H^{1/2} $-topology were already
present at that time.

When in the end of the eighties the first variant of this paper was
taking shape, it became clear that the conditions on curves/sheaves
simplify {\em enormously\/} if one restricts attention to gluings of boundary
circles which are fraction-linear, and gluing conditions for the bundles
are given by locally constant cocycles (relative to half-forms). To our
great surprise, with these restrictions the pants would glue together
into the Shottki model of the curve (see Section~\ref{s0.8}), and sheaves
glue into the bundle of half-forms on $ {\mathbb C}P^{1} $. This finished a complete
circle returning the settings back into context of the question of
Manin's of 1981!

Note that though the importance of Shottki model, of half-forms and
of $ H^{1/2} $-growth conditions was clear long time ago, it was the
formalization of the mollified restriction mapping which emphasized the
capacity of infinity and the parallelism with the usual study of Hardy
spaces $ {\mathcal H}^{2} $ and $ {\mathcal H}^{\infty} $. This formalization (together with the rest of this
paper) would not appear if not the fruitful discussions with A.~Tyurina
in spring of 1996.

During the last several years Feldman, Kn\"orrer and Trubowitz made a
major breakthrough using an unrelated approach (cf.~\cite{FelKnoTru96Inf}).

The author is most grateful to I.~M.~Gelfand, A.~Givental,
A.~Goncharov, D.~Kazhdan, M.~Kontsevich, Yu.~I.~Manin, H.~McKean,
V.~Serganova, A.~Tyurina and members of A.~Morozov's seminar for
discussions which directed many approaches applied here. It was the
patient work of V.~Serganova which improved readability of the most
obscure pieces of this paper.

This work was partially supported by NSF grant and Sloan
foundation fellowship.

\section{Geometry of half-forms }

\subsection{Half-forms and holomorphic half-forms } Consider an oriented real
manifold $ M $. For $ \alpha\in{\mathbb C} $ consider a (complex) linear bundle $ \Omega_{M}^{\alpha} $ with
transition functions $ \left(\det \frac{\partial Y}{\partial X}\right)^{\alpha} $ if $ Y=Y\left(X\right) $. We will write a section of
this bundle as $ \varphi\left(X\right)dX^{\alpha} $, so $ \varphi\left(X\right)dX^{\alpha}=\psi\left(Y\right)dY^{\alpha} $ if $ Y=Y\left(X\right) $ and
$ \varphi\left(X\right)\left(\det \frac{\partial Y}{\partial X}\right)^{\alpha}=\psi\left(Y\right) $. When $ M $ is clear from context, we will denote $ \Omega_{M}^{\alpha} $
by $ \Omega^{\alpha} $.

In what follows we will be most interested in $ \Omega^{1/2} $. Note that
$ \Omega^{1/2}\otimes\Omega^{1/2} $ is canonically isomorphic to the bundle of top de Rham forms.
For any two global sections $ \varphi $, $ \psi $ of $ \Omega^{1/2} $ the product $ \varphi\bar{\psi} $ is a top de Rham
form, thus the pairing
\begin{equation}
\left(\varphi,\psi\right) \mapsto \int\varphi\bar{\psi}
\notag\end{equation}
gives a structure of pre-Hilbert space on $ \Gamma\left(M,\Omega^{1/2}\right) $. In particular, a
notion of $ L_{2} $-section of $ \Omega^{1/2} $ is canonically defined.

If $ M $ is a complex-analytic manifold, the power $ \left(\det \frac{\partial Y}{\partial X}\right)^{\alpha} $ is not
uniquely defined (unless $ \alpha\in{\mathbb Z} $), thus $ \Omega^{\alpha} $ is not canonically defined. We use
an $ \operatorname{ad} $ hoc definition of an analogue of $ \Omega^{1/2} $ as follows:

Denote by $ \omega $ the top holomorphic de Rham bundle for $ M $. Consider a
linear bundle $ l $ with an isomorphism $ i\colon l\otimes l\simeq\omega $. We call sections of such a
bundle {\em holomorphic half-forms}, and will denote $ l = \omega^{1/2} $. In what follows
we will consider at most one bundle of half-forms on a given manifold.

Note that any local diffeomorphism $ M \to N $ induces a bundle of
half-forms on $ M $ given a bundle of half-forms on $ N $.

\subsection{Half-forms on complex curves }\label{s1.20}\myLabel{s1.20}\relax  In what follows we are most
interested
in complex curves. On a compact complex curve of genus $ g $ there are $ 2^{g} $
different bundles of half-forms. If $ g=0 $, then there is a uniquely defined
bundle of half-forms, isomorphic to $ {\mathcal O}\left(-1\right) $.

To pick up an isomorphism of $ {\mathcal O}\left(-1\right)^{\otimes2}={\mathcal O}\left(-2\right) $ with $ \omega $, one must fix an
element $ v\in\Lambda^{2}V^{*} $, here $ V $ is $ 2 $-dimensional space, projectivization of which
is $ {\mathbb C}P^{1} $. Indeed, let $ {\bold e} $ be the Euler vector field on $ V $, $ {\bold e} = z_{1}\frac{\partial}{\partial z_{1}} +
z_{2}\frac{\partial}{\partial z_{2}} $. Consider $ v $ as a constant $ 2 $-form on $ V $, then $ {\bold e}\lrcorner v $ is a $ 1 $-form on
$ V $ which induces a global section of $ \omega\otimes{\mathcal O}\left(2\right) $ on $ PV. $ Multiplication on this
section gives an isomorphism of $ {\mathcal O}\left(-2\right) $ and $ \omega $.

This shows that there should be an action of $ \operatorname{SL}\left(2,{\mathbb C}\right) $ on the bundle
of half-forms, compatible with the action of $ \operatorname{SL}\left(2,{\mathbb C}\right) $ on $ {\mathbb C}P^{1} $. Indeed, in
affine coordinate $ z $ on $ {\mathbb C}P^{1}\smallsetminus\left\{\infty\right\} $, a half-form may be written as $ \varphi\left(z\right)dz^{1/2} $.
Here if domain of $ \varphi $ contains $ \infty $, then $ \varphi\left(\infty\right)=0 $. To write an action of
$ m=\left(
\begin{matrix}
a & b \\ c & d
\end{matrix}
\right)\in\operatorname{SL}\left(2,{\mathbb C}\right) $, note that the derivative of the corresponding action
$ z'=\displaystyle\frac{az+b}{cz+d} $ on $ {\mathbb C}P^{1} $ has a canonical square root $ \frac{1}{cz+d} $. Thus
$ \varphi\left(z'\right)dz'{}^{1/2}\sim\varphi\left(z'\right)\frac{1}{cz+d}dz^{1/2} $.

Note that it is the restriction of this action on $ \operatorname{SL}\left(2,{\mathbb R}\right) $ what is
used to define automorphic forms on $ \left\{\operatorname{Im} z>0\right\} $. Note also that this action
cannot be pushed down to an action of $ \operatorname{PGL}\left(2,{\mathbb C}\right) $.

\subsection{Restriction onto $ S^{1} $ } Consider now a annulus $ B=\left\{r<|z|<R\right\} $ on a complex
plane. Obviously, there are only two different bundles of half-forms on
$ B $. The inclusion of $ B $ into $ {\mathbb C}P^{1} $ induces a bundle $ \omega^{1/2} $ of half-forms on $ B $.
Another bundle of half-forms is $ \omega^{1/2}\otimes\mu $, here $ \mu $ is a locally constant
sheaf with monodromy $ -1 $ on $ B $, obviously $ \mu^{2} $ is the constant sheaf $ \underline{{\mathbb C}}_{ } $. Note
that the holomorphic form $ \frac{dz}{iz} $ on $ B $ has no square root in $ \omega^{1/2} $, but
has two in $ \omega^{1/2}\otimes\mu $, one may write these square roots symbolically as
$ \pm\frac{\left(dz\right)^{1/2}}{i^{1/2}z^{1/2}} $.

Consider now a restriction of $ \omega^{1/2}\otimes\mu $ onto $ S^{1}=\left\{|z|=1\right\}\subset B $ (we assume
that $ r<1<R $). Consider a coordinate $ t = \operatorname{Im} \ln  z $ on $ S^{1} $, and identify
$ \frac{dz^{1/2}}{i^{1/2}z^{1/2}}\in\Gamma\left(S^{1},\omega^{1/2}\otimes\mu\right) $ with $ dt^{1/2}\in\Gamma\left(S^{1},\Omega_{S^{1}}^{1/2}\right) $. A simple
calculation shows that this identification gives an isomorphism of
$ \omega^{1/2}|_{S^{1}} $ with $ \Omega_{S^{1}}^{1/2}\otimes\mu $, which is preserved by real-analytic
diffeomorphisms of $ S^{1} $ with a given lift to $ \mu $.

In general, an immersion of $ S^{1} $ in a complex curve $ X $ with a bundle of
half-forms $ \omega $ induces an isomorphism of $ \omega^{1/2}|_{S^{1}} $ either with $ \Omega_{S^{1}}^{1/2}\otimes\mu $, or with
$ \Omega_{S^{1}}^{1/2} $. One may call immersions of the first type $ A $-cycles, of the second
one $ B $-cycles. Until Section~\ref{h9} the immersions we consider are going to be
$ A $-cycles, thus we will always have a factor $ \mu $.

Note that the above isomorphism is defined up to multiplication by
$ \pm1 $, and both $ \Omega_{S^{1}}^{1/2}\otimes\mu $ and $ \Omega_{S^{1}}^{1/2} $ have natural pre-Hilbert structures,
thus $ \omega^{1/2}|_{S^{1}} $ has a natural pre-Hilbert structure as well.

\section{Sobolev spaces }\label{h2}\myLabel{h2}\relax 

In what follows $ {\mathcal S} $ denotes the space of rapidly decreasing smooth
functions, $ {\mathcal D} $ denotes the space of smooth functions with compact support,
and $ {\mathcal S}' $ and $ {\mathcal D}' $ are the dual spaces. $ l_{2} $ denotes the Hilbert space of
square-integrable sequences, and $ L_{2} $ denotes the space of square-integrable
functions. The
symbol $ \bigoplus V_{i} $ denotes the space of sequences with finite number of non-zero
terms. If $ V_{i} $ are Hilbert spaces, $ \bigoplus_{l_{2}}V_{i} $ denotes the space of
square-integrable sequences $ \left(v_{i}\right) $ with $ v_{i}\in V_{i} $.

The material of this section is mostly standard, however, note that
we also discuss some topics which are not parts of standard curriculum of
Sobolev spaces: in Section~\ref{s2.40} we introduce {\em generalized\/} Sobolev
spaces, in Sections~\ref{s2.60} and~\ref{s2.70} we give descriptions of
supplementary series of representations of $ \operatorname{SL}\left(2,{\mathbb C}\right) $ and of Hardy space in
terms of Sobolev spaces. In addition to the results which are covered in
this section, in Section~\ref{s2.10} we introduce {\em mollifications\/} of
mappings of restriction to a submanifold (in the case $ s=1/2 $), and of
mapping of extension from submanifold (in the case $ s=0 $), both in the
cases when the non-mollified mappings are not continuous.

\subsection{Euclidean case }\label{s2.02}\myLabel{s2.02}\relax  Consider a vector space $ {\mathbb R}^{n} $ with the standard
Euclidean structure. Consider $ s\in{\mathbb R} $ and a norm $ \|\bullet\|_{s} $ on the set $ {\mathcal S} $ of rapidly
decreasing $ C^{\infty} $-functions:
\begin{equation}
\|f\left(x\right)\|_{s} = \int\left(1+|\xi|^{2}\right)^{s}|\widehat{f}\left(\xi\right)|^{2}d\xi,
\notag\end{equation}
here $ \widehat{f}\left(\xi\right) $ is the Fourier transform of $ f $. By definition, {\em the Sobolev space\/}
$ H^{s}\left({\mathbb R}^{n}\right) $ is the completion of $ {\mathcal S} $ w.r.t. this norm. It is naturally
isomorphic to the set of $ L_{2} $-functions $ g\left(\xi\right) $ w.r.t. the measure
$ \left(1+|\xi|^{2}\right)^{s/2}d\xi $. Since the inverse Fourier transform of $ g\left(\xi\right) $ is a
well-defined generalized function on $ {\mathbb R}^{n} $, there is a natural inclusion
$ H^{s}\left({\mathbb R}^{n}\right)\hookrightarrow{\mathcal S}' $, compatible with the inclusion $ {\mathcal S}\hookrightarrow H^{s}\left({\mathbb R}^{n}\right) $.

The pairing
\begin{equation}
\left(f,g\right) \mapsto \int f\bar{g}\,dx
\notag\end{equation}
extends to a pairing $ H^{s}\left({\mathbb R}^{n}\right)\otimes H^{-s}\left({\mathbb R}^{n}\right) \to {\mathbb C} $ which is a pairing of Hilbert
spaces. This pairing is compatible with the pairing between $ {\mathcal S} $ and $ {\mathcal S}' $.

Let $ D $ be a closed subset of $ {\mathbb R}^{n} $, $ U $ be an open subset of $ {\mathbb R}^{n} $. Define a
Hilbert subspace
\begin{equation}
\overset{\,\,{}_\circ}{H}^{s}\left(D\right)=H^{s}\left({\mathbb R}^{n}\right) \cap \left\{f\in{\mathcal S}' \mid \operatorname{Supp} f\subset D\right\},
\notag\end{equation}
and a quotient Hilbert space
\begin{equation}
H^{s}\left(U\right)=H^{s}\left({\mathbb R}^{n}\right)/\overset{\,\,{}_\circ}{H}\left({\mathbb R}^{n}\smallsetminus U\right).
\notag\end{equation}

The following properties of Sobolev spaces are most important in
analysis:
\begin{enumerate}
\item
If $ D $ is compact, $ V $, $ U $ are open, $ V\subset D\subset U\subset{\mathbb R}^{n} $, and $ \varphi\colon U \to \varphi\left(U\right) $ is a
diffeomorphism, then $ \varphi^{*} $ induces invertible bounded operators $ \overset{\,\,{}_\circ}{H}^{s}\left(\varphi\left(D\right)\right) \to
\overset{\,\,{}_\circ}{H}^{s}\left(D\right) $, $ H^{s}\left(V\right) \to H^{s}\left(V\right) $.
\item
If $ D $ is compact, then any differential operator $ A $ of degree $ d $ with
smooth coefficients gives a bounded operators $ \overset{\,\,{}_\circ}{H}^{s}\left(D\right) \to \overset{\,\,{}_\circ}{H}^{s-d}\left(D\right) $, $ H^{s}\left(V\right) \to
H^{s-d}\left(V\right) $.
\item
In the same way a pseudodifferential operator of degree $ d $
gives a bounded operators $ \overset{\,\,{}_\circ}{H}^{s}\left(D\right) \to H^{s-d}\left(U\right) $.
\item
The mapping of restriction of smooth functions onto $ {\mathbb R}^{n-k}\subset{\mathbb R}^{n} $ extends
to a bounded operator
\begin{equation}
r\colon H^{s}\left({\mathbb R}^{n}\right) \to H^{s-\frac{k}{2}}\left({\mathbb R}^{n-k}\right)\text{ if }s>\frac{k}{2}.
\notag\end{equation}
\item
Dually, extension-by-$ \delta $-function of generalized functions on $ {\mathbb R}^{n-k} $ to
generalized functions on $ {\mathbb R}^{n} $ gives a bounded operator $ e\colon H^{s}\left({\mathbb R}^{n-k}\right) \to
H^{s-\frac{k}{2}}\left({\mathbb R}^{n}\right) $ if $ s<0 $.
\item
If $ D $ is compact and has a smooth boundary, then $ \overset{\,\,{}_\circ}{H}^{s}\left(D\right)^{\perp}\subset H^{-s}\left({\mathbb R}^{n}\right) $
coincides with $ \overset{\,\,{}_\circ}{H}^{-s}\left(\overline{{\mathbb R}^{n}\smallsetminus D}\right) $.
\end{enumerate}
\subsection{Sobolev spaces on manifolds }\label{s2.20}\myLabel{s2.20}\relax  If $ M $ is a paracompact $ C^{\infty} $-manifold, then
$ H^{s}\left(M\right) $ can be defined as a subspace of generalized functions on $ M $
which consists of functions $ f $ such that for any $ U\subset M $ and $ \varphi\colon \widetilde{U} \overset{\sim}\to \widetilde{V}\subset{\mathbb R}^{n} $ the
restriction $ f|_{\widetilde{U}} $ satisfies $ \left(\varphi^{-1}\right)^{*}f\in H^{s}\left(V\right) $ (here $ \widetilde{U} $ is a neighborhood of $ \bar{U} $,
$ V=\varphi\left(U\right) $).

An alternative definition is that for an appropriate partition of
unity $ \left(U_{\alpha},\sigma_{\alpha}\right) $, $ \operatorname{Supp} \sigma_{\alpha}\subset\subset U_{\alpha} $, $ \sum_{\alpha}\sigma_{\alpha}=1 $, the products $ \sigma_{\alpha}f $ satisfy $ \left(\varphi_{\alpha}^{-1}\right)^{*}\left(\sigma_{\alpha}f\right)
\in \overset{\,\,{}_\circ}{H}\left(V_{\alpha}\right) $, here $ \varphi_{a}\colon U_{\alpha} \overset{\sim}\to V_{\alpha}\subset{\mathbb R}^{n} $.

In the same way one may define $ H^{s} $-sections of $ C^{\infty} $-vector bundles
over $ M $. We denote the space of $ H^{s} $-section of a bundle $ E $ by $ H^{s}\left(M,E\right) $.
Obviously, if $ M $ is compact, one may define a structure of Hilbert
space on $ H^{s}\left(M\right) $, but this structure is not uniquely defined. However, the
corresponding topology on $ H^{s}\left(M\right) $ is well-defined. For a general manifold $ M $
one defines a topology on $ H^{s}\left(M\right) $ as an inverse limit w.r.t. topologies on
$ H^{s}\left(U\right) $, $ U $ being open subsets of $ M $ diffeomorphic to bounded open subsets of
$ {\mathbb R}^{n} $.

If $ D $ is a closed subset of $ M $, and $ U $ is an open subset, define
$ \overset{\,\,{}_\circ}{H}^{s}\left(D\right)\subset H^{s}\left(M\right) $ as subspace of function with support in $ D $, and $ H^{s}\left(U\right) $ as
$ H^{s}\left(M\right)/\overset{\,\,{}_\circ}{H}^{s}\left(M\smallsetminus U\right) $.

Note that if $ M={\mathbb R}^{n} $, the above definition produces a different space
$ H_{\text{loc}}^{s}\left({\mathbb R}^{n}\right) $ of functions on $ {\mathbb R}^{n} $ than the space $ H^{s}\left({\mathbb R}^{n}\right) $ defined in Section
~\ref{s2.02}. A generalized function $ f $ on $ {\mathbb R}^{n} $ is in $ H_{\text{loc}}^{s}\left({\mathbb R}^{n}\right) $ if in any bounded
domain it is equal to a function from $ H^{s}\left({\mathbb R}^{n}\right) $.

The properties of Sobolev spaces $ H^{s}\left({\mathbb R}^{n}\right) $ have direct analogues for
$ H^{s}\left(M\right) $, thus diffeomorphisms of manifolds and differential operators act
on $ H^{s}\left(M\right) $, and one can restrict/extend Sobolev sections to/from
submanifolds. In particular, the existence of a parametrix for an elliptic
differential operator shows that
on a compact manifold $ M $ an elliptic operator $ A $ of degree $ d $ gives a Fredholm
operator $ H^{s}\left(M\right) \to H^{s-d}\left(M\right) $ (with suitable changes if $ A $ acts in vector
bundles).

Note that on compact manifolds one can define the Sobolev spaces
using arbitrary elliptic (pseudo)differential operators:

\begin{proposition} Consider a compact manifold $ M $ with a metric and a
positive self-adjoint elliptic operator $ A $ of degree $ d $. Then the $ s $-Sobolev
norm is equivalent to the norm
\begin{equation}
\|f\| = \int_{M}|\left(1+A\right)^{s/d}f|^{2}d\mu.
\notag\end{equation}
\end{proposition}

We will also use the following statement:

\begin{proposition} Consider two linear bundles $ {\mathcal L}_{1} $ and $ {\mathcal L}_{2} $ on a compact
manifold $ M $ such that $ {\mathcal L}_{1}\otimes{\mathcal L}_{2}\simeq\Omega^{\text{top}}\left(M\right) $. Then the pairing $ \int\alpha\beta $ between smooth
sections of $ {\mathcal L}_{1} $ and $ {\mathcal L}_{2} $ can be extended to non-degenerate pairing between
Hilbert spaces $ H^{s}\left(M,{\mathcal L}_{1}\right) $ and $ H^{-s}\left(M,{\mathcal L}_{2}\right) $ for an arbitrary $ s\in{\mathbb R} $.

\end{proposition}

\subsection{Capacity } Consider a manifold $ M $. We say that $ s $-{\em capacity\/} of a closed subset
$ S\subset M $ is
non-zero, if there is a non-zero function $ f\in H^{s}\left(M\right) $ such that $ \operatorname{Supp} f\subset S $. We
say that $ S $ has a {\em capacity dimension\/} $ \geq d $, if $ s $-capacity of $ S $ is non-zero
for $ s=-\frac{\dim  M - d}{2} $.

Since $ \delta $-function of a point belongs to $ H^{s}\left(M\right) $ with $ s<-\frac{\dim  M}{2} $, a
capacity dimension of a point is $ \geq d $ if $ d<0 $. It is clear that this
estimate cannot be improved. In the same
way the capacity dimension of a submanifold $ S $ is $ \geq d $ iff $ d<\dim  S $.

For a general subset $ S $ the same is true if one considers
Hausdorff dimension \cite{Hor83Dis}.

\subsection{Generalized Sobolev spaces }\label{s2.40}\myLabel{s2.40}\relax  In Section~\ref{s5.10} we consider the space
of $ H^{1/2} $-half-forms which are holomorphic outside of a closure $ \bar{U} $ of a union $ U $
of disks in $ {\mathbb C}P^{1} $. We will see that this space is too big for our
purposes if $ \bar{U} $
is ``much bigger'' than $ U $. We will need to
consider half-forms which ``are holomorphic'' in $ \bar{U}\smallsetminus U $ as well as in $ {\mathbb C}P^{1}\smallsetminus\bar{U} $.

To have this we need the half-forms to be {\em defined\/} in $ {\mathbb C}P^{1}\smallsetminus U $ instead
of $ {\mathbb C}P^{1}\smallsetminus\bar{U} $. The problem with this is that $ U $ is not closed, so the usual
definition of $ H^{1/2} $-half-forms as of a quotient-space does not work. This
shows the need for the following

\begin{definition} Consider a subset $ U $ of the manifold $ M $. Let $ \overset{\,\,{}_\circ} H^{s}\left(U\right) $ denotes
the closure of the subspace $ L\subset H^{s}\left(M\right) $
\begin{equation}
L=\left\{f\in H^{s}\left(M\right) \mid \operatorname{Supp} f\subset U\right\}.
\notag\end{equation}
Consider $ V\subset M $. Let $ H^{s}\left(V\right)=H^{s}\left(M\right)/\overset{\,\,{}_\circ} H^{s}\left(M\smallsetminus V\right) $. \end{definition}

Note that we do not require that the subset $ U $ is closed, and $ V $ is
open. If they are, then we get the standard definitions of $ \overset{\,\,{}_\circ}{H}^{s} $ and $ H^{s} $.
Obviously, $ \overset{\,\,{}_\circ} H^{s}\left(U\right)\subset\overset{\,\,{}_\circ} H^{s}\left(\bar{U}\right) $,
but this inclusion may be proper, as the following construction shows\footnote{The construction can be simplified, but in the current form it is an
example of domains we are going to deal with.}.

Consider a disjoint family of closed subsets $ V_{i}\subset M $, Let $ {\mathcal V}=\overline{\bigcup V_{i}}\smallsetminus
\bigcup V_{i} $. Suppose that
\begin{enumerate}
\item
The natural mapping $ \bigoplus_{l_{2}}\overset{\,\,{}_\circ} H^{s}\left(V_{i}\right) \xrightarrow[]{\iota} H^{s}\left(M\right) $ is a (continuous)
injection\footnote{I.e., the image is closed, and the mapping is an isomorphism on the
image.}.
\item
The $ s $-capacity of $ {\mathcal V} $ is positive.
\end{enumerate}
The first condition insures that the space $ \overset{\,\,{}_\circ} H^{s}\left(\bigcup V_{i}\right) $ is the image of
the mapping $ \iota $. Hence any non-zero function $ f\in\overset{\,\,{}_\circ} H^{s}\left(\bigcup V_{i}\right) $ satisfies the
condition $ \operatorname{Supp} f \cap \bigcup V_{i}\not=\varnothing $.

On the other hand, the second condition shows that there is a
non-zero function $ f\in H^{s}\left(M\right) $ such that $ \operatorname{Supp} f\subset{\mathcal V} $. Obviously,
$ f\in\overset{\,\,{}_\circ} H^{s}\left(\overline{\bigcup V_{i}}\right) $, but $ f\notin\overset{\,\,{}_\circ} H^{s}\left(\bigcup V_{i}\right) $.

In Section~\ref{s4.50} we show how to construct a family of disks
$ V_{i} $ which satisfy the first condition. The centers of these disks may be
an arbitrary prescribed locally discrete set. Moreover, one can easily find an
appropriate set of centers such that the corresponding set $ {\mathcal V} $ does not
depend on radii and coincides with an arbitrary given closed set $ {\mathcal V}_{0} $ with
empty interior.

In particular, if there exist a closed set $ {\mathcal V} $ with an empty interior
and non-zero $ s $-capacity, then one can construct a subset $ U $ of $ {\mathbb C}P^{1} $ such that
the inclusion $ \overset{\,\,{}_\circ} H^{s}\left(U\right)\subset\overset{\,\,{}_\circ} H^{s}\left(\bar{U}\right) $ is proper. Taking appropriate Cantor sets,
one gets that this is possible with any $ s<0 $.

On the other hand, if $ U $ has smooth boundary, then $ \overset{\,\,{}_\circ} H^{s}\left(U\right)=\overset{\,\,{}_\circ} H^{s}\left(\bar{U}\right) $.

\subsection{Rescaling on $ {\mathbb R}^{n} $ and $ {\mathbb C}^{n} $ }\label{s2.50}\myLabel{s2.50}\relax  Consider an $ H^{s} $-function $ f\left(x\right) $ on $ {\mathbb R}^{n} $ with
compact
support. The Fourier transform $ \widehat{f}\left(\xi\right) $ is real-analytic, thus the integral
\begin{equation}
\int|\xi|^{2s}|\widehat{f}\left(\xi\right)|^{2}d\xi,
\label{equ3.7}\end{equation}\myLabel{equ3.7,}\relax 
can be defined in the sense of generalized functions (i.e., analytic
continuation in $ s $) near $ \xi=0 $. Since $ f\in H^{s} $, the integral
converges near $ \infty $. Moreover, if $ s>-\frac{n}{2} $, the integral converges near $ \xi=0 $,
thus its value is a limit of Riemann sums, hence is non-negative.

In any case the integral defines a quadratic form on $ \overset{\,\,{}_\circ}{H}^{s}\left(D\right) $ for any
compact $ D $.

\begin{nwthrmi} If $ s>-\frac{n}{2} $, the integral defines a norm on $ \overset{\,\,{}_\circ}{H}^{s}\left(D\right) $, and this form
is equivalent to the Hilbert norm. \end{nwthrmi}

The advantage of the norm~\eqref{equ3.7} is the fact that it is covariant
with respect to dilatations. Define a mapping $ {\mathfrak D}_{a}^{s}\colon {\mathcal S} \to {\mathcal S} $ by
\begin{equation}
\left({\mathfrak D}_{a}^{s}f\right)\left(x\right) = a^{s}f\left(ax\right),\qquad a,s\in{\mathbb R},\quad a\not=0.
\notag\end{equation}
Then $ {\mathfrak D}_{a}^{-s+\frac{n}{2}} $ is an isometry w.r.t. the norm~\eqref{equ3.7}. Note also that
that $ {\mathfrak D}_{a}^{tn} $ is a natural action of a dilatation on sections of the linear
bundle $ \Omega^{t} $.

Combining this with the previous statement, we conclude that any
similarity transform\footnote{I.e., a composition of translations, rotations and dilatations.} $ {\bold T} $ acts as a uniformly bounded operator in the space
of $ H^{s} $-sections of $ \Omega^{\alpha} $
\begin{equation}
{\bold T}^{*}\colon \overset{\,\,{}_\circ}{H}^{s}\left({\bold T}D,\Omega^{\alpha}\right) \to \overset{\,\,{}_\circ}{H}^{s}\left(D,\Omega^{\alpha}\right),\qquad \alpha = \frac{1}{2}-\frac{s}{n},
\notag\end{equation}
as far as both $ D $ and $ {\bold T}D $ remain in the same disk $ \left\{|x|<R\right\} $, and $ R $ is fixed.

If we consider $ {\mathbb C}^{n}={\mathbb R}^{2n} $ and the linear bundle $ \omega^{1/2} $ over $ {\mathbb C}^{n} $, then the
action of dilatations\footnote{In fact one needs to consider $ 2 $-sheeted covering of the group of
similarity transforms.} on this bundle is the same\footnote{Strictly speaking, differs on a multiplication by a constant of
magnitude 1.} as on $ \Omega_{{\mathbb R}^{2n}}^{1/4} $. We
conclude that similarity transforms act as uniformly bounded transformations
\begin{equation}
\overset{\,\,{}_\circ}{H}^{n/2}\left({\bold T}D,\omega^{1/2}\right) \to \overset{\,\,{}_\circ}{H}^{n/2}\left(D,\omega^{1/2}\right)
\notag\end{equation}
as far as both $ D $ and $ {\bold T}D $ remain in the same disk $ \left\{|z|<R\right\} $, and $ R $ is fixed.

\subsection{$ \protect \operatorname{PGL}\left(2,{\mathbb C}\right) $-invariant realization of $ H^{1/2}\left({\mathbb C}P^{1}\right) $. }\label{s2.60}\myLabel{s2.60}\relax 

Consider double ratio
\begin{equation}
\left(a:b:c:d\right) = \frac{a-b}{b-d}:\frac{a-c}{c-d},\qquad a,b,c,d\in{\mathbb C}P^{1}.
\notag\end{equation}
Double ratio $ \left(z_{1}:z_{1}+\delta z_{1}:z_{2}:z_{2}+\delta z_{2}\right) $ gives a section $ \rho $ of the bundle
$ T^{*}{\mathbb C}P^{1}\boxtimes T^{*}{\mathbb C}P^{1} $ over $ {\mathbb C}P^{1}\times{\mathbb C}P^{1} $.
In affine coordinates it can be written as $ \frac{1}{\left(z_{1}-z_{2}\right)^{2}}dz_{1}\,dz_{2} $. By
construction this section is invariant with respect to the action of
$ \operatorname{PGL}\left(2,{\mathbb C}\right) $.

Consider $ K_{s}=\frac{\rho^{s} \bar{\rho}^{s}}{\Gamma\left(-2s+1\right)} $. It is a section of $ \Omega^{s}\left({\mathbb C}P_{{\mathbb R}}^{1}\right) \boxtimes \Omega^{s}\left({\mathbb C}P_{{\mathbb R}}^{1}\right) $,
which in local coordinates looks like
$ K_{s}\left(z_{1},z_{2}\right)=\frac{1}{|z_{1}-z_{2}|^{4s}\Gamma\left(-2s+1\right)}dZ_{1}^{s}\,dZ_{2}^{s} $ and has no singularity outside
of $ z_{1}=z_{2} $. Here $ dZ=dx\wedge dy $ if $ z=x+iy. $ The operator with kernel $ K_{s} $ defines a
pairing
\begin{equation}
\left(\alpha,\beta\right)_{s}=\int_{{\mathbb C}P^{1}\times{\mathbb C}P^{1}}K_{s}\left(z_{1},z_{2}\right)\alpha\left(z_{1}\right)\bar{\beta}\left(z_{1}\right)
\notag\end{equation}
on the space $ \Gamma\left(\Omega^{1-s}\left({\mathbb C}P_{{\mathbb R}}^{1}\right)\right) $ if $ s<\frac{1}{2} $, this pairing depends on $ s $
analytically \cite{GelShil58Gen}, and may be continued to an arbitrary $ s $.

\begin{proposition} $ \left(,\right)_{s} $ is a positive-definite Hermitian form on
$ \Gamma\left(\Omega^{1-s}\left({\mathbb C}P_{{\mathbb R}}^{1}\right)\right) $, if $ \left|s-\frac{1}{2}\right|<\frac{1}{2} $. It is equivalent to the $ H^{-1+2s} $-norm
on this space. \end{proposition}

\begin{proof} For $ s=\frac{1}{2} $ we get a standard pairing on $ \Omega^{1/2} $ (since $ K_{1/2} $ is a
$ \delta $-function), so it is sufficient to consider $ s\not=\frac{1}{2} $. Take an affine
coordinate system $ z $ on $ {\mathbb C}P^{1}\smallsetminus\left\{\infty\right\} $. A smooth section of $ \Omega^{t} $ near $ \infty $ is
represented by a smooth $ t $-form $ f\left(z\right)dZ^{t} $ with an asymptotic
\begin{equation}
|z|^{-4t}g\left(1/z\right)dZ^{t},
\notag\end{equation}
near $ \infty $. Here $ g\left(w\right) $ is a smooth function of $ x $ and $ y $, $ w=x+yi. $ Thus the
Fourier transform $ \widehat{f}\left(\zeta\right) $ of $ f $ is rapidly decreasing, smooth outside of $ \zeta=0 $,
and has an asymptotic expansion
\begin{equation}
h\left(\zeta\right)+|\zeta|^{2-4t}h_{1}\left(\zeta\right)+O\left(|\zeta|^{N}\right)
\notag\end{equation}
near 0, here $ N>0 $ is an arbitrary integer, $ h $ and $ h_{1} $ are smooth functions,
and $ t\not=\frac{1}{2} $.

The pairing $ \left(,\right)_{s} $ written in terms of $ \widehat{\alpha} $, $ \widehat{\beta} $ is the pairing
\begin{equation}
\operatorname{const}\cdot\int\widehat{\alpha}\left(\zeta\right)\bar{\widehat{\beta}}\left(\zeta\right)|\zeta|^{4s-2}d\zeta d\bar{\zeta}
\label{equ3.74}\end{equation}\myLabel{equ3.74,}\relax 
(in the sense of generalized functions) and we see
that for $ t=1-s $ and $ \left|s-\frac{1}{2}\right|<\frac{1}{2} $ the integral converges, thus the
pairing is positive.

On the other hand, the local equivalence of the norm~\eqref{equ3.74} with
$ H^{2s-1} $ is already proved in the previous section (cf. Equation~\eqref{equ3.7}). \end{proof}

\begin{remark} Note that what we got is a complementary series of unitary
representations of $ \operatorname{SL}\left(2,{\mathbb C}\right) $. \end{remark}

\begin{remark} \label{rem2.60}\myLabel{rem2.60}\relax  The defined pairings are compatible with duality $ \int\alpha\beta $
between $ \Omega^{s} $ and $ \Omega^{1-s} $. Note also that in the case $ s=0 $ the pairing has $ \operatorname{rk}=1 $,
dually, in the case $ s=1 $ the pairing becomes $ -\left(\Delta\alpha,\beta\right) $, which has constants
in the null-space. \end{remark}

Put $ s=\frac{1}{4} $. We get a $ \operatorname{PGL}\left(2,{\mathbb C}\right) $-invariant implementation\footnote{I.e., a Hilbert structure which is equivalent to the Sobolev Hilbert
structure which is in turn defined up to equivalence (compare with
Section~\ref{s2.20}).} of the space
$ H^{-1/2}\left(\Omega^{3/4}\right) $. Note now that $ \Omega^{3/4}=\omega^{3/4}\otimes\bar{\omega}^{3/4} $ differs not so much from the
space $ \omega^{1/2}\otimes\bar{\omega} $ we are most interested in. The difference $ \omega^{1/4}\otimes\bar{\omega}^{-1/4} $ is a
bundle with transitions functions $ \left(\frac{D}{\bar{D}}\right)^{1/4} $ of magnitude 1 (here
$ D=dz/dw $). In particular, a choice of affine coordinate $ z $ on $ {\mathbb C}P^{1} $ gives an
identification of $ \omega^{1/2}\otimes\bar{\omega} $ with $ \Omega^{3/4} $ (outside of infinity). This
identification commutes with natural actions of $ \operatorname{Aff}\left({\mathbb C}\right) $ on $ \Omega^{3/4} $ and
$ \omega^{1/2}\otimes\bar{\omega} $ (up to multiplication by a constant of magnitude 1). Here $ \operatorname{Aff}\left({\mathbb C}\right) $
is the group of affine transformations on $ {\mathbb C} $ (in fact to get an action on
$ \omega^{1/2}\otimes\bar{\omega} $ one needs to consider $ 2 $-covering of $ \operatorname{Aff}\left({\mathbb C}\right) $).

This explains the almost-invariance of $ H^{-1/2}\left({\mathbb C},\omega^{1/2}\otimes\bar{\omega}\right) $ with respect
to affine transformations, discussed in Section~\ref{s2.50}. Indeed, consider
a disk $ K\subset{\mathbb C}P^{1} $. On $ {\mathbb C}P^{1}\smallsetminus K $ multiplication by $ dz^{1/4}d\bar{z}^{-1/4} $ gives an isomorphism
of the sheaves $ \omega^{1/2}\otimes\bar{\omega} $ and $ \Omega^{3/4} $. Taking $ \infty $ as a center of $ K $, we obtain the
results of Section~\ref{s2.50} in the case $ n=1 $.

In what follows the following result is sufficient for us to show
conformal invariance of the objects we introduce:

\begin{amplification} \label{amp2.65}\myLabel{amp2.65}\relax  Fix a metric on $ {\mathbb C}P^{1} $ and $ \varepsilon>0 $. Consider $ \varphi\in\operatorname{SL}\left(2,{\mathbb C}\right) $, let
$ A\subset{\mathbb C}P^{1} $, $ B=\varphi\left(A\right) $. If both $ {\mathbb C}P^{1}\smallsetminus A $ and $ {\mathbb C}P^{1}\smallsetminus B $ contain disks of radius $ \varepsilon $ with
centers at $ c_{A} $, $ c_{B} $, then the norm of $ \varphi^{*}\colon \overset{\,\,{}_\circ}{H}^{-1/2}\left(B,\omega^{1/2}\otimes\bar{\omega}\right) \to
\overset{\,\,{}_\circ}{H}^{-1/2}\left(A,\omega^{1/2}\otimes\bar{\omega}\right) $ is bounded by
\begin{equation}
C\left(\varepsilon\right)\left(\frac{\operatorname{diam}\left(A\right)}{\operatorname{dist}\left(\varphi^{-1}\left(c_{B}\right),A\right)}\right)^{1/2}.
\notag\end{equation}
By duality, the same bound is valid for action of $ \varphi^{*} $ in $ H^{1/2}\left(\bullet,\omega^{1/2}\right) $.

\end{amplification}

\begin{proof} Indeed, $ SU\left(2\right) $ is compact, thus acts on $ H^{1/2}\left({\mathbb C}P^{1},\omega^{1/2}\otimes\bar{\omega}\right) $ by
uniformly bounded operators. Thus it is sufficient to consider the case when
a disk $ K $ of radius $ \varepsilon $ is fixed (say, if $ c_{A}=c_{B}=\infty $, thus $ |z| > 1/\varepsilon $), and
$ A\cap K=B\cap K=\varnothing $.

Let $ \alpha $, $ \beta $ be sections of $ \omega^{1/4}\otimes\bar{\omega}^{-1/4} $ and $ \omega^{-1/4}\otimes\bar{\omega}^{1/4} $ which coincide
with $ \bar{d}z^{1/4}d\bar{z}^{-1/4} $ and $ dz^{-1/4}d\bar{z}^{1/4} $ inside $ {\mathbb C}P^{1}\smallsetminus K $ correspondingly, and have a
compact support. Then $ \varphi^{*}|_{A}f $, $ f\in\omega^{1/2}\otimes\bar{\omega} $, is equal to
\begin{equation}
\mu_{\varphi}\beta\varphi_{3/4}^{*}|_{A}\left(\alpha f\right),\qquad \mu_{\varphi}=\left(\frac{\varphi'}{|\varphi'|}\right)^{1/2}
\notag\end{equation}
(here $ \varphi_{3/4}^{*} $ acts on $ \Omega^{3/4} $). On the other hand, multiplications by $ \alpha $ and $ \beta $
are bounded operators, and $ \varphi_{3/4}^{*} $ is unitary, so what remains to prove is
a bound on operator of multiplication by $ \mu_{\varphi} $ in $ H^{-1/2} $ or in $ H^{1/2} $. It is
sufficient to consider $ H^{1/2} $.

Let us apply interpolation theorem now: the norm of multiplication
by $ \mu_{\varphi} $ in $ H^{1/2} $ is no more than geometric mean of norms in $ H^{0}=L_{2} $ and in $ H^{1} $.
Since $ |\mu_{\varphi}|=1 $, it is unitary in $ H^{0} $, so it is sufficient to estimate the norm
in $ H^{1} $, which is bounded by $ \operatorname{const}\cdot\max  |\mu_{\varphi}'| $, i.e., by
$ \operatorname{const}\cdot\operatorname{dist}\left(\varphi^{-1}\left(\infty\right),A\right)^{-1} $.

To get the estimate in the theorem, note that we have a freedom of
rescaling $ A $, thus may assume that $ \operatorname{diam}\left(A\right)=1 $. \end{proof}

\subsection{Hardy space }\label{s2.70}\myLabel{s2.70}\relax  Let $ K = \left\{z \mid |z|\leq1\right\} $. The Hardy space $ {\mathcal H} $ is the
subspace of $ L_{2}\left(\partial K\right) $ consisting of functions with Fourier coefficients $ \left(a_{n}\right) $
which vanish for $ n<0 $. For a function $ f\in{\mathcal H} $, $ f=\sum a_{k}z^{k} $, $ |z|=1 $, let $ c_{f}=\sum a_{k}z^{k} $ be
defined for $ |z|<1 $. The latter series converges, and defines a holomorphic
function inside $ K $.

\begin{lemma} The mapping $ f \mapsto c_{f} $ is an injection $ L_{2}\left(\partial K\right) \to H^{1/2}\left(K\right) $. The
image of this injection coincides with
\begin{equation}
\operatorname{Ker}\left(H^{1/2}\left(K\right) \xrightarrow[]{\bar{\partial}} H^{-1/2}\left(K\right)\right).
\notag\end{equation}
\end{lemma}

\begin{proof} It is sufficient to prove that $ \|z^{k}\|_{H^{1/2}\left(K\right)} $ is bounded from
above and from below when $ k\in{\mathbb N} $. Consider a concentric disk $ K_{1} $ of radius
$ R<1 $. Since $ \|z^{k}\|_{H^{1/2}\left(K_{1}\right)}\leq\|z^{k}\|_{H^{1}\left(K_{1}\right)}=O\left(kR^{k}\right) $, it is $ o\left(1\right) $, thus it is
sufficient to consider the norm of $ z^{k} $ in a narrow annulus with external
boundary $ \partial K $.

In turn, taking coordinate $ \log  z=a+ib $, $ b\in{\mathbb R}/2\pi{\mathbb Z} $, it is sufficient to
consider $ H^{1/2} $-norm of $ e^{ka} $ on the half-line $ a\leq0 $, more precise, it is
sufficient to consider $ L_{2} $-norm of $ \left(\frac{d}{da}+k\right)^{1/2}e^{ka} $, $ a\leq0 $. Indeed, we need to
show that there is a continuation of $ e^{k\left(a+ib\right)} $ outside of $ a\leq0 $ such that
the $ H^{1/2} $-norm of this continuation is bounded, and that the $ H^{1/2} $-norm of
any such continuation is bounded from below. One can suppose that the
continuation is $ e^{ik b}\varphi\left(a\right) $, so one needs to estimate $ \int\left(1+k^{2}+\alpha^{2}\right)|\widehat{\varphi}\left(\alpha\right)|^{2}d\alpha $.
The operator $ \left(\frac{d}{da}+k\right)^{1/2} $ maps this norm to a norm equivalent to
$ L_{2} $-norm, moreover, it sends functions with support in $ a\geq0 $ into itself
(since $ k\geq0 $). Thus it defines an invertible operator $ H^{1/2}\left({\mathbb R}_{\leq0}\right) \to
H^{0}\left({\mathbb R}_{\leq0}\right)=L_{2}\left({\mathbb R}_{\leq0}\right) $.

This means that $ \left(\frac{d}{da}+k\right)^{1/2}e^{ka}|_{a\leq0}=\left(\frac{d}{da}+k\right)^{1/2}\varphi\left(a\right)|_{a\leq0} $, thus we
should not care about the choice of continuation $ \varphi $. Moreover, since
$ \left(\frac{d}{da}+k\right)^{1/2} $ does not move support to the left,
$ \left(\frac{d}{da}+k\right)^{1/2}e^{ka}=C\left(k\right)e^{ka} $, $ a\leq0 $. Since $ \left(\frac{d}{da}+k\right)^{1/2} $ is an operator of
convolution with $ \frac{a_{+}^{-3/2}e^{-ka}}{\Gamma\left(-1/2\right)} $ (in the sense of generalized
functions, $ x_{+}^{s} $ is 0 if $ x<0 $, $ x^{s} $ if $ x>0 $ and $ s>0 $, and is defined by analytic
continuation from the region $ s>0 $, where the above description defines it
completely \cite{GelShil58Gen}), to estimate $ C\left(k\right) $ one needs to calculate
\begin{align} -\int_{0}^{\infty}x_{+}^{-3/2}e^{-kx}e^{-kx}dx & \buildrel{\text{def}}\over{=} -\int_{0}^{\infty}x^{-3/2}\left(e^{-2kx}-1\right)dx
\notag\\
& = k^{1/2}\int_{0}^{\infty}\left(kx\right)^{-1/2}\frac{1-e^{-2kx}}{kx}d\left(kx\right) = k^{1/2}C.
\notag\end{align}
Thus we need to show that $ k^{1/2}\|e^{ka}|_{a\leq0}\|_{L_{2}} $ is bounded from above and from
below, which is obvious. \end{proof}

\begin{remark} It is clear that if $ f\in{\mathcal H}\cap C^{\infty} $, then $ f=c_{f}|_{\partial K} $. Moreover, if $ f\cdot\delta_{\partial K} $ is
the continuation of $ f $ to $ {\mathbb C} $ by $ \delta $-function, then $ c_{f}=\operatorname{const}\cdot\bar{\partial}^{-1}\left(f\cdot\delta_{\partial K}\right) $ (see
Section~\ref{s3.30} for description of $ \bar{\partial}^{-1} $). Several following sections are
dedicated to defining something similar in the case $ f\in{\mathcal H} $. \end{remark}

\begin{remark} \label{rem2.10}\myLabel{rem2.10}\relax  In Section~\ref{h4} we will construct a generalization of the
following modification of this result: instead of $ K $ we consider an
isomorphic domain $ {\mathbb C}P^{1}\smallsetminus K $. Instead of functions on $ {\mathbb C}P^{1} $ and $ \partial K $ we consider
sections of $ \omega^{1/2} $ and $ \Omega^{1/2}\otimes\mu $ correspondingly. \end{remark}

\section{Cauchy kernel for $ \omega^{1/2} $ }

\subsection{$ \bar{\partial} $ on piecewise-analytic functions }\label{s3.05}\myLabel{s3.05}\relax  Let $ M $ be a complex curve
with a
linear bundle $ {\mathcal L} $. Consider a domain $ D\subset M $ with a smooth boundary $ \partial D=\gamma $. Let
$ f\left(z\right) $, $ z\in\bar{D} $, be an analytic section of $ {\mathcal L} $ in $ D $ which continues as a smooth
section to $ \bar{D} $. Extend $ f $ outside of $ D $ as 0.

Consider $ \bar{\partial}f $. It is a section of $ {\mathcal L}\otimes\bar{\omega} $, which vanishes on $ M\smallsetminus\gamma $. Since $ f $
has a jump of the first kind along $ \gamma $, and $ \bar{\partial} $ is a differential operator of
the first order, $ \bar{\partial}f $ has at most a $ \delta $-function singularity along $ \gamma $. One may
locally write
\begin{equation}
\bar{\partial}f = g\left(z\right)\cdot\delta_{h\left(z\right)=0}.
\notag\end{equation}
Here $ g\left(z\right) $ is a smooth section of $ {\mathcal L}\otimes\bar{\omega} $, $ h $ is a local equation of $ \gamma $.

Let us calculate $ g|_{\gamma} $ in terms of $ f $. First of all, if one changes the
equation $ h=0 $ by a different equation $ hh_{0}=0 $ (here $ h_{0} $ has no zeros close to
the point in question), the coefficient at $ g\left(z\right) $ changes to
\begin{equation}
\delta_{hh_{0}=0}=\frac{1}{h_{0}}\delta_{h=0},
\notag\end{equation}
since $ \delta $-function has homogeneity degree $ -1 $. Thus $ g\left(z\right)|_{\gamma} $ should be
multiplied by $ h_{0}|_{\gamma} $ to preserve the same value of the product. This shows
that to get an invariant description of $ g\left(z\right)|_{\gamma} $ one needs to write
\begin{equation}
g\left(z\right)|_{\gamma} = G\left(z\right)\otimes dh|_{\gamma}.
\notag\end{equation}
Here $ dh|_{\gamma} $ is considered as a section of the conormal bundle $ N^{*}\gamma $ to $ \gamma $,
thus $ G\left(z\right) $ should be a section of $ \left({\mathcal L}\otimes\bar{\omega}\right)|_{\gamma}\otimes\left(N^{*}\gamma\right)^{*} = \left({\mathcal L}\otimes\bar{\omega}\right)|_{\gamma}\otimes N\gamma $.

Now $ G\left(z\right) $ does not depend on the parameterization of a neighborhood
of $ \gamma $, so it should be directly expressible in terms of $ f\left(z\right) $.

\begin{lemma} $ \bar{\omega}|_{\gamma} $ is canonically isomorphic to $ \left(N^{*}\gamma\right)\otimes{\mathbb C} $. \end{lemma}

\begin{proof} Indeed, $ \bar{\omega} $ is defined as a quotient of $ \left(T^{*}M_{{\mathbb R}}\right)\otimes{\mathbb C} $ by holomorphic
forms, moreover, $ \bar{\omega}_{{\mathbb R}} $ is an isomorphic
image of $ T^{*}M_{{\mathbb R}} $. Since $ N^{*}\gamma\subset T^{*}M_{{\mathbb R}} $, there is a natural non-zero mapping
$ N^{*}\gamma \to \bar{\omega}_{{\mathbb R}} $, which gives the required isomorphism after complexification. \end{proof}

As a corollary, $ \left({\mathcal L}\otimes\bar{\omega}\right)|_{\gamma}\otimes N\gamma \simeq {\mathcal L}|_{\gamma} $. Now the following fact is obvious
from calculations in local coordinates:

\begin{proposition} $ G\left(z\right) = -f\left(z\right)|_{\gamma} $. \end{proposition}

In particular, we see that one can write
\begin{equation}
\bar{\partial}f=-e\circ r\left(f\right)
\label{equ3.02}\end{equation}\myLabel{equ3.02,}\relax 
in terms of operators $ r $ of restriction and $ e $ of extension-by-$ \delta $-function.
Thus we can consider $ -e\circ r $ as an ``approximation'' to $ \bar{\partial} $ which is good on
functions which are analytic far from $ \gamma $. This approximation is exact on
functions with a jump of the first kind.

Yet another way to treat this identity is to write it as
\begin{equation}
-\bar{\partial}^{-1}\circ e\circ r\left(F\right) = \vartheta_{D}F.
\notag\end{equation}
Here $ F $ is a function which is holomorphic in a neighbourhood of $ \bar{D} $, $ \vartheta_{D} $ is
a function which is 1 on $ D $ and 0 outside of $ D $. We give a mollified
version of this statement in Section~\ref{s4.35}.

\subsection{Self-duality of $ \bar{\partial} $ } Consider a mapping $ \bar{\partial}\colon \omega^{1/2} \to \omega^{1/2}\otimes\bar{\omega} $. Note that the
spaces of sections of these bundles are dual w.r.t. the pairing
\begin{equation}
\left(\varphi,\psi\right) = \int_{M}\varphi\psi.
\label{equ3.20}\end{equation}\myLabel{equ3.20,}\relax 
\begin{lemma} $ \bar{\partial} $ is skew-symmetric w.r.t. the above pairing,
\begin{equation}
\left(\varphi_{1},\bar{\partial}\varphi_{2}\right)+\left(\varphi_{2},\bar{\partial}\varphi_{1}\right)=0.
\notag\end{equation}
\end{lemma}

\begin{proof} It is sufficient to show that $ \left(\varphi_{1},\bar{\partial}\varphi_{1}\right)=0 $, i.e., to study
\begin{equation}
\int_{M}\varphi\bar{\partial}\varphi = \frac{1}{2}\int_{M}\bar{\partial}\varphi^{2} = \frac{1}{2}\int_{M}d\varphi^{2}
\notag\end{equation}
which obviously vanishes. \end{proof}

Note that this supports the heuristic that the bundle $ \omega^{1/2} $ is ``the
best one'' of the powers of $ \omega $.

Used literally, the above considerations were applicable to the
space of smooth sections of the bundles in question. On the other hand,
the pairing~\eqref{equ3.20} extends to a continuous pairing between $ H^{1/2}\left(M,\omega^{1/2}\right) $
and $ H^{-1/2}\left(M,\omega^{1/2}\otimes\bar{\omega}\right) $, and the operator $ \bar{\partial} $ maps one of these spaces to
another. Using the facts from Section~\ref{s2.20} We get a

\begin{proposition} The operator $ \bar{\partial} $ gives a canonically defined
Fredholm symplectic form on $ H^{1/2}\left(M,\omega^{1/2}\right) $. This form is non-degenerate if
$ M={\mathbb C}P^{1} $. \end{proposition}

Here we call a bilinear form $ \alpha $ on a Hilbert space $ {\mathcal H} $ a {\em Fredholm form},
if the corresponding operator $ \alpha\colon H \to H^{*}=H $ is Fredholm.

Recall that in the case $ M={\mathbb C}P^{1} $ the space $ H^{1/2}\left(M,\omega^{1/2}\right) $ is the ``best''
of Sobolev spaces for $ \omega^{1/2} $, since it allows the action of $ \operatorname{SL}\left(2,{\mathbb C}\right) $ (at
least ``in small'') by bounded operators. Here we see another ``nice''
property of this particular Sobolev space: there is a canonically defined
invertible pairing on this space (invertible in the sense that the
corresponding operator is bounded).

\begin{remark} If we could construct a self-adjoint operator which enjoys the
above properties of the operator $ \bar{\partial} $, i.e., is $ \operatorname{SL}\left(2,{\mathbb C}\right) $-invariant and
Fredholm, one would be (almost) able
to define an invariant Sobolev space structure on $ H^{1/2}\left(M,\omega^{1/2}\right) $. The only
thing missing would be positive definiteness. \end{remark}

\subsection{The kernel of $ \bar{\partial}^{-1} $ }\label{s3.30}\myLabel{s3.30}\relax  Suppose $ M={\mathbb C}P^{1} $. In this case the operator $ \bar{\partial}:
{\mathcal D}'\left(\omega^{1/2}\right) \to {\mathcal D}'\left(\omega^{1/2}\otimes\bar{\omega}\right) $
is invertible, thus the operator $ \bar{\partial}^{-1} $ is canonically defined. Consider a
kernel $ K\left(x,y\right) $ of this operator,
\begin{equation}
\left(\bar{\partial}^{-1}f\right)\left(x\right) = \int K\left(x,y\right)f\left(y\right).
\notag\end{equation}
Obviously, $ K\left(x,y\right) $ is a section of $ \omega^{1/2} \boxtimes \omega^{1/2} $.

To calculate $ K\left(x,y\right) $ consider a $ \delta $-section $ \delta_{z_{0}} $ of $ \omega^{1/2} $ at $ z_{0} $ (defined
up to a scalar multiple). Then $ K\left(\bullet,z_{0}\right)=\operatorname{const}\cdot\bar{\partial}^{-1}\delta_{z_{0}} $ is a holomorphic
section of $ \omega^{1/2} $ on $ {\mathbb C}P^{1}\smallsetminus\left\{z_{0}\right\} $. It is easy to see that it has a simple pole
at $ z=z_{0} $, thus has no zeros (since $ \omega^{1/2}\simeq{\mathcal O}\left(-1\right) $). We conclude that $ K\left(x,y\right) $
has a simple pole at $ x=y $. There is only one such a section of $ \omega^{1/2} \boxtimes
\omega^{1/2} $. Now one may easily recognize the kernel for $ \bar{\partial}^{-1} $ in the Cauchy
formula
\begin{equation}
\bar{\partial}^{-1}\delta_{z=z_{0}}dz^{1/2}d\bar{z} = \frac{1}{2\pi i} \frac{1}{z-z_{0}}dz^{1/2}
\notag\end{equation}
The right-hand side $ f\left(z,z_{0}\right) $ has no pole at $ \infty $, and is uniquely determined
by this condition and the equation
\begin{equation}
\bar{\partial}f=\delta_{z=z_{0}}dz^{1/2}d\bar{z}
\notag\end{equation}
in $ {\mathbb C}P^{1}\smallsetminus\left\{\infty\right\} $. We conclude that
\begin{equation}
K\left(x,y\right) = \frac{1}{2\pi i} \frac{1}{x-y}dx^{1/2}dy^{1/2}.
\notag\end{equation}
\subsection{Cauchy kernel and $ L_{2} $ } Consider a smooth embedded curve $ \gamma\hookrightarrow{\mathbb C}P^{1} $ (we do
not suppose that $ \gamma $ is compact). Any smooth function $ f $ on $ \gamma $ with a
compact support may be extended as a generalized $ \delta $-function $ e\left(f\right) $ on $ {\mathbb C}P^{1} $
with a support on $ \gamma $, same for half-forms on $ \gamma $. Applying $ \bar{\partial}^{-1} $ to the
result, we obtain a form on $ {\mathbb C}P^{1} $ which is holomorphic outside of $ \gamma $.
Here we discuss when this mapping $ \bar{\partial}^{-1}\circ e $ allows $ f\in L_{2}\left(\gamma\right) $ instead of
$ f\in{\mathcal D}\left(\gamma\right) $.

\begin{remark} In what follows the main example of $ \gamma $ is a disjoint union of
infinite number of circles.

\end{remark}

Fix a point $ z\in{\mathbb C}P^{1} $ and suppose that some neighborhood of $ z $ does not
intersect $ \gamma $. Then for a half-form $ f $
\begin{align} \frac{\bar{\partial}^{-1}\circ e\left(f\right)}{dz^{1/2}d\bar{z}}|_{z_{0}} & \buildrel{\text{def}}\over{=} \left< \bar{\partial}^{-1}\circ e\left(f\right), \delta_{z_{0}}dz^{1/2}\right>
\notag\\
& = -\left< e\left(f\right), \bar{\partial}^{-1}\left(\delta_{z_{0}}dz^{1/2}\right)\right>
\notag\\
& =-\left< e\left(f\right), K\left(\bullet,z_{0}\right)dz^{-1/2} \right> = \int_{C}f\left(x\right)K\left(x,z_{0}\right)dz^{-1/2}.
\notag\end{align}
We conclude that the linear functional of calculating $ \bar{\partial}^{-1}\circ e\left(f\right) $ at $ z_{0} $ is
given by the kernel $ K\left(\bullet,z_{0}\right)|_{\bullet\in C} $. However, since $ K\left(\bullet,z_{0}\right) $ is bounded
outside of a neighborhood of $ z_{0} $, we see that

\begin{proposition} \label{prop4.12}\myLabel{prop4.12}\relax  The mapping $ \bar{\partial}^{-1}\circ e $ extends to a mapping from $ L_{2}\left(\gamma\right) $ to
the space $ \operatorname{Hol}\left({\mathbb C}P^{1}\smallsetminus\bar{\gamma}\right) $ of forms holomorphic outside of $ \gamma $ iff the length of
$ \gamma $ is finite. \end{proposition}

\section{Toy theory }\label{h35}\myLabel{h35}\relax 

In this section we study an example which shows the main framework
of our approach without any of the complications related to the necessity
to consider $ H^{1/2} $-spaces. This is an example related to dealing with
{\em partial period mapping\/} $ \Gamma\left(M,\omega\right) \to {\mathbb C}^{g} $ of taking periods of global
holomorphic forms along $ A $-cycles on the Riemann surface. In fact, we
start from studying a related space of global holomorphic functions which
are allowed to have constant jumps along $ A $-cycles.

\subsection{Toy global space } Consider a complex curve $ M $ and the mapping $ \partial\bar{\partial}\colon \Omega^{0} \to
\Omega^{1} $ (let us recall that $ \Omega^{\text{top}}=\Omega^{1} $ in our notations for fractional forms).
The integral $ -\int_{M}if\cdot\partial\bar{\partial} \bar{f} $ defines a sesquilinear form on functions with
compact support, it is Hermitian, and the only functions which are in the
null-space of this form are constants, since the integral is equal to
$ \int_{M}i\partial f\cdot\bar{\partial} \bar{f}\geq0 $.

This integral defines a pre-Hilbert structure on global sections of
$ \Omega^{0} $ modulo constants, this structure is compatible with the
Sobolev $ H^{1} $-topology.
Moreover, if $ M={\mathbb C}P^{1} $, this structure is $ \operatorname{PGL}\left(2,{\mathbb C}\right) $-invariant, and coincides
with the structure of supplementary series of representations of $ \operatorname{SL}\left(2,{\mathbb C}\right) $
(in notations of Section~\ref{s2.60} it is the case $ s=1 $, compare with Remark
~\ref{rem2.60}). As a result, we obtain a $ \operatorname{PGL}\left(2,{\mathbb C}\right) $-invariant realization of
Hilbert space $ H^{1}\left({\mathbb C}P^{1}\right)/\operatorname{const} $.

The dual space is the subspace $ H_{\int=0}^{-1}\left({\mathbb C}P^{1},\Omega^{1}\right) $ of $ H^{-1}\left({\mathbb C}P^{1},\Omega^{1}\right) $
consisting of forms with integral 0, thus the latter space also carries
an $ \operatorname{PGL}\left(2,{\mathbb C}\right) $-invariant Hilbert structure. To describe it, note that Remark
~\ref{rem2.60} is still applicable, thus in the limit $ s \to 0 $ the pre-Hilbert
structure of Section~\ref{s2.60} becomes a form of $ \operatorname{rk}=1 $ with $ H_{\int=0}^{-1}\left({\mathbb C}P^{1},\Omega^{1}\right) $ in
the null-space. As a corollary, $ \frac{d}{ds}|_{s=0}\left(\alpha,\beta\right)_{s} $ gives a correctly
defined positive form on $ H_{\int=0}^{-1}\left({\mathbb C}P^{1},\Omega^{1}\right) $ which should be dual to the
pairing on $ H^{1}\left({\mathbb C}P^{1}\right)/\operatorname{const} $.

Taking the $ s $-derivative of the kernel $ K_{s} $ from Section~\ref{s2.60}, we get
$ \log |z_{1}-z_{2}| $ (when restricted on $ {\mathbb C}\subset{\mathbb C}P^{1} $), which is indeed scaling-invariant
up to addition of a constant. (Note that the constant is irrelevant,
since the forms $ \alpha $ and $ \beta $ we are going to pair $ K $ with have vanishing
integral.) What the above argument shows is that the continuation of this
kernel to $ {\mathbb C}P^{1} $ is invariant w.r.t. fraction-linear mappings.

\subsection{Almost-perpendicularity } The description of the pairing on
$ H_{\int=0}^{-1}\left({\mathbb C}P^{1},\Omega^{1}\right) $ immediately implies

\begin{proposition} Consider two disjoint disks $ K_{1} $, $ K_{2} $ on $ {\mathbb C}P^{1} $ of conformal
distance $ l $. Let $ H_{1,2} = \overset{\,\,{}_\circ}{H}_{\int=0}^{-1}\left(K_{1,2},\Omega^{1}\right)\subset H_{\int=0}^{-1}\left({\mathbb C}P^{1},\Omega^{1}\right) $ . Let $ P_{l} $ be the
orthogonal projector from one subspace to another. Then
\begin{equation}
\|P_{l}\| \sim e^{-l}.
\notag\end{equation}
The equivalence means that the quotient of two sides remains bounded and
separated from 0 when $ l $ varies. \end{proposition}

\begin{proof} Since the natural norm on $ H_{\int=0}^{-1}\left(\Omega^{1}\right) $ is invariant with respect
to $ \operatorname{PGL} $, the angle between the subspaces $ H_{1,2} $ depends on $ l $ only. Note that
it is sufficient to prove the statement in the case $ l>\varepsilon $ for some fixed $ \varepsilon>0 $.
Now proceed as in the proof of Proposition~\ref{prop3.170}. The only change is
that the kernel is now $ \log |L+z_{1}-z_{2}| $, and on forms with integral 0 it is
equivalent to $ \left(\log |L+z_{1}-z_{2}|-\log  L-\operatorname{Re}\frac{z_{1}-z_{2}}{L}\right) $, which is $ L^{-2} $ times an
operator of rank 1, plus much smaller operator. \end{proof}

\begin{remark} Note that norm of this projector is much smaller than the
norm of the projector from Proposition~\ref{prop3.170} (if $ l $ is big enough). \end{remark}

\begin{corollary} \label{cor35.20}\myLabel{cor35.20}\relax  Consider a family of disjoint disks $ \left\{K_{i}\right\} $ in $ {\mathbb C}P^{1} $ with
pairwise conformal distances $ l_{ij} $, $ i\not=j $. Let $ l_{i i}=0 $. If the matrix $ \left(e^{-l_{ij}}\right) $
gives a bounded operator $ l_{2} \to l_{2} $, then the natural extension-by-0 mapping
\begin{equation}
\bigoplus_{l_{2}}\overset{\,\,{}_\circ}{H}_{\int=0}^{-1}\left(K_{i},\Omega^{1}\right) \to H_{\int=0}^{-1}\left({\mathbb C}P^{1},\Omega^{1}\right)
\notag\end{equation}
is a continuous injection, and, dually, the natural restriction mapping
\begin{equation}
H^{1}\left({\mathbb C}P^{1}\right)/\operatorname{const} \to \bigoplus_{l_{2}}H^{1}\left(K_{i}\right)/\operatorname{const}
\notag\end{equation}
is a continuous surjection. \end{corollary}

\subsection{Topology on the boundary }\label{s35.20}\myLabel{s35.20}\relax  Consider a disk $ K\subset{\mathbb C}P^{1} $. Given a
function $ f $ on $ \partial K $, one can consider its decomposition into a sum $ f_{+}+f_{-} $ of
functions which can be analytically extended into/outside of $ K $. Such a
decomposition
exists if $ f $ is in $ L_{2} $, and the summands are uniquely defined up to
addition of a constant.

If $ f\in H^{1/2}\left(\partial K\right) $, then $ f_{+} $ and $ f_{-} $ are $ H^{1} $-functions on $ K $ and on $ {\mathbb C}P^{1}\smallsetminus K $
(compare with Section~\ref{s2.70}), if $ f\in H^{1/2}\left(\partial K\right)/\operatorname{const} $, then $ f_{+}\in H^{1}\left(K\right)/\operatorname{const} $,
$ f_{-}\in H^{1}\left({\mathbb C}P^{1}\smallsetminus K\right)/\operatorname{const} $, and $ f_{\pm} $ are uniquely defined by $ f $. Thus $ H^{1} $-norms of $ f_{\pm} $
are correctly defined. In the same way as in Section~\ref{s2.70} one can prove

\begin{lemma} \label{lm35.30}\myLabel{lm35.30}\relax  The described above extension mapping
\begin{equation}
H^{1/2}\left(\partial K\right) \to H^{1}\left(K\right)/\operatorname{const} \oplus H^{1}\left({\mathbb C}P^{1}\smallsetminus K\right)/\operatorname{const}
\notag\end{equation}
is a continuous injection. \end{lemma}

Since the Hilbert structure on the right-hand side can be made
$ \operatorname{PGL}\left(2,{\mathbb C}\right) $-invariant, the extension mapping defines the $ \operatorname{PGL}\left(2,{\mathbb C}\right) $-invariant
realization of the Hilbert structure on $ H^{1/2}\left(\partial K\right) $. The property of
invariance can be described in the following way: let $ \varphi $ be a
fraction-linear mapping, $ \varphi\left(K\right)=K' $. Then $ \varphi^{*} $ defines a unitary operator
\begin{equation}
\varphi^{*}\colon H^{1/2}\left(\partial K'\right) \to H^{1/2}\left(\partial K\right).
\notag\end{equation}
\begin{remark} \label{rem35.35}\myLabel{rem35.35}\relax  In the same way as we defined a $ \operatorname{PGL}\left(2,{\mathbb C}\right) $-invariant
pairing on sections of $ \Omega^{1-s} $ on $ {\mathbb C}P^{1} $ (in Section~\ref{s2.60}), one can define a
$ \operatorname{PGL}\left(2,{\mathbb R}\right) $-invariant pairing on section of $ \Omega_{{\mathbb R}P^{1}}^{1-s} $ on $ {\mathbb R}P^{1}=S^{1} $. In this way
one gets so called supplementary series of representations of $ \operatorname{SL}\left(2,{\mathbb R}\right) $. By
expressing the pairing in appropriate coordinate systems on $ {\mathbb R}P^{1}=S^{1} $ one
can easily see that the described above $ \operatorname{SL}\left(2,{\mathbb R}\right) $-invariant Hilbert
structure coincides with one of these structures. \end{remark}

\begin{definition} The representation $ f=f_{+}+f_{-} $ decomposes $ H^{1/2}\left(\partial K\right)/\operatorname{const} $ into a
direct sum of two subspaces which we denote $ H_{+}^{1/2}\left(\partial K\right)/\operatorname{const} $ and
$ H_{-}^{1/2}\left(\partial K\right)/\operatorname{const} $. \end{definition}

Combining the above description of $ H^{1/2}\left(\partial K\right) $ with results of the
previous section, one obtains

\begin{corollary} In the conditions of Corollary~\ref{cor35.20} the natural
restriction mapping
\begin{equation}
H^{1}\left({\mathbb C}P^{1}\right)/\operatorname{const} \to \bigoplus_{l_{2}}H^{1/2}\left(K_{i}\right)/\operatorname{const}
\notag\end{equation}
is a continuous surjection. \end{corollary}

\subsection{Space of holomorphic functions }\label{s35.30}\myLabel{s35.30}\relax  Consider a family $ \left\{K_{i}\right\} $ of
disjoint disks in $ {\mathbb C}P^{1} $.

\begin{definition} We say that a generalized function $ f $ on $ {\mathbb C}P^{1} $ is
$ H^{1} $-{\em holomorphic in\/} $ {\mathbb C}P^{1}\smallsetminus\bigcup K_{i} $ if $ f\in H^{1}\left({\mathbb C}P^{1}\smallsetminus\bigcup K_{i}\right) $, and
$ \bar{\partial}f=0\in H^{0}\left({\mathbb C}P^{1}\smallsetminus\bigcup K_{i},\bar{\omega}\right)= L_{2}\left({\mathbb C}P^{1}\smallsetminus\bigcup K_{i},\bar{\omega}\right) $. Denote the
the space of $ H^{1} $-holomorphic functions in $ {\mathbb C}P^{1}\smallsetminus\bigcup K_{i} $ by $ {\mathcal H}^{\left(1\right)} $.
\end{definition}

Note that the Sobolev spaces in this definition are generalized
ones. Note also that the Hilbert norm on $ L_{2}\left({\mathbb C}P^{1}\smallsetminus\bigcup K_{i},\bar{\omega}\right) $ is canonically
defined by $ \|\alpha\|^{2}=-i\int\bar{\alpha}\alpha $.

\begin{theorem} \label{th35.15}\myLabel{th35.15}\relax  Suppose that the conformal distance between $ \partial K_{i} $ and $ \partial K_{j} $
is $ l_{ij} $, and the matrix $ \left(e^{-l_{ij}}\right) $ gives a bounded operator $ l_{2} \to l_{2} $.
Consider the mapping of taking the boundary value:
\begin{equation}
b\colon {\mathcal H}^{\left(1\right)}/\operatorname{const} \to \bigoplus_{l_{2}}H^{1/2}\left(\partial K_{i}\right)/\operatorname{const}.
\notag\end{equation}
Let $ b_{-} $ be the component of this mapping going into $ \bigoplus_{l_{2}}H_{-}^{1/2}\left(\partial K_{i}\right)/\operatorname{const} $.
Then $ b_{-} $ is a bounded invertible mapping. \end{theorem}

\begin{proof} We already know that $ b_{-} $ is bounded. To show that one can
reproduce a function by the component $ b_{-} $ of its restriction to a boundary
consider a function $ f\in{\mathcal H}^{\left(1\right)} $. By definition, it is a restriction of some
function $ g\in H^{1}\left({\mathbb C}P^{1}\right) $, and this function $ g $ is defined up to a function with
support in $ \bigcup K_{i} $. By the condition of $ f $ being holomorphic, $ \bar{\partial}g $ is an
element of $ L_{2}\left(\bigcup K_{i},\bar{\omega}\right)=\bigoplus_{l_{2}}L_{2}\left(K_{i},\bar{\omega}\right) $, and $ \bar{\partial}g $ is defined up to addition of
$ \bar{\partial}\overset{\,\,{}_\circ}{H}^{1}\left(\bigcup K_{i}\right) $. On the other hand, the topology on $ H^{1}\left(X\right) $ (here $ X $ is a manifold
with coordinates $ x_{k} $) can be defined by the norm $ \|f\|_{L_{2}}^{2}+\sum\|\partial_{x_{k}}f\|_{L_{2}}^{2} $, thus
$ \overset{\,\,{}_\circ}{H}^{1}\left(\bigcup K_{i}\right) = \bigoplus_{l_{2}}\overset{\,\,{}_\circ}{H}^{1}\left(K_{i}\right) $.

We conclude that $ \bar{\partial}g $ is a canonically defined element of
\begin{equation}
\bigoplus_{l_{2}}L_{2}\left(K_{i},\bar{\omega}\right)/\bar{\partial}\overset{\,\,{}_\circ}{H}^{1}\left(K_{i}\right).
\notag\end{equation}
\begin{lemma} Consider a disk $ K\subset{\mathbb C}P^{1} $. Then there exists a canonical
isomorphism
\begin{equation}
L_{2}\left(K,\bar{\omega}\right)/\bar{\partial}\overset{\,\,{}_\circ}{H}^{1}\left(K\right) \simeq H_{-}^{1/2}\left(\partial K\right)/\operatorname{const}.
\notag\end{equation}
\end{lemma}

\begin{proof} Let us start with a left-to-right mapping. Let
$ \alpha\in L_{2}\left(K,\bar{\omega}\right)/\bar{\partial}\overset{\,\,{}_\circ}{H}^{1}\left(K\right) $. Since $ \bar{\partial} $ is an elliptic operator without cokernel and
with $ 1 $-dimensional null-space, $ \bar{\partial}^{-1}\alpha $ is an $ H^{1} $-function defined up to
addition of a constant and addition of an element of $ \overset{\,\,{}_\circ}{H}^{1}\left(K\right) $. Thus the
restriction of $ f $ to $ \partial K $ is defined up to a constant, so it is an element
of $ H^{1/2}\left(\partial K\right)/\operatorname{const} $. Moreover, it is in $ H_{-}^{1/2}\left(\partial K\right)/\operatorname{const} $ since $ \bar{\partial}^{-1}\alpha $ is
holomorphic outside $ K $.

To get right-to-left mapping note that any element of
$ H_{-}^{1/2}\left(\partial K\right)/\operatorname{const} $ can be (by definition) continued to a holomorphic outside
of $ K $ function $ f $, and this function is of class $ H^{1} $. Thus to this
continuation one can apply the arguments which precede the lemma.
\end{proof}

To finish the proof of the theorem note that the knowledge of $ r_{-}\left(f\right) $
allows one to construct $ \bar{\partial}f $ (by lemma), so the application of $ \bar{\partial}^{-1} $
reconstructs $ f $.\end{proof}

\subsection{Toy Riemann--Roch theorem }\label{s35.40}\myLabel{s35.40}\relax  Consider a family of disjoint disks
$ K_{i} $, $ i\in I $, and the corresponding space of $ H^{1} $-holomorphic functions. Suppose
that $ I $ has an involution ' which interchanges two halves of $ I=I_{+}\amalg I_{+}' $. Fix
automorphisms $ \varphi_{i} $ of $ {\mathbb C}P^{1} $, $ i\in I $, such that $ \varphi_{i'}=\varphi_{i}^{-1} $, and $ \varphi_{i}\left(\partial K_{i}\right) $ is $ \partial K_{i'} $
with inverted orientation. Note that if the set $ I $ is finite, then after
gluing $ \partial K_{i} $ via $ \varphi_{i} $ one gets a curve of genus $ \operatorname{card}\left(I_{+}\right)=\operatorname{card}\left(I\right)/2 $.

Associate to a function $ f\in{\mathcal H}^{\left(1\right)} $ the jump of its boundary value after
such a gluing:
\begin{equation}
{\mathcal J}\colon f \mapsto \left(f|_{\partial K_{j}}-\varphi_{j}^{*}\left(f|_{\partial K_{j'}}\right)\right),\quad j\in I_{+}.
\label{equ4.101}\end{equation}\myLabel{equ4.101,}\relax 
\begin{theorem} \label{th35.45}\myLabel{th35.45}\relax  Suppose that the conformal distance between $ \partial K_{i} $ and $ \partial K_{j} $
is $ l_{ij} $, and the matrix $ \left(e^{-l_{ij}}-\delta_{ij}\right) $ gives a compact operator $ l_{2} \to l_{2} $.
Then the mapping
\begin{equation}
{\mathcal J}\colon {\mathcal H}^{\left(1\right)}/\operatorname{const} \to \bigoplus\Sb l_{2} \\ j\in I_{+}\endSb H^{1/2}\left(\partial K_{j}\right)/\operatorname{const}
\notag\end{equation}
is a Fredholm operator of index 0. \end{theorem}

\begin{proof} Consider the composition $ \widetilde{{\mathcal J}}={\mathcal J}\circ b_{-}^{-1} $. Let us show that $ \widetilde{{\mathcal J}} $ is
Fredholm of index 0, this would immediately imply the statement of the
theorem. The mapping $ \varphi_{j}^{*} $ interchanges $ \pm $-components of $ H^{1/2}\left(\partial K_{j}\right)/\operatorname{const} $ and
$ H^{1/2}\left(\partial K_{j'}\right)/\operatorname{const} $, thus one can identify
\begin{equation}
H^{1/2}\left(\partial K_{j}\right)/\operatorname{const} =H_{-}^{1/2}\left(\partial K_{j}\right)/\operatorname{const} \oplus H_{+}^{1/2}\left(\partial K_{j}\right)/\operatorname{const}
\notag\end{equation}
with $ H_{-}^{1/2}\left(\partial K_{j}\right)/\operatorname{const} \oplus H_{-}^{1/2}\left(\partial K_{j'}\right)/\operatorname{const} $. Denote the composition of $ \widetilde{{\mathcal J}} $ and
this identification by $ {\mathcal K}\colon \bigoplus_{l_{2}}H_{-}^{1/2}\left(\partial K_{j}\right)/\operatorname{const} \to \bigoplus_{l_{2}}H_{-}^{1/2}\left(\partial K_{j}\right)/\operatorname{const} $.
For $ f\in{\mathcal H}^{\left(1\right)} $ denote by $ f_{j\pm} $ the $ \pm $-components of $ f|_{\partial K_{j}} $. One can easily see
that $ {\mathcal K}\left(\left(f_{j-}\right)_{j\in I}\right)=\left(f_{j-}-\varphi_{j}^{*}\left(f_{j'+}\right)\right)_{j\in I} $, in other words, $ {\mathcal K}=\operatorname{id}-\varphi^{*}\circ b_{+}\circ b_{-}^{-1} $,
here $ b_{+}=b-b_{-} $. Since $ \varphi^{*}=\bigoplus_{l_{2}}\varphi_{j}^{*} $ is an isometry, to prove the theorem it is
enough to show that $ b_{+}\circ b_{-}^{-1} $ is compact. Now we investigate how to
reconstruct $ H_{+}^{1/2} $-components of restriction to a boundary via
$ H_{-}^{1/2} $-components.

Consider the operator with the Cauchy kernel $ \left(y-x\right)^{-1}dy $, $ x,y\in{\mathbb C}P^{1} $.
Restrict this kernel to $ \bigcup\partial K_{i} $. One gets an operator which sends functions
on $ \bigcup\partial K_{i} $ to functions on $ \bigcup\partial K_{i} $. Put zeros instead of the diagonal terms
(which send functions on $ \partial K_{i} $ to functions on $ \partial K_{i} $), and call the resulting
operator $ {\bold K} $. Obviously, the functions in the image can be holomorphically
extended to $ K_{i} $. Moreover, if a function was non-zero on $ \partial K_{i} $ only, and
could be analytically continued inside $ K_{i} $, then this function is in the
zero-space of $ {\bold K} $.

Thus the operator $ {\bold K} $ maps $ \bigoplus H_{-}^{1/2}\left(\partial K_{i}\right)/\operatorname{const} $ to $ \bigoplus H_{+}^{1/2}\left(\partial K_{i}\right)/\operatorname{const} $.
Moreover, as the proof of Theorem~\ref{th35.15} shows, $ b_{+}\circ b_{-}^{-1}= \frac{{\bold K}}{2\pi i} $, since
the pseudodifferential operator $ \bar{\partial}^{-1} $ has $ \frac{dy}{2\pi i\left(x-y\right)} $ as the null-space.

We conclude that the image $ V_{1}=b\left({\mathcal H}^{\left(1\right)}/\operatorname{const}\right) $ in $ \bigoplus H^{1/2}\left(\partial K_{i}\right)/\operatorname{const} $ is
a graph of operator $ \frac{{\bold K}}{2\pi i}\colon \bigoplus H_{-}^{1/2}\left(\partial K_{i}\right)/\operatorname{const} \to \bigoplus H_{+}^{1/2}\left(\partial K_{i}\right)/\operatorname{const} $. On
the other hand, consider the subspace $ V_{2} $ of $ \bigoplus H^{1/2}\left(\partial K_{i}\right)/\operatorname{const} $ consisting
of functions which are in the null-space of $ {\mathcal J} $, i.e., the subspace given
by the condition
\begin{equation}
f|_{\partial K_{j}}=\varphi_{j}^{*}\left(f|_{\partial K_{j'}}\right),\qquad j\in I_{+}.
\notag\end{equation}
Since $ \varphi_{j}^{*} $ interchanges $ H_{-}^{1/2}\left(\partial K_{j}\right)/\operatorname{const} $ and $ H_{+}^{1/2}\left(\partial K_{j'}\right)/\operatorname{const} $ and is
unitary, we conclude that the subspace $ V_{2} $ is a graph of a unitary
mapping $ \bigoplus\varphi_{i}^{*}|_{H_{+}^{1/2}\left(\partial K_{i}\right)} $.

Thus we are in conditions of the abstract Riemann--Roch theorem
(Theorem~\ref{th6.50}), which finishes the proof. \end{proof}

\begin{remark} This theorem is an infinite-genus variant of Riemann--Roch
theorem for the case of the bundle $ {\mathcal O} $ on an algebraic curve $ M $. Indeed, the
latter theorem says that the mapping $ \bar{\partial}\colon H^{s}\left(M,{\mathcal O}\right) \to H^{s-1}\left(M,\bar{\omega}\right) $
between Sobolev spaces is a Fredholm mapping of index $ 1-g\left(M\right) $. The theory
of elliptic operators says that in a case of smooth compact $ M $ the value
of $ s $ is irrelevant, but in our context we are forced to use the value
$ s=1 $.

The relation of the mapping $ \bar{\partial} $ to the mapping $ {\mathcal J} $ is the standard (in
mathematical physics) trick of {\em reduction to boundary}. We do cuts in the
curve $ M $ to obtain a region $ S={\mathbb C}P^{1}\smallsetminus\bigcup K_{i} $. Instead of $ H^{s}\left(M,{\mathcal O}\right) $ we consider
subspace $ V $ of $ H^{s}\left(S,{\mathcal O}\right) $ which consists of functions which satisfy gluing
conditions on the boundary. These gluing conditions are similar to $ {\mathcal J}f=0 $,
but are applied to functions (instead of functions modulo constants),
thus there is an extra condition per cut, total $ g $ extra conditions. Let
\begin{equation}
{\mathcal J}_{H^{1}}\colon H^{1}\left(S,{\mathcal O}\right)/\operatorname{const} \to \bigoplus_{l_{2}}H^{1/2}\left(\partial K_{j}\right)/\operatorname{const}
\notag\end{equation}
be defined by the formula~\eqref{equ4.101}. Thus $ V/\operatorname{const} $ is a subspace of $ \operatorname{Ker}{\mathcal J}_{H^{1}} $
of codimension $ g-1 $. Now the translation of the classical Riemann--Roch
theorem is that the mapping $ \bar{\partial}\colon \operatorname{Ker}{\mathcal J}_{H^{1}} \to H^{s-1}\left(S,\bar{\omega}\right) $ is of index 0. (Note
that the choice $ s=1 $ insures that there should be no gluing conditions for
elements of $ H^{s-1}\left(M,\bar{\omega}\right) $.)

Now to solve the equation $ \bar{\partial}f=h $, $ {\mathcal J}_{H^{1}}f=0 $ we reduce it to a boundary
value problem: we use the fact that $ \bar{\partial} $ is a surjection on $ {\mathbb C}P^{1} $, thus one
can immediately find a function $ f_{0} $ such that $ \bar{\partial}f_{0}=h $. Thus $ f_{1}=f-f_{0} $ should
satisfy $ \bar{\partial}f_{1}=0 $, $ {\mathcal J}_{H^{1}}f_{1}=-{\mathcal J}_{H^{1}}f_{0} $. This shows that the index of the mapping $ {\mathcal J} $
is 0 in the algebraic case $ g<\infty $. We see that Theorem~\ref{th35.45} is indeed an
infinite-dimensional analogue of Riemann--Roch theorem for the bundle $ {\mathcal O} $.
\end{remark}

\begin{remark} It is possible to modify Theorem~\ref{th35.45} to make it applicable
to deformations of the bundle $ {\mathcal O} $ as well (see Section~\ref{s5.30} for details).
Since to describe Jacobian we need the special case $ {\mathcal O} $ only, we do not
pursue this venue here. \end{remark}

\begin{remark} For the particular case of the bundle $ {\mathcal O} $ it is possible to prove
a much stronger result than Riemann--Roch theorem: that $ {\mathcal J} $ is an
isomorphism. However, since this result is not true out of context of toy
theory, we postpone its proof until Section~\ref{s7.90}, when it is needed for
description of Jacobian. \end{remark}

\section{$ H^{1/2} $-theory }

As we have seen in Section~\ref{h35}, Riemann--Roch theorem can be
easily proven given appropriate ingredients, such as {\em almost
perpendicularity\/} of functions with far-separated supports, ability to
{\em restrict\/} a global section to the boundary, ability to {\em invert\/} $ \bar{\partial} $ basing on
boundary values, and an ability to glue a function provided we know the
values inside and outside the given curve with appropriate compatibility
conditions along the curve. However, when one tries to apply the same
technique to the topology we are most interested in, i.e.,
$ H^{1/2} $-topology, significant difficulties arise.

\subsection{In an ideal world } In what follows we use the (complex) case $ n=1 $ of
Section~\ref{s2.50}. We have seen that the space $ H^{1/2}\left(M,\omega^{1/2}\right) $ should have a
special significance in studying the holomorphic half-forms. Moreover,
the results of Section~\ref{s2.70} suggest that one would be able to describe
holomorphic elements of $ H^{1/2}\left(D,\omega^{1/2}\right) $, $ D\subset M $, by their restrictions on $ \partial D $.
If our world were the ideal world, then the following properties would be
satisfied:
\begin{enumerate}
\item
The action of $ \operatorname{SL}\left(2,{\mathbb C}\right) $ on $ H^{1/2}\left({\mathbb C}P^{1},\omega^{1/2}\right) $ would be an action by
uniformly bounded operators;
\item
For an embedded curve $ \gamma \to {\mathbb C}P^{1} $ the restriction mappings
$ r\colon H^{s}\left({\mathbb C}P^{1},\omega^{1/2}\right) \to H^{s-1/2}\left(\gamma,\Omega^{1/2}\otimes\mu\right) $ (defined for $ s>\frac{1}{2} $) would be defined
for $ s=\frac{1}{2} $ as well;
\item
For the same curve the mapping $ e\colon H^{s}\left(\gamma,\Omega^{1/2}\otimes\mu\right) \to H^{s-1/2}\left({\mathbb C}P^{1},\omega^{1/2}\otimes\bar{\omega}\right) $ of
continuation by $ \delta $-function (defined for $ s<0 $) would be defined for $ s=0 $ as
well.
\item
If a curve $ \gamma $ divides a domain $ S $ into two parts $ S_{1} $, $ S_{2} $, then the
natural restriction mapping $ H^{1/2}\left(S\right) \to H^{1/2}\left(S_{1}\right)\oplus H^{1/2}\left(S_{2}\right) $ would be an
isomorphism.
\end{enumerate}

To make a long story short, in the ideal world the main technical
tools of this paper would behave in a civilized manner, which would
simplify the exposition a lot.

Identification of the Hardy space with the subspace of holomorphic
functions would be provided by the restriction $ r $. The composition
\begin{equation}
H^{1/2}\left({\mathbb C}P^{1},\omega^{1/2}\right) \xrightarrow[]{r} H^{0}\left(\gamma,\Omega^{1/2}\otimes\mu\right) \xrightarrow[]{e} H^{-1/2}\left({\mathbb C}P^{1},\omega^{1/2}\otimes\bar{\omega}\right)
\notag\end{equation}
would be a ``model'' of $ \bar{\partial} $-operator\footnote{In the sense of Section~\ref{s3.05}, i.e., it would coincide with $ -\bar{\partial} $ on
piecewise holomorphic functions.}
\begin{equation}
H^{1/2}\left({\mathbb C}P^{1},\omega^{1/2}\right) \xrightarrow[]{\bar{\partial}} H^{-1/2}\left({\mathbb C}P^{1},\omega^{1/2}\otimes\bar{\omega}\right),
\notag\end{equation}
so that the composition (Cauchy formula)
\begin{equation}
H^{1/2}\left({\mathbb C}P^{1},\omega^{1/2}\right) \xrightarrow[]{r} H^{0}\left(\gamma,\Omega^{1/2}\otimes\mu\right) \xrightarrow[]{e} H^{-1/2}\left({\mathbb C}P^{1},\omega^{1/2}\otimes\bar{\omega}\right)\xrightarrow[]{\bar{\partial}^{-1}} H^{1/2}\left({\mathbb C}P^{1},\omega^{1/2}\right)
\notag\end{equation}
would be continuous. By construction the functions in the image of this
operator are holomorphic outside of $ \gamma $, and the operator $ r $ gives
an injection
of this subspace into $ H^{0}\left(\gamma,\Omega^{1/2}\otimes\mu\right) $. Now one may consider
this as a Cauchy formula, since the composition is an identity operator
on the subspace of holomorphic functions in the domain $ D $ bounded by $ \gamma $.

All these operators would be canonically defined by $ D $, and
compatible with projective mappings $ {\mathbb C}P^{1} \to {\mathbb C}P^{1} $, $ D \to D_{1} $. Thus one would
be able to consider the image of $ r $ on subspace of holomorphic forms
in $ H^{1/2}\left(D,\omega^{1/2}\right) $ as a ``model'' of this space.

The last property in the list would allow us to glue together the
forms which are provided by different means on pieces $ S_{1} $ and $ S_{2} $.

Since we are confined to the current world, the above program will
not work, thus we need some workarounds against above three non-facts. We
will consider a Riemannian structure on $ {\mathbb C}P^{1} $, and will need to control the
size of domains in question, in order for results of Section~\ref{s2.50} to be
applicable. We will also need mollifications of operators $ r $ and $ e $, and
will need some restrictions on what we can glue together.

\subsection{Mollification }\label{s2.10}\myLabel{s2.10}\relax  Here we introduce a mapping
\begin{equation}
H^{0}\left(S^{1},\Omega^{1/2}\otimes\mu\right) \xrightarrow[]{\widetilde{e}} H^{-1/2}\left(S^{1}\times\left(-\varepsilon,\varepsilon\right),\omega^{1/2}\right)
\notag\end{equation}
which a closest existing analogue for the (non-existing) mapping $ e $ from
the previous section. It will be an invertible mapping onto its image,
and it will depend on additional parameters $ a_{n} $, $ n\in{\mathbb Z} $. (Later we will need
some particular choice of parameters $ \left(a_{n}\right) $, appropriate for the elliptic
operator $ \bar{\partial} $ we study.)

Fix $ \varepsilon>0 $ and a sequence $ \left(a_{n}\right)_{n\in{\mathbb Z}} $ such that $ 0<a<a_{n}<A $ for fixed
constants $ a $ and $ A $. Fix a smooth function $ \sigma\left(y\right) $, $ y\in\left(-\varepsilon,\varepsilon\right) $, with a
compact support and integral 1. Now map a half-form $ e^{2\pi i kx}dx^{1/2} $,
$ x\in S^{1}={\mathbb R}/{\mathbb Z} $, $ k\in{\mathbb Z}+\frac{1}{2} $, into
\begin{equation}
a_{k}e^{2\pi i kx}\frac{\sigma\left(|k|y\right)}{|k|}dz^{1/2},\qquad x\in S^{1}\text{, }y\in\left(-\varepsilon,\varepsilon\right),\quad z=x+iy,
\notag\end{equation}
and continue this mapping linearly to $ L_{2}\left(S^{1},\Omega^{1/2}\otimes\mu\right)=H^{0}\left(S^{1},\Omega^{1/2}\otimes\mu\right) $.

\begin{lemma} This mapping is an injection, i.e., an invertible mapping onto its
image. \end{lemma}

The dual mapping $ H^{1/2}\left(S^{1}\times\left(-\varepsilon,\varepsilon\right),\omega^{1/2}\right) \xrightarrow[]{\widetilde{r}} H^{0}\left(S^{1},\Omega^{1/2}\otimes\mu\right) $ given by
\begin{align} f\left(x,y\right)dz^{1/2} \mapsto g\left(x\right) & = \sum_{k}g_{k}e^{2\pi kx}dx^{1/2},
\notag\\
g_{k} & = \int f\left(x,y\right)e^{-2\pi kx}\frac{\sigma\left(|k|y\right)}{|k|}dx\,dy 
\notag\end{align}
has a similar property: it is a surjection, i.e., an invertible mapping
from the quotient by its null-space. It is a close analogue of the
(non-existing) mapping $ r $ from the previous section.

Suppose $ a_{k}=1 $ for any $ k $. If $ \varepsilon \to $ 0, then in weak topology for
mappings $ {\mathcal D} \to {\mathcal D}' $ the constructed mappings converge to the mappings of
extension as $ \delta $-function and restriction, but on Sobolev spaces the norms
of these mappings go to $ \infty $,

\begin{remark} In what follows we will use this mapping with following
modifications: we suppose that $ \operatorname{Supp}\sigma\subset\left[0,\varepsilon\right] $, thus we map a half-form on $ S^{1} $
to a half-form concentrated on a small collar to the {\em right\/} of $ S^{1} $ in
$ S^{1}\times\left[-\varepsilon,\varepsilon\right] $. In fact $ S^{1}\times\left[-\varepsilon,\varepsilon\right] $ will be identified with a
annulus
$ \left\{R\,e^{-2\pi\varepsilon}<|z-z_{0}|<R\,e^{2\pi\varepsilon}\right\} $ in $ {\mathbb C} $, and the mapping would send half-forms on
$ \left\{|z|=R\right\} $ to half-forms concentrated on the {\em outside\/} collar.

Similarly, the dual mapping
\begin{equation}
H^{1/2}\left(R\,e^{-2\pi\varepsilon}<|z-z_{0}|<R\,e^{2\pi\varepsilon}\right) \xrightarrow[]{\widetilde{r}} H^{0}\left(\left\{|z|=R\right\}\right)
\notag\end{equation}
will depend only on value of $ f\left(z\right) $ on the outside collar. \end{remark}

\subsection{Mollification suitable for $ \bar{\partial} $ }\label{s5.31}\myLabel{s5.31}\relax  Here we are going to fix the
values for the coefficients $ \left(a_{n}\right) $ from Section~\ref{s2.10} which are most
suitable for the operator $ \bar{\partial} $.

\begin{proposition} \label{prop4.15}\myLabel{prop4.15}\relax  Fix $ \varepsilon>0 $. Let $ \gamma= \left\{z \mid |z|=1\right\} $, $ \widetilde{K} = \left\{z \mid |z|<e^{2\pi\varepsilon}\right\} $. For
appropriate $ A $ and $ a $ there exists a sequence $ \left(a_{n}\right) $, $ 0<a<a_{n}<A $, such that
the corresponding restriction mapping $ \widetilde{r}\colon H^{1/2}\left(\widetilde{K},\omega^{1/2}\right) \to L_{2}\left(\gamma,\Omega^{1/2}\otimes\mu\right) $
coincides with the usual restriction on half-forms which are holomorphic
between $ \gamma $ and $ \widetilde{K} $. The numbers $ a_{n} $, thus the operator $ \widetilde{r} $, are uniquely
determined.

\end{proposition}

Dually,

\begin{proposition} \label{prop4.16}\myLabel{prop4.16}\relax  Let $ U=K\cup\left({\mathbb C}P^{1}\smallsetminus\widetilde{K}\right) $. Consider the Cauchy kernel
restricted to $ \widetilde{K} $. It gives two operators: an operator $ \overset{\,\,{}_\circ} H^{-1/2}\left(\widetilde{K}\right) \xrightarrow[]{\bar{\partial}^{-1}}
H^{1/2}\left(U\right) $, and an operator $ C^{\infty}\left(\gamma\right) \xrightarrow[]{\bar{\partial}^{-1}} {\mathcal D}'\left(U\right) $. Consider the mollification $ \widetilde{e} $
of the extension mapping $ e $ corresponding to a sequence $ \left(a_{n}\right) $. With
appropriate choice of the sequence $ \left(a_{n}\right) $ the composition
\begin{equation}
C^{\infty}\left(\gamma\right)\hookrightarrow L_{2}\left(\gamma\right) \xrightarrow[]{\widetilde{e}} \overset{\,\,{}_\circ} H^{-1/2}\left(\widetilde{K}\right) \xrightarrow[]{\bar{\partial}^{-1}} H^{1/2}\left(U\right) \hookrightarrow {\mathcal D}'\left(U\right)
\notag\end{equation}
coincides with the mapping $ C^{\infty}\left(\gamma\right) \xrightarrow[]{\bar{\partial}^{-1}} {\mathcal D}'\left(U\right) $. The numbers $ a_{n} $, thus the
operator $ \widetilde{e} $, are uniquely determined by the above condition.

\end{proposition}

\begin{remark} Note that the mapping $ \widetilde{r} $ is a left inverse to the mapping $ f \mapsto
c_{f} $ from Section~\ref{s2.70}. \end{remark}

\subsection{$ \protect \widetilde{e}\circ\protect \widetilde{r} $ as an approximation to $ \bar{\partial} $ }\label{s2.25}\myLabel{s2.25}\relax  Consider a circle $ \gamma $ in $ {\mathbb C} $. In what
follows we consider particular mollifications of the $ \delta $-inclusion $ C^{\infty}\left(\gamma\right)
\xrightarrow[]{e} {\mathcal D}'\left({\mathbb C}\right) $ and restriction $ {\mathcal D}\left({\mathbb C}\right) \xrightarrow[]{r} C^{\infty}\left(S^{1}\right) $, which correspond to the only
sequences $ \left(a_{n}\right) $ which satisfy the conditions of Propositions~\ref{prop4.15},
~\ref{prop4.16}. The notations $ \widetilde{e}_{\gamma} $ and $ \widetilde{r}_{\gamma} $ are reserved for these two mappings,
we may denote them $ \widetilde{e} $ and $ \widetilde{r} $ if the circle to apply them for is clear from
context.

A central tool in the following discussion is the mollification of
identity~\eqref{equ3.02}:

\begin{proposition} Consider a circle $ \gamma $ or radius $ \rho $ which bounds a disk $ K $.
Let $ \widetilde{K} $ be a concentric disk of radius $ \rho e^{2\pi\varepsilon} $. Let $ F $ be a holomorphic
half-form in $ \widetilde{K}\smallsetminus K $. Let
\begin{equation}
G=\bar{\partial}^{-1}\circ\widetilde{e}_{\gamma}\circ\widetilde{r}_{\gamma}\left(F\right).
\notag\end{equation}
If $ F $ is a holomorphic half-form in $ {\mathbb C}P^{1}\smallsetminus K $, then $ G $ coincides with $ F $ outside
of $ \widetilde{K} $, and is 0 in $ U $. If $ F $ is a holomorphic half-form in $ \widetilde{K} $, then $ G $
coincides with $ -F $ inside of $ K $, and is 0 outside of $ \widetilde{K} $. \end{proposition}

\begin{amplification} \label{amp4.21}\myLabel{amp4.21}\relax  Fix $ \varepsilon,R>0 $. Consider disks $ K_{\rho}=\left\{|z|<\rho\right\} $,
$ K'_{\rho}=\left\{|z|<\rho e^{2\pi\varepsilon}\right\} $, $ \widetilde{K}_{\rho}=\left\{|z|<\rho e^{4\pi\varepsilon}\right\} $ in $ {\mathbb C} $. Let $ \gamma=\partial K_{\rho} $. There exists a
mapping
\begin{equation}
\lambda_{\rho}\colon H^{1/2}\left(\widetilde{K}_{\rho},\omega^{1/2}\right) \to H^{1/2}\left(\widetilde{K}_{\rho},\omega^{1/2}\right)
\notag\end{equation}
such that
\begin{enumerate}
\item
For $ 0<\rho<R $ the mapping $ \lambda_{\rho} $ is continuous with the norm uniformly
bounded by a constant depending on $ \varepsilon $ and $ R $ only;
\item
$ \lambda_{\rho}\left(f\right)|_{K_{\rho}}=0 $ for any $ f\in H^{1/2}\left(\widetilde{K}_{\rho}\right) $;
\item
$ \left(f-\lambda_{\rho}\left(f\right)\right)|_{\widetilde{K}_{\rho}\smallsetminus K'_{\rho}}=0 $ for any $ f\in H^{1/2}\left(\widetilde{K}_{\rho}\right) $;
\item
if $ f $ is holomorphic in $ \widetilde{K}_{\rho}\smallsetminus K_{\rho} $, then $ \bar{\partial}\lambda_{\rho}\left(f\right) = \left(\widetilde{e}_{\gamma}\circ\widetilde{r}_{\gamma}\right)\left(f\right) $.
\end{enumerate}
\end{amplification}

\begin{proof} During the proof we abuse notations and do not mention the
bundle $ \omega^{1/2} $ in notations for Sobolev spaces.

Note that for a half-form $ f $ which is holomorphic outside of $ K $ all
the statements of the amplification are true if we take
$ \lambda_{\rho}\left(f\right)=\left(\bar{\partial}^{-1}\circ\widetilde{e}\circ\widetilde{r}\right)\left(f\right) $. If $ f $ is holomorphic inside of $ \widetilde{K} $, then
$ \lambda_{\rho}\left(f\right)=f-\left(\bar{\partial}^{-1}\circ\widetilde{e}\circ\widetilde{r}\right)\left(f\right) $ works. This uniquely defines $ \lambda_{\rho}\left(f\right) $ if $ f $ is
holomorphic inside of $ \widetilde{K}_{\rho}\smallsetminus K_{\rho} $. What we need is to adjust this formula to
the case of non-holomorphic half-forms.

Since $ R $ is fixed, we can consider the norm from Section~\ref{s2.50}
instead of the equivalent Sobolev $ H^{1/2} $-norm. Since for the former norm
the Sobolev spaces $ H^{1/2}\left(\widetilde{K}_{\rho}\right) $ with different $ \rho $ are naturally isomorphic, it
is enough to consider $ \rho=1 $.

Note that for a holomorphic $ f $ the image $ \lambda_{\rho}\left(f\right) $ depends on $ f|_{\widetilde{K}_{\rho}\smallsetminus K_{\rho}} $
only. We are going to define $ \lambda_{\rho} $ in general case such that it satisfies
the same condition. Thus $ \lambda_{\rho} $ is a mapping $ H^{1/2}\left(\widetilde{K}_{\rho}\smallsetminus K_{\rho}\right) \to H^{1/2}\left(\widetilde{K}_{\rho}\right) $ such
that a half-form in the image is 0 inside $ K_{\rho} $. We may substitute
conformally equivalent domain $ S^{1}\times\left(-2\varepsilon,2\varepsilon\right) $ instead of $ \widetilde{K}_{\rho}\smallsetminus K_{\rho} $.

Now $ H^{1/2}\left(S^{1}\times\left(-2\varepsilon,2\varepsilon\right)\right) $ is an orthogonal sum of subspaces $ L_{n} $,
$ n\in{\mathbb Z}+\frac{1}{2} $, spanned by half-forms of the form $ \varphi\left(y\right)e^{2\pi i\quad nx}dz^{1/2} $,
$ \left(x,y\right)\in S^{1}\times\left(-2\varepsilon,2\varepsilon\right) $, $ z=x+iy $, thus it is enough to construct uniformly
bounded mappings in these subspaces. For a given $ \lambda $ let $ \lambda^{\left(n\right)} $ be defined as
\begin{equation}
\lambda\left(\varphi\left(y\right)e^{2\pi i nx}dz^{1/2}\right)=\lambda^{\left(n\right)}\left(\varphi\left(y\right)\right)e^{2\pi i nx}dz^{1/2}.
\notag\end{equation}
Note that the conditions on $ \lambda^{\left(n\right)} $ are: $ \lambda^{\left(n\right)}\left(\varphi\left(y\right)\right) $ vanishes if $ y<0 $,
$ \lambda^{\left(n\right)}\left(\varphi\left(y\right)\right)=\varphi\left(y\right) $ if $ y>\varepsilon $, and $ \lambda^{\left(n\right)}\left(e^{ny}\right) $ has a prescribed value
(obtained basing on the last condition of the amplification).

Define $ \widetilde{\sigma}_{n}\left(y\right) $ via
\begin{equation}
\lambda^{\left(n\right)}\left(e^{ny}\right) = \widetilde{\sigma}_{n}\left(y\right)e^{ny}.
\notag\end{equation}
Obviously, $ \widetilde{\sigma}_{n}\left(y\right)=0 $ if $ y<0 $, $ \widetilde{\sigma}_{n}\left(y\right)=1 $ if $ y>\varepsilon $. The function $ \widetilde{\sigma}_{n}\left(y\right) $ is
uniquely determined by the last condition of the theorem.

After all these remarks {\em define\/} $ \lambda|_{L_{n}} $ by
\begin{equation}
\varphi\left(y\right)e^{2\pi i nx}dz^{1/2} \to \sigma_{n}\left(y\right)\varphi\left(y\right)e^{2\pi i nx}dz^{1/2}.
\notag\end{equation}
This mapping satisfies all the conditions of the amplification, with a
possible exception of uniform boundness of these mappings for different
$ n $. To prove the boundness, note that
\begin{equation}
\widetilde{\sigma}_{n}\left(y\right)=\widetilde{\sigma}_{1}\left(ny\right),
\notag\end{equation}
and the mappings
\begin{equation}
L_{1} \xrightarrow[]{m_{n}} L_{n}\colon \varphi\left(y\right)e^{2\pi i x} \mapsto n^{-1/2}\varphi\left(ny\right)e^{2\pi i nx}
\notag\end{equation}
are uniformly bounded together with their inverse mappings. \end{proof}

\subsection{Mollification for $ L_{2}\left(\Omega^{1/2}\otimes\mu\right) $ }\label{s4.35}\myLabel{s4.35}\relax  Consider a family of circles $ \gamma_{i} $ on
$ {\mathbb C}P^{1} $.
Fix a projective isomorphism of $ \gamma_{i} $ with $ \left\{|z|=1\right\} $ for every $ i $. This
isomorphism identifies an annulus $ \left\{e^{-2\pi\varepsilon}<|z|<e^{2\pi\varepsilon}\right\} $ with some neighborhood $ U_{i} $
of $ \gamma_{i} $.

The $ \bar{\partial} $-adjusted mollification of $ e $ gives a mapping
\begin{equation}
\widetilde{e}_{i}\colon L_{2}\left(\gamma_{i},\Omega^{1/2}\otimes\mu\right) \to H^{-1/2}\left({\mathbb C}P^{1},\omega^{1/2}\otimes\bar{\omega}\right).
\notag\end{equation}
This mapping is a injection of topological vector spaces. In the
following sections we discuss under which conditions $ \sum_{i}\widetilde{e}_{i} $ is an
injection. Here $ \sum_{i}\widetilde{e}_{i} $ is considered as a mapping from $ L_{2}\left(\coprod\gamma_{i},\Omega^{1/2}\otimes\mu\right) $ to
$ H^{-1/2}\left({\mathbb C}P^{1},\omega^{1/2}\otimes\bar{\omega}\right) $.

This injection is going to be a central tool in Section~\ref{s5.10},
which allows one to identify the space of ``admissible'' holomorphic
half-forms with the space of their boundary values. To provide the
criterion for this mapping to be an injection, we need to introduce a
notion of ``almost perpendicular'' subspaces, to study relative position of
Sobolev subspaces corresponding to different domains, and to investigate
the relative position of domains inside $ {\mathbb C}P^{1} $.

First, we suppose that $ \varepsilon $ is small enough and the subsets $ U_{i} $ do not
intersect. Then it is clear that $ \sum_{i}\widetilde{e}_{i} $ has no null-vectors. Obviously, $ \sum_{i}\widetilde{e}_{i} $
were an injection if the images of $ \widetilde{e}_{i} $ were perpendicular for
different $ i $.
We are going to investigate when these images are {\em almost perpendicular\/}
for big $ i $. In other words, when the sum of images is direct in the sense
of Hilbert topology (i.e., it is an orthogonal sum after an appropriate
change of the Hilbert norm to an equivalent one).

\subsection{Almost perpendicular subspaces }\label{s5.61}\myLabel{s5.61}\relax 

\begin{lemma} \label{lm4.41}\myLabel{lm4.41}\relax  Consider a Hilbert space $ H $ and a family of subspaces $ H_{i} $, $ i\in{\mathbb N} $.
Denote by $ a_{ij} $ the orthogonal projector $ H_{i} \to H_{j} $. Let
$ A=\left(a_{ij}\right) $ be
the matrix of an operator $ {\bold A}\colon \bigoplus H_{i} \to \coprod H_{i} $, $ \coprod H_{i} $ being the space of
arbitrary sequences $ \left(h_{i}\right) $, $ h_{i}\in H_{i} $. Then the natural mapping $ \bigoplus H_{i} \to H $
extends to a Fredholm mapping
\begin{equation}
\bigoplus_{l_{2}}H_{i} \xrightarrow[]{i} \operatorname{Im} i\subset H,\qquad \left(h_{i}\right) \mapsto \sum h_{i},
\notag\end{equation}
onto its image iff {\bf A }induces a Fredholm mapping $ \bigoplus_{l_{2}}H_{i} \to \bigoplus_{l_{2}}H_{i} $. \end{lemma}

\begin{proof} Consider $ {\bold A} $ as $ i^{*}\circ i $. \end{proof}

\begin{corollary} Suppose that the conditions of the previous lemma the matrix $ \alpha =
\left(\|a_{ij}\| - \delta_{ij}\right) $ corresponds to a compact mapping $ l_{2} \to l_{2} $. Then the mapping
$ \oplus_{l_{2}}H_{i} \xrightarrow[]{i} H $ is a Fredholm operator onto its image. In particular, this is
true if
\begin{equation}
\sum_{i\not=j}\|a_{ij}\|^{2}<\infty.
\notag\end{equation}
Moreover, if there is a number $ C $ and a family of unit vectors $ v_{i}\in V_{i} $ such
that
\begin{equation}
\|a_{ij}\| < C|\left(v_{i},v_{j}\right)|,
\notag\end{equation}
then the mapping $ i $ is continuous iff $ \alpha $ provides a bounded mapping $ l_{2} \to
l_{2} $.

\end{corollary}

\subsection{Conformal distance } In the applications we consider the curve $ \gamma $ is a
union of countably many connected components $ \gamma=\coprod\gamma_{i} $, and all the
components $ \gamma_{i} $ but a finite number are circles. Let $ \widetilde{e}_{i}=\widetilde{e}_{\gamma_{i}} $ be the mollified
extension mapping, and
\begin{equation}
H_{i}=\operatorname{Im} \widetilde{e}_{i}.
\notag\end{equation}
If $ H_{i} $ satisfy the conditions of Lemma~\ref{lm4.41}, then the mapping $ \sum_{i}\widetilde{e}_{i} $ from
Section~\ref{s4.35} is a Fredholm mapping onto its image. So
the next thing to study is how to calculate the norm of the projection
\begin{equation}
\operatorname{Im} \widetilde{e}_{i} \to \operatorname{Im} \widetilde{e}_{j}.
\label{equ4.33}\end{equation}\myLabel{equ4.33,}\relax 
The subspace $ \operatorname{Im} \widetilde{e}_{i} $ lies inside a subspace $ \overset{\,\,{}_\circ}{H}\left(\widetilde{K}_{i}\right) $, $ \widetilde{K}_{i} $ being
any disk which contains an appropriate neighborhood of $ \gamma_{i} $. Thus to
majorate the norm of the projection~\eqref{equ4.33} we start with
the case of two disjoint disks $ \widetilde{K}_{i} $, $ \widetilde{K}_{j} $ on $ {\mathbb C}P^{1} $, and consider the
orthogonal projection from $ \overset{\,\,{}_\circ}{H}^{-1/2}\left(\widetilde{K}_{i}\right) $ to $ \overset{\,\,{}_\circ}{H}^{-1/2}\left(\widetilde{K}_{j}\right) $ inside $ H^{-1/2}\left({\mathbb C}P^{1}\right) $. As
we will see, the norm of this projection can be majorated by a number
which depends only on some kind of {\em distance\/} between $ \widetilde{K}_{i} $ and $ \widetilde{K}_{j} $.

Two disjoint simple curves $ \gamma_{1} $, $ \gamma_{2} $ on $ {\mathbb C}P^{1} $ bound a tube $ U $,
which is
conformally equivalent to exactly one tube $ S^{1}\times\left[0,l\right] $, $ l>0 $. Here $ S^{1} $ is
a circle of circumference $ 2\pi $.

\begin{definition} We call the number $ l $ the {\em conformal distance\/} between $ \gamma_{1} $ and
$ \gamma_{2} $. \end{definition}

\begin{proposition} Consider three simple curves $ \gamma_{1} $, $ \gamma_{2} $, $ \gamma_{3} $ such that $ \gamma_{2} $
separates $ \gamma_{1} $ and $ \gamma_{3} $. Then the conformal distance $ l\left(\gamma_{1},\gamma_{3}\right) \geq
l\left(\gamma_{1},\gamma_{2}\right) +l\left(\gamma_{2},\gamma_{3}\right) $. \end{proposition}

\begin{proof} Conformal distance between $ \gamma $ and $ \gamma' $ is $ \geq l $ iff there exists a
function $ \varphi $ defined between $ \gamma $ and $ \gamma' $ such that $ \varphi|_{\gamma}=0 $, $ \varphi|_{\gamma'}=l $, and
the ``energy'' $ \int\partial\varphi\bar{\partial}\varphi \leq l $. Combining two such functions, one defined between
$ \gamma_{1} $ and $ \gamma_{2} $, another between $ \gamma_{2} $ and $ \gamma_{3} $, we obtain the statement. \end{proof}

\begin{lemma} The only $ \operatorname{PGL}_{2}\left({\mathbb C}\right) $-invariant of a couple of disjoint disks on $ {\mathbb C}P^{1} $
is the conformal distance. \end{lemma}

\subsection{Subspaces of $ \Omega^{3/4} $ }\label{s3.8}\myLabel{s3.8}\relax  Consider two disjoint disks $ K_{1} $, $ K_{2} $ on $ {\mathbb C}P^{1} $ of
conformal
distance $ l $. Since the natural norm on $ H^{-1/2}\left(\Omega^{3/4}\right) $ (see Section~\ref{s2.60}) is
invariant with respect
to $ \operatorname{PGL} $, the angle between the subspaces $ \overset{\,\,{}_\circ}{H}^{-1/2}\left(K_{1},\Omega^{3/4}\right) $ and
$ \overset{\,\,{}_\circ}{H}^{-1/2}\left(K_{2},\Omega^{3/4}\right) $ depends on $ l $ only. Let $ P_{l} $ be the orthogonal projector from
one subspace to another.

\begin{proposition} \label{prop3.170}\myLabel{prop3.170}\relax  Fix $ \varepsilon>0 $. If $ l\geq\varepsilon $, then
\begin{equation}
\|P_{l}\| \sim e^{-l/2}.
\notag\end{equation}
The equivalence means that the quotient of two sides remains bounded and
separated from 0. \end{proposition}

\begin{proof} The norm $ \|P_{l}\| $ is a smooth function of $ l $, thus one needs to
prove only the asymptotic when $ l \to \infty $. One may represent two disks of
conformal distance $ \approx l\gg 1 $ as $ |z|\leq1 $ and $ |z-e^{l/2}|\leq1 $. Let $ L=e^{l/2} $. Consider the
coordinate $ z_{1}=z $ in the first disk, and $ z_{2}=z-L $ in the second one. What we
need to prove is that the kernel $ K_{1}\left(z_{1},z_{2}\right)=\frac{1}{|z_{1}-z_{2}+L|} $ in $ \left\{\left(z_{1},z_{2}\right) |
|z_{1}|,|z_{2}|\leq1\right\} $ gives an operator of the norm $ \sim1/L $.

Since the radii are fixed now, we can consider functions instead of
$ 3/4 $-forms. Since the norm of operator with the kernel $ K_{0}=\frac{1}{L} $ is
$ O\left(\frac{1}{L}\right) $, it is enough to prove that the kernel $ K_{2}\left(z_{1},z_{2}\right)=\frac{1}{|z_{1}-z_{2}+L|} -
\frac{1}{L} $ corresponds to an operator of norm $ o\left(1/L\right) $. Let us estimate
Hilbert--Schmidt norm of this operator. It is equal to the
$ H^{1/2}\otimes_{l_{2}}H^{1/2} $-norm of $ K_{2} $. On the other hand, the last norm is bounded by
$ H^{1} $-norm, which is obviously $ O\left(L^{-2}\right) $.\end{proof}

\begin{remark} \label{rem3.75}\myLabel{rem3.75}\relax  Note that if $ 0<s<1 $, then a similar statement is true in
$ H^{1-2s}\left(\Omega^{s}\right) $. We will need only the following statement:
\begin{equation}
\|P_{l}\| = O\left(e^{-l/2}\right)\text{ if }0<s<\frac{3}{4},
\notag\end{equation}
and we will use it in the case $ s=\frac{1}{4} $ only. \end{remark}

\begin{corollary} \label{cor3.80}\myLabel{cor3.80}\relax  Consider a family of disjoint closed disks $ K_{i}\subset{\mathbb C}P^{1} $,
$ i\in{\mathbb N} $, with conformal distance between $ K_{i} $ and $ K_{j} $ being $ l_{ij} $. Put $ l_{i i}=0 $. Let
$ {\mathcal L} $ be a linear bundle on $ {\mathbb C}P^{1} $. Suppose that the closure
of $ \bigcup K_{i} $ does not coincide with $ {\mathbb C}P^{1} $. Suppose that $ \inf _{i\not=j}l_{ij}>0 $, and that the
matrix $ \left(e^{-l_{ij}/2}\right) $ corresponds to a bounded operator $ l_{2} \to l_{2} $. Then the
inclusions $ \overset{\,\,{}_\circ}{H}^{-1/2}\left(K_{i},{\mathcal L}\right)\hookrightarrow H^{-1/2}\left({\mathbb C}P^{1},{\mathcal L}\right) $ extend to a continuous injection
\begin{equation}
\bigoplus_{l_{2}}\overset{\,\,{}_\circ}{H}^{-1/2}\left(K_{i},{\mathcal L}\right)\hookrightarrow H^{-1/2}\left({\mathbb C}P^{1},{\mathcal L}\right).
\notag\end{equation}
Dually, the restrictions $ H^{1/2}\left({\mathbb C}P^{1},{\mathcal L}\right) \to H^{1/2}\left(K_{i},{\mathcal L}\right) $ extend to a continuous
surjection
\begin{equation}
H^{1/2}\left({\mathbb C}P^{1},{\mathcal L}\right) \to \bigoplus_{l_{2}}H^{1/2}\left(K_{i},{\mathcal L}\right).
\notag\end{equation}
\end{corollary}

\begin{proof} If $ {\mathcal L}=\Omega^{3/4} $ in the first part of the theorem, or $ {\mathcal L}=\Omega^{1/4} $ in the
second one, then everything is proved. Otherwise let $ U $ be an open subset
of $ {\mathbb C}P^{1}\smallsetminus\overline{\bigcup K_{i}} $. An isomorphism of $ {\mathcal L} $ and $ \Omega^{3/4} $ (or $ \Omega^{1/4} $) on $ {\mathbb C}P^{1}\smallsetminus U $ proves
the rest. \end{proof}

Using the results of Remark~\ref{rem3.75} we obtain the following
statement:

\begin{proposition} In the conditions of Corollary~\ref{cor3.80}
\begin{equation}
\bigoplus_{l_{2}}\overset{\,\,{}_\circ}{H}^{1/2}\left(K_{i},{\mathcal L}\right) \to H^{1/2}\left({\mathbb C}P^{1},{\mathcal L}\right).
\notag\end{equation}
(components being inclusions) is an isomorphism to its image. The image
of this mapping is $ \overset{\,\,{}_\circ}{H}^{1/2}\left(\bigcup_{i}K_{i},{\mathcal L}\right) $. \end{proposition}

\section{Generalized Hardy space }\label{h4}\myLabel{h4}\relax 

\subsection{Hilbert space of holomorphic half-forms }\label{s5.10}\myLabel{s5.10}\relax 

In this section we are going to construct a Hilbert space (or at
least a space with Hilbert topology) which models global holomorphic
half-forms on a curve which is not necessarily compact.

Consider a family of disjoint closed disks $ K_{i}\subset{\mathbb C}P^{1} $, $ i\in I $.
Let $ \varepsilon>0 $, $ \widetilde{K}_{i} $ be the concentric disk to $ K_{i} $ of radius $ e^{2\varepsilon}\cdot\operatorname{radius}\left(K_{i}\right) $.
Suppose that the disks $ \widetilde{K}_{i} $ do not intersect, and denote the conformal distance
between $ \widetilde{K}_{i} $ and $ \widetilde{K}_{j} $ by $ l_{ij} $. Put $ l_{i i}=0 $. Suppose that

\begin{nwthrmii} Say that the collection of disks is {\em well-separated\/} if
\begin{enumerate}
\item
$ \inf _{i\not=j}l_{ij}>0 $;
\item
$ \overline{\bigcup\widetilde{K}_{i}}\not={\mathbb C}P^{1} $; and
\item
the matrix $ \left(e^{-l_{ij}/2}\right) $ gives a bounded operator $ l_{2} \to l_{2} $;
\end{enumerate}
\end{nwthrmii}

Note that the interior of the closure $ \overline{{\mathbb C}P^{1}\smallsetminus\bigcup K_{i}} $ is non-empty by
the above conditions, and

\begin{proposition} \label{prop5.16}\myLabel{prop5.16}\relax  If disks $ K_{i} $ are well-separated, then $ \sum_{i}\operatorname{radius}\left(K_{i}\right)<\infty $,
thus Cauchy kernel defines a mapping from $ L_{2}\left(\bigcup\partial K_{i},\Omega^{1/2}\otimes\mu\right) $ to the space
of analytic functions in $ {\mathbb C}P^{1}\smallsetminus\overline{\bigcup K_{i}} $. \end{proposition}

\begin{proof} Indeed, since rows of $ \left(e^{-l_{ij}/2}\right) $ are in $ l_{2} $, $ \sum_{j} e^{-l_{ij}} <\infty $. On the
other hand, $ e^{-l_{ij}} $ is approximately proportional to radius of $ K_{j} $ when $ i $ is
fixed. \end{proof}

Moreover, if we slightly decrease the value of $ \varepsilon $, then the first two
conditions on the family $ \left\{K_{i}\right\} $ will be automatically satisfied, so only
the last condition is the new one.

\begin{definition} The {\em generalized Hardy space\/} $ {\mathcal H}\left({\mathbb C}P^{1}, \left\{K_{i}\right\}\right) $ is the space
\begin{equation}
\left\{f\left(z\right)\in H^{1/2}\left({\mathbb C}P^{1}\smallsetminus\bigcup_{i}K_{i}, \omega^{1/2}\right) \mid \bar{\partial}f=0\in H^{-1/2}\left({\mathbb C}P^{1}\smallsetminus\bigcup_{i}K_{i},\omega^{1/2}\otimes\bar{\omega}\right)\right\}.
\notag\end{equation}
\end{definition}

Note that elements of $ {\mathcal H}\left({\mathbb C}P^{1}, \left\{K_{i}\right\}\right) $ are holomorphic half-forms in
$ {\mathbb C}P^{1}\smallsetminus\overline{\bigcup\widetilde{K}_{i}} $. Note also that the Sobolev space in the definition is a
generalized Sobolev space, since $ {\mathbb C}P^{1}\smallsetminus\bigcup_{i}K_{i} $ is not open if the set $ I $ is
infinite.

\begin{remark} Note that if $ I $ consists of one element, then $ {\mathcal H} $ is the usual
Hardy space from Section~\ref{s2.70} with modifications outlined in Remark
~\ref{rem2.10}. \end{remark}

In Section~\ref{s4.90} we introduce a slightly weaker condition on
domains $ K_{i} $, and will use it instead of Condition A. Note that we are going
to mention the space $ {\mathcal H} $ only in cases when the Condition A is satisfied.

\subsection{Hilbert operator }\label{s5.20}\myLabel{s5.20}\relax  A section of $ \Omega^{1/2}\otimes\mu $ on a circle can be
decomposed into the Fourier series
\begin{equation}
f\left(t\right)dt^{1/2}=\sum_{n\in{\mathbb Z}+1/2}a_{n}e^{2\pi i nt}dt^{1/2}
\notag\end{equation}
(half-integers appear because of the factor $ \mu $), hence it can be written
as a sum of a component $ f_{+} $ which can be holomorphically extended to a
section of $ \omega^{1/2} $ inside the circle, and a component $ f_{-} $ which may be
extended outside of the circle. Moreover, norms of $ f_{\pm} $ are bounded by the
norm of $ f $, and the decomposition is unique.

We denote this decomposition $ L_{2}^{+}\left(\gamma,\Omega^{1/2}\otimes\mu\right) \oplus L_{2}^{-}\left(\gamma,\Omega^{1/2}\otimes\mu\right) $.

\begin{proposition} Consider the operator $ {\bold K}_{+} $ with Cauchy kernel acting in the
space $ \bigoplus_{l_{2}}L_{2}\left(\partial K_{i},\Omega^{1/2}\otimes\mu\right) $. Let $ {\bold K} $ be the operator $ \widetilde{{\bold K}} $ with diagonal blocks
$ L_{2}\left(\partial K_{i},\Omega^{1/2}\otimes\mu\right) \to L_{2}\left(\partial K_{i},\Omega^{1/2}\otimes\mu\right) $ removed. If the disks $ K_{i} $ are
well-separated, then $ {\bold K} $ is bounded. If we decompose
\begin{equation}
L_{2}\left(\partial K_{i},\Omega^{1/2}\otimes\mu\right) = L_{2}^{+}\left(\partial K_{i},\Omega^{1/2}\otimes\mu\right) \oplus L_{2}^{-}\left(\partial K_{i},\Omega^{1/2}\otimes\mu\right),
\notag\end{equation}
then the only nonzero blocks of $ {\bold K} $ act from $ L_{2}^{-}\left(\partial K_{i},\Omega^{1/2}\otimes\mu\right) $ to
$ L_{2}^{+}\left(\partial K_{j},\Omega^{1/2}\otimes\mu\right) $, $ i\not=j $. \end{proposition}

This statement is an immediate corollary of the first part of

\begin{lemma} \label{lm5.05}\myLabel{lm5.05}\relax  Consider a Hilbert space $ H=\bigoplus_{l_{2}}H_{i} $ and a linear operator
$ A\colon H \to H $ with bounded blocks $ A_{ij}\colon H_{j} \to H_{i} $. Let $ a_{ij}=\|A_{ij}\| $. Suppose that
the matrix $ \left(a_{ij}\right) $ gives a bounded operator $ \alpha\colon l_{2} \to l_{2} $. Then
\begin{enumerate}
\item
$ A $ is bounded, $ \|A\|\leq\|\alpha\| $;
\item
$ A $ is compact if $ \alpha $ is compact and each of operators $ A_{ij} $ is compact.
\end{enumerate}
\end{lemma}

\begin{proof} To bound the operator $ A $ it is enough to bound $ |\left(x, Ay\right)| $ for
$ |x|=|y|=1 $, $ x,y\in H $. On the other hand, decomposition of $ H $ gives $ x=\left(x_{i}\right) $,
$ y=\left(y_{i}\right) $, and $ |x|^{2}=\sum|x_{i}|^{2} $, $ |y|^{2}=\sum|y_{i}|^{2} $. This gives
\begin{equation}
|\left(x,Ay\right)| = \left|\sum_{ij}\left(x_{i},A_{ij}y_{j}\right)\right| \leq \sum_{ij}|\left(x_{i},A_{ij}y_{j}\right)| \leq \sum_{ij}\|A_{ij}\||x_{i}||y_{j}|
\notag\end{equation}
and the latter sum is just $ \left(\xi,\alpha\eta\right) $ with $ \xi=\left(|x_{i}|\right)\in l_{2} $, similarly for $ \eta $.

To prove the second part it is enough to show that we can
approximate $ A $ (in norm) by compact operators. Since $ \alpha $ is
compact, we can approximate it (in norm) by a matrix with finite number
of non-zero elements. The relation of norm of $ A $ and norm of $ \alpha $ shows that
$ A $ can be approximated by an operator $ A' $ which has the same blocks as $ A $
in finite number of places, all the rest is 0.

Since blocks of $ A $ are compact, $ A' $ is compact as well. We conclude
that $ A $ can be approximated with arbitrary precision by a compact
operator, thus $ A $ is compact.\end{proof}

\begin{corollary} \label{cor4.35}\myLabel{cor4.35}\relax  Any block-row and block-column of $ {\bold K} $ gives a compact
operator. \end{corollary}

\subsection{Boundary map } Here we consider which of the facts from Section~\ref{s2.70}
have sense for the generalized Hardy space as well.

Fix $ i\in I $, let $ S_{i}=\partial K_{i} $. Consider the mollified restriction
mapping
\begin{equation}
H^{1/2}\left(\widetilde{K}_{i}\smallsetminus K_{i},\omega^{1/2}\right) \xrightarrow[]{\widetilde{r}_{i}} L_{2}\left(S_{i},\Omega^{1/2}\otimes\mu\right).
\notag\end{equation}
Since any element $ f $ of the generalized Hardy space is $ H^{1/2} $ in $ \widetilde{K}_{i}\smallsetminus K_{i} $ (and
holomorphic inside this annulus), the restriction $ \widetilde{r}_{i}\left(f\right) =f|_{S_{i}} $ of this
half-form on $ S_{i} $ is an $ L_{2} $-section of $ \Omega^{1/2}\left(S_{i}\right)\otimes\mu $. In the rest of this
section we are going to abuse notations and denote $ \widetilde{r}_{i}\left(f\right) $ as $ f|_{S_{i}} $.

\begin{theorem} \label{th4.40}\myLabel{th4.40}\relax  Suppose that the disks $ K_{i} $ are well-separated. Then the
mappings $ \widetilde{r}_{i} $ taken together provide a mapping
\begin{equation}
{\mathcal H}\left({\mathbb C}P^{1},\left\{K_{i}\right\}\right) \xrightarrow[]{\widetilde{r}} \displaystyle\coprod L_{2}\left(S_{i},\Omega^{1/2}\otimes\mu\right),
\notag\end{equation}
which is fact is a continuous mapping
\begin{equation}
{\mathcal H}\left({\mathbb C}P^{1},\left\{K_{i}\right\}\right) \xrightarrow[]{\widetilde{r}} \bigoplus_{l_{2}}L_{2}\left(S_{i},\Omega^{1/2}\otimes\mu\right).
\notag\end{equation}
The image of $ \widetilde{r} $ is closed. Moreover, the mapping $ \widetilde{r} $ is invertible
onto its image.

\end{theorem}

\begin{proof} In fact, all the main ingredients for the proof of this theorem
are already here. An element $ f\in{\mathcal H}\left({\mathbb C}P^{1},\left\{K_{i}\right\}\right) $ is by definition a
restriction of some
$ H^{1/2} $-section $ g $ of $ \omega^{1/2} $ on the whole sphere $ {\mathbb C}P^{1} $. The restrictions of $ g $ to
$ \widetilde{K}_{i} $ give an element of $ \bigoplus_{l_{2}}H^{1/2}\left(\widetilde{K}_{i},\omega^{1/2}\right) $ by Corollary~\ref{cor3.80}. On the
other hand, given $ f $, the restriction of $ g $ on $ \widetilde{K}_{i} $ is defined up to a
section with support on $ K_{i} $. Hence the mollified restriction on $ S_{i} $ is
correctly defined, and it has a norm majorated by some multiple of the
norm of the restriction on $ \widetilde{K}_{i} $. Hence the restriction mapping $ \widetilde{r} $ is bounded
indeed.

To show that $ \operatorname{Im} \widetilde{r} $ is closed, let us construct a left inverse $ l $ this
operator. Then $ \widetilde{r}l $ is going to be a projection on $ \operatorname{Im} \widetilde{r} $, which will prove
the closeness.

Consider the mapping $ \lambda_{i} $ from Section~\ref{s2.25}, associated to the
disks $ K_{i} $, $ \widetilde{K}_{i} $. Then $ g - \lambda_{i}\left(g\right) $
\begin{enumerate}
\item
vanishes outside of $ \widetilde{K}_{i} $;
\item
coincides with $ g $ inside of $ K_{i} $;
\item
has a norm bounded by $ C\cdot\|g|_{\widetilde{K}_{i}}\| $;
\item
depends on values of $ g $ outside of $ K $ only.
\end{enumerate}
Combining all this together, we get
\begin{equation}
F = g+\sum_{i}\left(\lambda_{i}\left(g\right)-g\right)
\notag\end{equation}
which is a half-form of the norm bounded by
\begin{equation}
\|g\|+C'\cdot\left(\sum\|g|_{\widetilde{K}_{i}}\|^{2}\right)^{1/2}= O\left(\|g\|\right).
\notag\end{equation}
The half-form $ F $ is equal to $ g $ (thus to $ f $) outside of $ \bigcup\widetilde{K}_{i} $, and is equal to
0 inside all $ K_{i} $. Moreover, since $ f $ is holomorphic, one can calculate $ \bar{\partial}F|_{\widetilde{K}_{i}} $ as
$ \widetilde{e}_{i}\left(f|_{\partial K_{i}}\right) $.

\begin{proposition} Let $ f $ be a half-form from the generalized Hardy space.
There exists a half-form $ F\in H^{1/2}\left({\mathbb C}P^{1}\right) $ such that:
\begin{enumerate}
\item
$ \|F\| < C\cdot\|f\| $;
\item
$ f=F $ outside of $ \overline{\bigcup\widetilde{K}_{i}} $;
\item
$ F=0 $ inside of $ K_{i} $ for any $ i $;
\item
$ \bar{\partial}F = \sum\widetilde{e}_{i}\left(f|_{\partial K_{i}}\right) $.
\end{enumerate}
\end{proposition}

\begin{proof} The only part which needs proof is the last one. First of all,
the sum in the right-hand side obviously converges in $ H^{-1/2} $, since the
disks $ \widetilde{K}_{i} $ are separated far enough, and norms of the restrictions of $ f $
onto $ \partial K_{i} $ form a sequence in $ l_{2} $. The same arguments show that the
right-hand side is contained in the generalized Sobolev subspace
$ \overset{\,\,{}_\circ}{H}^{-1/2}\left(\bigcup\widetilde{K}_{i}\right) $. Moreover, $ \bar{\partial}\left(F-g\right) $ is contained in the same space, and $ \bar{\partial}f $ is
there by the definition of the generalized Hardy space. Since $ g $ coincides
with $ f $ outside of $ \bigcup K_{i} $, $ \bar{\partial}F\in\overset{\,\,{}_\circ}{H}^{-1/2}\left(\bigcup\widetilde{K}_{i}\right) $.

Thus we know that the difference of the right-hand side and
left-hand side is contained in $ \overset{\,\,{}_\circ}{H}^{-1/2}\left(\bigcup\widetilde{K}_{i}\right) $ and is 0 inside any $ \widetilde{K}_{i} $. On the
other hand, since $ \widetilde{K}_{i} $ are well-separated, $ \overset{\,\,{}_\circ}{H}^{-1/2}\left(\bigcup\widetilde{K}_{i}\right) = \bigoplus_{l_{2}}\overset{\,\,{}_\circ}{H}^{-1/2}\left(\widetilde{K}_{i}\right) $,
what finishes the proof.\end{proof}

\begin{remark} It is obvious that $ f \mapsto F $ is a continuous linear mapping
$ {\mathcal H}\left({\mathbb C}P^{1},\left\{K_{i}\right\}\right) \to H^{1/2}\left({\mathbb C}P^{1},\omega^{1/2}\right) $. Moreover, it is an injection. To show
this, one needs only to prove that $ f|_{\bigcup\widetilde{K}_{i}} $ is determined by $ F $ (and bounded
by the norm of $ F $). Since the disks $ \widetilde{K}_{i} $ are well-separated, it is enough
to show this for one particular disk $ \widetilde{K}_{i} $.

Since $ \widetilde{e}_{i} $ is an injection, $ \bar{\partial}F $ determines $ \widetilde{r}_{i}\left(f\right) $, thus $ \lambda_{i}\left(g\right) $ (by
construction of $ \lambda $), thus $ \left(F-f\right)|_{\widetilde{K}_{i}} $. We obtained

\end{remark}

\begin{corollary} The relation $ f \mapsto F $ is an injection $ {\mathcal H}\left({\mathbb C}P^{1}, \left\{K_{i}\right\}\right) \to
H^{1/2}\left({\mathbb C}P^{1},\omega^{1/2}\right) $. Since $ \bar{\partial}\colon H^{1/2}\left({\mathbb C}P^{1},\omega^{1/2}\right) \to H^{-1/2}\left({\mathbb C}P^{1},\omega^{1/2}\otimes\bar{\omega}\right) $ is an
isomorphism, the relation $ f \mapsto \bar{\partial}F $ gives an injection
$ {\mathcal H}\left({\mathbb C}P^{1}, \left\{K_{i}\right\}\right) \to H^{-1/2}\left({\mathbb C}P^{1},\omega^{1/2}\otimes\bar{\omega}\right) $. Moreover, the last mapping may be
pushed through $ {\mathcal H}\left({\mathbb C}P^{1}, \left\{K_{i}\right\}\right) \to \bigoplus_{l_{2}}H^{-1/2}\left(\widetilde{K}_{i},\omega^{1/2}\otimes\bar{\omega}\right) $. \end{corollary}

Since the operator $ \bar{\partial} $ has no null-space on $ H^{1/2}\left({\mathbb C}P^{1}, \omega^{1/2}\right) $, we obtain

\begin{corollary} $ \bar{\partial}^{-1}\left(\sum\widetilde{e}_{i}\left(f|_{\partial K_{i}}\right)\right) $ equals $ f $ outside of $ \bigcup\widetilde{K}_{i} $, i.e.,
modulo $ \overset{\,\,{}_\circ}{H}^{1/2}\left(\bigcup\widetilde{K}_{i}\right) $. \end{corollary}

Since $ \bar{\partial}F $ is determined by $ \widetilde{r}\left(f\right) $, we found a left inverse to the mapping
$ \widetilde{r} $, thus it is an injection and the image is closed. This finishes the
proof of Theorem~\ref{th4.40}. {}\end{proof}

The next step is to describe the image of the operator in question.
Let $ \gamma_{i}=\partial K_{i} $, $ \gamma=\bigcup\gamma_{i} $. We claim that an element of the image of the mapping
$ \widetilde{r} $ is uniquely determined by the minus-components, and any collection of
minus-components of bounded norm is possible.

\begin{proposition} \label{prop5.28}\myLabel{prop5.28}\relax  Consider a decomposition of the space $ L_{2}\left(\gamma_{i}, \Omega^{1/2}\right) =
L_{2}^{+}\left(\gamma_{i}\right) \oplus L_{2}^{-}\left(\gamma_{i}\right) $ into subspaces of half-forms which may be holomorphically
extended
inside the circle and outside the circle. The image of $ \widetilde{r} $ consists of
sequences $ \left(f_{i}^{\pm}\right) $ such that
\begin{equation}
f_{i}^{+}=\sum_{j\not=i}{\bold K}f_{j}^{-}.
\notag\end{equation}
Here $ {\bold K} $ is the Hilbert operator, i.e., the operator with Cauchy kernel. \end{proposition}

\begin{proof} First of all, note that the proof of Theorem~\ref{th4.40}
together with Propositions~\ref{prop4.12} and~\ref{prop5.16} shows that knowing
$ \widetilde{r}\left(f\right) $ one can write $ f $ by an explicit formula outside of $ \overline{\bigcup\widetilde{K}_{i}} $. On the
other hand, it is easy to see that if we throw away the restriction that
we want our operators to be continuous in $ H^{s} $-topology, one can
reconstruct $ f $ outside of $ \bigcup\bar{K}_{i} $.

Indeed, since the curve $ \gamma $ has a finite length, restriction of any
smooth section $ \alpha $ of $ \omega^{1/2} $ to $ \gamma $ is in $ L_{2}\left(\gamma\right) $, thus $ \int\alpha|_{S}\widetilde{r}\left(f\right) $ is correctly
defined. By duality, this means that extension $ e\left(\widetilde{r}\left(f\right)\right) $ of $ \widetilde{r}\left(f\right) $ to $ {\mathbb C}P^{1} $ by
$ \delta $-function is a correctly defined generalized section of $ \omega^{1/2}\otimes\bar{\omega} $. (Note
that we consider $ e $, not $ \widetilde{e}! $)

\begin{lemma} Half-form $ f $ coincides with $ \bar{\partial}^{-1}e\left(\widetilde{r}\left(f\right)\right) $ outside of $ \bigcup\bar{K}_{i} $. \end{lemma}

\begin{proof} Indeed, by the construction of $ \widetilde{e}_{i} $, $ \bar{\partial}^{-1} \widetilde{e}_{i}\left(\widetilde{r}_{i}\left(f\right)\right) $ coincides with
$ \bar{\partial}^{-1}e\left(\widetilde{r}_{i}\left(f\right)\right) $ outside of $ \widetilde{K}_{i}\smallsetminus K_{i} $, thus $ \bar{\partial}^{-1}e\left(\widetilde{r}\left(f\right)\right) $ coincides with $ f $ outside of
$ \bigcup\widetilde{K}_{i} $. On the other hand, $ \bar{\partial}^{-1}e\left(\widetilde{r}\left(f\right)\right) $ is holomorphic outside of $ \bar{\gamma} $, thus
this equality can be extended up to $ \gamma $. \end{proof}

Consider now $ f=\bar{\partial}^{-1}e\left(\widetilde{r}\left(f\right)\right) $ near $ \gamma_{i} $. Breaking $ \widetilde{r}\left(f\right) $ into two
components, $ \left(\widetilde{r}\left(f\right)-\widetilde{r}_{i}\left(f\right)\right) $ and $ \widetilde{r}_{i}\left(f\right) $, we conclude that
\begin{equation}
f=\bar{\partial}^{-1}e\left(\widetilde{r}\left(f\right)-\widetilde{r}_{i}\left(f\right)\right) + \bar{\partial}^{-1}e\left(\widetilde{r}_{i}\left(f\right)\right).
\notag\end{equation}
The first summand is holomorphic inside $ \widetilde{K}_{i} $ and coincides with $ \sum_{j\not=i}{\bold K}f_{j}^{-} $,
since $ {\bold K} $ kills $ f_{j}^{+} $. The second summand is holomorphic outside of $ \partial K_{i} $, so
its $ + $-part vanishes on $ \partial K_{i} $, which finishes the proof of Proposition
~\ref{prop5.28}. {}\end{proof}

\subsection{Gluing conditions }\label{s5.30}\myLabel{s5.30}\relax  We continue using notations of Section
~\ref{s5.10}. Suppose that the set $ I $ has an involution $ '\colon I\to I $ which
interchanges two subsets $ I_{+} $ and $ I_{+}' $, $ I=I_{+}\amalg I_{+}' $. Thus all the disks $ K_{i} $ are
divided into pairs $ \left(K_{i},K_{i'}\right) $, $ i\in I_{+} $. Fix a fraction-linear identification $ \varphi_{i} $
of $ \partial K_{i'} $ and $ \partial K_{i} $ which reverses the orientation of the circles, $ \varphi_{i'}=\varphi_{i}^{-1} $.
Let $ \varphi_{i} $ identifies the boundary of the disk $ \widetilde{K}_{i'} $ with the boundary of $ \overset{\,\,{}_\circ}{K}_{i} $, and
$ \partial\widetilde{K}_{i} $ with $ \partial\overset{\,\,{}_\circ}{K}_{i'} $, thus $ \overset{\,\,{}_\circ}{K}_{i}\subset K_{i}\subset\widetilde{K}_{i} $.

Let $ R_{i} $ be the annulus between $ \widetilde{K}_{i} $ and $ \overset{\,\,{}_\circ}{K}_{i} $. The mapping $ \varphi_{i} $ identifies
$ R_{i} $ with $ R_{i'} $. Let $ S $ be the part of $ {\mathbb C}P^{1} $ which lies outside of all the disks
$ \overset{\,\,{}_\circ}{K}_{i} $. Glue the annuli $ R_{i}\subset S $ with $ R_{i'}\subset S $ using $ \varphi_{i} $.

\begin{definition} A {\em model space\/} $ \bar{M} $ is the set obtained from $ S $ by identifying the
annuli $ R_{i} $ and $ R_{i'} $ using $ \varphi_{i} $. \end{definition}

Note that $ \bar{M} $ consists of two parts: a smooth manifold $ M $ which is the
image of $ S\smallsetminus\overline{\bigcup\overset{\,\,{}_\circ}{K}_{\bullet}} $, and the rest, which one should consider as ``infinity''
$ M_{\infty} $ of the manifold $ M $ (compare with Section~\ref{s0.10}). Unfortunately, the
topology on $ \bar{M} $ in neighborhood of infinity is not suitable for studying
the Riemann--Roch theorem, so we will not consider it in this paper.

Note that Section~\ref{s4.95} suggests a different topology on $ \bar{M} $. We will
use some features of this topology when we discuss a mapping into
universal Grassmannian.

\subsection{Strong sections and duality }\label{s6.50}\myLabel{s6.50}\relax  In the notations of the previous
section
consider now holomorphic functions $ \psi_{i} $ defined in $ R_{i} $. Suppose that $ \psi_{i} $ are
nowhere 0, and $ \psi_{i'}\circ\varphi_{i}=\psi_{i}^{-1} $.

Define $ {\mathcal L} $ to be a sheaf on $ \bar{M} $ associated with gluing conditions $ \psi_{i} $,
i.e., for $ U\subset\bar{M} $ the section of $ {\mathcal L} $ on $ U $ is a function $ f $ on $ \widetilde{U}\subset{\mathbb C}P^{1} $ such that
$ f\left(\varphi_{i}\left(x\right)\right)=\psi_{i}\left(x\right)f\left(x\right) $ whenever both sides have sense (here $ \widetilde{U} $ is an
appropriate covering subset of $ {\mathbb C}P^{1} $).

Let $ {\mathcal L}^{-1} $ be the sheaf associated with gluing conditions $ \psi_{i}^{-1} $.
Similarly define the tensor product of two sheaves defined via gluing
conditions.

The identifications $ \varphi_{i} $ fixes an identification $ \varphi_{i}^{*} $ of half-forms on
$ R_{i} $ and on $ R_{i'} $ up to a sign. Choose this sign for all $ i\in I_{+}. $\footnote{This corresponds to picking a representative of $ \varphi_{i} $ in $ 2 $-sheet cover
$ \operatorname{SL}\left(2,{\mathbb C}\right) \to \operatorname{PGL}\left(2,{\mathbb C}\right) $.} Let $ \omega^{1/2}\otimes{\mathcal L} $ be
the sheaf on $ \bar{M} $ consisting of half-forms on $ S $ such that
\begin{equation}
\varphi_{i}^{*}\left(\alpha|_{R_{i'}}\right) = \psi_{i}\cdot\alpha|_{R_{i}},\qquad i\in I.
\label{equ5.31}\end{equation}\myLabel{equ5.31,}\relax 
Similarly define $ \omega^{1/2}\otimes{\mathcal L}\otimes\bar{\omega} $.

To write an analogue of~\eqref{equ5.31} for sections in Sobolev spaces,
consider a subannulus $ \overset{\,\,{}_\circ}{R}_{i}\subset R_{i} $ such that $ \overset{\,\,{}_\circ}{R}_{i'}=\varphi_{i}\left(\overset{\,\,{}_\circ}{R}_{i}\right) $. We suppose that one can
find numbers $ C,D>0 $ and a sequence $ \left(s_{i}\right) $, $ i\in I $, $ |s_{i}|+C+D<1 $, such that the
annulus $ \overset{\,\,{}_\circ}{R}_{i} $ is described as $ \left(\left(s_{i}-C\right)\varepsilon,\left(s_{i}+C\right)\varepsilon\right)\times S^{1} $ in the conformal
coordinate system such that $ R_{i} $ is $ \left(-\varepsilon,\varepsilon\right)\times S^{1} $. Let $ K'_{i} $ be the disk bounded
by the inner boundary of $ \overset{\,\,{}_\circ}{R}_{i} $.

Consider a subset $ \overset{\,\,{}_\circ}{S}\subset S $, $ \overset{\,\,{}_\circ}{S}={\mathbb C}P^{1}\smallsetminus\bigcup_{i}K'_{i} $. One obtains the same manifold $ \bar{M} $
by gluing $ \overset{\,\,{}_\circ}{R}_{i}\subset\overset{\,\,{}_\circ}{S} $ as by gluing $ R\subset S $, but in what follows it will be more
convenient to have a choice of annuli $ \overset{\,\,{}_\circ}{R}_{i} $, and have them separated from
boundary of $ R_{i} $.

\begin{definition} Let $ H^{s}\left(\bar{M},\omega^{1/2}\otimes{\mathcal L}\right) $ be the subspace of $ H^{s}\left({\mathbb C}P^{1}\smallsetminus\bigcup_{i}K'_{i},\omega^{1/2}\right) $
consisting of sections which satisfy the gluing conditions~\eqref{equ5.31},
similarly define $ H^{s}\left(\bar{M},\omega^{1/2}\otimes{\mathcal L}\otimes\bar{\omega}\right) $. \end{definition}

Dually,

\begin{definition} Let $ \overset{\,\,{}_\circ}{H}^{s}\left(\bar{M},\omega^{1/2}\otimes{\mathcal L}\right) $ be the quotient of $ \overset{\,\,{}_\circ}{H}^{s}\left({\mathbb C}P^{1}\smallsetminus\bigcup_{i}K'_{i}, \omega^{1/2}\right) $ by
the subspace of $ \overset{\,\,{}_\circ}{H}^{s}\left(\bigcup_{i}\overset{\,\,{}_\circ}{R}_{i}, \omega^{1/2}\right) $ consisting of sections which satisfy the
gluing conditions
\begin{equation}
\varphi_{i}^{*}\left(\alpha|_{\overset{\,\,{}_\circ}{R}_{i'}}\right) = -\psi_{i}\cdot\alpha|_{\overset{\,\,{}_\circ}{R}_{i}},\qquad i\in I,
\label{equ5.32}\end{equation}\myLabel{equ5.32,}\relax 
similarly define $ \overset{\,\,{}_\circ}{H}^{s}\left(\bar{M},\omega^{1/2}\otimes{\mathcal L}\otimes\bar{\omega}\right) $. \end{definition}

Indeed, these definitions are dual due to

\begin{lemma} \label{lm5.22}\myLabel{lm5.22}\relax  If the disks $ K_{i} $ are well-separated, the spaces
$ H^{1/2}\left(\bar{M},\omega^{1/2}\otimes{\mathcal L}\right) $ and $ \overset{\,\,{}_\circ}{H}^{-1/2}\left(\bar{M},\omega^{1/2}\otimes{\mathcal L}^{-1}\otimes\bar{\omega}\right) $ are mutually dual w.r.t. the
pairing $ \int\alpha\beta $. \end{lemma}

\begin{proof} The spaces $ H^{s}\left({\mathbb C}P^{1}\smallsetminus\bigcup_{i}K'_{i},\omega^{1/2}\otimes{\mathcal L}\right) $ and $ \overset{\,\,{}_\circ}{H}^{-s}\left({\mathbb C}P^{1}\smallsetminus\bigcup_{i}K'_{i},\omega^{1/2}\otimes{\mathcal L}^{-1}\otimes\bar{\omega}\right) $
are mutually dual w.r.t. this pairing by definition. What remains to
prove is that the orthogonal complement to the subspace of
$ H^{s}\left({\mathbb C}P^{1}\smallsetminus\bigcup_{i}K'_{i},\omega^{1/2}\otimes{\mathcal L}\otimes\bar{\omega}\right) $ given by~\eqref{equ5.31} is given by
\begin{equation}
\varphi_{i}^{*}\left(\alpha|_{R_{i'}}\right) = -\psi_{i}^{-1}\cdot\alpha|_{R_{i}},\qquad i\in I,
\label{equ5.33}\end{equation}\myLabel{equ5.33,}\relax 
in $ \overset{\,\,{}_\circ}{H}^{-s}\left({\mathbb C}P^{1}\smallsetminus\bigcup_{i}K'_{i},\omega^{1/2}\otimes{\mathcal L}^{-1}\otimes\bar{\omega}\right) $, provided $ s=1/2 $. This statement if obvious
for any fixed $ i\in I $.

On the other hand, solutions to~\eqref{equ5.33} form a direct sum over $ i $, and
by the second part of Corollary~\ref{cor3.80}, solutions to~\eqref{equ5.31} form a
``direct'' intersection (i.e., an intersection of subspaces with almost
orthogonal complements).\end{proof}

\begin{lemma} \label{lm5.25}\myLabel{lm5.25}\relax  Suppose that the disks $ K_{i} $ are well-separated. The natural
mapping
\begin{equation}
\overset{\,\,{}_\circ}{H}^{-1/2}\left(\bar{M},\omega^{1/2}\otimes{\mathcal L}\otimes\bar{\omega}\right) \to H^{-1/2}\left(\bar{M},\omega^{1/2}\otimes{\mathcal L}\otimes\bar{\omega}\right)
\label{equ5.36}\end{equation}\myLabel{equ5.36,}\relax 
is an isomorphism. Dually,
\begin{equation}
\overset{\,\,{}_\circ}{H}^{1/2}\left(\bar{M},\omega^{1/2}\otimes{\mathcal L}\right) \to H^{1/2}\left(\bar{M},\omega^{1/2}\otimes{\mathcal L}\right)
\notag\end{equation}
is an isomorphism. \end{lemma}

\begin{proof} We may suppose that all the rings $ \overset{\,\,{}_\circ}{R}_{i} $ have conformal
distance between boundaries greater than $ 2C\varepsilon $, $ C>0 $. Consider a function
$ \sigma\left(x\right) $, $ x\in\left(-C\varepsilon,C\varepsilon\right) $, such that $ \sigma\left(x\right)=0 $ near the left end, $ \sigma\left(x\right)+\sigma\left(1-x\right)=1 $. This
gives a cut-off function in all the rings $ \overset{\,\,{}_\circ}{R}_{i} $, and we may extend it by 0
into all $ K_{i}' $, and by 1 into $ S $. We obtain a function on $ {\mathbb C}P^{1} $.

By the results of Section~\ref{s3.8} and Amplification~\ref{amp2.65} the
multiplication by this function is a bounded operator in $ H^{1/2} $. Clearly,
this operator provides an inverse to~\eqref{equ5.36}. \end{proof}

The operator $ \bar{\partial} $ gives mappings
\begin{align} \overset{\,\,{}_\circ}{H}^{s}\left(\bar{M},\omega^{1/2}\otimes{\mathcal L}\right) & \to \overset{\,\,{}_\circ}{H}^{s-1}\left(\bar{M},\omega^{1/2}\otimes{\mathcal L}\otimes\bar{\omega}\right),
\notag\\
H^{s}\left(\bar{M},\omega^{1/2}\otimes{\mathcal L}\right) & \to H^{s-1}\left(\bar{M},\omega^{1/2}\otimes{\mathcal L}\otimes\bar{\omega}\right),
\notag\end{align}
which by the previous lemma induce mappings
\begin{align} \overset{\,\,{}_\circ}{H}^{1/2}\left(\bar{M},\omega^{1/2}\otimes{\mathcal L}\right) & \xrightarrow[]{\bar{\partial}} H^{-1/2}\left(\bar{M},\omega^{1/2}\otimes{\mathcal L}\otimes\bar{\omega}\right),
\notag\\
H^{1/2}\left(\bar{M},\omega^{1/2}\otimes{\mathcal L}\right) & \xrightarrow[]{\bar{\partial}} \overset{\,\,{}_\circ}{H}^{-1/2}\left(\bar{M},\omega^{1/2}\otimes{\mathcal L}\otimes\bar{\omega}\right) 
\notag\end{align}
(obviously, it does not matter which of isomorphisms of Lemma~\ref{lm5.25} we
use to obtain these mappings).

\begin{definition} Define the space $ \Gamma_{\text{strong}}\left(\bar{M},\omega^{1/2}\otimes{\mathcal L}\right) $ of {\em strong global
holomorphic sections\/} of $ \omega^{1/2}\otimes{\mathcal L} $ as
\begin{equation}
\Gamma_{\text{strong}}\left(\bar{M},\omega^{1/2}\otimes{\mathcal L}\right) = \operatorname{Ker} \left(H^{1/2}\left(\bar{M},\omega^{1/2}\otimes{\mathcal L}\right) \xrightarrow[]{\bar{\partial}} H^{-1/2}\left(\bar{M},\omega^{1/2}\otimes{\mathcal L}\otimes\bar{\omega}\right)\right).
\notag\end{equation}
\end{definition}

Note that the conditions on strong global sections $ \alpha $ are: they are
holomorphic sections of $ \omega^{1/2}\otimes{\mathcal L}|_{M} $ (since $ \bar{\partial}\alpha $ vanishes on $ M $), they do not
grow very quick near $ M_{\infty} $ (since they belong to $ H^{1/2} $), and they have no
residue on $ M_{\infty} $ (since $ \bar{\partial}\alpha $ vanishes on $ M_{\infty} $).

\subsection{Weak sections }\label{s6.60}\myLabel{s6.60}\relax  The space $ \Gamma_{\text{strong}} $ from the last section is
good
for studying the duality conditions, but it is not suitable for for
description of global sections via boundary conditions.

\begin{definition} Define the space $ \Gamma_{\text{weak}}\left(\bar{M},\omega^{1/2}\otimes{\mathcal L}\right) $ of {\em weak global holomorphic
sections\/} of $ \omega^{1/2}\otimes{\mathcal L} $ as forms from the generalized Hardy space
$ \alpha\in{\mathcal H}\left({\mathbb C}P^{1},\left\{K_{i}\right\}\right) $ such that the (mollified) restrictions on the circles
$ \partial K_{i} $ satisfy the gluing conditions. In other words,
\begin{equation}
\Gamma_{\text{weak}}\left(\bar{M},\omega^{1/2}\otimes{\mathcal L}\right) = \left\{\alpha\in{\mathcal H}\left({\mathbb C}P^{1},\left\{K_{i}\right\}\right) \mid \psi_{i}\widetilde{r}_{i}\left(\alpha\right)=\varphi_{i}^{*}\left(\widetilde{r}_{i'}\left(\alpha\right)\right) \right\}.
\notag\end{equation}
\end{definition}

There is a natural mapping
\begin{equation}
\Gamma_{\text{strong}} \to \Gamma_{\text{weak}}.
\notag\end{equation}
\begin{theorem} \label{th5.31}\myLabel{th5.31}\relax  Suppose that the disks $ K_{i} $ are well-separated, and for
some $ A>1 $ either $ |\psi_{i}\left(z\right)|<A $, or $ |\psi_{i}\left(z\right)|>1/A $ for any $ i $ and $ z\in R_{i} $. Then the
above mapping is an isomorphism for an appropriate choice of annuli $ \overset{\,\,{}_\circ}{R}_{i} $. \end{theorem}

\begin{proof} Indeed, any weak section $ \alpha $ is a holomorphic form inside
$ {\mathbb C}P^{1}\smallsetminus\overline{\bigcup_{i}K_{i}} $. Consider $ \psi_{i}^{-1}\cdot\varphi_{i}^{*}\left(\alpha\right) $, it is a holomorphic form inside
$ K_{i}\smallsetminus\overset{\,\,{}_\circ}{K}_{i} $. Since $ \alpha $ satisfies gluing conditions, $ \alpha|_{\widetilde{K}_{i}\smallsetminus K_{i}} $ has the same Laurent
coefficients as $ \psi_{i}^{-1}\cdot\varphi_{i}^{*}\left(\alpha\right) $ (we used compatibility of $ \widetilde{r} $ with $ r $ on
holomorphic forms), thus these two forms are restrictions of the same
holomorphic form defined inside $ \widetilde{K}_{i}\smallsetminus\overset{\,\,{}_\circ}{K}_{i} $. We see that $ \alpha $ can be extended (as
a holomorphic form) into $ {\mathbb C}P^{1}\smallsetminus\overline{\bigcup_{i}\overset{\,\,{}_\circ}{K}_{i}} $, and this holomorphic form
satisfies the gluing conditions~\eqref{equ5.31}.

What remains to prove is that we can bound the norm of this
extension. On the other hand, this is a local statement, since one can
represent $ \alpha=\alpha_{1}+\alpha_{2} $, and $ \alpha=0 $ inside an appropriate circle concentric with
$ \widetilde{K}_{i} $, $ \alpha_{2}=0 $ outside of $ \bigcup\widetilde{K}_{i} $. Thus the only thing we need to prove is that
inside $ R_{i} $ the form $ \alpha $ can be extended across $ \partial K_{i} $ without increasing its
norm too much.

Take $ C=\frac{1}{4} $, let $ \varepsilon_{i}=C $ if $ |\psi_{i}\left(z\right)|<A $, $ \varepsilon_{i}=-C $ otherwise. With this
choice of $ \overset{\,\,{}_\circ}{R}_{i} $ we know that the $ H^{1/2} $-norms of $ \alpha|_{\widetilde{K}_{i}\smallsetminus K_{i}} $ and $ \alpha|_{K_{i}\smallsetminus K'_{i}} $ are
bounded by the norms of $ \alpha|_{\widetilde{K}_{i}\smallsetminus K_{i}} $ and $ \alpha|_{\widetilde{K}_{i'}\smallsetminus K_{i'}} $. Now the theorem becomes a
corollary of the following

\begin{lemma} Fix two numbers $ A>a>0 $. Consider an annulus $ R $ with concentric
boundaries and conformal distance between boundaries between $ A $ and $ a $.
Let $ V_{R} $ be the space of holomorphic half-forms in $ R $ which belong to
$ H^{1/2}\left(R\right) $. Then the mapping of taking boundary value
\begin{equation}
b\colon H^{1/2}\left(R,\omega^{1/2}\right) \to L_{2}\left(\partial R,\Omega^{1/2}\right)\colon \alpha \mapsto \alpha|_{\partial R}
\notag\end{equation}
is an invertible mapping to its (closed) image, and the norms of this
mapping and its inverse are bounded by numbers depending on $ A $ and $ a $
only.

\end{lemma}

This lemma is a variation of what we did in Section~\ref{s2.70}, with a
disk substituted with an annulus. It can be proven in the same way as the
case of a disk.

This finishes the proof of the theorem. \end{proof}

\begin{remark} To complement the notion of weak holomorphic sections, let us
define the spaces of ``weak'' $ H^{\pm1/2} $-section in such a way that
\begin{equation}
\Gamma_{\text{weak}}\left(\bar{M}, \omega^{1/2}\otimes{\mathcal L}\right) = \operatorname{Ker}\left(H_{\text{weak}}^{1/2}\left(\bar{M},\omega^{1/2}\otimes{\mathcal L}\right) \xrightarrow[]{\bar{\partial}} H_{\text{weak}}^{-1/2}\left(\bar{M},\omega^{1/2}\otimes{\mathcal L}\otimes\bar{\omega}\right) \right):
\notag\end{equation}
\end{remark}

\begin{definition} Let $ H_{\text{weak}}^{1/2}\left(\bar{M},\omega^{1/2}\otimes{\mathcal L}\right) $ be the subspace of $ H^{1/2}\left({\mathbb C}P^{1}\smallsetminus\bigcup_{i}K_{i},\omega^{1/2}\right) $
consisting of sections $ \alpha $ which satisfy the gluing conditions:
\begin{equation}
\psi_{i}\widetilde{r}_{i}\left(\alpha\right)=\varphi_{i}^{*}\left(\widetilde{r}_{i'}\left(\alpha\right)\right),\qquad i\in I.
\notag\end{equation}
Let $ H_{\text{weak}}^{-1/2}\left(\bar{M},\omega^{1/2}\otimes{\mathcal L}\otimes\bar{\omega}\right) = H^{-1/2}\left({\mathbb C}P^{1}\smallsetminus\bigcup_{i}K_{i},\omega^{1/2}\otimes\bar{\omega}\right) $.

\end{definition}

\subsection{Finite-degree bundles }\label{s4.90}\myLabel{s4.90}\relax  Let $ d_{i} $ be the {\em index\/} of $ \psi_{i} $, i.e., the
degree of the mapping $ \arg  \psi_{i}\colon K_{i} \to S^{1} $. We say that the collection of
gluing data $ \left\{\psi_{i}\right\} $ is of {\em finite degree\/} if $ d_{i}=0 $ for all but the finite
number of indices $ i\in I $. The collection $ \left\{\psi_{i}\right\} $ is {\em semibounded\/} if for some
fixed number $ C $ and any $ i $ either $ |\psi_{i}\left(z\right)| $ or $ |\psi\left(z\right)|_{i}^{-1} $ is bounded by $ C $ if
$ z\in R_{i} $. The {\em degree\/} of the finite-degree collection $ \left\{\psi_{i}\right\} $ is the sum
\begin{equation}
\sum_{i\in I_{+}}d_{i} = \frac{1}{2}\sum_{i\in I}d_{i}.
\notag\end{equation}
\begin{definition} The {\em degree\/} of a bundle $ {\mathcal L} $ defined by gluing conditions $ \psi_{i} $ is
the degree of the collection $ \left\{\psi_{i}\right\} $. \end{definition}

\begin{amplification} In what follows we are going to use the following
generalization of these constructions: we allow a substitution of a
finite number of simply-connected domains with smooth boundaries instead
of disks $ K_{i} $. (For such an $ i $ one should substitute any bigger domain
instead of $ \widetilde{K}_{i} $.) However, we still require that the identifications $ \varphi_{i} $ are
fraction-linear. \end{amplification}

\subsection{Stratification of infinity }\label{s4.95}\myLabel{s4.95}\relax  Consider a smooth curve $ \gamma $ on $ {\mathbb C}P^{1} $ of
finite length. Suppose that a metric on $ {\mathbb C}P^{1} $ is fixed, and $ z\notin\gamma $. Let
$ \rho_{k}\left(z,\gamma\right)=\|\operatorname{dist}\left(z,y\right)^{-k}\|_{L_{2}\left(\gamma\right)} $, here $ y\in\gamma $.

\begin{lemma} $ \rho_{k}\left(x,\gamma\right) $ is a semicontinuous function of $ x $, thus
\begin{equation}
D_{k,R}=\left\{z\in{\mathbb C}P^{1} \mid \rho_{k}\left(x,\gamma\right)\leq R \right\}
\notag\end{equation}
is a compact subset of $ {\mathbb C}P^{1}\smallsetminus\gamma $. \end{lemma}

Let $ D_{k}=\bigcup_{R\in{\mathbb R}}D_{k,R} $. It is a subset of $ {\mathbb C}P^{1}\smallsetminus\gamma $, moreover, $ {\mathbb C}P^{1}\smallsetminus\bar{\gamma}\subset D_{k} $.
Let $ \overset{\,\,{}_\circ}{D}_{k,R}=D_{k,R}\cap\left({\mathbb C}P^{1}\smallsetminus\bar{\gamma}\right) $, $ \overset{\,\,{}_\circ}{D}_{k}=\bigcup_{R\in{\mathbb R}}\overline{\overset{\,\,{}_\circ}{D}_{k,R}} $. Since length of $ \gamma $ is finite,
$ \overset{\,\,{}_\circ}{D}_{0}={\mathbb C}P^{1} $.

\begin{definition} Define a filtration of $ {\mathcal K}=\bar{\gamma}\smallsetminus\gamma $ by $ {\mathcal K}^{\left(k\right)}={\mathcal K}\cap\overset{\,\,{}_\circ}{D}_{k+1} $, $ k\geq-1 $. Let
$ {\mathcal K}^{\left(\infty\right)}=\bigcap{\mathcal K}^{\left(k\right)} $. \end{definition}

Suppose that $ \gamma $ is the boundary of a well-separated family of
circles, $ \gamma=\bigcup\partial K_{i} $.

\begin{theorem} Let $ z_{0}\in{\mathcal K}^{\left(k\right)} $, $ z $ be a coordinate system near $ z_{0} $, $ f\in{\mathcal H}\left({\mathbb C}P^{1},\left\{K_{i}\right\}\right) $.
Then $ f $ has an asymptotic decomposition
\begin{equation}
f\left(z\right)dz^{-1/2} = f_{0}+f_{1}\left(z-z_{0}\right)+f_{2}\left(z-z_{0}\right)^{2}+\dots +f_{k}\left(z-z_{0}\right)^{k}+o\left(z^{k}\right)
\notag\end{equation}
when $ z \to z_{0} $ along $ \overset{\,\,{}_\circ}{D}_{k,R} $ for an appropriate $ R\gg0 $. \end{theorem}

\begin{proof} The half-form $ f $ is holomorphic inside $ {\mathbb C}P^{1}\smallsetminus\overline{\bigcup K_{i}} $. The Cauchy
formula show that inside $ \overset{\,\,{}_\circ}{D}_{k,R} $ the derivatives $ f^{\left(l\right)} $, $ l\leq k $, are bounded.
Moreover, these derivatives are given by some integrals along $ \gamma $, and
these integrals remain well-defined in $ D_{k,R} $ as well. Let $ f_{l} $ be the values
of these integrals in $ z_{0}\in D_{k,R_{1}} $ (here $ R_{1}\gg0 $).

Since $ z_{0}\in\overset{\,\,{}_\circ}{D}_{k+1} $, it is in the closure of $ \overset{\,\,{}_\circ}{D}_{k+1,R} $ for an appropriate $ R $. It
is easy to see that $ f^{\left(k\right)}\left(z\right) $ has a limit $ f_{k} $ when $ z \to z_{0} $ along $ \overset{\,\,{}_\circ}{D}_{k+1,R} $. Same
is true for $ f^{\left(l\right)}\left(z\right) $, $ l\leq k $.

Consider the integral for
\begin{equation}
\frac{f\left(z\right)dz^{-1/2} - \left(f_{0}+f_{1}\left(z-z_{0}\right)+f_{2}\left(z-z_{0}\right)^{2}+\dots +f_{k}\left(z-z_{0}\right)^{k}\right)}{\left(z-z_{0}\right)^{k}}.
\notag\end{equation}
It is
\begin{equation}
\int_{\gamma}K_{k}\left(z,z_{0},\zeta\right)f\left(\zeta\right),
\notag\end{equation}
here
\begin{align} K\left(z,\zeta\right) & =\frac{d\zeta^{1/2}}{\zeta-z},
\notag\\
K_{k}\left(z,z_{0},\zeta\right) & =\frac{K\left(z,\zeta\right)-\sum_{l=0}^{k}\frac{d^{k}K}{dz^{k}}|_{z=z_{0}}\frac{\left(z-z_{0}\right)^{k}}{k!}}{\left(z-z_{0}\right)^{k}}.
\notag\\
& = \frac{z-z_{0}}{\left(\zeta-z_{0}\right)^{k+1}}K\left(z,\zeta\right) = \frac{1}{\left(\zeta-z_{0}\right)^{k}}\left(K\left(z,\zeta\right)-K\left(z_{0},\zeta\right)\right).
\notag\end{align}
We need to show that this integral goes to 0 when $ z \to z_{0} $ along $ \overset{\,\,{}_\circ}{D}_{k+1,R} $.
Consider a function $ \rho\left(r\right) $ such that $ \lim _{r\to0}\rho\left(r\right)=0 $, $ \lim _{r\to0}\rho\left(r\right)/r=\infty $. Here
$ \gamma_{r} $ is the intersection of $ \gamma $ with the disk of radius $ r $ about $ z_{0} $, instead
of $ \widetilde{r}_{\gamma}\left(f\right)\in L_{2}\left(\gamma\right) $ we write just $ f $.

Break the integral into three parts: two (which we do not want to
separate yet) along $ \gamma_{\rho\left(|z-z_{0}|\right)} $, the other along $ \gamma\smallsetminus\gamma_{\rho\left(|z-z_{0}|\right)} $. Since
\begin{equation}
K_{k}= \frac{z-z_{0}}{\zeta-z_{0}}\frac{K\left(z,\zeta\right)}{\left(\zeta-z_{0}\right)^{k}} = o\left(\frac{K\left(z_{0},\zeta\right)}{|\zeta-z_{0}|^{k}}\right)
\notag\end{equation}
along the second part of $ \gamma $ if $ |z-z_{0}| \to $ 0, and since
\begin{equation}
\frac{K\left(z_{0},\zeta\right)}{|z_{0}-\zeta|^{k}}\in L_{2}\left(\gamma\right)
\notag\end{equation}
as a function of $ \zeta $, it is enough to show that the first two part of the
integral go to 0. Since $ \lim _{r\to0} \int_{\gamma_{\rho\left(r\right)}}\left|f\left(\zeta\right)\right|^{2}=0 $, we need to show only
that the $ L_{2} $-norm of $ K_{k}\left(z,z_{0},\zeta\right) $ is bounded when $ \zeta\in\gamma_{\rho\left(|z-z_{0}|\right)} $.

Subdivide $ \gamma_{\rho\left(|z-z_{0}|\right)} $ once more: into $ \gamma_{\varepsilon|z-z_{0}|} $ and
$ \gamma_{\rho\left(|z-z_{0}|\right)}\smallsetminus\gamma_{\varepsilon|z-z_{0}|} $. Here $ \varepsilon\ll1 $. Since on the second part
$ K\left(z_{0},\zeta\right)=O\left(K\left(z,\zeta\right)\right) $, we see that $ K_{k}\left(z,z_{0},\zeta\right) = O\left(\frac{K\left(z,\zeta\right)}{\left(\zeta-z\right)^{k}}\right) $, thus $ K_{k} $
has a bounded $ L_{2} $-norm. We conclude that the integral along the second
part goes to 0 when $ z \to z_{0} $. On the first part $ K\left(z,\zeta\right)=O\left(K\left(z_{0},\zeta\right)\right) $, thus
\begin{equation}
K_{k}\left(z,z_{0},\zeta\right) = \frac{1}{\left(\zeta-z_{0}\right)^{k}}\left(K\left(z,\zeta\right)-K\left(z_{0},\zeta\right)\right) =
O\left(\frac{K\left(z_{0},\zeta\right)}{\left(\zeta-z_{0}\right)^{k}}\right)=O\left(\frac{1}{\left(z-z_{0}\right)^{k+1}}\right),
\notag\end{equation}
so it has a bounded $ L_{2} $-norm as well. \end{proof}

\begin{definition} Consider a model $ \bar{M} $ of a curve. Say that a point $ z\in M_{\infty} $ is of
smoothness $ C^{k} $, if $ z\in{\mathcal K}^{\left(k\right)} $, here $ {\mathcal K}=M_{\infty} $. \end{definition}

\begin{remark} Note that the stratification we used is related to the
following inclusion of $ M $ into $ L_{2}\left(\gamma,\Omega^{1/2}\otimes\mu\right) $:
\begin{equation}
z \mapsto \widetilde{r}_{\gamma}\left(K_{z}\right),
\notag\end{equation}
here $ K_{z}\left(\zeta\right)=\frac{d\zeta^{1/2}}{\zeta-z} $ is the Cauchy kernel. In other words,
\begin{equation}
z \mapsto \frac{d\zeta^{1/2}}{\zeta-z}|_{\gamma}
\notag\end{equation}
after a choice of coordinate $ \zeta $ on $ {\mathbb C}P^{1} $. The points of smoothness $ C^{0} $
correspond to limit points of this inclusion, the points of
smoothness $ C^{k} $ correspond to limit points of $ k $-jets continuation of this
mapping. \end{remark}

\section{Riemann--Roch theorems }

\subsection{Abstract Riemann--Roch theorem } We say that two vector subspaces $ V_{1} $,
$ V_{2} $ of a topological vector space $ H $ {\em satisfy the Riemann\/}--{\em Roch theorem\/} if
$ \dim  V_{1}\cap V_{2}<\infty $, and $ \operatorname{codim}\left(\overline{V_{1}+V_{2}}\right)<\infty $. We call the number
\begin{equation}
\dim  V_{1}\cap V_{2} - \operatorname{codim}\left(\overline{V_{1}+V_{2}}\right)
\notag\end{equation}
the {\em index\/} of two subspaces. If $ V_{1}+V_{2}=\overline{V_{1}+V_{2}} $, we say that $ V_{1},V_{2} $ satisfy
the {\em strong form\/} of the theorem.

\begin{remark} Note that if $ V_{1} $, $ V_{2} $ satisfy the strong form of the theorem,
then the natural mapping $ V_{1} \to V/V_{2} $ is a Fredholm mapping with the index
being the index of $ V_{1} $, $ V_{2} $. If $ V_{1} $, $ V_{2} $ satisfy the weak form of the
theorem, then this mapping is a continuous mapping $ p $ with $ \dim  \operatorname{Ker} p - \dim 
\operatorname{Coker} p $ being the index of $ V_{1},V_{2} $. Here $ \operatorname{Coker} $ is the quotient by the
closure of the image. \end{remark}

Consider a direct sum of two Hilbert spaces $ H=H_{1}\oplus H_{2} $. Consider two
closed vector subspaces $ V_{1},V_{2}\subset H $ such that the projection of $ V_{i} $ on $ H_{i} $
has no
null-space and a dense image. This means that one can consider $ V_{1} $ as a graph
of a mapping $ A_{1}\colon H_{1} \to H_{2} $, similarly $ V_{2} $ is a graph of $ A_{2}\colon H_{2} \to H_{1} $.
Mappings $ A_{1,2} $ are closed, but not necessarily bounded.

\begin{lemma}[abstract finiteness] If $ A_{2} $ is bounded, and $ A_{1}\circ A_{2} $ is compact,
then $ V_{1}\cap V_{2} $ is finite dimensional. \end{lemma}

\begin{proof} The projection of $ V_{1}\cap V_{2} $ to $ H_{2} $ is a subspace of $ \operatorname{Ker}\left(A_{1}\circ A_{2}-\boldsymbol1\right) $,
thus is finite-dimensional. \end{proof}

\begin{proposition}[strong form] If $ A_{1} $, $ A_{2} $ are bounded, and $ A_{1}\circ A_{2} $
is compact, then $ V_{1} $ and $ V_{2} $ satisfy the strong form of Riemann--Roch
theorem with index 0. \end{proposition}

\begin{proof} The projection of $ V_{1}\cap V_{2} $ to $ H_{2} $ is $ \operatorname{Ker}\left(A_{1}\circ A_{2}-\boldsymbol1\right) $, which implies
the statement about $ \dim  V_{1}\cap V_{2} $. The statement
about $ \operatorname{codim}\overline{V_{1}+V_{2}} $ follows from the fact that the orthogonal complements
to $ V_{1} $ and $ V_{2} $ satisfy the same conditions as $ V_{2} $ and $ V_{1} $ with linear
mappings being $ -A_{1}^{*} $, $ -A_{2}^{*} $.

To show that $ \dim  V_{1}\cap V_{2} = \operatorname{codim}\left(\overline{V_{1}+V_{2}}\right) $ note that $ \operatorname{Ker}\left(A_{1}\circ A_{2}-\boldsymbol1\right) $ is
dual to $ \operatorname{Coker}\left(A_{2}^{*}\circ A_{1}^{*}-\boldsymbol1\right) $. \end{proof}

\begin{proposition}[weak form] If $ A_{2} $ is bounded, and both
$ A_{1}\circ A_{2} $ and $ A_{1}^{*}\circ A_{2}^{*} $ are compact,
then $ V_{1} $ and $ V_{2} $ satisfy the Riemann--Roch theorem with
index 0. \end{proposition}

\begin{proof} We already know that $ \dim  V_{1}\cap V_{2} $ and $ \operatorname{codim} \overline{V_{1}+V_{2}} $ are finite.
The only thing to prove is that
\begin{equation}
\dim  \operatorname{Ker}\left(A_{1}\circ A_{2}-\boldsymbol1\right) = \dim  \operatorname{Coker}\left(A_{1}^{*}\circ A_{2}^{*}-\boldsymbol1\right).
\notag\end{equation}
Obviously, $ \left(A_{1}^{*}\circ A_{2}^{*}\right)^{*} $ is the closure of $ A_{2}\circ A_{1} $, thus
\begin{equation}
A_{2}\left(\operatorname{Ker}\left(A_{1}\circ A_{2}-\boldsymbol1\right)\right) \subset \operatorname{Ker}\left(\left(A_{1}^{*}\circ A_{2}^{*}\right)^{*}-\boldsymbol1\right),
\notag\end{equation}
hence
\begin{equation}
\dim  \operatorname{Ker}\left(A_{1}\circ A_{2}-\boldsymbol1\right) \leq \dim  \operatorname{Ker}\left(\left(A_{1}^{*}\circ A_{2}^{*}\right)^{*}-\boldsymbol1\right) = \dim  \operatorname{Coker}\left(A_{1}^{*}\circ A_{2}^{*}-\boldsymbol1\right).
\notag\end{equation}
Application of the same argument to the dual operators shows the opposite
unequality. \end{proof}

We say that two vector subspaces $ V $, $ V' $ of a vector space $ H $ are
{\em comparable}, if $ V\cap V' $ is of finite codimension in both $ V $ and $ V' $. The
{\em relative dimension\/} $ \operatorname{reldim}\left(V,V'\right) $ is $ \operatorname{codim}\left(V\cap V'\subset V\right)-\operatorname{codim}\left(V\cap V'\subset V'\right) $. The
following theorem is a direct corollary of the above statement:

\begin{theorem}[Riemann--Roch theorem]  \label{th6.50}\myLabel{th6.50}\relax  Consider two
vector subspaces $ V_{1,2}\subset H $ of a Hilbert space $ H=H_{1}\oplus H_{2} $. Suppose that $ V_{1} $ is
comparable with the graph of a closed mapping $ A_{1}\colon H_{1} \to H_{2} $ and the
relative dimension of $ V_{1} $ and this graph is $ d_{1} $. Suppose $ V_{2} $ is comparable
with the graph of a closed mapping $ A_{2}\colon H_{2} \to H_{1} $ and the relative
dimension of $ V_{2} $ and this graph is $ d_{2} $.
\begin{enumerate}
\item
{\bf(weak form) }If $ A_{2} $ is bounded, and both $ A_{1}\circ A_{2} $ and $ A_{1}^{*}\circ A_{2}^{*} $ are
compact, then $ V_{1} $ and $ V_{2} $ satisfy the Riemann--Roch theorem with the index
being $ d_{1}+d_{2} $.
\item
{\bf(strong form) }If both $ A_{1} $ and $ A_{2} $ are bounded, and $ A_{1}\circ A_{2} $ is compact,
then $ V_{1} $ and $ V_{2} $ satisfy the strong form of Riemann--Roch theorem with the
index being $ d_{1}+d_{2} $.
\end{enumerate}
\end{theorem}

\subsection{Riemann problem }\label{s5.5}\myLabel{s5.5}\relax  Consider a holomorphic function $ \psi\left(z\right) $ defined
in a annulus $ U=\left\{z \mid 1-\varepsilon\leq|z|\leq1+\varepsilon\right\} $, let $ S^{1}=\left\{z \mid |z|=1\right\} $. Suppose that $ \psi\left(z\right) $ is
nowhere 0, and consider the subspace $ V_{\psi} $ in $ L_{2}\left(S^{1},\omega^{1/2}\right)\oplus L_{2}\left(S^{1},\omega^{1/2}\right) $ consisting
of pairs of the form $ \left(f\left(z\right),\psi\left(z\right)f\left(z\right)\right) $. The Hilbert space $ L_{2}\left(S^{1},\omega^{1/2}\right) $ is a
direct sum of subspaces $ L_{2}^{\pm}\left(S^{1},\omega^{1/2}\right) $ consisting of forms which can be
holomorphically continued into two regions $ S^{1} $ divides $ {\mathbb C}P^{1} $ into.

Let $ p $ be the projection of $ V_{\psi} $ to
\begin{equation}
L_{2}^{+}\left(S^{1},\omega^{1/2}\right)\oplus L_{2}^{-}\left(S^{1},\omega^{1/2}\right)\subset L_{2}\left(S^{1},\omega^{1/2}\right)\oplus L_{2}\left(S^{1},\omega^{1/2}\right).
\notag\end{equation}
\begin{lemma} If $ \operatorname{ind}\psi=0 $, then $ p $ is invertible. \end{lemma}

\begin{proof} Consider a linear bundle $ {\mathcal L} $ over $ {\mathbb C}P^{1} $ with isomorphisms to
$ \omega^{1/2} $ over $ U^{+}=\left\{|z|<1+\varepsilon\right\} $ and $ U^{-}=\left\{|z|>1-\varepsilon\right\} $, and the gluing data being $ l^{+}=\psi l^{-} $.
Since $ \deg {\mathcal L}=\deg \omega^{1/2}+\operatorname{ind}\psi $, and a linear bundle over $ {\mathbb C}P^{1} $ is determined by
its degree up to an isomorphism, we see that $ {\mathcal L}\simeq\omega^{1/2} $ if $ \operatorname{ind}\psi=0 $.

Since $ \operatorname{Ker} p $ consists of global sections of $ {\mathcal L} $, $ \operatorname{Ker} p=\left\{0\right\} $. Similarly,
consideration of orthogonal complement to $ V_{\psi} $ shows that $ \operatorname{Coker} p=\left\{0\right\} $. This
finishes the proof, since the operator is obviously Fredholm. \end{proof}

Consider now another linear bundle $ {\mathcal L}' $ with trivializations over $ U^{+} $
and over $ U^{-} $ with gluing data $ l^{+}=\psi l^{-} $. The same arguments as above show
that $ {\mathcal L}' $ is trivial, thus it has a (unique up to multiplication by a
constant) global section. This means that there are functions $ l^{\pm} $ defined
on $ U^{\pm} $ such that $ l^{+}=\psi l^{-} $. Since this global section has no zeros, $ l^{\pm} $ have
no zeros inside the domain of definition, thus $ \psi=\left(l^{+}\right)^{-1}l^{-} $. We see that
any function $ \psi $ such that $ \operatorname{ind}\psi=0 $ can be represented as a product $ \psi=\psi_{+}\psi_{-} $ of
a parts $ \psi_{+} $, $ \psi_{-} $ which can be holomorphically extended inside/outside a
circle without zeros.

Let as write the mapping $ p^{-1} $ in terms of $ \psi_{\pm} $. For any $ L_{2} $-section $ \omega $ of
$ \Omega^{1/2}\otimes\mu $ on $ \left\{|z|=1\right\} $ let $ \omega=\omega_{+}+\omega_{-} $ be the (unique) decomposition of $ \omega $ into a
sum of forms which can holomorphically extended inside/outside of the
unit circle. Given $ \omega_{+} $ and $ \omega_{-}' $ we want to find $ \omega_{-} $ and $ \omega'_{+} $ from the
equality
\begin{equation}
\psi\left(\omega_{+}+\omega_{-}\right) = \omega'_{+}+\omega'_{-},\qquad \text{or\qquad }\psi_{-}\left(\omega_{+}+\omega_{-}\right) = \psi_{+}^{-1}\left(\omega'_{+}+\omega'_{-}\right).
\notag\end{equation}
Taking the $ + $-part we see that $ \left(\psi_{-}\omega_{+}\right)_{+} = \psi_{+}^{-1}\omega'_{+}+\left(\psi_{+}^{-1}\omega'_{-}\right)_{+} $, thus
\begin{equation}
\omega'_{+} = \psi_{+}\left(\psi_{-}\omega_{+}\right)_{+} - \psi_{+}\left(\psi_{+}^{-1}\omega'_{-}\right)_{+} = \psi\omega_{+} - \psi_{+}\left(\psi_{-}\omega_{+}\right)_{-} - \psi_{+}\left(\psi_{+}^{-1}\omega'_{-}\right)_{+},
\notag\end{equation}
similarly $ \left(\psi_{-}\omega_{+}\right)_{-}+\psi_{-}\omega_{-} = \left(\psi_{+}^{-1}\omega'_{-}\right)_{-} $, thus
\begin{equation}
\omega_{-} = \psi_{-}^{-1}\left(\psi_{+}^{-1}\omega'_{-}\right)_{-}-\psi_{-}^{-1}\left(\psi_{-}\omega_{+}\right)_{-} = \psi^{-1}\omega'_{-} - \psi_{-}^{-1}\left(\psi_{+}^{-1}\omega'_{-}\right)_{+} - \psi_{-}^{-1}\left(\psi_{-}\omega_{+}\right)_{-}.
\notag\end{equation}
These two formulae express $ \omega_{-} $ and $ \omega'_{+} $ in terms of $ \omega'_{-} $ and $ \omega_{+} $, thus give
an inverse mapping to $ p $.

Let $ |\psi|_{++}=\frac{\max  |\psi_{+}|}{\min |\psi_{+}|} $, $ |\psi|_{--}=\frac{\max  |\psi_{-}|}{\min |\psi_{-}|} $, $ |\psi|_{+-}=\max 
|\psi_{+}| \max  |\psi_{-}| $, $ |\psi|_{-+}=\max  |\psi_{+}|^{-1} \max  |\psi_{-}|^{-1} $,
$ |\psi|_{0}=\max \left(|\psi|_{++},|\psi|_{--},|\psi|_{+-},|\psi|_{-+}\right) $. Then the norm of $ p^{-1} $ is bounded by
$ C\cdot|\psi|_{0} $.

\begin{definition} \label{def5.155}\myLabel{def5.155}\relax  The {\em Riemann norm\/} $ \|\psi\|_{{\bold R}} $ of $ \psi $ is the norm of the
operator $ p^{-1} $. \end{definition}

The following lemma is a corollary of the fact that one can find
factorization $ \psi=\psi_{+}\psi_{-} $ using integral operators applied to $ \log \psi $:

\begin{lemma} \label{lm5.160}\myLabel{lm5.160}\relax  Let $ \log \psi\left(z\right) $ is defined using any branch of logarithm.
Then
\begin{equation}
\|\psi\|_{{\bold R}}< C \exp  C \max _{z\in U} |\log \psi\left(z\right)|
\notag\end{equation}
for an appropriate $ C $ (which depends on $ \varepsilon $ only). \end{lemma}

Let $ \Pi_{\psi} $ be the composition of $ p^{-1} $ with the projection of
$ L_{2}\left(S^{1},\omega^{1/2}\right)\oplus L_{2}\left(S^{1},\omega^{1/2}\right) $ to $ L_{2}^{-}\left(S^{1},\omega^{1/2}\right)\oplus L_{2}^{+}\left(S^{1},\omega^{1/2}\right) $,
\begin{equation}
\Pi_{\psi}\colon L_{2}^{+}\left(S^{1},\omega^{1/2}\right)\oplus L_{2}^{-}\left(S^{1},\omega^{1/2}\right) \to L_{2}^{-}\left(S^{1},\omega^{1/2}\right)\oplus L_{2}^{+}\left(S^{1},\omega^{1/2}\right).
\notag\end{equation}
In other words, $ \left(\left(f_{+},g_{+}\right),\left(f_{-},g_{-}\right)\right) $ lies on the graph of $ \Pi_{\psi} $ if
$ g_{+}+g_{-}=\psi\left(f_{+}+f_{-}\right) $.

We see that the norm of $ \Pi_{\psi} $ is bounded by $ C\cdot|\psi|_{0} $. Moreover, $ \Pi_{\psi} $ can
be written as a sum
\begin{equation}
\Pi_{\psi} = \left(
\begin{matrix}
0 & \psi^{-1} \\ \psi & 0
\end{matrix}
\right) + k.
\label{equ5.52}\end{equation}\myLabel{equ5.52,}\relax 
Obviously, $ k $ is compact (as any Hankel operator). Indeed, components of $ k $
look like $ f_{+} \mapsto \left(mf_{+}\right)_{-} $. The Hilbert operator $ f=f_{+}+f_{-} \buildrel{H}\over{\mapsto} f_{+}-f_{-} $
is a pseudodifferential operator of degree 0, thus components of $ K $ may
be written as $ \left[m,H\right] $, thus are pseudodifferential operators of degree $ -1 $,
thus compact.

\begin{lemma} \label{lm5.60}\myLabel{lm5.60}\relax  Consider an invertible function $ \psi $ defined in a
neighborhood of $ |z|=1 $ and having $ \operatorname{ind}=k $. Let $ V_{\psi}\subset H=L_{2}\left(S^{1}\right)\oplus L_{2}\left(S^{1}\right) $ be
\begin{equation}
V_{\psi}=\left\{\left(f_{1},f_{2}\right) \mid f_{2}=\psi f_{1}\right\}.
\notag\end{equation}
Define $ H^{\pm}=L_{2}\left(S^{1}\right)^{\pm}\oplus L_{2}\left(S^{1}\right)^{\mp} $. Then for an appropriate bounded operator $ \pi\colon H^{+}
\to H^{-} $ the $ \operatorname{graph}\left(\pi\right) $ is compatible with $ V_{\psi} $, and $ \operatorname{reldim}\left(H_{2},\operatorname{graph}\left(\pi\right)\right)=k $. \end{lemma}

\begin{proof} Let $ \psi\left(z\right)=z^{k}\psi_{0}\left(z\right) $. The function $ \psi_{0} $ has $ \operatorname{ind}=0 $, thus the
corresponding subspace $ V_{\psi_{0}} $ is the graph of $ \Pi_{\psi_{0}} $. Let
$ H_{0}^{\pm}=L_{2}\left(S^{1}\right)^{\pm}\oplus z^{k}L_{2}\left(S^{1}\right)^{\mp} $. We see that $ V_{\psi} $ is a graph of a bounded mapping $ H_{0}^{+}
\to H_{0}^{-} $.
Since $ H_{0}^{+} $ is compatible with $ H_{0} $ of relative dimension $ k $, we momentarily
obtain the required statement about $ V_{\psi} $. \end{proof}

\begin{remark} Note that $ \|\pi\| $ may be bounded in the same way as in Lemma
~\ref{lm5.160}, but neither this result, nor~\eqref{equ5.52} are going to be needed in
what follows. \end{remark}

\subsection{Finiteness theorem }\label{s7.30}\myLabel{s7.30}\relax  Consider a family of disks and
gluing conditions $ K_{\bullet} $, $ \varphi_{\bullet} $, $ \psi_{\bullet} $ from Section~\ref{s5.30}. Let
$ {\mathcal H}^{\pm}=\bigoplus_{l_{2}}L_{2}^{\pm}\left(\partial K_{i},\Omega^{1/2}\otimes\mu\right) $, $ i\in I $. Then the operator with Cauchy kernel defines
a Hilbert mapping $ {\bold K}\colon {\mathcal H}^{-} \to {\mathcal H}^{+} $ (see Section~\ref{s5.20}). This mapping depends
on the circles $ \partial K_{i} $ only, not on the gluing conditions $ \varphi_{\bullet} $, $ \psi_{\bullet} $.

Since the Hilbert structure on $ L_{2}\left(\partial K_{i},\Omega^{1/2}\otimes\mu\right) $ is invariant w.r.t.
fraction-linear mappings (as is decomposition into $ \pm $-parts), the
identification $ \varphi_{i}\colon \partial K_{i'} \to \partial K_{i} $ gives an isomorphism of Hilbert spaces
$ L_{2}\left(\partial K_{i},\Omega^{1/2}\otimes\mu\right) $ and $ L_{2}\left(\partial K_{i'},\Omega^{1/2}\otimes\mu\right) $ which interchanges $ + $-part and $ - $-part.

Suppose that $ \operatorname{ind}\psi_{i}=0 $. Then the operator $ \Pi_{\psi_{i}} $ from Section~\ref{s5.5}
together with an identification given by $ \varphi_{i} $ gives a mapping $ \Pi_{\varphi_{i},\psi_{i}} $
\begin{equation}
L_{2}^{+}\left(\partial K_{i},\Omega^{1/2}\otimes\mu\right)\oplus L_{2}^{+}\left(\partial K_{i'},\Omega^{1/2}\otimes\mu\right) \to L_{2}^{-}\left(\partial K_{i},\Omega^{1/2}\otimes\mu\right)\oplus L_{2}^{-}\left(\partial K_{i'},\Omega^{1/2}\otimes\mu\right),
\notag\end{equation}
the graph of this mapping consists of half-forms on $ \partial K_{i} $ and
$ \partial K_{i'} $ which differ by multiplication by $ \psi_{i} $. In other words, $ \left(\left(f_{+},g_{+}\right),\left(f_{-},g_{-}\right)\right) $
lies on this graph if
\begin{equation}
\left(g_{+}+g_{-}\right)\left(\varphi_{i}^{-1}\left(z\right)\right)=\psi_{i}\left(z\right)\left(f_{+}+f_{-}\right)\left(z\right),\qquad z\in\partial K_{i}.
\notag\end{equation}

If $ \operatorname{ind}\psi_{i}\not=0 $, instead of $ \Pi_{\varphi_{i},\psi_{i}} $ consider an arbitrary operator
between the same spaces such that the graph of $ P $ is comparable with the
set of half-forms which differ by multiplication by $ \psi_{i} $ (see Lemma
~\ref{lm5.60}). The
only fact important in what follows is that this is a bounded operator
(see Corollary~\ref{cor4.35}).

Let $ \Pi=\Pi_{\left\{\psi\right\}}=\bigoplus_{i\in I_{+}}\Pi_{\varphi_{i},\psi_{i}}\colon {\mathcal H}^{+} \to {\mathcal H}^{-} $. Similarly, let $ \Pi_{\left\{\psi^{-1}\right\}} $ corresponds
to gluing data $ \varphi_{i} $, $ \psi_{i}^{-1} $. Since $ \Pi $ consists of bounded diagonal blocks, it
is a closed operator.

Let $ \bar{M} $ be a curve determined by gluing conditions $ \varphi_{\bullet} $, $ {\mathcal L} $ be a bundle
on $ \bar{M} $ determined by gluing conditions $ \psi_{\bullet} $.

\begin{theorem} Suppose that the disks $ K_{i} $ are well separated, $ {\mathcal L} $ is a
finite-degree bundle, and $ \Pi_{\left\{\psi\right\}}\circ{\bold K} $ is compact. Then both
$ \Gamma_{\text{strong}}\left(\bar{M},\omega^{1/2}\otimes{\mathcal L}\right) $ and $ \Gamma_{\text{weak}}\left(\bar{M},\omega^{1/2}\otimes{\mathcal L}\right) $ are finite-dimensional. \end{theorem}

\begin{proof} We give only a sketch of a proof, since the details are the
same as in the case of Riemann--Roch theorem (see Section~\ref{s5.60}). Since
$ \Gamma_{\text{strong}} $ is identified with a subspace of $ \Gamma_{\text{weak}} $, it is enough to show that
$ \dim \Gamma_{\text{weak}}<\infty $. On the other hand, $ \Gamma_{\text{weak}} $ is defined as an intersection of the
generalized Hardy space $ {\mathcal H}\left({\mathbb C}P^{1},\left\{K_{i}\right\}\right) $ with the subspace of forms which
satisfy the
gluing conditions. We are going to reduce the statement of the theorem
to the abstract finiteness theorem.

To do this, note that Proposition~\ref{prop5.28} identifies $ {\mathcal H}\left({\mathbb C}P^{1},\left\{K_{i}\right\}\right) $
with the graph of $ {\bold K}\colon {\mathcal H}^{-} \to {\mathcal H}^{+} $, thus the only thing to note is that fact
that the subspace of forms which satisfy the gluing conditions is
compatible with the graph of $ \Pi_{\left\{\psi\right\}} $. \end{proof}

\subsection{Duality }\label{s7.40}\myLabel{s7.40}\relax  Consider the curve $ \bar{M} $ and a sheaf $ {\mathcal L} $ from the previous
section and the dual sheaf $ {\mathcal L}^{-1} $ with inverse gluing conditions $ \left\{\psi_{i}^{-1}\right\} $.

\begin{theorem} \label{th5.70}\myLabel{th5.70}\relax  Suppose that the disks $ K_{i} $ are well-separated, and $ {\mathcal L} $
is of finite degree. The mappings
\begin{equation}
H^{1/2}\left(\bar{M},\omega^{1/2}\otimes{\mathcal L}\right) \xrightarrow[]{\bar{\partial}} H^{-1/2}\left(\bar{M},\omega^{1/2}\otimes{\mathcal L}\otimes\bar{\omega}\right)
\notag\end{equation}
and
\begin{equation}
H^{1/2}\left(\bar{M},\omega^{1/2}\otimes{\mathcal L}^{-1}\right) \xrightarrow[]{-\bar{\partial}} H^{-1/2}\left(\bar{M},\omega^{1/2}\otimes{\mathcal L}^{-1}\otimes\bar{\omega}\right)
\notag\end{equation}
are mutually dual, thus dimension of null-space of one mapping is equal
to the dimension of cokernel of another one. \end{theorem}

\begin{proof} This is a direct corollary of Lemmas~\ref{lm5.22} and~\ref{lm5.25}. \end{proof}

Using Theorems~\ref{th5.70} and~\ref{th5.31} together with selfduality of $ {\bold K} $
and the fact that $ \Pi_{\psi}^{t}=\Pi_{\psi^{-1}} $, we obtain

\begin{corollary} Suppose that the disks $ K_{i} $ are well separated, and $ {\mathcal L} $ is
finite-degree and semibounded, and both $ \Pi_{\left\{\psi\right\}}\circ{\bold K} $ and $ \Pi_{\left\{\psi^{-1}\right\}}\circ{\bold K} $ are
compact. Then the bounded operator
\begin{equation}
H_{\text{weak}}^{1/2}\left(\bar{M},\omega^{1/2}\otimes{\mathcal L}\right) \xrightarrow[]{\bar{\partial}} H_{\text{weak}}^{-1/2}\left(\bar{M},\omega^{1/2}\otimes{\mathcal L}\otimes\bar{\omega}\right)
\notag\end{equation}
has finite-dimensional null-space and cokernel. \end{corollary}

\subsection{Riemann--Roch for curves }\label{s5.60}\myLabel{s5.60}\relax  Now we are ready to state

\begin{theorem}[Riemann--Roch]  \label{th5.15}\myLabel{th5.15}\relax  Suppose that the disks $ K_{i} $ are
well-separated, and $ \psi_{i} $ is a finite-degree family.
\begin{enumerate}
\item
If both $ \Pi_{\left\{\psi\right\}}\circ{\bold K} $ and $ \Pi_{\left\{\psi^{-1}\right\}}\circ{\bold K} $ are compact operators then the mapping
\begin{equation}
H_{\text{weak}}^{1/2}\left(\bar{M},\omega^{1/2}\otimes{\mathcal L}\right) \xrightarrow[]{\bar{\partial}} H_{\text{weak}}^{-1/2}\left(\bar{M},\omega^{1/2}\otimes{\mathcal L}\otimes\bar{\omega}\right)
\notag\end{equation}
has finite-dimensional null-space and cokernel, and
\begin{equation}
\dim  \operatorname{Ker} \bar{\partial} - \dim  \operatorname{Coker} \bar{\partial} = \sum_{i\in I_{+}} \operatorname{ind} \psi_{i}.
\notag\end{equation}
\item
If $ \Pi_{\left\{\psi\right\}} $ is bounded, and $ {\bold K} $ is compact, then the mapping is Fredholm
of index $ \sum_{i\in I_{+}} \operatorname{ind} \psi_{i} $.
\end{enumerate}
\end{theorem}

\begin{proof} The only thing we need to do is to describe the null-space and
the image of $ \bar{\partial} $. We are going to reduce this description to the abstract
Riemann--Roch theorem. First of all, the null-space consists of
half-forms in the generalized Hardy space, so $ \widetilde{r} $ maps it injectively to a
subspace of $ \bigoplus_{l_{2}}L_{2}\left(\partial K_{i},\Omega^{1/2}\otimes\mu\right) $. On the other hand, the image of $ \operatorname{Ker}\bar{\partial} $ in
$ \bigoplus_{l_{2}}L_{2}\left(\partial K_{i},\Omega^{1/2}\otimes\mu\right) $ is described as intersection of the image of the
generalized Hardy space and a subspace in $ \bigoplus_{l_{2}}L_{2}\left(\partial K_{i},\Omega^{1/2}\otimes\mu\right) $ consisting
of half-forms which satisfy the gluing conditions.

The first condition can be written as
\begin{equation}
f_{i}^{+}=\sum_{j\not=i}{\bold K}_{ij}f_{j}^{-},
\notag\end{equation}
here $ f_{i}^{\pm} $ are $ \pm $-parts of the component of $ f $ in $ L_{2}\left(\partial K_{i},\Omega^{1/2}\otimes\mu\right) $. The second
condition is
\begin{equation}
f'_{i'}=\psi_{i}f_{i}.
\notag\end{equation}
Here $ f\left(t\right)'=f\left(-t\right) $, and we suppose that we use compatible parameterizations
of $ \partial K_{i} $ and of $ \partial K_{i'} $, i.e., such that $ \varphi_{i} $ send parameter $ t $ on $ \partial K_{i} $ to
parameter $ -t $ on $ \partial K_{i'} $.

Consider the decomposition $ {\mathcal H}={\mathcal H}^{+}\oplus{\mathcal H}^{-} $. Vectors which satisfy the
first condition are in the graph of operator $ {\bold K}\colon {\mathcal H}^{-} \to {\mathcal H}^{+} $. Consider the
vector space of vectors $ {\mathcal H}_{2} $ which satisfy the second conditions. It is a
direct sum over $ i\in I_{+} $ of subspaces of
\begin{equation}
L_{2}\left(\partial K_{i},\Omega^{1/2}\otimes\mu\right)\oplus L_{2}\left(\partial K_{i'},\Omega^{1/2}\otimes\mu\right),
\notag\end{equation}
and all the components but a finite number have $ \psi $ with $ \operatorname{ind}=0 $, thus are
described as graphs of $ \Pi_{\varphi_{i},\psi_{i}} $. As a corollary, we conclude that this
subspace is compatible with the graph of the mapping $ \Pi $.

To finish the description of the null-space the only thing which remains
to prove is to show that the relative dimension of $ {\mathcal H}_{2} $ and $ \operatorname{graph}\left(\Pi\right) $ is
$ \sum_{i\in I_{+}}\operatorname{ind}\psi_{i} $. However, because of decomposition of $ {\mathcal H}_{2} $ into a direct sum it
follows from the corresponding fact for each component, i.e., from
Lemma~\ref{lm5.60}.

To describe the image, consider
\begin{equation}
\alpha\in H_{\text{weak}}^{-1/2}\left(\bar{M},\omega^{1/2}\otimes{\mathcal L}\otimes\bar{\omega}\right) = H^{-1/2}\left({\mathbb C}P^{1}\smallsetminus\bigcup K_{i},\omega^{1/2}\otimes\bar{\omega}\right).
\notag\end{equation}
Since $ \bar{\partial} $ is an isomorphism
\begin{equation}
H^{1/2}\left({\mathbb C}P^{1},\omega^{1/2}\right) \to H^{-1/2}\left({\mathbb C}P^{1},\omega^{1/2}\otimes\bar{\omega}\right),
\notag\end{equation}
the element $ \bar{\partial}^{-1}\alpha\in H^{1/2}\left({\mathbb C}P^{1},\omega^{1/2}\right) $ is defined up to addition of an
element of $ \bar{\partial}^{-1} \overset{\,\,{}_\circ}{H}^{-1/2}\left(\bigcup K_{i},\omega^{1/2}\otimes\bar{\omega}\right) $. On the other hand, the latter space is
$ {\mathcal H}\left({\mathbb C}P^{1},\left\{K_{i}\right\}\right) $, i.e., $ \bar{\partial}^{-1} $ gives a correctly defined isomorphism
\begin{equation}
H^{-1/2}\left({\mathbb C}P^{1}\smallsetminus\bigcup K_{i},\omega^{1/2}\otimes\bar{\omega}\right) \to H^{1/2}\left({\mathbb C}P^{1},\omega^{1/2}\right)/{\mathcal H}\left({\mathbb C}P^{1},\left\{K_{i}\right\}\right),
\notag\end{equation}
hence $ \bar{\partial}^{-1}\alpha $ is an element of the latter space. Thus
\begin{equation}
\alpha\in\operatorname{Im}\left(H_{\text{weak}}^{1/2}\left(\bar{M},\omega^{1/2}\otimes{\mathcal L}\right) \xrightarrow[]{\bar{\partial}} H_{\text{weak}}^{-1/2}\left(\bar{M},\omega^{1/2}\otimes{\mathcal L}\otimes\bar{\omega}\right)\right)
\notag\end{equation}
is equivalent to
\begin{equation}
\partial^{-1}\alpha\in {\mathcal H}_{2}/{\mathcal H}\left({\mathbb C}P^{1},\left\{K_{i}\right\}\right) \buildrel{\text{def}}\over{=}\left({\mathcal H}\left({\mathbb C}P^{1},\left\{K_{i}\right\}\right)+{\mathcal H}_{2}\right)/{\mathcal H}\left({\mathbb C}P^{1},\left\{K_{i}\right\}\right)\text{ .}
\notag\end{equation}
In particular, if $ {\mathcal H}\left({\mathbb C}P^{1},\left\{K_{i}\right\}\right)+{\mathcal H}_{2} $ is closed, then $ \operatorname{Im} \bar{\partial} $ is closed, and in
any case the codimension of the closure of $ \operatorname{Im}\bar{\partial} $ is equal to codimension of
the closure of $ {\mathcal H}\left({\mathbb C}P^{1},\left\{K_{i}\right\}\right)\oplus{\mathcal H}_{2} $. \end{proof}

From Lemma~\ref{lm5.160} we momentarily obtain

\begin{corollary} Suppose that $ {\bold K} $ is compact, and $ \psi_{i}=\exp \Phi_{i} $ for all $ i $ but a
finite number. If for an appropriate $ C $ and any $ i $
\begin{equation}
|\Phi_{i}| < C,
\notag\end{equation}
then the operator
\begin{equation}
H_{\text{weak}}^{1/2}\left(\bar{M},\omega^{1/2}\otimes{\mathcal L}\right) \xrightarrow[]{\bar{\partial}} H_{\text{weak}}^{-1/2}\left(\bar{M},\omega^{1/2}\otimes{\mathcal L}\otimes\bar{\omega}\right)
\notag\end{equation}
is Fredholm. \end{corollary}

\begin{definition} Call a pair $ \left(\bar{M},{\mathcal L}\right) $ {\em admissible\/} if $ \bar{M} $ is given by a family of
well-separated disks, $ {\mathcal L} $ is of finite degree, and both $ \Pi_{\left\{\psi\right\}}\circ{\bold K} $ and $ \Pi_{\left\{\psi^{-1}\right\}}\circ{\bold K} $
are compact. Call $ \bar{M} $ {\em admissible\/} if $ \left(\bar{M},\boldsymbol1\right) $ is admissible. Here {\bf1 }is a bundle
over $ \bar{M} $ with $ \psi_{i}\equiv 1 $. \end{definition}

\subsection{Criterion of admissibility of a curve }\label{s4.50}\myLabel{s4.50}\relax  Recall that in Section
~\ref{s5.10} we considered a matrix $ \left(e^{-l_{ij}/2}\right) $ which was supposed to give a
bounded operator in $ l_{2} $.

\begin{proposition} If the matrix $ \left(e^{-l_{ij}/2}-\delta_{ij}\right) $ gives a compact operator $ l_{2} \to
l_{2} $, then $ \bar{M} $ is admissible. \end{proposition}

\begin{proof} Indeed, for the bundle $ {\mathcal L}=\boldsymbol1 $ the operator $ \Pi $ is an isometry, thus
we need to show that $ {\bold K} $ is compact, which is a corollary of Lemma~\ref{lm5.05}. \end{proof}

\begin{corollary} If $ \sum_{i\not=j}e^{-l_{ij}}<\infty $, then $ \bar{M} $ is admissible. \end{corollary}

\begin{corollary} Fix a locally discrete subset $ I\subset{\mathbb C}P^{1} $ (i.e., for any point $ i $
of $ I $ there is a punctured neighborhood of $ i $ which does not intersect
$ I $). Then there exists a family of disks $ K_{i} $, $ i\in I $ with $ \operatorname{center}\left(K_{i}\right)=i $ such
that for any involution ' and any gluing data $ \varphi_{i} $ the corresponding curve $ \bar{M} $
is admissible. Moreover, one can chose $ K_{i} $ in such a way that $ M_{\infty} $ consists
of points of smoothness $ C^{\infty} $ (see Section~\ref{s4.95}). \end{corollary}

\begin{corollary} For an arbitrary nowhere dense subset $ N $ of $ {\mathbb C}P^{1} $ there exists
an admissible curve $ \bar{M} $ such that $ M_{\infty}=N $. Moreover, it is possible to make
every point of $ N $ to be of smoothness $ C^{\infty} $. \end{corollary}

\begin{remark} While the above statements are obvious, note that construction
of examples and counterexamples may be simplified a lot by an additional
restriction:
\begin{equation}
\operatorname{dist}\left(i,j\right) \geq \varepsilon\cdot\operatorname{dist}\left(i,N\right)\qquad \text{for any }i,j\in I\text{, }i\not=j.
\notag\end{equation}
Here $ \varepsilon\ll1 $. To construct such a family $ I $ for a given nowhere dense set $ N $,
let $ N_{k}=\left\{z \mid 2^{-k-1}\leq\operatorname{dist}\left(z,N\right)\leq2^{-k}\right\} $. Fix $ k $, and let $ \delta=2^{-k} $. Consider a $ \delta/8 $-net
for $ N_{k} $.
By removing some points from this net one can obtain $ \delta/4 $-net such
that it does not have two points closer than $ \delta/8 $. This net is necessarily
finite.

Now consider the union of these finite sets over $ k\in{\mathbb N} $, and again
remove net points from $ N_{k+1} $ which are closer than $ \delta/8 $ to net points
in $ N_{k} $. One obviously obtains a set $ I $ with required properties.

Now to chose the radius of $ K_{i} $ denote by $ n_{k} $ the number of points in
$ I\cap N_{k} $. Let $ \operatorname{radius}\left(K_{i}\right)=\frac{f\left(k\right)}{n_{k}} $ if $ i\in N_{k} $, here $ f\left(k\right) $ is rapidly decreasing
function. Picking appropriate $ f\left(k\right) $, one obtains disks which satisfy the
given above requirements. \end{remark}

The following property of admissibility is obvious:

\begin{nwthrmiii} If we change a finite number of contours $ \partial K_{i} $ and/or a finite
numbers of identifications $ \varphi_{i} $, this does not change the admissibility of
the resulting curve. \end{nwthrmiii}

\subsection{Filling the gap }\label{s7.70}\myLabel{s7.70}\relax  This section contains heuristics only, so
anyone
interested exceptionally in exact results should proceed directly to
Section~\ref{h8}.

Consider once more the criterion of admissibility of a curve. It
says that if the operator $ {\bold K} $ is compact, the curve is admissible (together
with any bundle which is defined by gluing functions $ \psi_{i} $ with bounded
Riemann norm $ |\psi|_{{\bold R}} $). On the other hand, to obtain this result we study the
generalized Hardy space, which is correctly defined if $ {\bold K} $ is bounded.

Thus the tools we use leave a gap between the objects for which the
analysis is applicable (i.e., $ {\bold K} $ is bounded), and objects for which we get
the admissibility (i.e., $ {\bold K} $ is compact). How can we use the existence of
this gap?

We propose to consider this gap as a confirmation that the
Riemann--Roch theorem we obtained is {\em almost unimprovable}. Indeed, we know
that the finiteness condition can be written as $ \dim  \operatorname{Ker}\left(\Pi\circ{\bold K}-1\right)<\infty $. On the
other hand, the invertible operator $ \Pi $ depends on the family $ \left(\psi_{i}\right) $ which
has very high degree of freedom (even if we consider the strong form of
Riemann--Roch theorem, so $ \Pi $ is required to be bounded), and one should
expect that the condition
\begin{equation}
\dim  \operatorname{Ker}\left(\Pi\circ{\bold K}-1\right)<\infty\text{ for every choice of }\psi_{i}\text{ with }|\psi_{\bullet}|_{{\bold R}}<\infty
\notag\end{equation}
should be very close to the condition
\begin{equation}
{\bold K}\text{ has discrete spectrum near }|\lambda|=1.
\label{equ4.88}\end{equation}\myLabel{equ4.88,}\relax 
In other words, the pairs $ \left(\Pi,{\bold K}\right) $ with an invertible $ \Pi $, bounded $ {\bold K} $, and
{\em infinite-dimensional\/} null-space of $ \Pi{\bold K}-1 $ are plentiful (at least if we
drop ``geometric'' conditions on $ \Pi $ and $ {\bold K} $, and consider abstract operators),
and they form a ``natural boundary'' of the set of curves for which
Riemann--Roch theorem has a chance to be true. This natural boundary is
quite close to the boundary of the set of compact operators, which is
another confirmation of our thesis.

Note that~\eqref{equ4.88} may lead to a different description of possible
$ {\bold K} $, like $ \frac{{\bold K}}{1+{\bold K}^{2}} $ being compact. It is unclear, however, whether the above
boundary separates the set of compact operators as a connected component
of the set of operators $ {\bold K} $ with compact $ \frac{{\bold K}}{1+{\bold K}^{2}} $.

\subsection{Moduli space } The description of a complex curve by disks $ \left\{K_{i}\right\} $,
involution ' and gluings $ \varphi_{i} $ leaves a feeling of being incomplete, since
in the case of finite genus it is enough to provide just gluings $ \varphi_{i} $ which
generate a subgroup of $ \operatorname{SL}\left(2,{\mathbb C}\right) $. To describe the quotient by this subgroup
one can take any fundamental domain for the subgroup, and different
choices of the fundamental domain result in the same geometric data.

To get a similar description in the case of infinite genus, note
that in Section~\ref{s5.10} instead of the restriction that $ \widetilde{K}_{i} $ {\em is\/} a concentric
with $ K_{i} $ disk of radius $ e^{2\varepsilon}\cdot\operatorname{radius}\left(K_{i}\right) $ one can require that $ \widetilde{K}_{i} $ {\em contains\/}
such a circle. This leads to the construction of strong sections of
$ \omega^{1/2}\otimes{\mathcal L} $ in the same way as in Section~\ref{s6.50}, and Theorem~\ref{th5.31} can be
refined as

\begin{amplification} Fix a metric on $ {\mathbb C}P^{1} $ and $ \varepsilon>0 $. Consider a family of
elements $ \varphi_{i} $ of $ \operatorname{SL}\left(2,{\mathbb C}\right) $. Suppose that there exists a family of disjoint
domains $ \widetilde{K}_{i}\subset{\mathbb C}P^{1} $, $ i\in I $, which satisfy the following properties:
\begin{enumerate}
\item
All $ \widetilde{K}_{i} $ but a finite number are disks;
\item
Let $ K'_{i} $ be a concentric with $ \widetilde{K}_{i} $ disk of radius $ \left(1-\varepsilon\right)\operatorname{radius}\left(\widetilde{K}_{i}\right) $ (or
any domain in $ \widetilde{K}_{i} $ if $ \widetilde{K}_{i} $ is not a disk). Suppose that for an involution ':
$ I\to $I one has $ \varphi_{i'}=\varphi_{i}^{-1} $, $ \varphi_{i}\left(K'_{i}\right)\cup K'_{i'}={\mathbb C}P^{1} $;
\item
Let $ \psi_{i} $ be holomorphic functions defined inside $ \widetilde{K}_{i}\cap\varphi_{i}^{-1}\left(\widetilde{K}_{i'}\right) $ which
satisfy $ \psi_{i'}\circ\varphi_{i}=\psi_{i}^{-1} $, and $ \operatorname{ind}\psi_{i}=0 $ for all but a finite number of $ i\in I $.
\item
Suppose that the pairwise conformal distances $ l_{ij} $ between $ \widetilde{K}_{i} $, $ \widetilde{K}_{j} $
satisfy the condition that the matrix $ \left(e^{-l_{ij}/2}-\delta_{ij}\right) $ gives a compact
mapping $ l_{2} \to l_{2} $, and all the functions $ \psi_{i} $ are bounded taken together.
\end{enumerate}

Let $ \bar{M} $ be a curve given by gluing $ S=\bigcap\varphi_{i}^{-1}\left(K'_{i}\right) $ together via $ \varphi_{i} $, $ {\mathcal L} $ be
a bundle on $ \bar{M} $ given by cocycle $ \psi_{i} $. Define strong sections of $ \omega^{1/2}\otimes{\mathcal L} $
associated with family $ \left\{\widetilde{K}_{i}\right\} $ as forms $ \alpha $ in $ H^{1/2}\left(\bigcap\varphi_{i}^{-1}\left(K'_{i}\right),\omega^{1/2}\right) $ which
satisfy $ \varphi_{i}^{*}\left(\alpha\right)=\psi_{i}\cdot\alpha $ whenever both sides have sense.

Let $ I=I_{+}\amalg I_{+}' $. For any choice of circles $ \gamma_{i} $ in $ K_{i}'\cap\varphi_{i}^{-1}\left(K'_{i'}\right) $, $ i\in I_{+} $, let
$ K_{i} $ be the disk bounded by $ \gamma_{i} $ inside $ K'_{i} $ (with appropriate modifications
if $ K_{i} $ is not a disk). Let $ K_{i'}={\mathbb C}P^{1}\smallsetminus\varphi_{i}\left(K_{i}\right) $, $ i\in I_{+} $. Then the space of strong
sections of $ \omega^{1/2}\otimes{\mathcal L} $ coincides with the set of weak sections associated to
the family $ \left\{K_{i}\right\} $.

In particular, the Riemann--Roch theorem (in the strong form) is
valid for strong sections associated to the family $ \left\{\widetilde{K}_{i}\right\} $. \end{amplification}

The proof of this statement is a corollary of the proof of
Theorem~\ref{th5.31}.

We see that at least for bundles described by bounded cocycles one
does not need to specify {\em precisely\/} the circles which cut the curve, one
can vary them in wide ranges (which depend on the Kleinian group) without
any change to the geometric data. One possible objection to usability of
the above theorem is that for different choices $ \left\{\gamma_{i}^{\left(1\right)}\right\} $, $ \left\{\gamma_{i}^{\left(2\right)}\right\} $ of the
circles the cocycles $ \psi_{i} $ which describe the bundle need to be defined
everywhere between the circles for the theorem to be applicable. However,
in Section~\ref{s7.90} we will see that any bounded bundle of degree 0 may be
described by constant cocycles, thus this objection becomes void. (If the
degree is not 0, one can take all the cocycles to be constant except one
pair described by rational functions.)

Call a Kleinian group {\em admissible\/} if it has generators $ \varphi_{i} $, $ i\in I_{+} $, which
satisfy the conditions of the amplification. The following conjecture
would show that the curve (with a fixed family of $ A $-cycles) is
completely described by the corresponding Kleinian group, and the
restriction on distances between $ \widetilde{K}_{i} $ is in fact the restriction on the
Kleinian group.

Consider an admissible Kleinian group. A {\em fundamental family\/} is
a family of domains $ \widetilde{K}_{\bullet} $ which satisfies the conditions of amplification. If
$ \widetilde{K}_{i}^{\left(1\right)}\subset\widetilde{K}_{i}^{\left(2\right)} $, $ i\in I $, then call these families {\em equivalent}, and continue this
relation by transitivity. The amplification shows that equivalent
families lead to the same spaces of sections of bundles.

\begin{conjecture} Any two fundamental families for an admissible
Kleinian group are equivalent. \end{conjecture}

Let the {\em fine moduli space\/} be the set of admissible Kleinian groups
up to conjugation in $ \operatorname{PGL}\left(2,{\mathbb C}\right) $.

\begin{conjecture} Consider a complex curve $ M $. Let $ \gamma_{i}^{\left(1\right)} $, $ i\in I^{\left(1\right)} $, be a
disjoint family of smooth embedded cycles in $ M $ such that $ M\smallsetminus\bigcup\gamma_{i} $ is
conformally equivalent to a fundamental domain of an admissible Kleinian
group. Let $ \gamma_{i}^{\left(2\right)} $ be a different family which satisfies same conditions.
Then all the cycles $ \gamma_{i}^{\left(1\right)} $ except a finite number are homotopic to cycles
in $ \gamma_{i}^{\left(2\right)} $. \end{conjecture}

\begin{conjecture} Two families of cuts from the previous conjecture lead to
the same set of $ C^{n} $-points at infinity (see Section~\ref{s4.95}) for any $ n\geq0 $. \end{conjecture}

\begin{definition} Suppose that families $ \left\{K_{i}^{\left(1\right)},',\varphi_{i}^{\left(1\right)}\right\} $ and $ \left\{K_{i}^{\left(2\right)},',\varphi_{i}^{\left(2\right)}\right\} $
have all the accumulation points of class $ C^{0} $. Say that these families
{\em describe the same curve\/} if the set of finite points of the corresponding
curves are equivalent as complex manifolds, and the equivalence continues
by continuity to the points at infinity. Let the {\em moduli space\/} be the set
of equivalence classes of such families up to relationship that they
describe the same curve. \end{definition}

Note that it is not reasonable to drop the consideration of points
at infinity.

\begin{example} Indeed, consider a curve $ \bar{M} $ with the Serpinsky carpet as the
set $ M_{\infty} $ of accumulation points of disks. The Serpinsky carpet breaks $ {\mathbb C}P^{1} $
into a union of triangles. Suppose that each triangle contains exactly
one disk $ K_{i} $. Then the smooth points on $ M $ form a disjoint union of tubes,
so the only invariant is the conformal lengths of these tubes. This is
one parameter per handle, much smaller than three parameters per handle
as one would expect to have from finite-genus theory. The remaining
parameters must be contained in the data for gluing the boundary of each
smooth component to the set of points at infinity. \end{example}

\section{Set of admissible bundles }\label{h8}\myLabel{h8}\relax 

Here we investigate the structure of the set $ {\mathfrak L} $ of admissible bundles
$ {\mathcal L} $ over the given curve (described by the model space $ \bar{M} $, as in Section
~\ref{s5.30}). The mapping $ {\mathcal L} \mapsto \left(\Pi_{\left\{\psi\right\}}\circ{\bold K},\Pi_{\left\{\psi^{-1}\right\}}\circ{\bold K} \right) $ into a pair of compact
operators allows one to consider the topology on $ {\mathfrak L} $ induced by the
operator topology on the space of compact operators.

Thus $ {\mathfrak L} $ is a topological space. We will describe some remarkable
subsets of $ {\mathfrak L} $ and algebraic structures on these subsets.

\subsection{Exceptional indices } Since operators $ \Pi_{\varphi_{i},\psi_{i}} $ are bounded, and row-blocks
and column-blocks of the operator $ {\bold K} $ are compact, we can conclude that
instead of compactness of $ \Pi{\bold K} $ one can require compactness of restriction
of $ \Pi{\bold K} $ to the direct sum of all but a finite number of $ L_{2}\left(\partial K_{i},\Omega^{1}\otimes\mu\right) $. We
call the indices of excluded contours {\em exceptional indices}.

In particular, we can include all the non-circular contours and all
contours $ \partial K_{i} $ such that $ \operatorname{ind}\psi_{i}\not=0 $ in the set of exceptional indices. Thus
whenever we discuss admissibility conditions we can suppose that $ \operatorname{ind} \psi_{i}=0 $
and all the $ K_{i} $ are disks.

The following facts are obvious:

\begin{proposition} If we change a finite number of functions $ \psi_{i} $, this does
not change the admissibility of the resulting bundle. If we multiply
functions $ \psi_{i} $ by constants $ c_{i} $ with $ |c_{i}| $, $ |c_{i}^{-1}| $ being bounded,
this does not change the admissibility of the bundle. \end{proposition}

\subsection{Hilbert--Schmidt bundles } In practice the condition of admissibility
is very hard to use directly, since there is no practically useful
criterion of compactness which is necessary and sufficient. The closest
simple-to-check approximation is the Hilbert--Schmidt condition.

\begin{definition} We say that a curve with a bundle $ \left(\bar{M},{\mathcal L}\right) $ is {\em Hilbert\/}--{\em Schmidt\/}
if both operators $ \Pi_{\left\{\psi\right\}}{\bold K} $ and $ \Pi_{\left\{\psi^{-1}\right\}}{\bold K} $ are Hilbert--Schmidt operators. \end{definition}

\begin{remark} Note that since the Hilbert structure on the sections of
the bundle of half-forms on a curve is canonically defined, so is the
notion of Hilbert--Schmidt operator. \end{remark}

\subsection{Involution } Consider an admissible bundle $ {\mathcal L} $ and the dual bundle $ {\mathcal L}^{-1} $
(defined by inverse gluing conditions). The following statement is a
direct corollary of definitions:

\begin{proposition} If $ {\mathcal L} $ is admissible, so is $ {\mathcal L}^{-1} $. If $ {\mathcal L} $ is Hilbert--Schmidt,
so is $ {\mathcal L}^{-1} $. \end{proposition}

\subsection{Hilbert--Schmidt criterion }\label{s8.40}\myLabel{s8.40}\relax  From the description of a
solution of the Riemann problem one can easily get the
following criterion:

\begin{theorem} Consider a family of disks $ K_{i} $, $ i\in I $, with an involution ' on I
and gluing data $ \varphi_{i} $, $ \psi_{i} $.
\begin{enumerate}
\item
If
\begin{equation}
\sum_{i\not=j}\left(\left|\psi_{i}\right|_{{\bold R}}^{2}+\left|\psi_{i}^{-1}\right|_{{\bold R}}^{2}\right)e^{-l_{ij}} < \infty
\notag\end{equation}
then $ {\mathcal L} $ is Hilbert--Schmidt (here $ ||_{{\bold R}} $ is the Riemann norm, see Definition
~\ref{def5.155}).
\item
Suppose that all but a finite
number of functions $ \psi_{i} $ are constant. Then the corresponding curve $ M $ with a
bundle $ {\mathcal L} $ is Hilbert--Schmidt iff
\begin{equation}
\sum_{i\not=j}\left(\left|\psi_{i}\right|^{2}+\left|\psi_{i}\right|^{-2}\right)e^{-l_{ij}} < \infty
\notag\end{equation}
(here the indices $ i $ for which $ \psi_{i} $ is not constant are excluded from
summation).
\end{enumerate}

Here $ l_{ij} $ is the conformal distance between $ \partial K_{i} $ and $ \partial K_{j} $.

\end{theorem}

This theorem follows immediately from the following

\begin{lemma} \label{lm6.30}\myLabel{lm6.30}\relax  Let $ k_{ij} $, $ i,j\in I $ be the Hilbert--Schmidt norm of the block of $ {\bold K} $
which maps $ L_{2}^{-}\left(\partial K_{i},\omega^{1/2}\otimes\mu\right) \to L_{2}^{+}\left(\partial K_{j},\omega^{1/2}\otimes\mu\right) $. Then
\begin{equation}
k_{ij} = \sum\Sb s>0 \\ s\in{\mathbb Z}+1/2\endSb e^{-sl_{ij}} = O\left(e^{-l_{ij}/2}\right).
\notag\end{equation}
\end{lemma}

\begin{proof} It is enough to prove that the characteristic numbers of the
block of $ {\bold K} $ which is the mapping
\begin{equation}
L_{2}^{-}\left(\partial K_{i},\omega^{1/2}\otimes\mu\right) \to L_{2}^{+}\left(\partial K_{j},\omega^{1/2}\otimes\mu\right)
\notag\end{equation}
are $ e^{-sl_{ij}} $ for $ s>0 $, $ s\in{\mathbb Z}+\frac{1}{2} $. Since $ {\bold K} $ is invariant with respect to
fraction-linear transformations, we can assume that $ \partial K_{i} $ and $ \partial K_{j} $ bound a
tube $ S^{1}\times\left(0,l_{ij}\right) $ of circumference $ 2\pi $ and length $ l_{ij} $. Then characteristic
vectors of $ {\bold K} $ correspond to holomorphic $ 1/2 $-forms
\begin{equation}
e^{i sx}e^{-sy}dz^{1/2},\qquad z=x+iy,\quad \left(x,y\right)\in S^{1}\times\left(0,l_{ij}\right),
\notag\end{equation}
the restrictions on $ y=0 $ being the characteristic vectors themselves, the
restrictions on $ y=l_{ij} $ being their images, which are $ e^{-sl_{ij}} $ times smaller.
Now the condition that the restriction on $ y=0 $ is a section of $ \Omega\otimes\mu $ gives
the condition that $ s\in{\mathbb Z}+\frac{1}{2} $, the condition that this restriction is in
$ - $-part of $ L_{2} $ gives $ s>0 $.\end{proof}

\begin{corollary}
\begin{enumerate}
\item
A curve $ \bar{M} $ is Hilbert--Schmidt iff
\begin{equation}
\sum_{i\not=j}e^{-l_{ij}} < \infty;
\notag\end{equation}
\item
If a curve $ \bar{M} $ allows a Hilbert--Schmidt bundle, it is
Hilbert--Schmidt itself.
\end{enumerate}
\end{corollary}

\begin{proof} Indeed, only the second part requires proof, and it follows
from
\begin{equation}
\sum_{ij}k_{ij}^{2} \leq \frac{1}{2}\sum_{ij}\left(|\psi_{i}|^{2}+|\psi_{i}|^{-2}\right)k_{ij}^{2}.
\notag\end{equation}
\end{proof}

\subsection{$ \protect \log  $-convexity }\label{s8.50}\myLabel{s8.50}\relax  We were not able to prove the following

\begin{conjecture} Let $ {\mathcal L} $ and $ {\mathcal M} $ are two bundles on $ \bar{M} $ given by gluing
conditions. Suppose that both $ {\mathcal M}\otimes{\mathcal L} $ and $ {\mathcal M}\otimes{\mathcal L}^{-1} $ are admissible. Then $ {\mathcal M} $ is
admissible as well. \end{conjecture}

However, the following statement is true:

\begin{proposition} \label{prop6.70}\myLabel{prop6.70}\relax  Let $ {\mathcal L} $ and $ {\mathcal M} $ are two bundles on $ \bar{M} $ given by
gluing
conditions with functions $ \psi_{i} $ and $ \gamma_{i} $, and all but a finite number of these
functions are constant. Suppose that both $ {\mathcal M}\otimes{\mathcal L} $ and $ {\mathcal M}\otimes{\mathcal L}^{-1} $ are
Hilbert--Schmidt. Then $ {\mathcal M} $ is Hilbert--Schmidt as well. \end{proposition}

\begin{proof} We may suppose that all $ \psi_{i} $ and $ \gamma_{i} $ are constant. Then in
notations of Lemma~\ref{lm6.30} we know that
\begin{gather} \sum_{ij}\left(\left|\psi_{i}\gamma_{i}\right|^{2}+\left|\psi_{i}\gamma_{i}\right|^{-2}\right)k_{ij}^{2} < \infty,
\notag\\
\sum_{ij}\left(\left|\psi_{i}\gamma_{i}^{-1}\right|^{2}+\left|\psi_{i}\gamma_{i}^{-1}\right|^{-2}\right)k_{ij}^{2} < \infty,
\notag\end{gather}
and want to prove that
\begin{equation}
\sum_{ij}\left(\left|\psi_{i}\right|^{2}+\left|\psi_{i}\right|^{-2}\right)k_{ij}^{2} < \infty.
\notag\end{equation}
However, this is an obvious corollary of relation between geometric
mean and arithmetic mean. \end{proof}

\subsection{Types of admissible bundles } In what follows we are going to study
Hilbert--Schmidt bundles, thus we may assume that the curve is
Hilbert--Schmidt itself.

\begin{definition}
\begin{enumerate}
\item
Call the bundle $ \omega\otimes{\mathcal L} $ on a curve $ \bar{M} $ {\em real\/} if all the gluing
functions $ \psi_{i} $ for this bundle are constants of magnitude 1.
\item
Call the bundle $ \omega\otimes{\mathcal L} $ on a curve $ \bar{M} $ {\em bounded\/} if all the gluing
functions $ \psi_{i} $ for this bundle taken together are bounded in $ ||_{{\bold R}} $-norm.
\item
Call the bundle $ \omega\otimes{\mathcal L} $ on a curve $ \bar{M} $ {\em strongly Hilbert\/}--{\em Schmidt\/} if
the bundle $ \omega\otimes{\mathcal L}^{n} $ is Hilbert--Schmidt for any $ n\in{\mathbb Z} $.
\end{enumerate}
\end{definition}

\begin{lemma} If the bundles $ \omega\otimes{\mathcal L} $ and $ \omega\otimes{\mathcal M} $ on a curve $ \bar{M} $ are strongly
Hilbert--Schmidt, then $ \omega\otimes{\mathcal L}\otimes{\mathcal M} $ is also strongly Hilbert--Schmidt. \end{lemma}

\begin{proof} Since $ \omega\otimes{\mathcal L}^{2n} $ and $ \omega\otimes{\mathcal M}^{2n} $ are Hilbert--Schmidt, such is $ \omega\otimes{\mathcal L}^{n}\otimes{\mathcal M}^{n} $
by the $ \log  $-convexity. \end{proof}

\begin{lemma} Any real bundle is bounded. Any bounded bundle is strongly
Hilbert--Schmidt. \end{lemma}

\begin{proposition} Consider an admissible curve $ \bar{M} $ and a bundle $ {\mathcal L} $ defined by
gluing conditions $ \psi_{i} $. Suppose that for any $ N>0 $
\begin{equation}
\sum_{i\not=j}\left(\max _{R_{i}} \left|\psi_{i}\right|^{N}+\max _{R_{i}}\left|\psi_{i}\right|^{-N}\right)e^{-l_{ij}} < \infty.
\notag\end{equation}
Then $ {\mathcal L} $ is strongly Hilbert--Schmidt. \end{proposition}

\begin{proof} Indeed, this is a direct corollary of Lemma~\ref{lm5.160}. \end{proof}

Topology on $ {\mathfrak L} $ defines a topology on the set of bounded bundles.
Note that the latter topology is very easy to describe. Fix a number $ M>1 $,
and consider the subset of bundles with $ |\psi_{i}|_{{\bold R}}<M $, $ i\in I $. Then the topology
on this subset is the topology of direct product. This topology is
important for the description of Jacobian in Section~\ref{s9.70}.

\subsection{Multiplicators, equivalence and Jacobians }\label{s8.7}\myLabel{s8.7}\relax  Consider what can play
a r\^ole of a mapping $ {\mathfrak m}\colon {\mathcal L}_{1} \to {\mathcal L}_{2} $ between two bundles on $ \bar{M} $ defined by
gluing conditions. Inside the ``smooth'' part of the curve such mapping
should be a multiplication by a section $ a $ of a holomorphic bundle, thus
$ {\mathfrak m} $ is determined by a function $ a $ which is holomorphic inside $ {\mathbb C}P^{1}\smallsetminus\overline{\bigcup K_{i}} $.

We start with discussing heuristics for the properties of the
function $ a $. On one hand, $ {\mathfrak m} $ should send holomorphic sections to
holomorphic one. If $ \bar{\partial}a $ is defined in a neighborhood of the infinity $ M_{\infty} $,
and is not 0, then by Leibniz rule $ {\mathfrak m} $ will not send holomorphic sections
of $ {\mathcal L}_{1} $ to holomorphic section of $ {\mathcal L}_{2} $.

On the other hand, $ {\mathfrak m} $ should send $ H^{1/2} $-sections to $ H^{1/2} $-sections.
Since $ 1\in H^{1/2} $, $ a\in H^{1/2} $. Moreover, if $ {\mathfrak m} $ is bounded, multiplication by $ a $
should send $ L_{2} $-sections of $ \Omega^{1/2}\left(\bigcup_{i}\partial K_{i}\right)\otimes\mu $ to $ L_{2} $-sections. Since the only
multiplicators in $ L_{2} $ are essentially bounded functions, we conclude that
the restriction that $ a $ is bounded on the domain of definition and is in
$ H^{1/2} $ looks like a particularly good candidate.

\begin{proposition} \label{prop8.140}\myLabel{prop8.140}\relax  Consider a function $ a $ which is holomorphic
inside
$ {\mathbb C}P^{1}\smallsetminus\overline{\bigcup K_{i}} $. Identify $ {\mathcal H}\left({\mathbb C}P^{1},\left\{K_{i}\right\}\right) $ with a subspace of the space of
analytic functions on $ {\mathbb C}P^{1}\smallsetminus\overline{\bigcup K_{i}} $. Suppose that the disks $ K_{i} $ are well
separated. Multiplication by $ a $ preserves $ {\mathcal H}\left({\mathbb C}P^{1},\left\{K_{i}\right\}\right) $ if and only if $ a $
is a restriction of an element $ \widetilde{a}\in H^{1/2}\left({\mathbb C}P^{1}\smallsetminus\bigcup K_{i}\right) $ such that
$ \bar{\partial}\widetilde{a}=0\in H^{-1/2}\left({\mathbb C}P^{1}\smallsetminus\bigcup K_{i}\right) $ and the restriction $ \widetilde{r}\left(\widetilde{a}\right)\in L_{2}\left(\bigcup\partial K_{i}\right) $ of $ \widetilde{a} $ to the
boundary is essentially bounded. \end{proposition}

\begin{proof} Indeed, we may suppose that $ \infty\in K_{i} $ for some $ i\in I $. Let us show the
``only if'' part first. Considering $ a\cdot dz^{1/2} $ we see that if $ a $ is a
multiplicator in $ {\mathcal H}\left({\mathbb C}P^{1},\left\{K_{i}\right\}\right) $, then $ \widetilde{a} $ satisfying the first two conditions
of the proposition exists. Since the operator of multiplication by $ a $ is
automatically bounded in $ {\mathcal H}\left({\mathbb C}P^{1},\left\{K_{i}\right\}\right) $, we see that the boundary value of
$ \widetilde{a} $ is essentially bounded by the norm $ M $ of this operator.

Indeed, if $ \varepsilon>0 $ and the essential supremum of $ \widetilde{r}\left(\widetilde{a}\right) $ is bigger than
$ M/\left(1-\varepsilon\right) $, then there is an arc in $ \bigcup\partial K_{i} $ such that $ |\widetilde{r}\left(\widetilde{a}\right)|>M/\left(1-\varepsilon\right) $ on a
subset of this arc of relative measure greater than $ 1-\varepsilon $. Taking
$ \|a\cdot\frac{dz^{1/2}}{z-z_{0}}\|_{L_{2}\left(\bigcup\partial K_{i}\right)} $ with $ z_{0} $ close to this arc and inside $ \bigcup K_{i} $, we
obtain a contradiction.

The proof of the ``if'' part consists of three parts. First, let us
show that if $ \alpha\in{\mathcal H}\left({\mathbb C}P^{1},\left\{K_{i}\right\}\right) $, then $ a\alpha $ is a generalized half-form correctly
defined up to addition of a half-form with support in $ \bigcup\partial K_{i} $ (here we say
that the support of a generalized function $ \beta $ is in a set $ U $---not
necessarily closed---if it is a weak limit of generalized functions $ \beta_{n} $
such that $ \operatorname{Supp}\beta_{n}\subset U $). Indeed, for any manifold $ M $ the formula
$ \left<\alpha\beta,\varphi\right>\buildrel{\text{def}}\over{=}\left<\alpha,\varphi\beta\right> $, $ \varphi\in{\mathcal D}\left(M\right) $, $ \alpha\in H^{s}\left(M\right) $, $ \beta\in H^{-s}\left(M,\Omega^{\text{top}}\right) $, shows that there is
a natural pairing $ \left(\alpha,\beta\right) \mapsto \alpha\beta $ of $ H^{s}\left(M\right) $ with $ H^{-s}\left(M,\Omega^{\text{top}}\right) $ with values in
$ {\mathcal D}'\left(M\right) $.

This pairing is weakly bicontinuous, moreover, for any
smooth vector field $ v\in\operatorname{Vect}\left(M\right) $ the Leibniz identity
\begin{equation}
v\left(\alpha\beta\right)=\left(v\alpha\right)\beta+\alpha\left(v\beta\right)
\notag\end{equation}
holds if $ \alpha\in H^{s}\left(M\right) $, $ \beta\in H^{1-s}\left(M,\Omega^{\text{top}}\right) $. Since both $ a $ and $ \alpha $ are of smoothness
$ H^{1/2} $, we see that $ a\alpha $ is indeed a generalized function defined with the
described above ambiguity (take $ s=1/2 $, and make appropriate changes to
adjust the above discussion to half-forms). Moreover, $ \bar{\partial}\left(a\alpha\right) =\left(\bar{\partial}a\right)\alpha +
a\left(\bar{\partial}\alpha\right) $, thus $ \bar{\partial}\left(a\alpha\right)=0 $ as a generalized function defined up to addition
of a function with support in $ \bigcup\partial K_{i} $. What remains to prove is that
$ \bar{\partial}\left(a\alpha\right)\in H^{-1/2}\left({\mathbb C}P^{1},\omega^{1/2}\otimes\bar{\omega}\right) $ after an appropriate choice of continuation of $ a\alpha $
to $ {\mathbb C}P^{1} $.

As a second step, fix $ i\in I $ and show that $ a\alpha $ has an appropriate
extension into $ K_{i} $. Consider restriction of $ a\alpha $ to a small collar outside
$ \partial K_{i} $.

\begin{lemma} Consider two concentric circles $ K\subset\widetilde{K}\subset{\mathbb C} $, and a holomorphic
half-form $ \beta $ in $ \widetilde{K}\smallsetminus K $. Consider a $ L_{2} $-half-form $ B $ on $ \partial K $, and suppose that
Laurent coefficients of $ \beta $ coincide with Fourier coefficients of $ B $. Then
$ \beta $ has an $ H^{1/2} $-continuation into $ K $ with the $ H^{1/2} $-norm being $ O\left(\|B\|_{L_{2}}\right) $. \end{lemma}

\begin{proof} Indeed, we know that positive Fourier coefficients of $ B $ are
$ O\left(e^{-\varepsilon k}\right) $, and negative are in $ l_{2} $. Consider positive and
negative parts of $ \beta $ separately. The positive part automatically continues
into $ K $, and the bound is a corollary of results of Section~\ref{s2.70}.

Consider now the negative part. We can suppose that $ K=\left\{z \mid |z|<1 \right\} $,
and $ \beta $ is $ \sum B_{k}z^{-k}dz^{1/2} $, $ \left(B_{k}\right)\in l_{2} $. Extend $ \beta $ into $ K $ as $ \sum B_{k}\bar{z}^{k}dz^{1/2} $. Let us
show that this extension satisfies the conditions of the lemma.

It is sufficient to show that the $ H^{-1/2} $-norm of $ \bar{\partial}\beta=\sum kB_{k}\bar{z}^{k-1}\vartheta_{K_{1}}dz^{1/2}d\bar{z} $ is
bounded (here $ \vartheta_{K_{1}} $ is the characteristic function of $ K_{1} $). Since different
components of this form are perpendicular in $ H^{-1/2} $, it is sufficient to bound
the norm of one component. To do this it is sufficient to apply the methods
of Section~\ref{s2.70}. {}\end{proof}

Application of this lemma to $ a\alpha $ shows that $ a\alpha $ may be extended into
every disk $ K_{i} $, moreover, that after these extensions the norms
$ \|\bar{\partial}\left(a\alpha\right)|_{K_{i}}\|_{H^{-1/2}} $ form a sequence in $ l_{2} $. Indeed, $ \widetilde{r}_{i}\left(a\right) $ is a bounded
function on $ \partial K_{i} $, thus $ \widetilde{r}_{i}\left(a\right)\widetilde{r}_{i}\left(\alpha\right) $ has its $ L_{2} $-norm bounded by $ L_{2} $-norm of
$ \widetilde{r}_{i}\left(\alpha\right) $, and the latter norms (for different $ i $) form a sequence in $ l_{2} $.

Third, we need to show that the above extensions can be glued
together to an extension to $ {\mathbb C}P^{1} $. Consider an arbitrary
generalized-function-extension of $ a\alpha $ to $ {\mathbb C}P^{1} $ and a generalized function
$ \bar{\partial}\left(a\alpha\right) $. It is a generalized function with support in the disjoint union
$ \bigcup K_{i} $, thus the components $ \bar{\partial}\left(a\alpha\right)|_{K_{i}} $ are well defined generalized forms $ b_{i} $,
$ \operatorname{Supp} b_{i}\subset K_{i} $. On the other hand, above we constructed an $ H^{1/2} $-extension of
$ a\alpha $ to $ K_{i} $, let $ \overset{\,\,{}_\circ}{b}_{i} $ be $ \bar{\partial}\left(a\alpha\right) $ obtained from this extension.

\begin{lemma} $ b_{i}-\overset{\,\,{}_\circ}{b}_{i}=\bar{\partial}\beta_{i} $, $ \beta_{i} $ being a generalized function with support in $ K_{i} $. \end{lemma}

This lemma is an obvious corollary of the fact that
restriction of $ \bar{\partial}^{-1}\left(b_{i}-\overset{\,\,{}_\circ}{b}_{i}\right) $ to a collar around $ \partial K_{i} $ can be holomorphically
extended inside of $ K_{i} $. Since both series $ \sum b_{i} $ and $ \sum\overset{\,\,{}_\circ}{b}_{i} $ converge
(one in $ {\mathcal D}' $, another in $ H^{-1/2} $), we conclude that $ \sum\left(b_{i}-\overset{\,\,{}_\circ}{b}_{i}\right) $ converges in $ {\mathcal D}' $,
thus $ \sum\beta_{i} $ converges in $ {\mathcal D}' $.

We conclude that $ \beta=a\alpha-\sum\beta_{i} $ is a generalized function such that
\begin{enumerate}
\item
$ \bar{\partial}\beta $ has support in $ \bigcup K_{i} $;
\item
$ \beta $ coincides with $ a\alpha $ inside $ {\mathbb C}P^{1}\smallsetminus\overline{\bigcup K_{i}} $;
\item
$ \bar{\partial}\beta $ is in $ H^{-1/2} $.
\end{enumerate}
Last condition implies $ \beta\in H^{1/2} $, which shows that $ a\alpha $ can be extended to
an element of $ H^{1/2}\left({\mathbb C}P^{1},\omega^{1/2}\right) $, thus $ a\alpha\in{\mathcal H}\left({\mathbb C}P^{1},\left\{K_{i}\right\}\right) $. This finishes the
proof of Proposition~\ref{prop8.140}. {}\end{proof}

\begin{definition} Call a bounded operator $ {\mathcal M}\colon {\mathcal H}\left({\mathbb C}P^{1},\left\{K_{i}\right\}\right) \to {\mathcal H}\left({\mathbb C}P^{1},\left\{K_{i}\right\}\right) $ a
{\em multiplicator}, if for some point $ z_{0}\in{\mathbb C}P^{1}\smallsetminus\overline{\bigcup K_{i}} $ there is a formal power
series $ \nu $ at $ z_{0} $ such that for any $ f\in{\mathcal H}\left({\mathbb C}P^{1},\left\{K_{i}\right\}\right) $ formal power series $ {\mathcal M}f $
and $ \nu f $ coincide. \end{definition}

(This is just a formal way to say that $ {\mathcal M} $ is a bounded operator of
multiplication by a holomorphic function.)

\begin{definition} Let $ {\mathcal H}^{\infty}\left({\mathbb C}P^{1},\left\{K_{i}\right\}\right)\subset H^{1/2}\left({\mathbb C}P^{1}\smallsetminus\bigcup K_{i}\right) $ consists of functions $ f $ such
that the restriction $ \widetilde{r}\left(f\right) $ to $ \bigcup\partial K_{i} $ (which is automatically in $ L_{2}\left(\bigcup\partial K_{i}\right) $)
is essentially bounded. Define a norm on this space by taking the
essential maximum of $ |\widetilde{r}\left(f\right)| $. \end{definition}

\begin{amplification} Multiplicators form an algebra. The set of multiplicators
coincides with $ {\mathcal H}^{\infty}\left({\mathbb C}P^{1},\left\{K_{i}\right\}\right) $. One can choose a norm in $ {\mathcal H}\left({\mathbb C}P^{1},\left\{K_{i}\right\}\right) $ in such
a way that the operator norm of any multiplicator coincides with the
$ {\mathcal H}^{\infty} $-norm. \end{amplification}

\begin{proof} The only statement which needs a proof is the last one, and the
norm in question is the norm induced by the inclusion into
$ \bigoplus_{l_{2}}L_{2}\left(\partial K_{i},\Omega^{1/2}\otimes\mu\right) $. \end{proof}

\begin{corollary} Let $ z_{0}\in{\mathbb C}P^{1}\smallsetminus\overline{\bigcup K_{i}} $, $ f\in{\mathcal H}^{\infty}\left({\mathbb C}P^{1},\left\{K_{i}\right\}\right) $. Then $ |f\left(z_{0}\right)|\leq \|f\|_{{\mathcal H}^{\infty}} $. \end{corollary}

\begin{proof} One can suppose that $ \infty\in K_{i} $ for some $ i\in I $, so that $ dz^{1/2}\in{\mathcal H} $. Since
$ L_{2} $-norm on $ {\mathbb C}P^{1}\smallsetminus\overline{\bigcup K_{i}} $ is majorated by $ H^{1/2} $-norm, we see that
$ \|f^{n}dz^{1/2}\|_{L_{2}}\leq C\cdot\|f\|_{{\mathcal H}^{\infty}}^{n} $, $ n\in{\mathbb N} $. Hence $ |f\left(z_{0}\right)|^{n}\leq C\cdot\|f\|_{{\mathcal H}^{\infty}}^{n} $, thus $ |f\left(z_{0}\right)|\leq\|f\|_{{\mathcal H}^{\infty}} $.
\end{proof}

Last conditions on $ {\mathfrak m} $ is that it should preserve the gluing conditions.
In particular, if $ {\mathcal L}_{1} $ is defined by gluing conditions $ \psi_{i} $, and $ {\mathcal L}_{2} $ by $ \xi_{i} $,
then
\begin{equation}
\xi_{i}=\frac{a\circ\varphi_{i}}{a}\psi_{i}.
\label{equ7.65}\end{equation}\myLabel{equ7.65,}\relax 
\begin{definition} Say that linear bundles $ {\mathcal L}_{1} $, $ {\mathcal L}_{2} $ defined by gluing conditions
$ \left\{\psi_{i}\right\} $ and $ \left\{\xi_{i}\right\} $ are {\em bounded-equivalent\/} if there exists a function
$ a\in{\mathcal H}^{\infty}\left({\mathbb C}P^{1},\left\{K_{i}\right\}\right) $ such that $ \xi_{i}=\frac{a\circ\varphi_{i}}{a}\psi_{i} $ and $ a^{-1}\in{\mathcal H}^{\infty}\left({\mathbb C}P^{1},\left\{K_{i}\right\}\right) $. \end{definition}

\begin{remark} It is obvious that a bundle which is bounded-equivalent to a
bounded bundle is bounded itself. Moreover, if a bundle is
bounded-equivalent to a (strongly) Hilbert--Schmidt bundle, it is
(strongly) Hilbert--Schmidt itself. \end{remark}

\begin{definition} The {\em Jacobian\/} is the set of equivalence classes of admissible
bundles of degree 0. The {\em bounded Jacobian\/} is the subset of Jacobian which
consists of classes of bounded bundles, similarly for {\em(strongly)
Hilbert\/}--{\em Schmidt Jacobian}, and {\em real Jacobian}. {\em Constant Jacobian\/} is formed
from classes of bundles defined by constant gluing functions, similarly
one can define different flavors of constant Jacobians. \end{definition}

Multiplication by an appropriate rational function with zeros and
poles inside $ \bigcup\overset{\,\,{}_\circ}{K}_{i} $ shows that

\begin{proposition} Any admissible bundle of degree 0 is bounded-equivalent
to a bundle with all the gluing functions $ \psi_{i} $ of index 0. \end{proposition}

\subsection{Divisors }\label{s8.80}\myLabel{s8.80}\relax  Consider a model $ \left({\mathbb C}P^{1}, \left\{K_{\bullet}\right\}, \left\{\varphi_{\bullet}\right\}\right) $ of a curve $ \bar{M} $,
and a
rational function $ a $ on $ {\mathbb C}P^{1} $ with a divisor $ D $. Suppose that the part of $ D $
inside $ \bigcup R_{i} $ is invariant w.r.t. $ \varphi_{\bullet} $, and $ D $ does not intersect with the
infinity $ M_{\infty}\subset{\mathbb C}P^{1} $ of the curve $ M $. If the bundle $ {\mathcal L} $ with gluing data $ \left\{\psi_{i}\right\} $ is
admissible, so is the bundle $ L\left(D\right) $ given by the gluing data
\begin{equation}
\xi_{i}=\frac{a\circ\varphi_{i}}{a}\psi_{i}.
\notag\end{equation}
Obviously, $ \deg {\mathcal L}\left(D\right)=\deg {\mathcal L}+\deg ' D $, if $ \deg ' D $ is the degree of the part of $ D $
inside $ M\smallsetminus\bigcup R_{i} $ plus half the degree of the part of $ D $ inside $ \bigcup R_{i} $.
Moreover, if we change the part of $ D $ inside $ \bigcup\overset{\,\,{}_\circ}{K}_{i} $, the bundle $ {\mathcal L}\left(D\right) $ will
change to a bounded-equivalent one.

In particular, to any finite subset $ D_{0} $ of $ M $ (with integer
multiplicities fixed for any point of $ D_{0} $) we associate a transformation $ {\mathcal L}
\mapsto {\mathcal L}\left(D_{0}\right) $ where the right-hand side is defined up to equivalence. Note
that if we fix a point $ Z\in\bigcup\overset{\,\,{}_\circ}{K}_{i} $, then we can complete any divisor $ D_{0} $ on $ M $
to a divisor on $ {\mathbb C}P^{1} $ of degree 0 by adding some multiple of $ Z $, thus one
can define $ {\mathcal L}\left(D_{0}\right) $ uniquely. Moreover, if $ D_{0} $ is of degree 0, then $ {\mathcal L}\left(D_{0}\right) $
does not depend on the choice of $ Z $.

Consider now two sequences of points $ \left(x_{k}\right),\left(y_{k}\right)\subset{\mathbb C}P^{1} $. Let
$ a_{k}\left(z\right)=\frac{z-x_{k}}{z-y_{k}} $, and $ a\left(z\right)=\prod_{k}a_{k}\left(z\right) $. If $ x_{k} $ is sufficiently close to $ y_{k} $,
and both these points are in appropriate neighborhood of $ K_{k} $,
then the infinite product converges and defines a bounded function on
$ {\mathbb C}P^{1}\smallsetminus\overline{\bigcup\overset{\,\,{}_\circ}{K}_{i}} $. We conclude that it is possible to consider also some
``infinite'' divisors (of finite degree) on $ M $.

Consider now a finite divisor $ D\subset{\mathbb C}P^{1} $ such that $ D\cap M_{\infty} $ consists of one
point $ z_{0}\in M_{\infty} $ (with some multiplicity). Changing $ D $ to $ D-\left(\deg  D\right)\cdot Z $, we can
consider $ D $ as a divisor of a rational function $ a $, thus $ {\mathcal L}\left(D\right) $ is correctly
defined by gluing conditions conditions~\eqref{equ7.65}. For $ {\mathcal L}\left(D\right) $ to be
admissible for any such $ D $ and any admissible $ {\mathcal L} $ it is sufficient that for
some constant $ C $
\begin{equation}
\operatorname{dist}\left(z_{0},\operatorname{center}\left(K_{i'}\right)\right) \leq C\cdot\operatorname{dist}\left(z_{0},\operatorname{center}\left(K_{i}\right)\right),\qquad i\in I
\notag\end{equation}
(since the disks $ K_{i} $ are well-separated). On the other hand, suppose that
all the bundles with gluing functions
\begin{equation}
\Psi_{i}=\left(\frac{\operatorname{dist}\left(z_{0},\operatorname{center}\left(K_{i'}\right)\right)}{\operatorname{dist}\left(z_{0},\operatorname{center}\left(K_{i}\right)\right)}\right)^{N},\qquad n\in{\mathbb Z},
\notag\end{equation}
are Hilbert--Schmidt. Then by $ \log  $-convexity the bundle $ L\left(D\right) $ is strongly
Hilbert--Schmidt if $ {\mathcal L} $ is strongly Hilbert--Schmidt.

\begin{definition} Say that the point $ z_{0}\in M_{\infty} $ is {\em bounded\/} if for any $ N\in{\mathbb Z} $ the sheaf
with constant gluing functions
\begin{equation}
\Psi_{i}=\left(\frac{\operatorname{dist}\left(z_{0},\operatorname{center}\left(K_{i'}\right)\right)}{\operatorname{dist}\left(z_{0},\operatorname{center}\left(K_{i}\right)\right)}\right)^{N},\qquad n\in{\mathbb Z}\text{, }i\in I,
\notag\end{equation}
is Hilbert--Schmidt. \end{definition}

We see that if the point $ z_{0}\in M_{\infty} $ is bounded, then the function
$ a_{N}\left(z\right)=\left(\frac{z-z_{0}}{z-Z}\right)^{N} $ defines a strongly Hilbert--Schmidt sheaf of degree
1. In particular, one can consider finite divisors on $ \bar{M} $ which consist of
points of $ M $ and bounded points on $ M_{\infty} $. Similarly, one can also consider
some infinite divisors on $ \bar{M} $.

\subsection{Universal Grassmannian and bundles }\label{s8.90}\myLabel{s8.90}\relax 

\begin{definition} Let $ V $ is the vector space of sequences $ \left(a_{k}\right) $, $ k\in{\mathbb Z} $, such that
$ a_{k}=0 $ for $ k\ll0 $. The vector space $ V $ carries a natural topology of inductive
limit of projective limits. Let
\begin{equation}
V_{+}=\left\{\left(a_{k}\right) \mid a_{k}=0\text{ if }k<0\right\},\qquad V_{-}=\left\{\left(a_{k}\right)\in V \mid a_{k}=0\text{ if }k\geq0\right\}.
\notag\end{equation}
Say that the vector subspace $ W\subset V $ is {\em admissible\/} if $ W\cap V_{-} $ is
finite-dimensional, and $ W+V_{-} $ is
closed and has a finite codimension in $ V $. Let {\em universal Grassmannian\/} $ {\mathcal G} $ be
the set of admissible vector subspaces with natural topology. \end{definition}

Consider an admissible pair $ \left(\bar{M},{\mathcal L}\right) $ and a smooth point $ z_{0}\in M $. Pick up a
coordinate system $ z $ in neighborhood of $ z_{0} $, and a half-form $ f $ defined in
the same neighborhood. Now any global section of $ {\mathcal L} $ may be written as
$ g\left(z\right)f\left(z\right) $ with $ g\left(z\right) $ being a holomorphic function which is correctly
defined in the neighborhood of $ z_{0} $.

Call $ f $ a {\em local section\/} of $ {\mathcal L} $ near $ z_{0} $. If $ {\mathcal L}' $ is an equivalent to $ {\mathcal L} $
bundle, and the equivalence is given by multiplication by $ a $, we obtain a
local section $ af $ of $ {\mathcal L}' $. Say that local sections $ f $ and $ af $ are {\em equivalent}.

Consider now a bundle $ {\mathcal L}\left(k\cdot z_{0}\right) $, defined using a point $ Z\in\overset{\,\,{}_\circ}{K}_{i} $ (as in the
previous section). A section $ h\left(\zeta\right) $, $ \zeta\in{\mathbb C}P^{1} $, of this bundle can be
identified (via multiplication by $ \left(\frac{\zeta-Z}{\zeta-z_{0}}\right)^{k} $) with a ``meromorphic
section'' of $ {\mathcal L} $, i.e., one can write it as $ g\left(z\right)f\left(z\right) $, $ g\left(z\right) $ being a
meromorphic function correctly defined in the neighborhood of $ z_{0} $.
Moreover, $ g\left(z\right) $ has poles only at $ z_{0} $. One can momentarily see that $ g\left(z\right) $
does not depend on the choice of the point $ Z $, more precise, for a
different choice of $ Z $ there is a different choice of $ h $ which gives the
same $ g\left(z\right) $.

Associate to any such function $ g\left(z\right) $ the sequence of its Laurent
coefficients. Consider this sequence as an element of $ V $. Let $ W $ be the
vector subspace of $ V $ spanned by all possible functions $ g\left(z\right) $ for
bundles $ {\mathcal L}\left(k\cdot z_{0}\right) $, $ k\in{\mathbb Z} $. The Riemann--Roch theorem momentarily implies that
$ W $ is admissible. Indeed, the condition that $ W\cap V_{-} $ is
finite-dimensional means that $ {\mathcal L} $ has a finite-dimensional space of
global sections. The condition on $ W+V_{-} $ is implied by the following fact:

\begin{lemma} Consider a admissible semibounded bundle $ {\mathcal L} $ on a curve $ \bar{M} $, and a
point $ z_{0}\in M $. Then
\begin{enumerate}
\item
for big enough $ k $ the bundle $ {\mathcal L}\left(-k\cdot z_{0}\right) $ has no sections.
\item
for big enough $ k $ the bundle $ {\mathcal L}\left(\left(k+1\right)\cdot z_{0}\right) $ has one more section than
$ {\mathcal L}\left(k\cdot z_{0}\right) $.
\end{enumerate}
\end{lemma}

\begin{proof} The second statement is a corollary of the first one, of
Riemann--Roch theorem and duality. The first one is obvious, since the
global sections of $ {\mathcal L}\left(-k\cdot z_{0}\right) $ are naturally identified with global
sections of $ {\mathcal L} $ which have a zero of $ k $-th order at $ z_{0} $. \end{proof}

We obtained

\begin{proposition} To each admissible pair $ \left(\bar{M},{\mathcal L}\right) $ with a fixed smooth point
$ z_{0}\in M $, a coordinate system near $ z_{0} $, and a local section $ f $ of $ {\mathcal L} $ one can
associate a point $ W\in{\mathcal G} $. If we change $ {\mathcal L} $ and $ f $ to an equivalent bundle with
a local section, $ W $ does not change. \end{proposition}

One can generalize this proposition to some points at infinity. If
$ z_{0}\in M_{\infty} $, define a {\em local section\/} of $ {\mathcal L} $ near $ z_{0} $ in the same way as above, i.e.,
as a non-zero section of $ \omega^{1/2} $ near $ z_{0} $. In what follows we use only the
$ \infty $-jet of this section.

\begin{amplification} Consider an admissible pair $ \left(\bar{M},{\mathcal L}\right) $ with a fixed point at
infinity $ z_{0}\in M_{\infty} $, a coordinate system near $ z_{0} $, and a local section $ f $ of $ {\mathcal L} $.
Suppose that $ {\mathcal L} $ is Hilbert--Schmidt, and $ z_{0} $ is bounded and of class $ C^{\infty} $. To
this data we can associate a point $ W\in{\mathcal G} $. If we change $ {\mathcal L} $ and $ f $ to an
equivalent bundle with a local section, $ W $ does not change. \end{amplification}

\begin{proof} Since $ z_{0} $ is bounded, and $ {\mathcal L} $ is Hilbert--Schmidt, $ {\mathcal L}\left(k\cdot z_{0}\right) $ is
Hilbert--Schmidt too, thus the Riemann--Roch theorem is applicable. Since
$ z_{0} $ is of class $ C^{\infty} $, any section of $ {\mathcal L}\left(k\cdot z_{0}\right) $ has a (formal) Taylor series
near $ z_{0} $, thus the corresponding section of $ {\mathcal L} $ has a (formal) Laurent
series.

What remains to be proved is the fact that equivalence between
bundles can be pushed to infinity points of class $ C^{\infty} $. Note that one can
associate an element of the generalized Hardy space to any bounded
function $ a $. Indeed, suppose that $ \zeta=\infty\in K_{i} $, then $ a\cdot d\zeta^{1/2} $ is an element of
$ {\mathcal H} $. In particular, $ a $ has an asymptotic expansion near $ z_{0} $, thus
multiplication by $ a $ maps (formal) Laurent series at $ z_{0} $ to themselves. \end{proof}

\section{Structure of Jacobian }\label{h9}\myLabel{h9}\relax 

\subsection{Constant Jacobian }\label{s7.90}\myLabel{s7.90}\relax  Consider the involution $ '\colon I\to I $. It defines a
transposition matrix $ t=\left(t_{ij}\right) $, $ i,j\in I $, $ t_{ij}=\delta_{ij'} $. We use the results of toy
theory (see Section~\ref{h35}) to obtain the following

\begin{theorem} \label{th7.25}\myLabel{th7.25}\relax  Let $ l_{ij} $ be pairwise conformal distances between disks
$ K_{i} $. Suppose that the matrix $ {\mathcal R}=\left(e^{-l_{ij}}-\delta_{ij}\right) $ defines a compact mapping $ l_{2} \to
l_{2} $ and a bounded mapping $ l_{\infty} \to l_{\infty} $, and that for some $ N>0 $ the matrix $ \left(t{\mathcal R}\right)^{N} $
defines a mapping $ l_{\infty} \to l_{2} $. Then for any admissible curve obtained by
gluing circles $ \partial K_{i} $ the constant bounded Jacobian coincides with the
bounded Jacobian. \end{theorem}

\begin{proof} Indeed, to show this we need to show that for any cycle $ \left\{\psi_{i}\right\}_{i\in I} $
with $ \deg \psi_{i}=0 $ and $ |\psi_{i}|<C $ one can find a bounded collection of constants $ c_{i} $
such that $ c_{i}\psi_{i} =\frac{a\circ\varphi_{i}}{a} $ for some $ a $ such that $ a,a^{-1}\in{\mathcal H}^{\infty}\left({\mathbb C}P^{1},\left\{K_{i}\right\}\right) $.

Taking logarithms, we see that it is sufficient to show that the mapping
\begin{equation}
{\mathcal J}_{b}\colon f+C \mapsto \left(f|_{\partial K_{j}}-\varphi^{*}\left(f|_{\partial K_{j'}}\right)+C_{j}\right)_{j\in I_{+}}
\notag\end{equation}
from bounded analytic functions (modulo constants) to functions on
boundary (modulo
constants) with bounded $ \pm $-parts is surjective. In Section~\ref{s35.40} we have
seen that (given the first condition of the theorem) a similar mapping
\begin{equation}
{\mathcal J}\colon {\mathcal H}^{\left(1\right)}/\operatorname{const} \to \bigoplus\Sb l_{2} \\ j\in I_{+}\endSb H^{1/2}\left(\partial K_{j}\right)/\operatorname{const}
\notag\end{equation}
is a bounded mapping of index 0. We are going to prove the
surjectivity by using a combination of following lemmas:

\begin{lemma} \label{lm8.31}\myLabel{lm8.31}\relax  If $ \left(e^{-l_{ij}}-\delta_{ij}\right) $ gives a compact operator $ l_{2} \to l_{2} $, then the
mapping $ {\mathcal J} $ is a bijection. \end{lemma}

\begin{proof} Since we know the index of $ {\mathcal J} $, it is sufficient to show that
$ \operatorname{Ker}{\mathcal J}=0 $. Let $ f\in{\mathcal H}^{\left(1\right)}/\operatorname{const} $ and $ {\mathcal J}f=0 $. We are going to show that $ \|f\|_{H^{1}}=0 $.
Since $ \partial $ is elliptic, and $ \bar{\partial}f=0 $, it is sufficient to show that $ \|\partial f\|_{L_{2}}=0 $,
i.e., that $ \int\partial f\wedge\bar{\partial}\bar{f}=0 $, the integral is taken over $ {\mathbb C}P^{1}\smallsetminus\bigcup K_{i} $. Denote the
domain of integration by $ S $. Note that the above integral is well defined,
since $ \partial f\in H^{0}\left(S\right)=L_{2}\left(S\right) $.

First of all, since $ \wedge $-product of any two sections of $ \omega $ is 0, one can
add $ \partial\bar{f} $ to $ \bar{\partial}\bar{f} $ without changing the integral, thus it is sufficient to show
that $ \int_{S}\partial f\wedge d\bar{f} = $ 0. Additionally, since $ \bar{\partial}f=0 $ on $ S $, one can change $ \partial f $ to $ df $
without changing the integral, thus it is sufficient to show that
$ \int_{S}df\wedge d\bar{f}=0 $. Second, use the duality identity $ \int_{S}d\alpha=\int_{\partial S}\alpha $. Applying it
(at first formally), we get
\begin{align} \int_{S}df\wedge d\bar{f} & = \int_{S}d\left(f\,d\bar{f}\right) = \int_{\partial S}f\,d\bar{f} = -\sum_{i\in I}\int_{\partial K_{i}}f\,d\bar{f}
\notag\\
& = -\sum_{i\in I_{+}}\left(\int_{\partial K_{i}}f\,d\bar{f} + \int_{\partial K_{i'}}f\,d\bar{f}\right).
\notag\end{align}
We want to show that $ \int_{\partial K_{i}}f\,d\bar{f} + \int_{\partial K_{i'}}f\,d\bar{f}=0 $. Indeed,
$ \int_{\partial K_{i'}}f\,d\bar{f}=-\int_{\partial K_{i}}\varphi_{i'}^{*}\left(f\,d\bar{f}\right) $, the sign appears since $ \varphi_{i} $ changes the
orientations of $ \partial K_{\bullet} $. Now the above equality becomes obvious, since
$ f\circ\varphi_{i}=f+C_{i} $ when both sides are defined, thus the sum of the integrals is
$ -C_{i}\int_{\partial K_{i}}d\bar{f}=0 $.

What remains to be proved is that one can indeed apply the formula
$ \int_{S}d\alpha=\int_{\partial S}\alpha $ in our situation, when $ \alpha $ is not smooth, and $ S $ has non-smooth
boundary. Consider an $ H^{1} $-extension $ g $ of $ f $ to $ {\mathbb C}P^{1} $. Since $ dg\wedge d\bar{g}\in L_{1}\left(\Omega^{\text{top}}\right) $,
the integral $ \int_{S}df\wedge d\bar{f} $ can be represented as
\begin{equation}
\int_{{\mathbb C}P^{1}}dg\wedge d\bar{g} - \int_{\bigcup K_{i}}dg\wedge d\bar{g} = \int_{{\mathbb C}P^{1}}dg\wedge d\bar{g} - \sum\int_{K_{i}}dg\wedge d\bar{g}.
\notag\end{equation}
It remains to prove that $ \int_{{\mathbb C}P^{1}}dg\wedge d\bar{g} = $ 0, and $ \int_{K_{i}}dg\wedge d\bar{g} = \int_{\partial K_{i}}g\,d\bar{g}. $ Note
that in both identities the boundary is already smooth (it is empty in
the first one!).

Let us prove that if $ M $ is a two-dimension manifold, $ S $ is a
compact subset of $ M $ with a smooth boundary, and $ g_{1},g_{2}\in H^{1}\left(M\right) $, then
\begin{equation}
\int_{S}dg_{1}\wedge dg_{2}=\int_{\partial S}g_{1}\,dg_{2}.
\notag\end{equation}
Here we understand the right-hand side as a natural pairing between
$ g_{1}|_{\partial S}\in H^{1/2}\left(\partial S\right) $ and $ d\left(g_{2}|_{\partial S}\right)\in H^{-1/2}\left(\partial S,\Omega_{\partial S}^{\text{top}}\right) $. However, both sides define
bounded bilinear functionals on $ H^{1}\left(M\right) $, thus it is sufficient to check
them on a dense subset $ C^{\infty}\left(M\right) $, where they are true due to de Rham theory.
\end{proof}

To continue the proof of the theorem, note that the operator $ {\mathcal J} $ is
related to the operator $ {\bold K} $ in the following way: identify $ {\mathcal H}^{\left(1\right)}/\operatorname{const} $ with
$ \bigoplus\Sb l_{2} \\ i\in I\endSb H_{-}^{1/2}\left(\partial K_{i}\right)/\operatorname{const} $. Let $ \Pi_{i} $ be the identification of $ H_{+}^{1/2}\left(\partial K_{i}\right) $ with
$ H_{-}^{1/2}\left(\partial K_{i'}\right) $ via $ \varphi^{*} $, let $ \Pi^{\left(1\right)} $ be the direct sum of such identifications.
As we have seen it in Section~\ref{s35.40}, the image of $ {\mathcal J} $ coincides with the
image of $ \Pi^{\left(1\right)}\circ{\bold K}-1 $ (which is a mapping
\begin{equation}
\bigoplus_{l_{2}}H_{-}^{1/2}\left(\partial K_{i}\right)/\operatorname{const} \to \bigoplus_{l_{2}}H_{-}^{1/2}\left(\partial K_{i}\right)/\operatorname{const},
\notag\end{equation}
but we identify $ \bigoplus_{l_{2}}H_{-}^{1/2}\left(\partial K_{i}\right)/\operatorname{const} $ with $ \bigoplus\Sb l_{2} \\ j\in I_{+}\endSb H^{1/2}\left(\partial K_{j}\right)/\operatorname{const} $ via $ \varphi_{\bullet}^{*} $).

\begin{lemma} \label{lm8.33}\myLabel{lm8.33}\relax 
\begin{enumerate}
\item
Consider two disjoint disks $ K_{1} $, $ K_{2} $ on $ {\mathbb C}P^{1} $ of conformal distance $ l $.
The operator $ {\bold K} $ with smooth kernel $ \frac{dy}{y-x} $ defines bounded operators
$ H^{1/2}\left(\partial K_{1}\right) \to L_{\infty}\left(\partial K_{2}\right) $, $ L_{\infty}\left(\partial K_{1}\right) \to L_{\infty}\left(\partial K_{2}\right) $, $ H^{1/2}\left(\partial K_{1}\right) \to H^{1/2}\left(\partial K_{2}\right) $, $ L_{\infty}\left(\partial K_{1}\right)
\to H^{1/2}\left(\partial K_{2}\right) $ with the norms being $ O\left(e^{-l}\right) $.
\item
Consider a family of disjoint disks $ \left\{K_{i}\right\} $ on $ {\mathbb C}P^{1} $ with conformal
distances $ l_{ij} $, and suppose that the matrix $ {\mathcal R}=\left(e^{-l_{ij}}-\delta_{ij}\right) $ defines bounded
mappings $ l_{2} \to l_{2} $ and $ l_{\infty} \to l_{\infty} $, and for some $ N>0 $ the matrix $ \left(t{\mathcal R}\right)^{N} $ defines
a mapping $ l_{\infty} \to l_{2} $. Then the operator $ {\bold K} $ gives bounded mappings
$ \bigoplus_{l_{2}}H^{1/2}\left(\partial K_{i}\right) \to \bigoplus_{l_{2}}H^{1/2}\left(\partial K_{i}\right) $ and $ L_{\infty}\left(\bigcup\partial K_{i}\right) \to L_{\infty}\left(\bigcup\partial K_{i}\right) $, moreover, the
operator $ \left(\Pi^{\left(1\right)}{\bold K}\right)^{N} $ gives a bounded mapping $ L_{\infty}\left(\bigcup\partial K_{i}\right) \to \bigoplus_{l_{2}}H^{1/2}\left(\partial K_{i}\right) $.
\item
In the conditions of the previous part of the lemma let
$ f_{-}\in L_{\infty}\left(\bigcup\partial K_{i}\right) $. Suppose that $ + $-parts of $ f_{-} $ on all the $ \partial K_{i} $ vanish. Then $ {\bold K}f_{-} $
is a bounded analytic function.
\end{enumerate}
\end{lemma}

Now we are able to invert the mapping $ {\mathcal J}_{b} $. Consider a fixed function
$ f $ on $ \bigcup\partial K_{i} $ with bounded $ + $-part and bounded $ - $-part. If $ f-f\circ\varphi_{i} $ is a
constant for any $ i\in I $, then $ + $-part of $ f|_{\partial K_{i}} $ can be reconstructed basing on
-part of $ f|_{\partial K_{i'}} $, thus one can identify $ f $ (modulo constants) with the
collection of $ - $-parts of $ f $ on $ \bigcup_{i\in I}\partial K_{i} $ (modulo constants).

The restrictions on $ l_{ij} $ show that the radii of disks $ K_{i} $ for a
sequence from $ l_{1} $. Hence the Cauchy kernel allows one to construct an
analytic function $ \widetilde{f} $ on $ S $ such that the $ - $-parts of $ \widetilde{f}|_{\partial K_{i}} $ and $ f|_{\partial K_{i}} $
coincide (modulo constants). Moreover, $ {\mathcal J}_{b}\widetilde{f} = \left(\Pi^{\left(1\right)}\circ{\bold K}-1\right)\left(f\right) $. Since $ f $ is
bounded, it is in $ L_{2} $, thus $ \widetilde{f}\in H^{1/2} $. Since $ \widetilde{f}|_{\bigcup\partial K_{i}} $ is bounded, $ \widetilde{f} $ is a
multiplicator.

We see that $ \Pi^{\left(1\right)}\circ{\bold K}\left(f\right) $ is a well-defined bounded function on
$ \bigcup_{j\in I_{+}}\partial K_{j} $, moreover, $ \left(\Pi^{\left(1\right)}\circ{\bold K}\right)^{N}\left(f\right)\in\bigoplus_{l_{2}}H^{1/2}\left(\partial K_{j}\right) $. In particular,
$ \left(\Pi^{\left(1\right)}\circ{\bold K}\right)^{N}\left(f\right)={\mathcal J}F $ for some function $ F\in{\mathcal H}^{\left(1\right)} $. To show that $ f $ is in the image
of $ \Pi^{\left(1\right)}\circ{\bold K}-1 $, it remains to prove is that $ f-\left(\Pi^{\left(1\right)}\circ{\bold K}\right)^{N}\left(f\right) $ is in the image
of $ \Pi^{\left(1\right)}\circ{\bold K}-1 $, what is obvious.

This finishes the proof of Theorem~\ref{th7.25}. {}\end{proof}

To demonstrate a particular case in which the assumptions of Theorem
~\ref{th7.25} hold, consider

\begin{lemma} Suppose that the matrix $ \left(a_{ij}\right) $ gives a Hilbert--Schmidt operator
$ l_{2} \to l_{2} $, i.e., $ \sum|a_{ij}|^{2}<\infty $. Then the matrix $ \left(a_{ij}^{2}\right) $ gives a bounded
operator $ l_{\infty} \to l_{2} $.

If $ a_{ij}>0 $ for any $ i,j\in I $, and the matrix $ \left(a_{ij}\right) $ gives a compact mapping
$ l_{2} \to l_{2} $, then $ \left(a_{ij}^{2}\right) $ also gives a compact mapping $ l_{2} \to l_{2} $. \end{lemma}

\begin{corollary} If $ \bar{M} $ is a Hilbert--Schmidt curve, then the bounded Jacobian
of $ \bar{M} $ coincides with the constant bounded Jacobian. \end{corollary}

In Section~\ref{s9.70} we show that under suitable restrictions on $ \bar{M} $ the
above Jacobians are in $ 1 $-to-1 correspondence with quotients of
topological vector spaces by $ {\mathbb Z} $-lattices. On the other hand, this
statement is obvious when applied to real Jacobian, since it is a direct
product of circles $ |\psi_{i}|=1 $. To construct an isomorphism of this product
with a quotient by a lattice, one can take the Hilbert space $ l_{2} $ with the
standard basis $ e_{i} $, and the lattice spanned by $ \alpha_{i}e_{i} $, with $ \alpha_{i} \to $ 0.

The description of the real part of Jacobian was first obtained in
\cite{McKTru76Hil} in the case of a real hyperelliptic curve of a special form.

\subsection{The partial period mapping }\label{s9.20}\myLabel{s9.20}\relax  Consider a theory similar to one
discussed in Section~\ref{s35.30}, but related to holomorphic $ 1 $-forms instead
of holomorphic functions. Consider a family $ \left\{K_{i}\right\} $ of disjoint disks in
$ {\mathbb C}P^{1} $.

\begin{definition} We say that a generalized-function section $ \alpha $ of the linear
bundle $ \omega $ on $ {\mathbb C}P^{1} $ is $ H^{0} $-{\em holomorphic in\/} $ {\mathbb C}P^{1}\smallsetminus\bigcup K_{i} $ if $ \alpha\in H^{0}\left({\mathbb C}P^{1}\smallsetminus\bigcup K_{i},\omega\right) =
L_{2}\left({\mathbb C}P^{1}\smallsetminus\bigcup K_{i},\omega\right) $, and $ \bar{\partial}\alpha=0\in H^{-1}\left({\mathbb C}P^{1}\smallsetminus\bigcup K_{i},\omega\otimes\bar{\omega}\right) $. Denote the the space of
$ H^{0} $-holomorphic forms in $ {\mathbb C}P^{1}\smallsetminus\bigcup K_{i} $ by $ {\mathcal H}^{\left(0\right)} $. \end{definition}

Since one can extend a $ H^{0} $-holomorphic form $ \alpha $ into any disk $ K_{i} $ by 0
without increasing $ H^{0} $-norm,
$ \alpha $ has a canonical extension to $ {\mathbb C}P^{1} $, and $ \bar{\partial}\alpha $ is concentrated on $ \bigcup\partial K_{i} $.
Moreover, if $ K_{i} $ has a neighborhood which does not intersect with other
disks, one can define a {\em boundary value\/} $ \alpha|_{\partial K_{i}} $ of $ \alpha $ on $ \partial K_{i} $, which is a
element of $ H^{-1/2}\left(\partial K_{i},\Omega_{\partial K_{i}}^{1}\right) $ (note that if $ \gamma $ is a curve in $ {\mathbb C}P^{1} $, then
$ \omega|_{\gamma}\simeq\Omega_{\gamma}^{1} $; one can get the $ H^{-1/2} $-restriction on the smoothness of the
boundary value in the same way as in Section~\ref{s2.70}). Now $ \bar{\partial}\alpha $ can be
described as the extension of $ -\alpha|_{\partial K_{i}} $ from $ \partial K_{i} $ to $ {\mathbb C}P^{1} $ by $ \delta $-function (see
~\eqref{equ3.02}). (Such an extension is a correctly defined continuous mapping
$ H^{-1/2}\left(\partial K_{i},\Omega_{\partial K_{i}}^{1}\right) \to H^{-1}\left({\mathbb C}P^{1},\omega\otimes\bar{\omega}\right) $.)

Since the norm in $ H^{1}\left({\mathbb C}P^{1}\right) $ can be described by a local (non-invariant)
formula $ \|\alpha\|_{l_{2}}^{2}+\|\alpha_{,x}\|_{l_{2}}^{2}+\|\alpha_{,y}\|_{l_{2}}^{2} $ (here $ x+iy $ is the local coordinate on
$ {\mathbb C}P^{1} $), functions with non-intersecting support are orthogonal in $ H^{1} $.
Dually, if $ U_{i} $ are pairwise non-intersecting, then the natural mapping
\begin{equation}
H^{-1}\left({\mathbb C}P^{1},\omega\otimes\bar{\omega}\right) \to \bigoplus_{l_{2}}H^{-1}\left(U_{i},\omega\otimes\bar{\omega}\right)
\notag\end{equation}
is a continuous epimorphism.

\begin{definition} Fix $ \varepsilon>0 $. Say that disks $ \left\{K_{i}\right\} $ have a {\em thickening\/} $ \left\{U_{i}\right\} $ if $ K_{i} $
have non-intersecting neighborhoods $ U_{i} $ such that a pair $ K_{i}\subset U_{i} $ is
conformally equivalent to $ \left\{|z|<1\right\}\subset\left\{|z|<1+\varepsilon\right\} $. Say that $ \left\{K_{i}\right\} $ have a {\em uniform
thickening\/} if the neighborhood $ U_{i} $ can be picked up to be concentric to $ K_{i} $
circles (assume that a metric on $ {\mathbb C}P^{1} $ is fixed). \end{definition}

(Note that the existence of uniform thickening does not depend on
the metric on $ {\mathbb C}P^{1} $.)

In particular, if $ \left\{K_{i}\right\} $ have a thickening $ \left\{U_{i}\right\} $, then the above
arguments show that $ \left(\|\alpha|_{\partial K_{i}}\|_{H^{-1/2}}\right)_{i\in I}\in l_{2} $. Note also that $ H^{-1/2}\left(\partial K_{i},\Omega_{\partial K_{i}}^{1}\right) $
is dual to $ H^{1/2}\left(\partial K_{i}\right) $, thus the subspace $ H_{\int=0}^{-1/2}\left(\partial K_{i},\Omega_{\partial K_{i}}^{1}\right) $ of forms with
vanishing integral has an $ \operatorname{PGL}\left(2,{\mathbb C}\right) $-invariant Hilbert structure (in the
same sense as in Section~\ref{s35.20}).

\begin{definition} Call the mapping $ \alpha \mapsto \left(\int_{\partial K_{i}}\alpha\right)_{i\in I} $, $ H^{0}\left({\mathbb C}P^{1}\smallsetminus\bigcup K_{i},\omega\right) \to l_{2} $ the
{\em integration mapping}. Fix an involution $ '\colon I\to I $ which interchanges two
parts of $ I $, $ I=I_{+}\amalg I_{+}' $, and fraction-linear identifications $ \varphi_{i}\colon \partial K_{i} \to \partial K_{i} $
(with the same conditions as in Section~\ref{s35.40}), then a {\em global
holomorphic form\/} $ \alpha $ is an element of $ {\mathcal H}^{\left(0\right)} $ such that $ \varphi_{i}^{*}\left(\alpha|_{\partial K_{i'}}\right)=\alpha|_{\partial K_{i}} $.
Define the {\em partial period mapping\/} $ P $ by restriction of the integration
mapping to global holomorphic forms, and taking integrals only for $ i\in I_{+} $.
Let $ \bar{M} $ be the curve obtained by gluing $ \partial K_{i} $ together via $ \varphi_{i} $. Denote by
$ \Gamma\left(\bar{M},\omega\right) $ the space of global holomorphic forms on $ \bar{M} $, i.e., the space of
global holomorphic forms compatible with $ \varphi_{\bullet} $. \end{definition}

The significant difference of this case and the case of Section
~\ref{s35.30} is that the mapping $ \bar{\partial}\colon H^{0}\left({\mathbb C}P^{1},\omega\right) \to H^{-1}\left({\mathbb C}P^{1},\omega\otimes\bar{\omega}\right) $ has no
null-space, but has $ 1 $-dimensional cokernel, thus it is not easy to
reconstruct $ \alpha $ basing on $ \bar{\partial}\alpha $ by local formulae. However, it is easy to
prove

\begin{lemma} \label{lm9.90}\myLabel{lm9.90}\relax  Suppose that $ \left\{K_{i}\right\} $ has a thickening. Then the partial period
mapping $ P $ satisfies the identity
\begin{equation}
\dim  \operatorname{Ker} P = \dim  \operatorname{Ker} {\mathcal J}\colon {\mathcal H}^{\left(1\right)}/\operatorname{const} \to \bigoplus_{l_{2}}H^{1/2}\left(\partial K_{j}\right)/\operatorname{const}.
\notag\end{equation}
\end{lemma}

\begin{proof} We claim that $ \partial $ gives a mapping from one null-space to another
one. Let $ \alpha\in\operatorname{Ker} P $. Then $ \alpha|_{\partial K_{i}}\in H^{-1/2}\left(\partial K_{i},\Omega_{\partial K_{i}}^{1}\right) $ has an antiderivative
$ f_{i}\in H^{1/2}\left(\partial K_{i}\right) $, and the harmonic extension of $ f_{i} $ into $ K_{i} $ has a bounded
$ H^{1} $-norm. Since $ \alpha $ is closed in $ U_{i}\smallsetminus K_{i} $, it has an antiderivative. Since $ \partial $ is
elliptic, the antiderivative has smoothness $ H^{1} $. It is
possible to pick up the constant in such a way that two antiderivatives
coincide on $ \partial K_{i} $. Taking de Rham differential of resulting $ H^{1} $-function, we
see that one can extend $ \alpha $ into $ K_{i} $ preserving the closeness, and the
related increase of the norm of $ \alpha $ is bounded by $ \|\alpha|_{\partial K_{i}}\|_{H^{-1/2}} $.

Repeating this operation for all the $ K_{i} $, we obtain a closed
extension of $ \alpha $ to $ {\mathbb C}P^{1} $ with the $ H^{0} $-norm bounded by $ \|\alpha\|_{H^{0}} $ (here we use
locality of $ L_{2} $-norm). Thus $ \alpha $ has a $ H^{1} $-antiderivative $ f $ on $ {\mathbb C}P^{1}\smallsetminus\bigcup K_{i} $, thus
$ \alpha=df $, hence $ \alpha=\partial f $. Obviously, $ f\in\operatorname{Ker} {\mathcal J} $.

On the other hand, if $ f\in\operatorname{Ker} {\mathcal J} $, then $ \partial f\in\operatorname{Ker} P $.\end{proof}

\begin{remark} As we have seen in Section~\ref{s7.90}, in assumptions of Theorem
~\ref{th35.45} the dimension in the lemma is 0, thus $ P $ defines an inclusion of
global holomorphic forms into $ l_{2} $. Note that we do not claim that this
inclusion is a monomorphism---estimates in Section~\ref{s9.12} show that it is
not if $ \operatorname{card}\left(I\right)=\infty $. \end{remark}

\subsection{$ A $-Periods }\label{s9.12}\myLabel{s9.12}\relax  To describe the image of the integration mapping,
consider the mapping $ \bar{\partial}^{-1}\colon H_{\int=0}^{-1}\left({\mathbb C}P^{1},\omega\otimes\bar{\omega}\right) \to H^{0}\left({\mathbb C}P^{1},\omega\right) $. It is a
continuous elliptic operator. To make the formulae localizable, fix a top
form $ \beta_{0} $ with integral 1, and extend $ \bar{\partial}^{-1} $ to any form on $ {\mathbb C}P^{1} $ by $ \bar{\partial}^{-1}\beta\buildrel{\text{def}}\over{=}
\bar{\partial}^{-1}\left(\beta-\beta_{0}\int_{{\mathbb C}P^{1}}\beta\right) $. This operator is still elliptic, but is not an
isomorphism any more.

Let $ \alpha\in{\mathcal H}^{\left(0\right)} $, consider $ \alpha|_{\partial K_{i}} $. In contrast with the cases of Section
~\ref{s35.20} and Proposition~\ref{prop5.28} the $ - $-components of $ \left(\alpha|_{\partial K_{i}}\right)_{i\in I} $ are not
arbitrary. To describe possible values of $ \left(\alpha|_{\partial K_{i}}\right)_{i\in I} $, break $ \alpha|_{\partial K_{i}} $ into
two parts, one $ \alpha_{i}^{\left(0\right)} $ with integral 0, another one $ \alpha_{i}^{\left(1\right)} $ proportional to
some fixed $ 1 $-form $ \mu_{i} $ on $ \partial K_{i} $. Use the same letters for
$ \delta $-function-extensions of these forms to $ {\mathbb C}P^{1} $. Suppose that $ 1 $-forms $ \mu_{i} $ (one
per disk boundary) are normalized to have integral 1 and have uniformly
bounded $ H^{-1/2} $-norms when $ i $ varies, and that disks $ K_{i} $ have a thickening.
Then $ \sum\alpha_{i}^{\left(1\right)} $ converges in $ \bigoplus_{l_{2}}H^{-1/2}\left(\partial K_{i}\right) $, same for $ \sum\alpha_{i}^{\left(0\right)} $. Note that if
conditions of Theorem~\ref{th35.45} hold, then $ \alpha^{\left(0\right)}=\bar{\partial}^{-1}\left(\sum\alpha_{i}^{\left(0\right)}\right)\in H^{0}\left({\mathbb C}P^{1},\omega\right) $.
Since $ \alpha=\bar{\partial}^{-1}\sum_{i}\alpha|_{\partial K_{i}} $, $ \alpha^{\left(1\right)}=\bar{\partial}^{-1}\left(\sum\alpha_{i}^{\left(1 \right)}\right)=\alpha-\bar{\partial}^{-1}\left(\sum\alpha_{i}^{\left(0\right)}\right) $ is also in $ H^{0}\left({\mathbb C}P^{1},\omega\right) $.

We see that any form $ \alpha\in{\mathcal H}^{\left(0\right)} $ can be represented as a sum $ \alpha^{\left(0\right)}+\alpha^{\left(1\right)} $,
$ \alpha^{\left(0\right)},\alpha^{\left(1\right)}\in{\mathcal H}^{\left(0\right)} $ such that integrals of $ \alpha^{\left(0\right)} $ around each disk $ K_{i} $ vanishes,
and $ \alpha^{\left(1\right)} $ is a linear combination of $ 1 $-forms $ \bar{\partial}^{-1}\mu_{i} $, $ i\in I $. Consider the
subspace $ {\mathcal H}_{0}^{\left(0\right)}\subset{\mathcal H}^{\left(0\right)} $ consisting of $ 1 $-forms which satisfy $ \int_{\partial K_{i}}\alpha=0 $ for any
$ i\in I $, and the subspace $ {\mathcal H}_{1}^{\left(0\right)} $ of forms satisfying $ \int_{\partial K_{i}}\alpha=\int_{\partial K_{i}'}\alpha $ for any $ i\in I $.
Let $ {\mathfrak c}={\mathcal H}^{\left(0\right)}/{\mathcal H}_{0}^{\left(0\right)} $, $ {\mathfrak c}''={\mathcal H}_{1}^{\left(0\right)}/{\mathcal H}_{0}^{\left(0\right)} $. The decompositon $ \alpha=\alpha^{\left(0\right)}+\alpha^{\left(1\right)} $ gives a
splitting of $ {\mathcal H}^{\left(0\right)} $ into a direct sum of $ {\mathcal H}_{0}^{\left(0\right)} $ and the span of $ \bar{\partial}^{-1}\mu_{i} $, $ i\in I $.
As a corollary, the images of forms $ \bar{\partial}^{-1}\mu_{i} $ form a basis in $ {\mathfrak c} $, and the Gram
matrix of the pairing in $ {\mathfrak c} $ can be calculated as $ \left(\bar{\partial}^{-1}\mu_{i},\bar{\partial}^{-1}\mu_{j}\right)_{L_{2}\left({\mathbb C}P^{1}\right)} $.
Similarly, projections of $ \bar{\partial}^{-1}\mu_{i}-\bar{\partial}^{-1}\mu_{i'} $, $ i\in I_{+} $, form a basis in $ {\mathfrak c}'' $, and one
can easily calculate the Gram matrix of this basis. Let us estimate
elements of these Gram matrices.

We may assume that $ \infty\in K_{i_{0}} $ for some $ i_{0}\in I $, and that $ \operatorname{Supp}\beta_{0}\subset K_{i_{0}} $. Then
the operator $ \bar{\partial}^{-1} $ restricted to $ {\mathbb C}P^{1}\smallsetminus K_{i_{0}} $ has $ \frac{dx}{y-x} $ as a kernel, and we
may suppose that $ \mu_{i} $ is $ d\vartheta_{i}/2\pi $, $ \vartheta_{i} $ being the natural angle coordinate on
$ \partial K_{i}\subset{\mathbb C} $. Let us calculate the Gram matrix $ \left(G_{ij}\right) $ for $ \bar{\partial}^{-1}d\vartheta_{i}\in{\mathcal H}^{\left(0\right)} $. If the
radius of $ K_{i} $ is $ r_{i} $, and the distance between centers of $ K_{i} $ and $ K_{j} $ is $ d_{ij} $,
then $ G_{ij} $ can be estimated as
\begin{equation}
\int\Sb r_{i}<|z|<C \\ |z-d_{ij}|>r_{j}\endSb\frac{1}{z\overline{\left(z-d_{ij}\right)}}dx\,dy
\notag\end{equation}
here $ C\gg0 $ (and depends on $ \beta_{0} $), $ z=x+iy. $ Note that the integral over $ |z|>C/2 $
is greater than a constant which does not depends on $ d_{ij} $ and $ r_{i},r_{j} $, thus
this integral is bounded from below. If $ i=j $, then it behaves as $ \log \frac{1}{r} $.
If $ i\not=j $, then the part outside the disk of radius $ 2d_{ij} $ is $ \sim\log \frac{1}{d_{ij}} $, and
the part inside this disk is scaling-invariant. Thus to estimate the
second part one may assume $ d_{ij}=1 $. However, the integral converges
absolutely inside the whole disk $ |z|<2 $, thus this part is bounded.

Thus we have estimates for the elements of the Gram matrix\footnote{Note that these estimates may be not sufficient to describe the Hilbert
structure on $ {\mathfrak c} $ up to equivalence.}.
Knowledge of this Gram matrix gives a complete description of the
subspace $ {\mathcal B} $ of $ H^{-1}\left({\mathbb C}P^{1},\omega\otimes\bar{\omega}\right) $ spanned by extensions-by-$ \delta $-function of forms
on $ \partial K_{i} $, $ i\in I $.

Similarly one can estimate the elements of the Gram matrix for $ {\mathfrak c}'' $.
Since elements of $ {\mathfrak c}'' $ correspond to $ 2 $-forms with integral 0, this Gram
matrix does not change when we apply a conformal transformation to disks
$ K_{i} $. The diagonal elements are $ \frac{l_{i i'}}{\pi} $, and the off-diagonal ones are $ C
\frac{1}{\pi}\log |\lambda\left(c_{i},c_{i'},c_{j},c_{j'}\right)|+O\left(\sum\varepsilon_{k}r_{k}\right) $, here $ \lambda $ is the double ratio, $ c_{\bullet} $ and
$ r_{\bullet} $ are the center and radius of $ K_{\bullet} $. The constants $ \varepsilon_{\bullet} $ can be calculated as
$ c_{i}=\frac{c_{j}-c_{j'}}{\left(c_{i}-c_{j}\right)\left(c_{i}-c_{j'}\right)} $.

\begin{remark} Note that we have estimates for the elements of this Gram
matrix. It is easy to check that if the matrix $ \left(r_{i}/d_{ij}\right)_{i,j\in I} $ gives a
compact operator $ l_{2} \to l_{2} $, these estimates allow one to reconstruct the
Hilbert norm up to equivalence, thus to describe the space $ {\mathfrak c}'' $ completely.
\end{remark}

Note that the diagonal entries of the Gram matrix for $ {\mathfrak c} $ are not
bounded. Since any element in the image of the integration mapping has a
finite norm w.r.t. this matrix, this shows, in particular, that the image
of the integration mapping does not coincide with $ l_{2} $. Similarly, the Gram
matrix for $ {\mathfrak c}'' $ contains arbitrarily big elements, thus $ {\mathfrak c}'' $ also differs
from $ l_{2} $, thus image of partial period mapping differs from $ l_{2} $. As we will
see it in Section~\ref{s9.40}, this image coincides with the space $ {\mathfrak c}' $ defined
below.

So far the spaces we consider depended on the relative position of
disks $ K_{i} $ only. Now suppose that a fraction-linear orientation-changing
identification $ \varphi_{i} $ of $ \partial K_{i} $ and $ \partial K_{i'} $ is fixed, and consider the subspace $ {\mathcal B}' $
of $ {\mathcal B} $ which is spanned by extensions-by-$ \delta $-function of forms on $ \partial K_{i} $, $ i\in I $,
which are preserved by the identifications $ \varphi_{i} $. Let $ {\mathcal B}_{0} $ be the subspace of
$ {\mathcal B} $ spanned by extensions-by-$ \delta $-function of forms on $ \partial K_{i} $, $ i\in I $, with integral
0 along each $ \partial K_{i} $, let $ {\mathcal B}_{0}'={\mathcal B}'\cap{\mathcal B}_{0} $.

It is clear that $ {\mathcal B}/{\mathcal B}_{0}\simeq{\mathfrak c} $, Let us describe $ {\mathcal B}'/{\mathcal B}_{0}' $. Obviously, the
image of $ {\mathcal B}'/{\mathcal B}_{0}' $ in $ {\mathfrak c} $ is a subspace of $ {\mathfrak c}'' $, however, if the identifications
$ \varphi_{i} $ ``squeeze'' the circles $ \partial K_{i} $ too much, this image may be a proper
subspace of $ {\mathfrak c}'' $.

Indeed, the dual statement to Remark~\ref{rem35.35} shows that one can
define a $ \operatorname{PGL}\left(2,{\mathbb R}\right) $-invariant norm on $ H_{\int=0}^{-1/2}\left(S^{1},\Omega^{\text{top}}\right) $ which is equivalent
to Hilbert norm. However, there is no $ \operatorname{PGL}\left(2,{\mathbb R}\right) $-invariant norm on the
whole space $ H^{-1/2}\left(S^{1},\Omega^{\text{top}}\right) $. For $ g\in\operatorname{PGL}\left(2,{\mathbb R}\right) $ consider how $ g $ changes the
Sobolev norm on $ H^{-1/2}\left(S^{1},\Omega^{\text{top}}\right) $. It is enough to estimate $ N\left(g\right)=\|g^{*}\left(\beta\right)\| $, $ \beta $
being an arbitrary form with non-zero integral. It is clear that
different choices of Sobolev norm on $ H^{-1/2}\left(S^{1},\Omega^{\text{top}}\right) $\footnote{Recall that it is defined up to equivalence only.} and different
choices of $ \beta $ would change $ N\left(g\right) $ to something of the form $ CN\left(g\right)+O\left(1\right) $ only.
Since $ \operatorname{PGL}\left(2,{\mathbb R}\right) $ has a compact subgroup $ U\left(1\right) $ of dimension 1, it is enough
to estimate $ N\left(g\right) $ on double classes $ U\left(1\right)\backslash\operatorname{PGL}\left(2,{\mathbb R}\right)/U\left(1\right) $ only, thus one can
assume $ g=\operatorname{diag}\left(\lambda,1\right) $, $ \lambda\geq1 $.

In turn, it is sufficient to estimate $ H^{-1/2} $-norm of $ \frac{\lambda dx}{1+\lambda^{2}x^{2}} $ on
$ {\mathbb R} $, and one can do it explicitly since the Fourier transfrom can be easily
calculated, it is proportional to $ e^{-\lambda|\xi|} $. Thus $ N\left(g\right)=C \log \lambda + O\left(1\right) $. Using
this, one can easily obtain

\begin{proposition} Consider the image $ {\mathfrak c}' $ of $ {\mathcal B}'/{\mathcal B}_{0}' $ in $ {\mathfrak c}={\mathcal B}/{\mathcal B}_{0} $. Consider $ \varphi_{i} $ as
elements of $ \operatorname{SL}\left(2,{\mathbb R}\right)\subset{\mathbb R}^{4} $. The space $ {\mathfrak c}' $ consists of sequences $ \left(p_{i}\right)\in{\mathfrak c}'' $ which
satisfy $ \left(p_{i}\log \|\varphi_{i}\|\right)\in l_{2} $. \end{proposition}

We see that $ {\mathfrak c}' $ can be described by Gram matrix
$ G'_{ij}=\delta_{ij}\log \|\varphi_{i}\|+G_{ij}'' $, $ G'' $ being the Gram matrix for $ {\mathfrak c}'' $.

\subsection{$ B $-periods }\label{s9.41}\myLabel{s9.41}\relax  What we defined in Section~\ref{s9.20} was integration
of global holomorphic forms along $ A $-cycles. Finite-genus theory shows
that it is important to study additional integrals along $ B $-cycles. They
should depend on the choice
of cuts on the Riemann surface, and we are not in the conditions when one
can easily proceed with such cuts.

Since we have only $ H^{0} $-smoothness of global holomorphic forms, one
cannot invent a priori bounds on integrals of these forms along arbitrary
curves. Indeed, the ``infinity'' $ M_{\infty}\subset{\mathbb C}P^{1} $ can break the complex sphere into
infinitely many
connected components (see the example in Section~\ref{s0.40}), thus one cannot
assume that the cuts do not
intersect $ M_{\infty} $. Moreover, even if there is only one connected component, it
is not clear how to make infinitely many cuts in $ {\mathbb C}P^{1} $ which would not
intersect each other. In fact it {\em is\/} possible to make such cuts, but in
general the lengths of these cuts form a quickly increasing sequence.

However, if one has a strip $ \left(-\varepsilon,\varepsilon\right)\times\left(a,b\right) $ embedded into the complex
curve, then the value of $ \int_{a}^{b} $ averaged along $ \left(-\varepsilon,\varepsilon\right) $ is well defined and
may be bounded as $ O\left(\varepsilon^{-1/2}\right) $. Indeed, this average is $ L_{2} $-pairing
\begin{equation}
\frac{1}{2\varepsilon}\int\alpha\wedge dx
\notag\end{equation}
with $ dx/\varepsilon $, here $ \left(x,y\right) $ are coordinates on the strip $ -\varepsilon<x<\varepsilon $.

Thus we are not going to define the $ B $-periods as integrals over
curves, but as some averaged integrals. On the level of would-be homology
of the surface we will pair the form with cycles on the surface which are
not ``geometric'' cycles, but linear combinations of them.

Consider the involution ' and gluings $ \varphi_{i} $. Suppose that $ 0\in K_{i_{0}} $,
$ \infty\in K_{i_{0}'} $, and that the disks $ \left\{K_{i}\right\} $ have a uniform thickening.
Parameterize
the set of rays going from 0 to $ \infty $ by the angle $ \vartheta $, and deform each ray
slightly in such a way that it avoids all the circles $ K_{i} $. We require
that deformations of two rays which differ by the angle $ \Delta\vartheta $ do not become
closer than $ \varepsilon\cdot\Delta\vartheta $ (outside of $ K_{i_{0}} $). Say, let $ \widetilde{K}_{i} $ be the concentric with $ K_{i} $
disk of radius $ \left(1+2\varepsilon\right)\operatorname{radius}\left(K_{i}\right) $, and deform the ray $ {\mathcal R} $ inside $ \widetilde{K}_{i} $ so that
it moves along an appropriate arc going between $ K_{i} $ and $ \widetilde{K}_{i} $. If
the
$ \operatorname{dist}\left({\mathcal R},\operatorname{center}\left(K_{i}\right)\right)/\operatorname{radius}\left(\widetilde{K}_{i}\right)=\rho\leq1 $, one can take the radius of the arc to
be $ \left(1+\varepsilon+\varepsilon\rho\right)\operatorname{radius}\left(K_{i}\right) $.

In fact we need to choose whether the arcs are going to leave the
disk $ K_{i} $
on the right or on the left. We use the following algorithm: let $ O_{i} $ be
the fixed point of $ \varphi_{i} $ inside $ K_{i} $. If the ray $ {\mathcal R} $ leaves $ O_{i} $ on the right, let
the deformation leave $ K_{i} $ on the right, otherwise leave it on the left.
(This choice is going to be important in the proof of Proposition
~\ref{prop8.33}. Note that the choice of the direction of the turnout does not
depend on the metric on $ {\mathbb C}P^{1} $.)

Now we assume that 0 and $ \infty $ are fixed points of $ \varphi_{i_{0}} $. As a last
correction, in a neighborhood of $ K_{i_{0}'} $ change a ray $ {\mathcal R} $ to a logarithmic
spiral so that the intersection of $ {\mathcal R} $ with $ \partial K_{i_{0}} $ and intersection of $ {\mathcal R} $ with
$ \partial K_{i_{0}'} $ are glued together by $ \varphi_{i_{0}} $. Say, if $ \partial K_{i_{0}'} = \left\{z \mid |z|=R\right\} $, take the
part of a spiral $ d\vartheta = C\,dr/\varepsilon r $ inside $ \left\{z \mid e^{-\varepsilon}R <|z|<R\right\} $, here $ C=\operatorname{Arg} \varphi_{i_{0}} $ is
defined up to addition of a multiple of $ 2\pi $. Now each deformed ray
represents a closed curve after the gluing by $ \varphi_{i_{0}} $ is performed.

Under these conditions the averaged over $ \vartheta\in\left[0,2\pi\right] $ integral over the
deformed rays (between $ \partial K_{i_{0}} $ and $ \partial K_{i_{0}'} $) is correctly defined. It
represents a combination of cycles, thus is a cycle itself. Denote this
linear functional on $ {\mathcal H}^{\left(0\right)}\left({\mathbb C}P^{1},\left\{K_{i}\right\}\right) $ by $ {\mathcal Q}_{i_{0}} $.

Restricting this linear functional to global holomorphic forms, we
call the integral along this cycle $ B $-{\em period\/} of the global holomorphic
form {\em from\/} $ \partial K_{i_{0}} $ {\em to\/} $ \partial K_{i_{0}'} $. The above description shows that

\begin{lemma} $ B $-period of $ \alpha $ from $ \partial K_{i} $ to $ \partial K_{i'} $ is bounded by $ C\cdot l_{i i'}\|\alpha\|_{L_{2}} $, here
$ l_{i i'} $ is the conformal distance between $ \partial K_{i_{0}} $ and $ \partial K_{i_{0}'} $. The constant $ C $
depends on $ \varepsilon $ only. \end{lemma}

\begin{definition} Associate to $ \alpha\in\Gamma\left(\bar{M},\omega\right) $ the sequence $ \left(q_{j}\right)_{j\in I_{+}} $, $ q_{j} $ being the
$ B $-period of $ \alpha $ from $ \partial K_{j} $ to $ \partial K_{j'} $. Denote this mapping by $ Q $. \end{definition}

\begin{remark} Note that in the case of finite genus one chooses the
$ B $-cycles to be non-intersecting. This assures that the matrix of periods
is symmetric. The above construction takes average of
different cycles connecting given points, thus one can get an impression
that the resulting matrix of periods will have much worse properties than
in the standard settings.

However, the properties turn out to be exactly the same, due to our
choice of direction of turnouts, and the following surprising result: \end{remark}

\begin{proposition} Consider 4 different points $ x_{0},x_{1},y_{0},y_{1}\in{\mathbb C}P^{1} $. Parameterize
the set of (circular) arcs connecting $ x_{0} $ with $ x_{1} $ by the angle at $ x_{0} $, and
do the same with arcs connecting $ y_{0} $ with $ y_{1} $. This provides a measure on
the set of arcs connecting $ x_{0} $ with $ x_{1} $, same for $ y_{0} $ and $ y_{1} $. Then the
average index of intersection of an arc $ x_{0}x_{1} $ with an arc $ y_{0}y_{1} $ is 0. \end{proposition}

\subsection{Space of periods }\label{s9.40}\myLabel{s9.40}\relax  In Section~\ref{s9.12} we have shown that under
mild assumption the space $ {\mathcal H}^{\left(0\right)} $ (which, loosely speaking, consists of
holomorphic forms with jumps along the cuts on the curve) is a sum of
two components: forms with integrals 0 along $ A $-cycles (which can be
described by locale data on each cut), and some explicitely defined
Hilbert subspace $ {\mathfrak c} $ of ``small'' dimension (two basis vectors per each cut
on the curve). Only the elements of $ {\mathcal H}^{\left(0\right)} $ which have equal integrals along
two sides of the cut have a chance to correspond to a global holomorphic
form on $ \bar{M} $, thus only subspace $ {\mathfrak c}'\subset{\mathfrak c} $ is interesting for us.

\begin{theorem} Suppose that the matrix $ \left(e^{-l_{ij}}-\delta_{ij}\right) $ defines a compact
operator $ l_{2} \to l_{2} $, and the disks $ K_{i} $ have a thickening. Then the partial
period mapping $ P $ sends global holomorphic forms to elements of $ {\mathfrak c}' $.
Moreover, it is an isomorphism onto $ {\mathfrak c}' $. \end{theorem}

\begin{proof} The first statement is an immediate corollary of the
description of the elements of $ {\mathcal H}^{\left(0\right)} $ via the space $ {\mathcal B}' $ in Section~\ref{s9.12}.
We start the proof of the second one by showing that $ P $ is a component of
a Fredholm operator of index 0.

In the notations of Section~\ref{s9.12}, consider the mapping $ {\mathcal H}_{1}^{\left(0\right)} \xrightarrow[]{\pi}
{\mathfrak c}'' $. Let $ {\mathcal H}_{2}^{\left(0\right)}=\pi^{-1}\left({\mathfrak c}'\right) $. Then the mapping
\begin{equation}
{\mathcal J}_{\omega}\colon {\mathcal H}_{2}^{\left(0\right)} \to \bigoplus\Sb l_{2} \\ i\in I_{+}\endSb H_{\int=0}^{-1/2}\left(\partial K_{i},\Omega_{\partial K_{i}}^{1}\right)\colon \alpha \mapsto \left(\alpha|_{\partial K_{i}}-\varphi^{*}\left(\alpha|_{\partial K_{i'}}\right)\right).
\notag\end{equation}
is continuous. The same arguments as in the proof of Theorem~\ref{th35.45}
show that $ {\mathcal J}_{\omega}|_{{\mathcal H}_{0}^{\left(0\right)}} $ is a Fredholm operator of index 0. Combining this
mapping with the projection $ {\mathcal H}_{1}^{\left(0\right)} \xrightarrow[]{\pi} {\mathcal H}_{2}^{\left(0\right)}/{\mathcal H}_{0}^{\left(0\right)}={\mathfrak c}' $, we see that $ \alpha \mapsto
\left({\mathcal J}_{\omega}\left(\alpha\right),\pi\left(\alpha\right)\right) $ is a Fredholm mapping of index 0. If $ \alpha $ is in the null-space
of this mapping, then $ \alpha $ is a non-trivial global holomorphic form with
vanishing $ A $-periods, thus $ \alpha\in\operatorname{Ker} P $, and $ \alpha=0 $. Hence this mapping is an
isomorphism. Since $ P=\pi|_{\operatorname{Ker}{\mathcal J}_{\omega}} $, it is an isomorphism as well. \end{proof}

Let $ \alpha_{i}=\bar{\partial}^{-1}\beta_{i} $. The theorem says that for any element $ \left(c_{i}\right)\in{\mathfrak c}' $ one can
find an element $ \alpha\in{\mathcal H}^{\left(0\right)} $ such that integrals of $ \alpha $ along $ \partial K_{i} $ are 0 for any
$ i\in I $, and $ \sum c_{i}\alpha_{i}+\alpha $ is a global holomorphic form (which automatically has
$ A $-periods $ c_{i} $). Moreover, $ \|\alpha\|_{{\mathcal H}^{\left(0\right)}}\leq C\cdot\|\left(c_{i}\right)\|_{{\mathfrak c}'} $.

\begin{proposition} \label{prop9.42}\myLabel{prop9.42}\relax  Let $ \beta\in{\mathcal H}^{\left(0\right)} $, and integrals of $ \beta $ along $ \partial K_{i} $ are 0 for any
$ i\in I $, let $ \left(c_{i}\right)\in{\mathfrak c}' $. Then $ \left|\sum_{j\in I_{+}}c_{j}{\mathcal Q}_{j}\left(\beta\right)\right| \leq C\|\beta\|_{{\mathcal H}^{\left(0\right)}}\cdot\|\left(c_{i}\right)\|_{{\mathfrak c}'} $ for an
appropriate $ C $ which does not depend on $ \beta $ and $ \left(c_{i}\right) $. \end{proposition}

\begin{proof} Suppose that $ \beta $ can be holomorphically continued into all disks
$ K_{i} $ except $ K_{i_{0}} $. Then $ \beta=\partial f $, and $ f $ is holomorphic outside $ K_{i_{0}} $. Let
$ j,j'\not=i_{0} $. The
construction of $ {\mathcal Q}_{j} $ shows that one can calculate $ {\mathcal Q}_{j}\left(\beta\right) $ as $ f\left(y_{1}\right)-f\left(y_{0}\right) $, $ y_{0} $,
$ y_{1} $ being two fixed points of $ \varphi_{j} $. Let $ \Psi_{j} $ be the fraction-linear function
with a zero and a pole at $ y_{0} $ and $ y_{1} $. One can momentarily see that
$ {\mathcal Q}_{j}\left(\beta\right)=\int_{\partial K_{i_{0}}}\beta\log \Psi_{j} $. Similarly, $ {\mathcal Q}_{i_{0}}\left(\beta\right)=0 $.

Similar statements are true for forms $ \beta $ such that they can be
holomorphically continued into all the disks $ K_{i} $ except a finite number.
Since any form $ \beta $ which satisfies conditions of the proposition can be
approximated by such forms, we see that it is sufficient to show that
$ \sum_{j\not=i,i'}c_{j}\log \Psi_{j} $ converges in $ \bigoplus_{l_{2}}H^{1/2}\left(\partial K_{i}\right)/\operatorname{const} $, or that $ \sum_{j\not=i,i'}c_{j}d\Psi_{j}/\Psi_{j} $
converges in $ \bigoplus H^{-1/2}\left(\partial K_{i},\Omega^{1}\right) $. In turn, it is sufficient to show
convergence of $ \sum_{j\not=i,i'}c_{j}d\Psi_{j}/\Psi_{j} $ in $ \bigoplus_{l_{2}}H^{0}\left(K_{i},\omega\right)=L_{2}\left(\bigcup K_{i},\omega\right) $. On the other
hand, if $ \sigma_{j} $ is 1 in $ K_{j} $ and $ K_{j'} $, and 0 otherwise, then
$ \bar{\partial}\sum_{j\not=i,i'}c_{j}\sigma_{j}d\Psi_{j}/\Psi_{j} =\sum_{j\not=i,i'}c_{j}\beta_{j}\in H^{-1} $, hence $ \sum c_{j}\sigma_{j}d\Psi_{j}/\Psi_{j} $ converges in $ L_{2} $. \end{proof}

\begin{proposition} \label{prop9.25}\myLabel{prop9.25}\relax  Let $ \left(d_{i}\right),\left(c_{i}\right)\in{\mathfrak c}' $, and $ \alpha\in{\mathcal H}^{\left(0\right)} $ such that integrals of
$ \alpha $ along $ \partial K_{i} $ are 0 for any $ i\in I $, and $ \alpha_{c}=\sum c_{i}\alpha_{i}+\alpha $ is a global holomorphic
form on $ \bar{M} $ (automatically with $ A $-periods $ c_{i} $). Then $ \sum d_{j}{\mathcal Q}_{j}\left(\alpha_{c}\right)
\leq C\|\left(c_{i}\right)\|_{{\mathfrak c}'}\cdot\|\left(d_{i}\right)\|_{{\mathfrak c}'} $. \end{proposition}

\begin{proof} It is sufficient to show that $ \sum d_{j}{\mathcal Q}_{j}\left(\sum c_{i}\alpha_{i}\right)\leq C_{1}\|\left(c_{i}\right)\|_{{\mathfrak c}'}\cdot\|\left(d_{i}\right)\|_{{\mathfrak c}'} $. In
the notations of the previous proposition $ \alpha_{i}-\alpha_{i'}=d\Psi_{i}/\Psi_{i} $ outside of $ K_{i} $ and
$ K_{i'} $. One can assume that only a finite number of $ c_{i} $ and $ d_{i} $ is non-zero.
Then $ \sum_{i\in I}c_{i}\alpha_{i} $ is $ \sum_{i\in I_{+}}c_{i}d\Psi_{i}/\Psi_{i} $.

Now it should be obvious that $ {\mathcal Q}_{j}\left(\alpha_{i}\right)=\int\sigma_{i}\sigma_{j}\alpha_{i}\wedge\alpha_{j} $, which finishes the
proof. \end{proof}

\subsection{Period matrix }\label{s9.60}\myLabel{s9.60}\relax  Since $ P $ is an isomorphism $ \Gamma\left(\bar{M},\omega\right) \to {\mathfrak c}' $, one can
consider
the mapping $ \Omega=Q\circ P^{-1} $. Write this mapping using coordinate ``basis'' in $ {\mathfrak c}' $:

\begin{definition} Let $ \widetilde{\alpha}_{j} $, $ j\in I_{+} $, be the global holomorphic form on $ \bar{M} $ such that
$ P\left(\widetilde{\alpha}_{i}\right) $ has 1 on $ j $-th position, $ -1 $ at $ j' $-position, 0 at the other
positions. Let $ \Omega_{ij} $ be the $ B $-period of $ \widetilde{\alpha}_{j} $ from $ \partial K_{i} $ to $ \partial K_{i'} $, $ i\in I_{+} $. \end{definition}

\begin{proposition} \label{prop8.33}\myLabel{prop8.33}\relax  Let $ \alpha\in\Gamma\left(\bar{M},\omega\right) $ and the sequence $ \left(p_{j}\right)=P\left(\alpha\right) $ has only a
finite number of non-zero elements. Let $ Q\left(\alpha\right)=\left(q_{j}\right) $. Then $ \|\alpha\|_{L_{2}}^{2}=
i\sum_{I_{+}}\left(\bar{p}_{j}q_{j}-p_{j}\bar{q}_{j}\right) $. \end{proposition}

\begin{proof} Consider a representative of $ \alpha $ in $ {\mathcal H}^{\left(0\right)}\left({\mathbb C}P^{1},\left\{K_{i}\right\}\right) $. Since all the
$ A $-periods of $ \alpha $ but a finite number are 0, $ \alpha $ can be extended (without
changing the norm too much) into all the disks but a finite number
preserving the closeness. Consider the remaining disks. Since the
integral
of $ \alpha $ around $ K_{i} $ is opposite to the integral around $ K_{i'} $, we see that
if we
connect $ K_{i} $ and $ K_{i'} $ by a cut, then the integral of $ \alpha $ around the resulting
hole is 0. Make smooth cuts which connect the remaining disks pairwise
(according to $ ' $) and do not intersect. One can suppose that two ends
of the cut---one on $ \partial K_{i} $, another on $ \partial K_{i'} $---are identified by $ \varphi_{i} $.

After the cuts are performed, on the resulting domain $ 1 $-form $ \alpha $ is
closed, and the integral along any component of (piecewise-smooth)
boundary is 0. Thus one can write $ \alpha = df $, $ f $ being a function of
smoothness $ H^{1} $. The restriction of $ f $ to any smooth curve is well-defined,
and is of smoothness $ H^{1/2} $. When one goes from one side of the cut to
another one along $ \partial K_{i} $, $ f $ grows by $ p_{i} $. Let $ \overset{\,\,{}_\circ}{q}_{i} $ be the change of $ f $ when one
goes from $ \partial K_{i} $ to $ \partial K_{i'} $ along the cut (choosing the side of the cut so that
the direction is counterclockwise).

Note that $ f $ is holomorphic near $ \partial K_{i} $ and $ \partial K_{i'} $ thus the {\em value\/} of $ f $ at
points is well-defined, thus the change of $ f $ along the cut is
well-defined. (In generic point $ z $ of the cut $ f\left(z\right) $ is not correctly
defined, since $ f $ is only of smoothness $ H^{1/2} $.)

Let $ \gamma_{i} $ be the part of the boundary of the domain consisting of the
circles $ \partial K_{i} $, $ \partial K_{i'} $ and both sides of the cut which connects
them.

Now take into account that $ \|\alpha\|_{L_{2}}^{2}=i\int\alpha\wedge\bar{\alpha}=i\int\partial f\wedge\bar{\partial}\bar{f} $. Proceeding as in
Section~\ref{s7.90}, we see that the only change to the arguments is that
instead of taking integrals along $ \partial K_{i}\cup\partial K_{i'} $, one needs to take some
integrals along $ \gamma_{i} $. As there, the integral along $ \partial K_{i}\cup\partial K_{i'} $ vanishes, so
what remains is
\begin{equation}
\|\alpha\|_{L_{2}}^{2}=-i\sum_{i}\int_{\gamma_{i}}f\,d\bar{f}
\notag\end{equation}
summation being over $ i $ such that $ p_{i}\not=0 $. The cycle $ \gamma_{i} $ consists of 4 parts:
two going around $ \partial K_{i} $ and $ \partial K_{i'} $, another two going along sides of the cut.
On the first two parts $ d\bar{f} $ are identified via $ \varphi_{i} $ (with opposite signs),
and $ f $ differs by $ \overset{\,\,{}_\circ}{q}_{i} $, thus the total integral is $ -\bar{p}_{i}\overset{\,\,{}_\circ}{q}_{i} $. On the second two
$ d\bar{f} $ coincide (but the orientation is opposite), and $ f $ differs by $ p_{i} $, thus
the total integral is $ p_{i}\overline{\overset{\,\,{}_\circ}{q}_{i}} $.

What remains to prove is that we may substitute $ q_{i} $ instead of $ \overset{\,\,{}_\circ}{q}_{i} $. To
do this one needs to investigate the relationship between $ q_{i} $ and $ \overset{\,\,{}_\circ}{q} $.
First, one can describe $ \overset{\,\,{}_\circ}{q}_{i} $ as an integral of $ df $ along one side of the
cut. Indeed, though $ d\left(f|_{\gamma_{i}}\right) $ is of smoothness $ H^{-1/2} $ on each part of $ \gamma_{i} $, it
is actually analytic near circles $ \partial K_{i} $, $ \partial K_{i'} $, thus the pairing with the
fundamental cycle of the interval (which has jumps at the ends of the
interval!) is well-defined.

Second, one can suppose that the cuts are in fact piecewise-smooth
(as far as non-smooth points are in $ {\mathbb C}P^{1}\smallsetminus\overline{\bigcup K_{i}} $), and consist of arcs of
circles. Let us recall that $ q_{i} $ is the average integral of $ \alpha $ along ``rays''
which connect 0 and $ \infty $ (after an appropriate choice of coordinate system).
Here ``rays'' are curves with differ from rays inside circles $ \widetilde{K}_{j} $ only, and
consist of arcs of circles. What is more, one does not need to deform a
ray into a ``ray'' inside the circle $ K_{j} $ as far as $ p_{j}=0 $, since a closed
continuation of $ \alpha $ inside $ K_{j} $ is already fixed, thus the deformation does
not change the value of the integral along the ``ray''. Note that each
ray intersects cuts along a finite number of points (except a finite
number of rays which may contain whole pieces of cuts), thus the same is
true for ``rays''.

We see that the integral of $ \alpha $ along each ``ray'' is now well defined
(recall that in general setting only the average was well-defined),
moreover, it is easy to calculate this integral using representation
$ \alpha=df. $ The integral along a ``ray'' $ {\mathcal R} $ is equal to the change of $ f $ on
the ends minus the jumps of $ f $ at the finite number of points where ``ray''
intersects cuts, thus it is
\begin{equation}
f\left(\operatorname{end}\left({\mathcal R}\right)\right)-f\left(\operatorname{start}\left({\mathcal R}\right)\right)-\sum_{j}n_{ij}p_{j},\qquad n_{ij}\in{\mathbb Z}\text{, }i,j\in I.
\notag\end{equation}
The change of $ f $ is equal to $ \overset{\,\,{}_\circ}{q}_{i} $, thus the average value of the integral is
$ \overset{\,\,{}_\circ}{q}_{i} $ minus sum of some real multiples of $ p_{j} $:
\begin{equation}
q_{i}=\overset{\,\,{}_\circ}{q}_{i}-\sum_{j}\nu_{ij}p_{j},\qquad \nu_{ij}\in{\mathbb R}\text{, }i,j\in I,
\notag\end{equation}
and $ \nu_{ij} $ are averaged values of $ n_{ij} $.

Finally, use our choice of direction of turnout around the disk $ K_{i} $. It
insures that the following fact is true:

\begin{lemma} $ \nu_{ij}=\nu_{ji} $ if $ i\not=i $. \end{lemma}

\begin{proof} Extend the cut between $ \partial K_{j} $ and $ \partial K_{j'} $ to fixed points of $ \varphi_{j} $ (one
inside each of $ K_{j} $ and $ K_{j'} $) along
straight intervals. Call the resulting curve $ \gamma'_{j} $. Then $ 2\pi\nu_{ij} $ is the
change of $ \operatorname{Arg} z $ along the curve $ \gamma_{j}' $ (we again suppose that 0 and $ \infty $ are
fixed points of $ \varphi_{i} $). Similarly, $ \pi\nu_{ji} $ is the change of $ \operatorname{Arg}\Psi_{j}\left(z\right) $ along $ \gamma'_{i} $,
here $ \Psi_{j}\left(z\right) $ is a fraction-linear function with a pole and a zero at fixed
points of $ \varphi_{j} $.

We need to show that $ \operatorname{Im} \int_{\gamma_{j}'}\frac{dz}{z}+\int_{\gamma_{i}'}\frac{d\Psi_{j}\left(z\right)}{\Psi_{j}\left(z\right)}=0 $. However,
\begin{equation}
\int_{\gamma_{j}'}\frac{dz}{z}=\frac{1}{2\pi i} \int_{\gamma_{j}'}\frac{dz}{z}\operatorname{Jump}\left(\log \Psi_{j}\left(z\right)\right) =\frac{1}{2\pi i}\int_{\gamma_{j}''}\log  \Psi_{j}\left(z\right) d \log  z
\notag\end{equation}
here we take an arbitrary branch of $ \log \Psi\left(z\right) $ defined outside of $ \gamma_{j}' $, and
$ \gamma_{j}'' $ is a loop around $ \gamma_{j}' $. Now the identity is obvious, since $ \gamma_{j}'' $ is
homotopic to $ -\gamma_{i}'' $. \end{proof}

The lemma implies that plugging in $ \overset{\,\,{}_\circ}{q}_{i} $ instead of $ q_{i} $ into the formula
of the proposition gives the same value, which finishes the proof of
Proposition~\ref{prop8.33}. {} \end{proof}

\begin{proposition} \label{prop8.35}\myLabel{prop8.35}\relax  Let $ \alpha,\beta\in\Gamma\left(\bar{M},\omega\right) $ and the sequences $ \left(p_{j}\right)=P\left(\alpha\right) $, $ \left(p_{j}'\right)=P\left(\beta\right) $
have only a finite number of non-zero elements. Let $ Q\left(\alpha\right)=\left(q_{j}\right) $,
$ Q\left(\beta\right)=\left(q_{j}'\right) $. Then $ \sum_{I_{+}}\left(p_{j}q'_{j}-p'_{j}q_{j}\right) = $ 0. \end{proposition}

\begin{proof} The proof of Proposition~\ref{prop8.33} with minor changes is
applicable, the integral to consider is $ 0=\int\alpha\beta $. \end{proof}

The following statement is an immediate corollary of Proposition
~\ref{prop9.25}:

\begin{proposition} In the conditions of the previous proposition
$ \sum_{I_{+}}p_{j}q'_{j}\leq C\|\alpha\|\cdot\|\beta\| $. \end{proposition}

\begin{corollary} The period matrix $ \Omega_{ij} $ is symmetric, $ \operatorname{Im}\Omega_{ij} $ is a positive
real symmetric matrix, and the matrix $ \Omega_{ij} $ defines a bounded symmetric
form on the Hilbert space with the pairing given by the Gram matrix
$ \operatorname{Im}\Omega_{ij} $. \end{corollary}

\subsection{Bounded Jacobian as a torus }\label{s9.70}\myLabel{s9.70}\relax  In Section~\ref{s7.90} we have seen that
(under mild assumptions) the Jacobian coincides with the constant
Jacobian. Given a topological space $ S $ and a set $ I $, let $ S_{l_{\infty}}^{I}\buildrel{\text{def}}\over{=}\bigcup_{K}K^{I} $,
here $ K $ runs over compact subsets of $ S $. Obviously, there is a surjection
from $ \left({\mathbb C}^{*}\right)_{l_{\infty}}^{I_{+}} $ to the constant bounded Jacobian, which sends a sequence
$ \left(\psi_{i}\right)_{i\in I} $ such that $ \psi_{i'}=\psi_{i}^{-1} $ into a corresponding line bundle. Here we are
going to show that the kernel of this surjection is a lattice in
$ \left({\mathbb C}^{*}\right)_{l_{\infty}}^{I_{+}} $, as in finite-genus case. To avoid defining a lattice in
$ \left({\mathbb C}^{*}\right)_{l_{\infty}}^{I_{+}} $, consider $ \left({\mathbb C}^{*}\right)_{l_{\infty}}^{I_{+}} $ as the set of coordinate-wise exponents of
$ {\mathbb C}_{l_{\infty}}^{I_{+}} $. We obtain a mapping from $ {\mathbb C}_{l_{\infty}}^{I_{+}} $ to the bounded Jacobian which sends
$ \left(2\pi i{\mathbb Z}\right)_{l_{\infty}}^{I_{+}} $ to the origin in the Jacobian.

We are going to prove that there is a lattice $ L\subset{\mathbb C}^{I_{+}} $ which goes to
the origin. First, construct generators of this lattice:

\begin{proposition} \label{prop9.72}\myLabel{prop9.72}\relax  Suppose that the matrix $ \left(e^{-l_{ij}}-\delta_{ij}\right) $ gives a compact
operator $ l_{2} \to l_{2} $ and a compact operator $ l_{1} \to l_{1} $. Let $ j\in I_{+} $. There exists
a global holomorphic form $ \alpha^{\left(j\right)} $ on $ \bar{M} $ such that the $ A $-periods of $ \alpha^{\left(j\right)} $
vanish except for the $ j $-th one, which is equal to $ 2\pi i $. Then $ a_{j}=\exp \int\alpha^{\left(j\right)} $
is an element of $ {\mathcal H}^{\infty} $. The corresponding cocycle $ \psi_{kj}=a_{j}\circ\varphi_{k}/a_{j} $ is (locally)
constant, and coincides with $ \exp  2\pi i\Omega_{kj} $. \end{proposition}

\begin{proof} Existence of $ \alpha^{\left(j\right)} $ is a corollary of results of Section~\ref{s9.40}.
The only statement we need to prove is that $ a_{j}\in{\mathcal H}^{\infty} $. In turn, it is
sufficient to prove that $ \operatorname{Im}\Omega_{kj} $ is bounded (for a fixed $ j $). On the other
hand, $ \operatorname{Im}\Omega_{kj} $ consists of two parts which correspond to decomposition
$ \alpha^{\left(j\right)}=2\pi i\alpha_{j}+o^{\left(j\right)} $, here $ \alpha_{j} $ is defined as in Section~\ref{s9.40}, and $ o^{\left(j\right)}\in{\mathcal H}^{\left(0\right)} $
and has $ A $-periods 0. The description of the cycle for $ B $-period shows
that first part of $ \operatorname{Im}\Omega_{kj} $ is bounded by $ \log |l_{jj'}| $, thus it is sufficient
to estimate the second part $ \operatorname{Im}{\mathcal Q}_{k}\left(o^{\left(j\right)}\right) $.

This estimate follows from the following lemma:

\begin{lemma} Suppose that the matrix $ \left(e^{-l_{ij}}-\delta_{ij}\right) $ gives a compact operator
$ l_{2} \to l_{2} $ and a compact operator $ l_{1} \to l_{1} $.
\begin{enumerate}
\item
Let $ \alpha\in H^{-1/2}\left(\partial K_{i},\Omega^{\text{top}}\right) $. Denote by the same letter
extension-by-$ \delta $-function of $ \alpha $ to $ {\mathbb C}P^{1} $. Then $ \bar{\partial}^{-1}\alpha\in{\mathcal H}^{\left(0\right)} $, and
$ {\mathcal Q}_{k}\left(\bar{\partial}^{-1}\alpha\right)=O\left(\|\alpha\|_{H^{-1/2}}\right) $ uniformly in $ k\not=i $.
\item
The sequence $ \left(\|o^{\left(j\right)}|_{\partial K_{k}}\|_{H^{-1/2}}\right)_{k\in I}\in l_{1} $.
\end{enumerate}
\end{lemma}

\begin{proof} The first part follows from the explicit construction of cycles
for $ B $-periods. Prove the second part.

Since sum of radii of disks $ K_{i} $ is finite, one can estimate that the
sequence $ \left(\|\alpha_{j}|_{\partial K_{k}}\|_{H^{-1/2}}\right)_{k\in I}\in l_{1} $. Denote the space of such $ 1 $-forms on
$ \bigcup_{i}\partial K_{i} $ by $ {\mathcal H}_{1} $. On the other hand, $ o^{\left(j\right)} $ is uniquely determined by the
conditions $ {\mathcal J}_{\omega}o^{\left(j\right)}=-2\pi i{\mathcal J}_{\omega}\alpha_{j} $ and $ o^{\left(j\right)}\in{\mathcal H}_{0}^{\left(0\right)} $, here $ {\mathcal J}_{\omega} $ is the operator from
Section~\ref{s9.40}. What remains to prove is that $ {\mathcal J}_{\omega} $ sends the subspace
$ {\mathcal H}_{0}^{\left(0\right)}\cap{\mathcal H}_{1} $ onto itself. We know that $ {\mathcal J}_{\omega} $ sends $ {\mathcal H}_{0}^{\left(0\right)} $ onto itself.
However, from the
restrictions on $ l_{ij} $ one immediately obtains that $ {\mathcal J}_{\omega} $ gives a Fredholm
operatorar $ {\mathcal H}_{1} \to {\mathcal H}_{1} $ of index 0. Since $ {\mathcal J}_{\omega} $ has no kernel, this operator is
an isomorphism. \end{proof}

This finishes the proof of Proposition~\ref{prop9.72}. {}\end{proof}

\begin{definition} Call a multiplicator $ a\in{\mathcal H}^{\infty} $ {\em lattice-like}, if $ a^{-1}\in{\mathcal H}^{\infty} $
and $ \psi_{k}=a^{-1}\cdot\left(a\circ\varphi_{k}\right) $ is a constant function on $ \partial K_{k} $ for every $ k\in I $. \end{definition}

Now to each $ j\in I_{+} $ we associated a lattice-like multiplicator $ a_{j} $,
which induces a cocycle $ \left(\psi_{k}\right)_{k\in I_{+}} $, $ \psi_{k}=\exp  2\pi i\Omega_{kj} $. The bundle given by this
cocycle is isomorphic to a trivial one (via $ a_{j} $). Since $ \operatorname{Im}\Omega_{kj} $ gives a
positive Hermitian form (defined at least on real sequences of finite
length), columns of $ \operatorname{Im}\Omega_{kj} $ generate the subspace which is
dense in the set of real sequences (with topology of direct product),
which shows in conditions of Proposition~\ref{prop9.72}

\begin{proposition} Consider the space $ {\mathbb C}_{l_{\infty}}^{I_{+}} $ of bounded sequences $ \left(\Psi_{k}\right) $ with
topology induced from the direct product $ {\mathbb C}^{I_{+}} $. Consider the subgroup $ L $
generated by rows $ \left(\Omega_{kj}\right) $, $ j\in I_{+} $, and $ {\mathbb Z}_{l_{\infty}}^{I_{+}} $. This subgroup is a lattice,
i.e., its $ {\mathbb R} $-span is dense in $ {\mathbb C}_{l_{\infty}}^{I_{+}} $. \end{proposition}

To finish the description of bounded Jacobian, it remains to prove
that any lattice-like multiplicator $ a $ is a product of powers of $ a_{j} $,
$ j\in I_{+} $. Recall that the index of function $ a $ is the change of $ \frac{\operatorname{Arg} a}{2\pi i} $
along a contour. Note that the {\em index\/} of $ a_{j} $ around $ \partial K_{k} $ is $ \pm1 $ if $ k=j,j' $,
and is 0 otherwise. Thus one reconstruct the needed powers of $ a_{j} $ by
taking indices of $ a $ around $ \partial K_{j} $, $ j\in I_{+} $. The only things we need to prove
is that a lattice-like multiplicator $ a $ with all indices around $ \partial K_{j} $ being
0 is constant, and only a finite number of indices of a lattice-like
multiplicator is non-zero. The first statement is a direct corollary of
Lemma~\ref{lm8.31} and the following

\begin{lemma} \label{lm9.91}\myLabel{lm9.91}\relax  Suppose that the conformal distances $ l_{ij} $ satisfy the
condition that $ \left(e^{-l_{ij}}-\delta_{ij}\right) $ gives a bounded operator $ l_{\infty} \to l_{2} $. Then any
lattice-like multiplicator is of class $ {\mathcal H}^{\left(1\right)} $. \end{lemma}

Recall that the space $ {\mathcal H}^{\left(1\right)} $ was defined in Section~\ref{s35.30}.

\begin{proof} Since $ a\,dz^{1/2}\in{\mathcal H}\left({\mathbb C}P^{1},\left\{K_{j}\right\}\right) $, Theorem~\ref{th4.40} shows that $ a $ can be
reconstructed from its restriction to $ \bigcup_{j\in I}\partial K_{j} $ using Cauchy formula. In
particular, $ + $-part of $ \partial a $ on $ \partial K_{j} $ is given by an integral along $ \bigcup_{k\not=j}\partial K_{k} $.
The restrictions on $ l_{ij} $ guarantie that the $ + $-parts of $ \partial a $ on the circles
$ \partial K_{j} $ is in $ \oplus_{l_{2}}H^{-1/2}\left(\partial K_{j}\right) $ modulo constants. Since $ a^{-1}\cdot\left(a\circ\varphi_{k}\right) $ is constant
and bounded, $ - $-parts of $ \partial a $ are also in this space. Since the integral of
$ \partial a $ along each circle vanishes, arguments similar to ones in Section
~\ref{s35.30} show that $ \partial a\in L_{2}\left({\mathbb C}P^{1}\smallsetminus\bigcup K_{i},\omega\right) $, thus $ a\in{\mathcal H}^{\left(1\right)} $. \end{proof}

Iteration of the procedure used in the proof of the lemma shows that

\begin{amplification} \label{amp9.94}\myLabel{amp9.94}\relax  Conclusions of Lemma~\ref{lm9.91} remain true if
$ \left(e^{-l_{ij}}-\delta_{ij}\right) $ gives a bounded operator $ l_{\infty} \to l_{\infty} $, and some power of this
operator sends $ l_{\infty} $ to $ l_{2} $. \end{amplification}

To prove that only a finite number of indices of a lattice-like
multiplicator $ a $ is non-zero, it is enough to show that $ a^{-1}\partial a\in L_{2}\left(U,\omega\right) $, $ U $
being homotopic to $ {\mathbb C}P^{1}\smallsetminus\bigcup K_{i} $. Indeed, indices of $ a $ are proportional to
periods of $ a^{-1}\partial a $, and periods of $ L_{2} $-form form a sequence in $ l_{2} $.

In turn, since $ a^{-1} $ is a multiplicator, thus is bounded, this follows
from the fact that $ \partial a\in L_{2}\left(U,\omega\right) $. We obtain

\begin{theorem} Consider a family of non-intersecting disks $ K_{i}\subset{\mathbb C}P^{1} $ with
pairwise conformal distances. Let $ {\mathcal R} $ be the matrix $ \left(e^{-l_{ij}}-\delta_{ij}\right) $. Suppose
that $ {\mathcal R} $ gives a compact mapping $ l_{2} \to l_{2} $ and a bounded mapping $ l_{\infty} \to l_{\infty} $,
and that some power of $ {\mathcal R} $ gives a bounded mapping $ l_{\infty} \to l_{2} $. If disks $ K_{i} $
have a uniform thickening, then the bounded Jacobian coincides with the
quotient of $ {\mathbb C}_{l_{\infty}}^{I_{+}} $ by the lattice generated by $ {\mathbb Z}_{l_{\infty}}^{I_{+}} $ and rows of the
period matrix $ \left(\Omega_{ij}\right) $. \end{theorem}

\begin{remark} Formally speaking, we defined the Jacobians in the case when
$ \left(e^{-l_{ij}/2}-\delta_{ij}\right) $ gives a compact operator $ l_{2} \to l_{2} $, thus the statement of
theorem is abuse of notations. However, it is easy to define all the
ingredients needed for the definition of the bounded Jacobian as
far as the conditions of Amplification~\ref{amp9.94} are satisfied. \end{remark}

Let us recall that the conditions of the theorem are automatically
satisfied for any Hilbert--Schmidt curve, thus we get a complete
description of the bounded Jacobian in this case. It is similar to the
finite-genus case, where Jacobian is a quotient of a finite-dimensional
complex vector space by a lattice.

Moreover, note that the topology on the constant bounded Jacobian is
inherited from the topology of direct product on $ {\mathbb C}^{I_{+}} $, thus the above
description allows one to reconstruct the topology on the bounded
Jacobian as well.

\subsection{Rigged Hodge structure }\label{s9.80}\myLabel{s9.80}\relax  Let us wrap the results of the previous
section into the familiar form of Hodge structures. In fact the resulted
structure will be a hybrid of a Hodge structure and a structure of a
rigged topological vector space.

Let $ I_{+} $ be an arbitrary set. Let $ H_{{\mathbb Z}}^{A} $, $ H_{{\mathbb R}}^{A} $, $ H_{{\mathbb C}}^{A} $ be the spaces of
sequences $ \left(p_{i}\right)_{i\in I_{+}} $ with only a finite number of non-zero terms (with
integer/real/complex terms), $ H_{{\mathbb Z}}^{B} $, $ H_{{\mathbb R}}^{B} $, $ H_{{\mathbb C}}^{B} $ be the spaces of sequences
$ \left(q_{i}\right)_{i\in I_{+}} $ without any restriction on growth, and $ H_{{\mathbb Z}}= H_{{\mathbb Z}}^{A}\oplus H_{{\mathbb Z}}^{B} $, $ H_{{\mathbb R}}=H_{{\mathbb R}}^{A}\oplus H_{{\mathbb R}}^{B} $,
$ H_{{\mathbb C}}=H_{{\mathbb C}}^{A}\oplus H_{{\mathbb C}}^{B} $. One can define an operation of complex conjugation on the
spaces $ H_{{\mathbb C}}^{\bullet} $, this operation leaves $ H_{{\mathbb R}}^{\bullet} $ fixed. The spaces $ H_{{\mathbb Z}} $, $ H_{{\mathbb R}} $, $ H_{{\mathbb C}} $ have a
natural symplectic structure
\begin{equation}
\left[\left(\left(p_{i}\right),\left(q_{i}\right)\right),\left(\left(p'_{i}\right),\left(q_{i}'\right)\right)\right] = \sum_{i}\left(p_{i}q_{i}'-p'_{i}q_{i}\right)
\notag\end{equation}
such that the components $ H_{\bullet}^{A,B} $ are Lagrangian and mutually dual. Fix an
arbitrary mapping $ \Omega\colon H_{{\mathbb C}}^{A} \to H_{{\mathbb C}}^{B} $. Let $ \Omega_{1}\colon H_{{\mathbb C}}^{A} \to H_{{\mathbb C}}^{A}\oplus H_{{\mathbb C}}^{B}=H_{{\mathbb C}} $ be $ \operatorname{id}\oplus\Omega $.
Denote $ \operatorname{Im}\Omega_{1}\subset H_{{\mathbb C}} $ by $ H^{1,0} $. Suppose that (compare Proposition~\ref{prop8.35}) $ H^{1,0} $
is Lagrangian. Moreover, suppose (compare Proposition~\ref{prop8.33}) that if
$ \alpha\in H^{1,0} $, then $ \operatorname{Im} \left[\alpha,\bar{\alpha}\right]\gg0 $ in the following sense:
\begin{enumerate}
\item
$ \operatorname{Im} \left[\alpha,\bar{\alpha}\right]\geq0 $ (here $ \alpha \mapsto \bar{\alpha} $ is the complex conjugation on $ H_{{\mathbb C}} $);
\item
the equality is achieved only if $ \alpha=0 $;
\item
the completion $ {\mathfrak h}^{1,0} $ of $ H^{1,0} $ w.r.t. the norm $ \|\alpha\|^{2}=\operatorname{Im} \left[\alpha,\bar{\alpha}\right] $ has no
vectors of length 0.
\end{enumerate}
Let $ H^{0,1} = \left\{\alpha\in H_{{\mathbb C}} \mid \bar{\alpha}\in H^{1,0}\right\} $. The mapping $ \Omega_{1} $ gives an inclusion $ i $ of $ H_{{\mathbb C}}^{A} $
into $ {\mathfrak h}^{1,0} $.

\begin{definition} A {\em rigged Hodge structure\/} is a subspace $ H^{1,0}\subset H_{{\mathbb C}} $ which
satisfies the above conditions, and such that the projection $ p\colon H^{1,0} \to
H_{{\mathbb C}}^{B} $ can be continuously extended to a mapping $ p\colon {\mathfrak h}^{1,0} \to H_{{\mathbb C}}^{B} $. \end{definition}

Note the mappings $ H_{{\mathbb C}}^{A} \xrightarrow[]{i} {\mathfrak h}^{1,0} \xrightarrow[]{p} H_{{\mathbb C}}^{B} $ equip $ {\mathfrak h}^{1,0} $ with a structure
of a rigged topological vector space \cite{GelVil64Gen}. It is clear that $ H^{1,0} $
is a rigged Hodge structure iff $ \Omega $ is symmetric, $ \operatorname{Im}\Omega>0 $, $ -C\cdot\operatorname{Im}\Omega<\operatorname{Re}\Omega<C\cdot\operatorname{Im}\Omega $.

From now on suppose that $ H^{1,0} $ is a rigged Hodge structure. Let $ {\mathfrak h}^{0,1} $
be the complexly conjugated space to $ {\mathfrak h}^{1,0} $, $ {\mathfrak h}_{{\mathbb C}}={\mathfrak h}^{1,0}\oplus{\mathfrak h}^{0,1} $. The vector space
$ {\mathfrak h}_{{\mathbb C}} $ has a natural operation of complex conjugation $ \left(\alpha,\alpha'\right) \mapsto \left(\bar{\alpha}',\bar{\alpha}\right) $, let
$ {\mathfrak h}_{{\mathbb R}} $ be the subspace $ \left(\alpha,\bar{\alpha}\right) $ of fixed points of this complex conjugation.
There is a natural extension of the projection $ {\mathfrak h}^{1,0} \to H_{{\mathbb C}}^{B} $ to $ {\mathfrak h}_{{\mathbb C}} \to H_{{\mathbb C}}^{B} $.
This mapping is compatible with the complex conjugation, thus induces a
mapping $ {\mathfrak h}_{{\mathbb R}} \to H_{{\mathbb R}}^{B} $. Let $ {\mathfrak h}_{{\mathbb C}}^{A} $ be the kernel of the mapping $ {\mathfrak h}_{{\mathbb C}} \to H_{{\mathbb C}}^{B} $,
similarly for $ {\mathfrak h}_{{\mathbb R}}^{A} $.

Let $ \bar{\Omega}_{1}\colon H_{{\mathbb C}}^{A} \to {\mathfrak h}^{0,1} $ be the complex conjugate to the mapping $ \Omega_{1}\colon H_{{\mathbb C}}^{A}
\to {\mathfrak h}^{1,0} $, $ \bar{\Omega}_{1}v=\overline{\Omega_{1}\bar{v}} $. Define a pre-Hilbert structure on $ H_{{\mathbb C}}^{A} $ via
$ \|v\|_{A}^{2}=\|\Omega_{1}v\|_{{\mathfrak h}^{1,0}} $. Then $ \Omega $ extends to a mapping from Hilbert completion
$ \left(H_{{\mathbb C}}^{A}\right)^{\text{compl}} $ of $ H_{{\mathbb C}}^{A} $ to $ {\mathfrak h}^{1,0} $, same for $ \bar{\Omega}_{1} $ and $ {\mathfrak h}^{0,1} $. Let $ {\mathfrak h}_{{\mathbb C}}^{B} $ be the image of
this completion w.r.t. $ \Omega_{1}-\bar{\Omega}_{1} $. Identify $ {\mathfrak h}_{{\mathbb C}}^{B} $ with a subspace of $ H_{{\mathbb C}}^{B} $ via the
projection $ {\mathfrak h}_{{\mathbb C}} \to H_{{\mathbb C}}^{B} $.

The restriction on $ \Omega $ immediately imply

\begin{proposition} $ {\mathfrak h}_{{\mathbb C}} $ carries a natural symplectic structure, the
corresponding mapping $ {\mathfrak h}_{{\mathbb C}} \to {\mathfrak h}_{{\mathbb C}}^{*} $ is invertible, and $ {\mathfrak h}_{{\mathbb C}} $ is a direct sum of
Lagrangiann subspaces $ {\mathfrak h}_{{\mathbb C}}^{A} $ and $ {\mathfrak h}_{{\mathbb C}}^{B} $. \end{proposition}

Since the subspace $ {\mathfrak h}_{{\mathbb C}}^{B} $ is stable w.r.t. complex conjugation, one can
define $ {\mathfrak h}_{{\mathbb R}}^{B} $. Similarly, define $ {\mathfrak h}_{{\mathbb Z}}^{A} $ as the image of $ H_{{\mathbb Z}}^{A} $ in $ {\mathfrak h}_{{\mathbb C}}^{A} $. Note that
$ {\mathfrak h}_{{\mathbb C}}^{A} $ and $ {\mathfrak h}_{{\mathbb C}}^{B} $ are mutually dual, so $ {\mathfrak h}_{{\mathbb Z}}^{A} $ is identified with a lattice in
$ {\mathfrak h}_{{\mathbb C}}^{B} $. Say that $ v\in{\mathfrak h}_{{\mathbb Z}}^{B}\subset{\mathfrak h}_{{\mathbb C}}^{B} $ if $ \left(v,w\right)=0 $ for any $ w\in{\mathfrak h}_{{\mathbb Z}}^{A} $.

\begin{proposition} The subset $ {\mathfrak h}_{{\mathbb Z}}^{B} $ generates $ {\mathfrak h}_{{\mathbb C}}^{B} $. \end{proposition}

\begin{proof} Let $ {\bold e}_{j} $, $ j\in I_{+} $, be the natural basis of $ H_{{\mathbb Z}}^{A} $. Let $ V_{j} $ be a subspace
of $ {\mathfrak h}_{{\mathbb C}}^{A} $ generated by $ {\bold e}_{k} $, $ k\not=j $. If $ V_{j}={\mathfrak h}_{{\mathbb C}}^{A} $, then $ {\bold e}_{j} $ is a linear combination
of $ {\bold e}_{k} $, $ k\not=j $, which contradicts $ \operatorname{Im}\Omega \gg $ 0. Appropriate multiple of a normal
vector to $ V_{i} $ is in $ {\mathfrak h}_{{\mathbb Z}}^{B} $. It is obvious that these elements generate $ {\mathfrak h}_{{\mathbb C}}^{B} $. \end{proof}

\begin{definition} Say that the rigged Hopf structure is {\em integer\/} if $ {\mathfrak h}_{{\mathbb Z}}^{A} $ is
closed in $ {\mathfrak h}_{{\mathbb C}}^{A} $. \end{definition}

The condition of being integer is a restriction from below on $ \operatorname{Im}\Omega $,
say, it prohibits $ \operatorname{Im}\Omega=\operatorname{diag}\left(\lambda_{i}\right) $, $ \lambda_{i}\in l_{2} $.

We obtained 4 subspaces $ {\mathfrak h}_{{\mathbb C}}^{A} $, $ {\mathfrak h}_{{\mathbb C}}^{B} $, $ {\mathfrak h}^{1,0} $ and $ {\mathfrak h}^{0,1} $ of $ {\mathfrak h}_{{\mathbb C}} $. Note that
$ {\mathfrak h}_{{\mathbb R}}^{A} $, $ {\mathfrak h}_{{\mathbb R}}^{B} $ are generated by $ {\mathbb Z} $-lattices $ {\mathfrak h}_{{\mathbb Z}}^{A} $, $ {\mathfrak h}_{{\mathbb Z}}^{B} $ in them. If we consider $ {\mathfrak h}_{{\mathbb C}} $
and $ {\mathfrak h}^{0,1} $ as real vector spaces, then $ {\mathfrak h}_{{\mathbb C}} ={\mathfrak h}_{{\mathbb R}}\oplus{\mathfrak h}^{0,1} $. Note that $ {\mathfrak h}^{1,0} $ is
naturally identified with $ {\mathfrak h}_{{\mathbb C}}/{\mathfrak h}^{0,1} $, thus $ {\mathfrak h}^{1,0}\simeq{\mathfrak h}_{{\mathbb R}} $. Let $ {\mathfrak L} $ be the image of
$ {\mathfrak h}_{{\mathbb Z}}^{A}\oplus{\mathfrak h}_{{\mathbb Z}}^{B}\subset{\mathfrak h}_{{\mathbb R}}\subset{\mathfrak h}_{{\mathbb C}} $ in $ {\mathfrak h}_{{\mathbb C}}/{\mathfrak h}^{0,1} $, identify $ {\mathfrak L} $ with a subgroup of $ {\mathfrak h}^{1,0} $. It is clear
that $ {\mathfrak L} $ is a lattice in $ {\mathfrak h}^{1,0}\simeq{\mathfrak h}_{{\mathbb R}} $.

Start from these 4 subspaces of $ {\mathfrak h}_{{\mathbb C}} $, the complex conjugation and
symplectic structure on $ {\mathfrak h}_{{\mathbb C}} $, and lattices in the real parts of the first
two spaces. Try to reconstruct the initial rigged Hodge structure.
Suppose that $ {\mathfrak h}_{{\mathbb C}}^{A} $, $ {\mathfrak h}_{{\mathbb C}}^{B} $ are stable w.r.t. the complex conjugation, and that
$ {\mathfrak h}^{1,0} $, $ {\mathfrak h}^{0,1} $ are interchanged by complex conjugation. Suppose that the last
two subspaces project isomorphically on any subspace of the first two
(along the other one). Thus $ {\mathfrak h}^{1,0} $ and $ {\mathfrak h}^{0,1} $ are graphs of invertible
mappings $ \Omega_{1},\bar{\Omega}_{1}\colon {\mathfrak h}_{{\mathbb C}}^{A} \to {\mathfrak h}_{{\mathbb C}}^{B} $, these mappings are mutually complex
conjugate.

Now the symplectic structure on $ {\mathfrak h}_{{\mathbb C}} $ identifies $ {\mathfrak h}_{{\mathbb Z}}^{B} $ with a subset of
the dual lattice to $ {\mathfrak h}_{{\mathbb Z}}^{A} $. This reconstructs the groups $ H_{{\mathbb Z}}^{A} $, $ H_{{\mathbb Z}}^{B} $, together
with $ \Omega_{1} $ they allow to reconstruct the mapping $ \Omega $ from $ H_{{\mathbb C}}^{A}=H_{{\mathbb Z}}^{A}\otimes{\mathbb C} $ to
$ H_{{\mathbb C}}^{B}=H_{{\mathbb Z}}^{B}\otimes{\mathbb C} $. The positivity of $ \operatorname{Im}\Omega $ is translated into positivity of $ \left(\alpha,\bar{\alpha}\right) $,
$ \alpha\in{\mathfrak h}^{1,0} $, and the condition $ {\mathfrak h}^{1,0}\oplus{\mathfrak h}^{0,1}={\mathfrak h} $.

Now the results of this chapter imply

\begin{proposition} Consider a family of disks $ K_{i}\subset{\mathbb C}P^{1} $, $ i\in I $, with pairwise
conformal distances $ l_{ij} $. Let $ {\mathcal R}=\left(e^{-l_{ij}}-\delta_{ij}\right) $. Suppose that $ {\mathcal R} $ gives a
compact mapping $ l_{2} \to l_{2} $ and a bounded mapping $ l_{\infty} \to l_{\infty} $, and that some
power of $ {\mathcal R} $ gives a bounded mapping $ l_{\infty} \to l_{2} $. Consider an involution $ ':
I\to $I, $ I=I_{+}\coprod I'_{+} $, glue boundaries of disks pairwise using fraction-linear
mappings $ \varphi_{i}=\varphi_{i'}^{-1} $, let $ \left(\Omega_{jk}\right) $, $ j,k\in I_{+} $, be the matrix of periods of the
resulting curve $ \bar{M} $. Then $ \Omega_{jk} $ gives rise to a rigged Hodge structure. \end{proposition}

Note that it is natural to call the space $ {\mathfrak h}_{{\mathbb C}} $ of this rigged Hodge
structure {\em the first cohomology space\/} of the curve $ \bar{M} $. The space $ {\mathfrak h}^{1,0} $ can
be identified with $ \Gamma\left(\bar{M},\omega\right) $. As in finite-genus case, the quotient $ {\mathfrak h}^{1,0}/{\mathfrak L} $
is identified with the (bounded) Jacobian of the curve.

\bibliography{ref,outref,mathsci}

\providecommand{\bysame}{\leavevmode\hbox to3em{\hrulefill}\thinspace}
\begin{thebibliography}{10}

\bibitem{BatBloGuil95Sym}
D.~B{\"a}ttig, A.~M. Bloch, J.-C. Guillot, and T.~Kappeler, \emph{On the
  symplectic structure of the phase space for periodic {K}d{V}, {T}oda, and
  defocusing {N}{L}{S}}, Duke Math. J. \textbf{79} (1995), no.~3, 549--604.

\bibitem{BeiMan86Mum}
A.~A. Be{\u\i}linson and Yu.~I. Manin, \emph{The {M}umford form and the
  {P}olyakov measure in string theory}, Comm. Math. Phys. \textbf{107} (1986),
  no.~3, 359--376.

\bibitem{FelKnoTru96Inf}
J.~Feldman, H.~Kn\"orrer, and Trubowitz E., \emph{Infinite genus riemann
  surfaces}, Canadian {M}athematical {S}ociety 1945-1995 (James~B. Carrell and
  Ram Murty, eds.), no.~3, Canadian Mathematical Society, Ottawa, 1996,
  pp.~91--112.

\bibitem{GelShil58Gen}
I.~M. Gelfand and G.~E. Shilov, \emph{Generalized functions and operations},
  Generalized functions, vol.~1, Nauka, Moscow, 1958.

\bibitem{GelVil64Gen}
I.~M. Gelfand and N.~Ya. Vilenkin, \emph{Generalized functions. {V}ol. 4},
  Academic Press [Harcourt Brace Jovanovich Publishers], New York, 1964 [1977],
  Applications of harmonic analysis, Translated from the Russian by Amiel
  Feinstein.

\bibitem{Hor83Dis}
Lars H{\"o}rmander, \emph{The analysis of linear partial differential
  operators. {I}}, Grundlehren der Mathematischen Wissenschaften [Fundamental
  Principles of Mathematical Sciences], vol. 256, Springer-Verlag, Berlin,
  1983, Distribution theory and Fourier analysis.

\bibitem{McKTru76Hil}
H.~P. McKean and C.~Trubowitz, \emph{Hill's operator and hyperelliptic function
  theory in the presence of infinitely many branch points}, Comm. Pure Appl.
  Math. \textbf{29} (1976), no.~1, 143--226.

\bibitem{Shio86Char}
Takahiro Shiota, \emph{Characterization of {J}acobian varieties in terms of
  soliton equations}, Inventiones Mathematicae \textbf{83} (1986), 333--382.

\bibitem{ZakhFinGap2}
Ilya Zakharevich, \emph{The direct and inverse scattering transforms for a
  {S}turm - {L}iouville equation in finite interval as smooth maps. {T}he
  geometry of a set of finite gap potentials}, Inverse Methods in Action (P.~C.
  Sabatier, ed.), Springer, 1990, p.~246.

\bibitem{ZakhFinGap3}
\bysame, \emph{The spectral problems for a {S}turm - {L}iouville equation with
  a parameter. {T}he geometry of a set of finite gap potentials}, Soviet
  Mathematics---Doklady \textbf{43} (1991), no.~2, 439--444.

\end{thebibliography}
\end{document}